\documentclass[10pt,reqno,a4paper]{amsart}
\usepackage{amsthm}
\usepackage{mathtools}
\usepackage{amssymb}
\usepackage{hyperref}
\usepackage{enumitem}
\usepackage{graphicx}
\usepackage{comment}
\usepackage{xcolor}
\usepackage{bm}
\usepackage{mathrsfs} 
\usepackage[T1]{fontenc}
\usepackage[english]{babel}
\usepackage{microtype}
\numberwithin{equation}{section}
\numberwithin{figure}{section}
\theoremstyle{plain}
\newtheorem{thm}{\protect\theoremname}[section]
\theoremstyle{plain}
\newtheorem{lem}[thm]{\protect\lemmaname}
\theoremstyle{plain}
\newtheorem{prop}[thm]{\protect\propositionname}
\theoremstyle{plain}
\newtheorem{cor}[thm]{\protect\corollaryname}
\theoremstyle{remark}
\newtheorem{rem}[thm]{\protect\remarkname}
\theoremstyle{definition}
\newtheorem{definition}[thm]{\protect\definitionname}
\theoremstyle{definition}
\newtheorem{assum}[thm]{\protect\assumptionname}
\renewcommand{\Re}{\mathrm{Re}}
\renewcommand{\Im}{\mathrm{Im}}

\providecommand{\theoremname}{Theorem}
\providecommand{\corollaryname}{Corollary}
\providecommand{\definitionname}{Definition}
\providecommand{\lemmaname}{Lemma}
\providecommand{\remarkname}{Remark}
\providecommand{\propositionname}{Proposition}
\providecommand{\assumptionname}{Assumption}

\makeatother

\providecommand{\lemmaname}{Lemma}
\providecommand{\propositionname}{Proposition}
\providecommand{\remarkname}{Remark}
\providecommand{\theoremname}{Theorem}

\begin{document}
	\global\long\def\R{\mathbf{\mathbb{R}}}%
	\global\long\def\C{\mathbf{\mathbb{C}}}%
	\global\long\def\Z{\mathbf{\mathbb{Z}}}%
	\global\long\def\N{\mathbf{\mathbb{N}}}%
	\global\long\def\T{\mathbb{T}}%
	\global\long\def\Im{\mathrm{Im}}%
	\global\long\def\Re{\mathrm{Re}}%
	\global\long\def\Hc{\mathcal{H}}%
	\global\long\def\M{\mathbb{M}}%
	\global\long\def\P{\mathbb{P}}%
	\global\long\def\L{\mathcal{L}}%
	\global\long\def\F{\mathcal{\mathcal{F}}}%
	\global\long\def\s{\sigma}%
	\global\long\def\Rc{\mathcal{R}}%
	\global\long\def\W{\tilde{W}}%
	\global\long\def\G{\mathcal{G}}%
	\global\long\def\d{\partial}%
	\global\long\def\mc#1{\mathcal{\mathcal{#1}}}%
	\global\long\def\Right{\Rightarrow}%
	\global\long\def\Left{\Leftarrow}%
	\global\long\def\les{\lesssim}%
	\global\long\def\hook{\hookrightarrow}%
	\global\long\def\D{\mathbf{D}}%
	\global\long\def\rad{\mathrm{rad}}%
	\global\long\def\d{\partial}%
	\global\long\def\jp#1{\langle#1\rangle}%
	\global\long\def\norm#1{\|#1\|}%
	\global\long\def\ol#1{\overline{#1}}%
	\global\long\def\wt#1{\widehat{#1}}%
	\global\long\def\tilde#1{\widetilde{#1}}%
	\global\long\def\br#1{(#1)}%
	\global\long\def\Bb#1{\Big(#1\Big)}%
	\global\long\def\bb#1{\big(#1\big)}%
	\global\long\def\lr#1{\left(#1\right)}%
	\global\long\def\la{\lambda}%
	\global\long\def\al{\alpha}%
	\global\long\def\be{\beta}%
	\global\long\def\ga{\gamma}%
	\global\long\def\La{\Lambda}%
	\global\long\def\De{\Delta}%
	\global\long\def\na{\nabla}%
	\global\long\def\fl{\flat}%
	\global\long\def\sh{\sharp}%
	\global\long\def\calN{\mathcal{N}}%
	\global\long\def\avg{\mathrm{avg}}%
	\global\long\def\bbR{\mathbf{\mathbb{R}}}%
	\global\long\def\bbC{\mathbf{\mathbb{C}}}%
	\global\long\def\bbZ{\mathbf{\mathbb{Z}}}%
	\global\long\def\bbN{\mathbf{\mathbb{N}}}%
	\global\long\def\bbT{\mathbb{T}}%
	\global\long\def\bfD{\mathbf{D}}%
	\global\long\def\calF{\mathcal{\mathcal{F}}}%
	\global\long\def\calH{\mathcal{H}}%
	\global\long\def\calL{\mathcal{L}}%
	\global\long\def\calO{\mathcal{O}}%
	\global\long\def\calR{\mathcal{R}}%
	\global\long\def\calA{\mathcal{A}}%
    \global\long\def\calC{\mathcal{C}}%
	\global\long\def\calM{\mathcal{M}}%
	\global\long\def\calE{\mathcal{E}}%
    \global\long\def\calG{\mathcal{G}}%
	\global\long\def\calB{\mathcal{B}}%
	\global\long\def\calU{\mathcal{U}}%
	\global\long\def\calV{\mathcal{V}}%
	\global\long\def\calK{\mathcal{K}}%
	\global\long\def\calQ{\mathcal{Q}}%
    \global\long\def\calZ{\mathcal{Z}}%
    \global\long\def\calP{\mathcal{P}}%
    \global\long\def\frakc{\mathfrak{c}}%
    \global\long\def\frakC{\mathfrak{C}}%
    \global\long\def\fraka{\mathfrak{a}}%
    \global\long\def\frakA{\mathfrak{A}}%
	\global\long\def\Lmb{\Lambda}%
	\global\long\def\eps{\varepsilon}%
	\global\long\def\lmb{\lambda}%
	\global\long\def\gmm{\gamma}%
	\global\long\def\rd{\partial}%
	\global\long\def\aleq{\lesssim}%
	\global\long\def\chf{\mathbf{1}}%
	\global\long\def\td#1{\widetilde{#1}}%
	\global\long\def\sgn{\mathrm{sgn}}%
    \global\long\def\Red#1{\textcolor{red}{#1}}%
	\global\long\def\Mod{\textnormal{\textbf{Mod}}}%
	
\subjclass[2020]{35B44 (primary), 35Q55, 37K10} 

\keywords{Calogero--Moser derivative nonlinear Schrödinger equation, continuum Calogero--Moser model, classification, quantized blow-up, soliton resolution, hierarchy of conservation laws}

\title[Calogero--Moser derivative NLS]{Classification of single-bubble blow-up solutions for Calogero--Moser derivative nonlinear Schrödinger equation}

\begin{abstract}

We study the Calogero--Moser derivative nonlinear Schrödinger equation (CM-DNLS), 
a mass-critical and completely integrable dispersive model. Recent works established finite-time blow-up constructions and soliton resolution, describing the asymptotic behaviors of blow-up solutions.

In this paper, we go beyond soliton resolution and provide a sharp classification of finite-time blow-up dynamics in the \textit{single-bubble} regime. 
Assuming that a solution blows up at time $0<T<\infty$ with a single-soliton profile, we determine all possible blow-up rates.
For initial data in $H^{2L+1}(\bbR)$ with $L\ge1$, we prove a dichotomy: either the solution 
lies in a \emph{quantized regime}, where the scaling parameter satisfies
\[
    \lambda(t)\sim (T-t)^{2k},\qquad 1\le k\le L,
\]
with convergent phase and translation parameters, or it lies in an \emph{exotic regime}, where the blow-up rate satisfies $\lambda(t)\lesssim (T-t)^{2L+\frac 32}$.
To our knowledge, this is the first classification result for quantized blow-up dynamics
in the class of dispersive models. We provide a framework for identifying the quantized blow-up rates in classification problems.

The proof relies on a modulation analysis combined with the hierarchy of conservation laws provided by the complete integrability of (CM-DNLS).
However, it does not use \emph{more refined integrability-based techniques}, such as the inverse scattering method, the method of commuting flows, or the explicit formula. As a result, our analysis applies beyond the chiral solutions.

\end{abstract}

\author{Uihyeon Jeong}
\email{juih26@kaist.ac.kr}
\address{Department of Mathematical Sciences, Korea Advanced Institute of Science and Technology, 291 Daehak-ro, Yuseong-gu, Daejeon 34141, Korea}

\author{Kihyun Kim}
\email{kihyun.kim@snu.ac.kr}
\address{Department of Mathematical Sciences and Research Institute of Mathematics, Seoul National University, 1 Gwanak-ro, Gwanak-gu, Seoul 08826, Korea}

\author{Taegyu Kim}
\email{k1216300@kias.re.kr}
\address{School of Mathematics, Korea Institute for Advanced Study, 85 Hoegiro Dongdaemun-gu, Seoul 02455, Korea}

\author{Soonsik Kwon}
\email{soonsikk@kaist.edu}
\address{Department of Mathematical Sciences, Korea Advanced Institute of Science and Technology, 291 Daehak-ro, Yuseong-gu, Daejeon 34141, Korea}

\maketitle
\setcounter{tocdepth}{1}
\tableofcontents{}

% Adding sections...
\section{Introduction}
%\textcolor{red}{Edited: Comment 1, Comment 5 translation part, Comment 6 loglog, Comment 7 qunatized blow-up, Strategy of the proof}
%\textcolor{red}{Proof of Proposition~\ref{prop:radiation estimate 2}: some identities are removed}
%\textcolor{red}{The comments on the long version have been moved here.}
We consider the \textit{Calogero--Moser derivative nonlinear Schrödinger equation} (CM-DNLS)
\begin{align}\label{CMdnls}
	\begin{cases}
		i\partial_tu+\partial_{xx}u+2D_+(|u|^{2})u=0, \quad (t,x) \in \mathbb{R}\times\mathbb{R}
		\\
		u(0)=u_0
	\end{cases} \tag{CM-DNLS}
\end{align}
where $u:[-T,T]\times \mathbb{R}\to \mathbb{C}$, $D_+ = D \Pi_+$, $D=-i\partial_x$, and $\Pi_+$ is the Fourier projection onto positive frequencies. \eqref{CMdnls} is a newly introduced nonlinear Schrödinger-type equation that has attracted significant interest due to its diverse mathematical structures. Notably, \eqref{CMdnls} is mass-critical, (pseudo-)conformally invariant, and completely integrable \cite{GerardLenzmann2024CPAM}.
Despite being completely integrable, it exhibits finite-time blow-up as shown in \cite{KimKimKwon2024arxiv,JeongKim2024arXiv}.
Furthermore, \emph{soliton resolution} for \eqref{CMdnls} is also established in \cite{KimTKwon2024arxivSolResol}.

In this paper, we take a step beyond soliton resolution \textit{by classifying the dynamical behavior (including the blow-up speeds) of single-bubble blow-up solutions}. 
This serves as a first step towards a complete classification of multi-bubble blow-up dynamics.

\subsection{Calogero--Moser derivative nonlinear Schrödinger equation}

%\ \vspace{5bp}

Equation~\eqref{CMdnls} was derived in \cite{AbanovBettelheimWiegmann2009FormalContinuum} 
as a formal continuum limit of the classical Calogero--Moser system
\begin{equation*}
    \frac{d^2x_j}{dt^2}=\sum_{k\neq j}^N \frac{8}{(x_j-x_k)^3}, \qquad j=1,\dots,N,
\end{equation*}
a Hamiltonian model of inverse-square interacting particles that is completely integrable 
\cite{Calogero1971ClassicCalogerMoser,CalogeroMarchioro1974ClassicCalogerMoser,Moser1975ClassicCalogerMoser,OlshanetskyPerelomov1976Invent}.
In this sense, \eqref{CMdnls} is also known as the \emph{continuum Calogero--Moser model} (CCM). 
%Its defocusing variant, introduced earlier in \cite{PelinovskyGrimshaw1995IntermediateDefocusingCMDNLSLaxPair}, appears as a special case of the intermediate nonlinear Schrödinger equation.

Gérard and Lenzmann \cite{GerardLenzmann2024CPAM} initiated the rigorous mathematical study of \eqref{CMdnls}, establishing its complete integrability on the \emph{Hardy--Sobolev space} $H_{+}^{s}(\mathbb{R})$ defined by
\[
    H_{+}^{s}(\mathbb{R})\coloneqq H^{s}(\R)\cap L_{+}^{2}(\R),\qquad L_{+}^{2}(\R)\coloneqq\{f\in L^{2}(\mathbb{R}):\text{supp}\,\widehat{f}\subset[0,\infty)\},
\]
where $\widehat{f}$ denotes the Fourier transform of $f$.
The positive-frequency constraint $u(t)\in L^2_+$, referred to as the \emph{chiral condition} or \emph{chirality}, is preserved under the \eqref{CMdnls} flow. Within this chiral framework, \eqref{CMdnls} admits a Lax pair structure
\begin{equation}
    \partial_t\mathcal{L}_{\textnormal{Lax}}=[\mathcal{P}_{\textnormal{Lax}},\mathcal{L}_{\textnormal{Lax}}], \label{eq:LaxPair-evol}
\end{equation}
where $u$-dependent operators $\mathcal{L}_{\textnormal{Lax}}$ and $\mathcal{P}_{\textnormal{Lax}}$ are given by%\footnote{Formulations \eqref{eq:LaxPair-evol} and \eqref{eq:LaxPair} coincide with that in \cite{GerardLenzmann2024CPAM} on $H^s_+(\bbR)$ and was presented in this form in the introduction of \cite{KLV2025CAMS}.}
\begin{align}
    \mathcal{L}_{\textnormal{Lax}}=-i\partial_{x}-u\Pi_{+}\overline{u},\quad\text{and}\quad\mathcal{P}_{\textnormal{Lax}}=i\partial_{xx}+2iuD_{+}\overline{u}.\label{eq:LaxPair}
\end{align}
A similar Lax pair structure for the full intermediate NLS was found in \cite{PelinovskyGrimshaw1995IntermediateDefocusingCMDNLSLaxPair}. Beyond this, integrability of \eqref{CMdnls} further manifests through \emph{the explicit formula}, a concept originally introduced by Gérard and his collaborators in the context of other integrable models \cite{GerardGrellier2015ExplictFormulaCubicSzego1,Gerard2023BOequexplicitFormula,GerardPushnitski2024CMP}. For \eqref{CMdnls}, a rigorous derivation of the explicit formula was carried out in \cite{KLV2025CAMS} for sufficiently regular or subthreshold data.
Notably, although the explicit formula itself is restricted to the chiral framework,
the Lax pair structure \eqref{eq:LaxPair-evol} and \eqref{eq:LaxPair} holds in $H^s(\bbR)$ without chirality \cite{ChenPelinovsky2025,KimKimKwon2024arxiv}.

We briefly review symmetries and conservation laws. \eqref{CMdnls} enjoys (time and space) translation and phase rotation symmetries associated to conservation of energy, mass, and momentum:
\begin{gather*}
    \tilde{E}(u)=\frac{1}{2}\int_{\mathbb{R}}\left|\partial_x u-i\Pi_+(|u|^2)u\right|^2dx, %\label{eq:intro-E-tilde} 
    \\
    M(u)=\int_{\mathbb{R}}|u|^{2}dx, \quad \tilde{P}(u)=\Re\int_{\mathbb{R}}(\overline{u}Du-\frac{1}{2}|u|^{4})dx. \nonumber
\end{gather*}
In addition, it has Galilean invariance
\begin{align*}
	u(t,x)\mapsto e^{i\nu x-i\nu^2t}u(t,x-2\nu t), \quad \nu \in \mathbb{R},
\end{align*}
\textit{$L^{2}$-scaling symmetry} 
\begin{align*}
    u(t,x)\mapsto\lambda^{-\frac{1}{2}}u(\lambda^{-2}t,\lambda^{-1}x),\quad \lambda>0,
\end{align*}
and pseudo-conformal symmetry
\begin{equation}
    u(t,x)\mapsto \frac{1}{|t|^{1/2}}e^{ix^{2}/4t}u\left(-\frac{1}{t},\frac{x}{|t|}\right).\label{eq:pseudo-conf-transf}
\end{equation}

The \emph{ground state} (or \emph{soliton}) defined as the non-trivial zero energy solution (so it is static) plays a fundamental role in the description of global behavior of solutions. It is known (see \cite{GerardLenzmann2024CPAM}) that the ground state is unique up to scaling, phase rotation, and translation symmetries and is explicitly given by 
\begin{equation*}
    \mathcal{R}(x)=\frac{\sqrt{2}}{x+i}\in H_{+}^{1}(\mathbb{R})\quad\text{with}\quad M(\mathcal{R})=2\pi\quad\text{and}\quad\widetilde{E}(\mathcal{R})=0,%\label{eq:solition CMDNLS}
\end{equation*}
Applying the pseudo-conformal transform \eqref{eq:pseudo-conf-transf} to $\mathcal R$, one obtains an explicit finite-time blow-up solution
\begin{equation}
    S(t,x)\coloneqq \frac{1}{t^{1/2}} e^{i x^{2}/4t}\,
    \mathcal R\!\left(\frac{x}{t}\right)
    \in L^{2}(\bbR), \qquad \forall\, t>0.
    \label{eq:st def}
\end{equation}
We emphasize that this solution is \emph{non-chiral} and exhibits only \emph{limited Sobolev regularity}: $S(t)\in H^{s}(\bbR)$ for $0\le s<\tfrac12$, but not for $s\ge \tfrac12$.
In particular, $S(t)\notin H^{1}(\bbR)$ and therefore does not have finite energy.

\subsection{Long-time dynamics and soliton resolution}

%\ \vspace{5bp}

We first review the Cauchy problem for \eqref{CMdnls}. Local well-posedness for \eqref{CMdnls} in $H^s(\R)$ with $s>1/2$ was obtained by de Moura and Pilod~\cite{MouraPilod2010CMDNLSLocal}. Gérard and Lenzmann~\cite{GerardLenzmann2024CPAM} extended the local well-posedness theory to $H^s_+(\bbR)$ for $s>1/2$. More recently, local well-posedness in $H^s(\bbR)$ for $s>\frac14$ was established without chirality~\cite{CFL2025arXiv} by Chapouto, Forlano, and Laurens.

For the long-time dynamics, Gérard and Lenzmann~\cite{GerardLenzmann2024CPAM} derived several consequences of complete integrability for chiral solutions. They proved that if $M(u_0)<M(\calR)=2\pi$ and $u_0\in H^k_+ (\bbR)$ with integer $k\geq1$, then the solution is global and uniformly $H^k_+$ bounded. At the threshold $M(u_0)=M(\mathcal R)$, they showed that finite-time blow-up does not occur in $H^1$. Above the threshold, they constructed explicit $N$-soliton solutions described by a complexified Calogero--Moser system, whose $H^s$ norms blows up in infinite time, $\|u(t)\|_{H^s}\sim |t|^{2s}$ as $|t|\to\infty$.

The explicit formula based on the integrable structure was later obtained by Killip, Laurens, and Vișan~\cite{KLV2025CAMS}, who also established global well-posedness in the scaling-critical space $L^2_+(\bbR)$ for subthreshold mass $M(u_0)<M(\mathcal R)$. 
This formula was subsequently used to prove the existence of near-threshold (infinite-time) blow-up solutions~\cite{HoganKowalski2024PAA}.

In addition to these long time aspects, many authors studied other aspects of \eqref{CMdnls}. These works include the scattering theory via distorted Fourier transform~\cite{FrankRead2025arXiv}, the zero-dispersion limit~\cite{Badreddine2024SIAM}, and dynamics on nonzero backgrounds~\cite{Matsuno2023,ChenPelinovsky2025}.
We also mention several closely related models, such as periodic variants~\cite{Badreddine2024PAA,Badreddine2024AIHPC}, system extensions~\cite{BF2023PhysD,Sun2024LMPsystem}, numerical simulations~\cite{BCD2024arxivNumerical}, and the defocusing case~\cite{Chen2025arXivNonzero,Chen2025arXivScatter,ABIK2512arXiv}; see also the references therein.

Recent works by the present authors~\cite{KimKimKwon2024arxiv,JeongKim2024arXiv} constructed explicit families of smooth finite-time blow-up solutions for data near the threshold. Following these developments, soliton resolution was established in~\cite{KimTKwon2024arxivSolResol}. We recall statement in the finite-time blow-up case, which serves as a structural guideline for the present work. 

%I changed z* to u*.
\begin{thm}[Soliton resolution for finite-time blow-up solutions \cite{KimTKwon2024arxivSolResol}]\label{thm:Soliton resolution} 
    Let $u$ be a $H^1$-solution to \eqref{CMdnls} which blows up at a time $T<+\infty$. Then, there exist an integer $N\in\bbN$ with $1\leq N\leq\frac{M(u_{0})}{M( \mathcal{R})}$, a time $0<\tau<T$, modulation parameters $(\lambda_{j}(t),\gamma_{j}(t),x_{j}(t)):C^{1}([\tau,T))\to\bbR^{+}\times\bbR/2\pi\bbZ\times\bbR$ for $j=1,2,\cdots,N$, and asymptotic profile $u^{\ast}\in L^{2}$ so that $u(t)$ admits the decomposition 
    \begin{align}
        u(t)-\sum_{j=1}^{N}\frac{e^{i\gamma_j(t)}}{\sqrt{\lambda_{j}(t)}} \mathcal{R}\left(\frac{\cdot - x_{j}(t)}{\lambda_{j}(t)}\right)
        \to u^{\ast}\text{ in }L^{2}\text{ as }t\to T,\label{eq:soliton resol ungauge}
    \end{align}
    and satisfies the following properties: 
    \begin{itemize}
        \item (Asymptotic orthogonality and no bubble tree) For all $1\leq i\neq j\leq N$,
        \begin{align}
            \left|\frac{x_{i}(t)-x_{j}(t)}{\lambda_{i}(t)}\right|\to\infty\quad\text{as}\quad t\to T.%\label{eq:asymptotic orthogonality}
        \end{align}
        \item (Convergence of translation parameters) For all $j=1,2,\cdots,N$,
        \begin{equation}
            \lim_{t\to T}x_{j}(t)\eqqcolon x_{j}^{\ast}\quad \text{exists.} \label{eq:convergence translation}
        \end{equation}
        
        \item (Further information about $u^{\ast}$) We have $M(u^{\ast})=M(u_{0})-N\cdot M(\mathcal{R})$.
        In addition, if $u_{0}\in H^{1,1} \coloneqq H^1 \cap x^{-1}L^2$, then $xu^{\ast}\in L^{2}$. 
        \item (Bound on the blow-up speed) We have $\|u(t)\|_{\dot{H}^{1}}\sim\max_{j}(\la_{j}(t)^{-1})$,
        and 
        \begin{align}
            \lambda_{j}(t)\lesssim T-t \quad \text{ for all }1\le j\le N. \label{eq:lambda bdd}
        \end{align}
        
    \end{itemize}
\end{thm}
For an analogous soliton resolution statement for global solutions, see~\cite{KimTKwon2024arxivSolResol}.

%\Red{I revised a lot the following paragraphs, so please check.}

A natural but difficult question following soliton resolution is the \emph{classification of all possible asymptotic behaviors of parameters of solitons}, namely, $\lmb_j (t)$, $x_j(t)$, and $\gmm_j (t)$ in the resolution \eqref{eq:soliton resol ungauge}. Given that soliton resolution appears universally for many dispersive PDEs (and critical PDEs), the classification problem is also universal. There has been a huge literature on \emph{constructing} various scenarios, but \emph{classification} is still poorly understood. We will review some literature on this aspect in our comments on the main results. 

Returning to \eqref{CMdnls}, we may focus only on \emph{finite-time blow-up solutions}, thanks to pseudoconformal invariance and the fact that it preserves regularities $L^2$, $H^{1,1}$, and so on. As $x_j(t)$ are already known to be convergent \eqref{eq:convergence translation}, identifying the leading order behaviors of the scales $\lmb_j(t)$ and phases $\gmm_j(t)$ becomes the main goal.

We restrict our attention to the single-bubble regime ($N=1$). A complete treatment of general multi-bubble configurations remains as a major challenge. However, this single-bubble regime is still \emph{challenging} in the classification aspect because it already exhibits (countably-)\emph{infinitely many scenarios}
\begin{align*}
    \lambda(t)\sim (T-t)^{2k},\quad k\in\bbN,
\end{align*}
arising from \emph{smooth compactly supported} initial data (constructed in \cite{KimKimKwon2024arxiv,JeongKim2024arXiv}). These rates are called the \emph{quantized rates}. These constructions highlight rich possibilities of single-bubble blow-up solutions even for smooth solutions.
The purpose of this work is to formulate and obtain a complete classification of single-bubble regular blow-up solutions for \eqref{CMdnls}.

\subsection{Main result}

As mentioned before, we consider the simplest case among $N$-soliton configurations, i.e., the \emph{single-bubble case} ($N=1$).
\begin{assum}[Single-bubble blow-up]\label{assum:one soliton}
    We assume $u(t)$ is a $H^1$ finite-time blow-up solution to \eqref{CMdnls} with $N=1$ in Theorem~\ref{thm:Soliton resolution}, i.e.,
    \begin{align}
        u(t) - \frac{e^{i\gamma(t)}}{\sqrt{\lambda(t)}} \mathcal{R}\left(\frac{\cdot - x(t)}{\lambda(t)}\right) \to u^* \quad \text{in } L^2 \text{ as } t \to T<+\infty, \label{eq:one soliton cmdnls}
    \end{align}
    and $x(t)\to x^*$ as $t\to T$ (see \eqref{eq:convergence translation}). We call this scenario \textit{the single-bubble (finite-time) blow-up}.
    %\Red{I deleted ``for some $x^*$'' because $x^*$ was already mentioned in the soliton resolution theorem.}
\end{assum}

Estimate \eqref{eq:lambda bdd} provides a universal bound on the blow-up rate for $H^1$-solutions. Under Assumption~\ref{assum:one soliton}, this bound can be sharpened to
\begin{equation}\label{eq:frac 32 bound}
\lambda(t)\lesssim (T-t)^{\frac{3}{2}} .
\end{equation}
The proof of \eqref{eq:frac 32 bound} follows the same strategy as the argument leading to (1.14) in \cite[Theorem~1.1]{KimKwonOh2025AJM}. For the reader’s convenience, we present the proof in Appendix~\ref{appendix B}.

Depending on the regularity of the solution — from the energy solutions up to smooth solutions — we obtain a complete classification of possible blow-up rates.
\begin{thm}[Classification of single-bubble finite-time blow-up solutions]\label{thm:classfy CMDNLS}
    Let $L \ge 1$ be an integer, and let $u(t)$ be a solution to \eqref{CMdnls} satisfying Assumption~\ref{assum:one soliton}, with blow-up time $0<T<+\infty$ and initial data
    \begin{equation*}
        u(0)=u_0 \in H^{2L+1}(\bbR).
    \end{equation*}
    Then, the solution $u(t)$ falls into precisely one of the following scenarios:

    (1) \textbf{Quantized regime.} There exists an integer $1\leq k\leq L$ such that the modulation parameters $(\lambda(t),\gamma(t),x(t))$ in the decomposition \eqref{eq:one soliton cmdnls} satisfy
    \begin{align*}
        \lambda(t)&=(\ell+o_{t\to T}(1))(T-t)^{2k}, 
        \\
        \gamma(t)&=\gamma^*+o_{t\to T}(1),
        \\
        x(t)&=x^*+2\nu^*(T-t)+O((T-t)^2),
    \end{align*}
    for some $\ell \in (0,\infty)$, $\gamma^*\in \bbR/2\pi\bbZ$, and $\nu^* \in \bbR$.

    (2) \textbf{Exotic regime.} The scale and translation parameters $(\lambda(t),x(t))$ in the decomposition \eqref{eq:one soliton cmdnls}  satisfy
    \begin{align*}
        \lambda(t)&\lesssim (T-t)^{2L+\frac32},\\
        x(t)&=x^*+2\nu^*(T-t)+O((T-t)^2),
    \end{align*}
    for some $\nu^* \in \bbR$.
    
\end{thm}

Note that our result, in particular, applies to chiral solutions. In fact, to deal with chiral solutions, we need to work without symmetry assumption. 

%\Red{Regarding the statement:}

%\Red{- First, I deleted for some $x^*$ because this was already included in our assumption.}

%\Red{- I deleted the dependence $\ell=\ell (u)$ and $\gmm^*=\gmm^*(u)$. To mention the dependence, all $k, \ell, \gmm^*, \nu^*$ should depend on $u$.}

%\Red{- Corresponding change in Theorem 3.3?}

%\Red{- Can we have $L=0$ version with the bound $\lmb(t) \aleq (T-t)^{\frac32}$? If the proof is less than one or two pages, it seems better to be included, as a separate proposition or a remark.}

As a consequence of Theorem~\ref{thm:classfy CMDNLS}, we obtain a classification result for all \textit{smooth} single-bubble finite-time blow-up solutions.

\begin{cor}[Classification of smooth single-bubble blow-up solutions]\label{cor:smooth}
   Let $u(t)$ be a single-bubble finite-time blow-up solution to \eqref{CMdnls} as in Assumption~\ref{assum:one soliton}, with initial data 
    \begin{equation*}
        u(0)=u_0 \in H^\infty(\bbR)\coloneqq \bigcap_{k=1}^\infty H^k(\bbR).
    \end{equation*}  
   Then, either 
   
   (1) the solution $u(t)$ blows up in a quantized regime for some $k\in\bbN$, or
   
   (2) its blow-up rate is super-polynomial, i.e., $\lambda(t)\lesssim (T-t)^{s}$ for all $s>0$.
\end{cor}

In view of soliton resolution Theorem~\ref{thm:Soliton resolution}, any solution with mass $M(u) < 2M(\mathcal{R})$ cannot be decomposed into two or more solitons. This fact, together with our main Theorem~\ref{thm:classfy CMDNLS} provides a classification result for solutions of mass less than twice the mass of the soliton.

\begin{cor}[Classification of global dynamics below two soliton mass]\label{cor:two soliton under} Assume $M(u_0)<2M(\mathcal{R})$ and $u_0\in  H^{2L+1}(\bbR)$. Then, the following hold:
    \begin{enumerate}
        \item (Finite-time blow-up solutions) 
        If the solution $u(t)$ blows up in finite time, then it follows the blow-up scenario described in Theorem~\ref{thm:classfy CMDNLS}.
        
        \item (Global solutions)
        Assume further $\langle x \rangle^{\ell} u_0\in H^{2L+1-\ell}$ for all $0\le \ell \le 2L+1$. 
        If the solution $u(t)$ is global, then either $u(t)$ scatters, or there exist modulation parameters $\lambda(t)$, $\gamma(t)$, and $\nu(t)$ such that
        \begin{equation*}
            u(t) - \frac{e^{i(\gamma(t)+\nu(t)x-t\nu(t)^{2})}}{\sqrt{\lambda(t)}}
            \mathcal{R}\!\left(\frac{\cdot - 2t\nu(t)}{\lambda(t)}\right)
            - e^{it\partial_{xx}}u^{\ast}
            \to 0,
            \quad u^{\ast}\in L^{2},
        \end{equation*}
        in $L^{2}$ as $t\to\infty$, with $\nu(t)\to\nu^{\ast}$, and where the scaling
        parameter $\lambda(t)$ exhibits either a quantized rate $\lambda(t)\sim t^{-(2k-1)}$ for some integer $1\le k\le L$, in which case $\gamma(t)\to\gamma^{\ast}$, or an exotic rate $\lambda(t)\lesssim t^{-(2L+\frac12)}$ as $t\to\infty$.
    \end{enumerate}
\end{cor}

The proof of Corollary~\ref{cor:two soliton under} proceeds by applying first the pseudo-conformal transform. 
In the finite-time blow-up case, the solution reduces to the $N=1$ case of Theorem~\ref{thm:Soliton resolution} due to the mass constraint, so we can apply our main theorem.
For global solutions, one follows Step~3 of the proof of \cite[Theorem~1.2]{KimTKwon2024arxivSolResol}.

\subsection{Comments on Theorem~\ref{thm:classfy CMDNLS}}

%\ \vspace{5bp}

%\Red{Can you replace all items into Remark environments? This is easier to refer.}
%\emph{0.~Classification of quantized rates.}~
\begin{rem}[Classification of quantized rates]
To our best knowledge, Theorem~\ref{thm:classfy CMDNLS} establishes for the first time a classification result identifying infinitely many blow-up rates (the quantized ones) in dispersive models. Earlier classification results in dispersive equations, such as \cite{Kim2025JEMS} (for large symmetric data) and \cite{Raphael2005MathAnnalen, MerleRaphael2005CMP, MerleRaphael2006JAMS, MartelMerleRaphael2014Acta} (near-soliton dynamics), proved the rigidity of one or partially two blow-up rates in the single-bubble blow-up scenario. See our comments below for more discussions on classification results in dispersive equations and quantized blow-up dynamics.
\end{rem}

%\ \vspace{5bp}

%\Red{I'm wondering very much with this `Method and novelties' item. It is typically very similar to the beginning of Strategy of the proof. Perhaps it's better to delete the beginning of Strategy of the proof, say ``(We recommend the reader to see ...) before getting into details,'' and unify the current item (under the different name ``Method, difficulties, and novelties''). Anyway, I rewrote the beginning of Strategy of proof (starting the paragraph with difficulties), so one can unify this item and the revised beginning of Strategy of proof.}
%\emph{1.~Method and novelties.}~
\begin{rem}[Method, difficulties, and novelties]\label{rem:1.7 method}
Our analysis follows the general philosophy as the forward construction of blow-up dynamics developed in earlier works such as \cite{RodnianskiSterbenz2010,RaphaelRodnianski2012,RaphaelSchweyer2013CPAMHeat,RaphaelSchweyer2014AnalPDEHeatQuantized,MerleRaphaelRodnianski2013Invention,MerleRaphaelRodnianski2015CambJMath}. 
This approach has been further adapted to \eqref{CMdnls} in \cite{KimKimKwon2024arxiv,JeongKim2024arXiv}.
The same strategy has also been extended to various models; see, for example, \cite{HillairetRaphael2012AnalPDE,RaphaelSchweyer2014MathAnn,Collot2017AnalPDE,Collot2018MemAmer,KimKwon2020blowup,KimKwonOh2020blowup}. 
While relying on analytical tools developed in blow-up constructions, the present work takes a classification perspective, which is motivated by \cite{MerleRaphael2006JAMS,RS2011JAMSinhomo,MartelMerleRaphael2014Acta,Kim2025JEMS,KimMerle2025CPAM}. From the viewpoint of classification, the work~\cite{Kim2025JEMS} on the self-dual Chern--Simons--Schrödinger equation is particularly relevant to the present paper.

In addition, we interpret \eqref{CMdnls} as a system of higher-order equations and exploit this structure to obtain control of higher-order norms, as developed in Chern--Simons--Schrödinger equations \cite{KimKwonOh2020blowup,Kim2025JEMS}. %This strategy is particularly natural for \eqref{CMdnls}, since the equation is completely integrable.
One of our main ingredients is \textit{a hierarchy of associated equations and conservation laws}, as in earlier works \cite{KimKimKwon2024arxiv,JeongKim2024arXiv}.
However, we do not make use of \emph{more refined integrability-based techniques}, such as the inverse scattering method, the method of commuting flows, or the explicit formula.
As a result, our analysis does not rely on chirality and applies beyond the chiral setting.

Unlike in a constructive setting, the classification problem does not permit additional dynamical assumptions beyond the occurrence of blow-up. One must account for all possible blow-up scenarios that are compatible with the equation. Consequently, no a priori dynamical restrictions can be imposed, and in particular one cannot assume that the scaling parameter $\lambda(t)$ converges to zero monotonically.
This adds another difficulty to our problem (especially compared to  \cite{Kim2025JEMS} for a different model) because our setting already admits \emph{infinitely many possible scenarios}: the aforementioned quantized blow-up regimes. Therefore, a refined analysis, depending also on the regularity of initial data, is required to identify these quantized blow-up rates among all possible blow-up scenarios, without imposing monotonicity or other restrictive dynamical assumptions.

Our main novelty of the present work is the development of a systematic classification scheme in this setting where multiple blow-up rates are \emph{a priori} admissible. The identification of these quantized rates is a central part of the analysis. Our scheme uniquely determines the rate selected by the dynamics.
%In contrast to situations governed by situations where there is only one (or partially two) blow-up rate, the dynamics here admit a discrete family of possible rates, whose identification is a central part of the analysis.

A key difficulty in implementing this classification scheme is that \textit{nontrivial radiation effects} are present even in the one-dimensional setting.
In particular, when closing the modulation equations, radiation terms necessarily enter the resulting dynamical ODE system and must be tracked through their own time evolution.
This forces the radiation component to be incorporated into the dynamical description, yielding a closed system suitable for classification.   
\end{rem}

%\Red{It seems better not to mention rotational instability here, to be focused. So it seems better to delete:Although the assumption $\lambda(t)\to 0$ would rule out rotational instability toward global solutions,}

%\vspace{5bp} \emph{2.~Existence of scenarios in classification.}~
\begin{rem}[Existence of scenarios in classification]
To complete a full description of dynamics, it is necessary to address the existence of solutions realizing each of the scenarios allowed by the classification.

Quantized blow-up solutions with smooth (say $H^\infty$) data have been constructed in \cite{KimKimKwon2024arxiv,JeongKim2024arXiv}.
In \cite{KimKimKwon2024arxiv}, smooth \emph{chiral} blow-up solution was obtained in the case $L=1$.
By contrast, \cite{JeongKim2024arXiv} constructs smooth quantized blow-up solutions for the general case $L\ge 1$, but not in chiral class.
While a smooth chiral quantized blow-up construction for general $L$ is still missing, it is natural to expect that such solutions could be constructed by combining the ideas developed in previous two works \cite{KimKimKwon2024arxiv,JeongKim2024arXiv} and current analysis. In particular, the preparation of chiral initial data in \cite[Section~5.6]{KimKimKwon2024arxiv} should be revisited.

Regarding the exotic regime of Theorem~\ref{thm:classfy CMDNLS}, we expect the existence of exotic blow-up solutions (with initial data $u_0\in H^{2L+1}(\bbR)$) with the rate 
\[
    \lambda(t)\sim (T-t)^{2L+\frac32+\nu},
\]
where $\nu>0$ can be any real number, in contrast to being quantized. The limited regularity (not being smooth) concerns here. Although we have not constructed examples exhibiting this \emph{continuum} of rates (hence this construction remains open), the similarity between \eqref{CMdnls} and the 3D energy-critical NLS (at least in numerology as mentioned in \cite{KimKimKwon2024arxiv}) and the existence of blow-up solutions for the latter equation exhibiting a continuum of rates \cite{OrtolevaPerelman2013,Schmid2023arXiv} strongly suggest our expectation. These works moreover suggest that one \emph{cannot} expect the phase $\gmm(t)$ to converge. See also \cite{KriegerSchlagTataru2008Invent,KriegerSchlagTataru2009Duke,Perelman2014CMP} and \cite{JendrejLawrieRodriguez2022ASENS,KimKwonOh2024arXiv} for related models.

Finally, for smooth solutions, as suggested by Corollary~\ref{cor:smooth}, any exotic blow-up would necessarily occur at a super-polynomial rate. Whether such smooth exotic blow-up solutions exist at all remains an interesting open problem.
\end{rem}

%\vspace{5bp}\emph{3.~Low regularity blow-up solutions.}~
\begin{rem}[Low regularity blow-up solutions]
The explicit pseudo-conformal blow-up solution $S(t)$ defined in \eqref{eq:st def} provides a canonical example of finite-time blow-up outside the energy space $H^1(\bbR)$. Note that the bound $\lmb(t) \aleq (T-t)^{\frac 32}$ of \eqref{eq:frac 32 bound} proved in Appendix~\ref{appendix B}, which formally corresponds to the $L=0$ case of Theorem~\ref{thm:classfy CMDNLS}, \emph{rules out} the pseudo-conformal rate for any finite energy single-bubble blow-up solutions.

Formally replacing $2L+1$ by any $s>0$ in Theorem~\ref{thm:classfy CMDNLS} gives interesting statements (conjectures at this moment). The pseudo-conformal rate will be compatible when $s\leq\frac12$, and this almost matches with the regularity of $S(t)\in H^s(\bbR)$ for $s<\frac12$. Since $S(t)$ is not chiral, the construction of \emph{chiral} pseudo-conformal blow-up solutions in this low-regularity regime is an interesting open problem. Note that the local well-posedness is currently known in $H^s(\bbR)$ for $s>\frac14$ by \cite{CFL2025arXiv}.
\end{rem}

%\vspace{5bp}\emph{4.~Classification of single-soliton global solution.}

\begin{rem}[Classification of single-soliton global solution]
We assume suitable spatial decay in Corollary~\ref{cor:two soliton under} in order to apply the pseudo-conformal transform. However, we believe that spatial decay assumptions are necessary to pinpoint the quantized rates. By pseudo-conformal symmetry again, spatial decay of initial data for global solutions and regularity of data in finite-time blow-up solutions are intrinsically linked. From this viewpoint, decay assumptions are not merely technical, but play a structural role in the classification of global dynamics.

However, it should be noted that the weighted space $H^{1,1}$ does not contain the soliton $\mathcal R$ and pure multi-solitons constructed in \cite{GerardLenzmann2024CPAM} due to their slow decay. In fact, it turns out that the scale of a single soliton arising from $H^{1,1}$ initial data must \textit{contract} $\lambda(t)\lesssim t^{-\frac12}$ as $t\to\infty$ by \eqref{eq:frac 32 bound} and pseudo-conformal symmetry.
\end{rem}

%\vspace{5bp}\emph{5.~Further properties of asymptotic profile.}~
\begin{rem}[Further properties of asymptotic profile]
We expect that the asymptotic profile $u^*$ %\Red{$u^*$ or $z^*$? $u^*$ was not defined at this moment. I changed z* to u* in Thm 1.1 and Assum 1.2}
associated with quantized blow-up solutions admits a more refined description. 
In particular, one expects that $u^*\in H^{2k-\frac12-}$ and that, near $x^*$, the profile admits an expansion of the form
\begin{align*}
    u^*(x)=q|x-x^*|^{2k-2}(x-x^*)\cdot\chi(x-x^*) +o_{x\to x^*}(|x-x^*|^{2k-1}),
\end{align*}
for some $q\in \bbC$. At present, we do not attempt to justify these refined asymptotic properties; see \cite{JeongKim2024arXiv} for a discussion of the related difficulties.
\end{rem}

%\vspace{5bp}\emph{6.~Classification of blow-up dynamics in dispersive equations.}~
\begin{rem}[Classification of blow-up dynamics in dispersive equations]
We review known classification results for blow-up dynamics in other dispersive equations.
In the 1-dimensional nonlinear wave equation, the blow-up dynamics have been completely classified in \cite{MerleZaag2012AJM, CoteZaag2013CPAM}.
In higher-dimensional NLW, although the precise blow-up profiles remain unknown, the possible blow-up \emph{rates} have been classified in \cite{MerleZaag2003AJM, MerleZaag2005Math.Ann.}.

For the mass-critical NLS, a series of works by Merle and Rapha\"el \cite{MerleRaphael2005AnnMath,MerleRaphael2003GAFA,Raphael2005MathAnnalen,
MerleRaphael2004Invent,MerleRaphael2006JAMS,MerleRaphael2005CMP} established a detailed description of blow-up dynamics in the near-soliton (single-bubble) mass regime. In this setting, any finite-time blow-up is known to satisfy either $\lambda(t)\sim \sqrt{2\pi (T-t)/\log|\log(T-t)|}$ or $\lambda(t)\lesssim T-t$.
It is conjectured that the latter can be strengthened to the pseudo-conformal rate $\lambda(t)\sim T-t$, but such a rigidity statement remains open.
We also refer to \cite{RS2011JAMSinhomo} for an inhomogeneous NLS model, where a minimal-mass blow-up solution is constructed and shown to be unique.

%Another closely related example is the mass-critical generalized Korteweg--de Vries equation. In a series of works \cite{MartelMerleRaphael2014Acta,MartelMerleRaphael2015JEMS,MartelMerleRaphael2015pisa}, \Red{[cite Martel--Pilod finite-point blow-up]} the near-soliton (single-bubble) dynamics were investigated. Roughly speaking, for $H^1$ solutions with strong spatial decay on the right, one of the following scenarios occurs: the solution blows up in finite time with $\lambda(t)\sim T-t$, exits the soliton tube, or scatters to a modulated soliton with some fixed scale and velocity. Moreover, the former two regimes are shown to be stable.

%% I chose the old version (v1) of this paragraph.

Another closely related example is the mass-critical generalized Korteweg--de Vries equation.
The work of Martel--Merle--Raphaël \cite{MartelMerleRaphael2014Acta} classifies the near-soliton (single-bubble) dynamics for $H^1$ solutions with strong spatial decay to the right; the solution blows up in finite time with $\lmb(t)\sim T-t$, exits the soliton tube, or scatters to a modulated soliton with some fixed scale and velocity. Moreover, the former two regimes are shown to be stable. This classification result further complemented by the works \cite{MartelMerleRaphael2015JEMS, MartelMerleRaphael2015pisa} and also by \cite{MartelPilod2024pisa}, for example.

For the self-dual Chern--Simons--Schr\"odinger equation (CSS), it was shown under suitable symmetry assumption (equivariant with $m\ge 1$) that any blow-up must concentrate exactly one soliton, and that the blow-up rate is necessarily pseudo-conformal, namely $\lambda(t)\sim T-t$; see
\cite{KimKwon2019,KimKwon2020blowup,KimKwonOh2025AJM,Kim2025JEMS}.

Finally, although not dispersive, related classification questions have also been studied for parabolic equations using modulation analysis.
We refer to \cite{KimMerle2025CPAM} for the nonlinear heat equation (NLH) and the harmonic map heat flow. The second author and Merle classified the dynamical behaviors of multi-bubble towers.
\end{rem}

%\vspace{5bp}\emph{7.~Quantized blow-up.}~
\begin{rem}[Other quantized blow-up results]
Quantized (type-II) blow-up dynamics were first studied in the setting of energy-supercritical NLH.
The foundational mechanisms underlying quantized blow-up were first uncovered in \cite{HerreroVelazquez1992,HerreroVelazquez1994CRASPSI}, with subsequent works \cite{Mizoguchi2007Math.Ann.,Mizoguchi2011TranAMS} establishing a rigorous classification by using the maximum principle, which is not available in dispersive models. For the critical NLH, see \cite{FilippasHerreroVelazquez2000,delPinoMussoWeiZhang2020HeatQuantized,Harada2020AIHPC} for construction results.
%The analysis in this parabolic framework relies essentially on the maximum principle \cite{MatanoMerle2004CPAM}; see also \cite{FilippasHerreroVelazquez2000,delPinoMussoWeiZhang2020HeatQuantized,Harada2020AIHPC}. %MatanoMerle2004CPAM,MatanoMerle2009JFA
%\Red{I think \cite{delPinoMussoWeiZhang2020HeatQuantized,Harada2020AIHPC} should be separated.}
Starting from \cite{RaphaelSchweyer2014AnalPDEHeatQuantized,MerleRaphaelRodnianski2015CambJMath}, quantized blow-up dynamics have since been constructed via modulation analysis in a wide variety of models, including parabolic equations, dispersive equations, and compressible fluid models; see, for instance,
\cite{Collot2018MemAmer,GhoulIbrahimNguyen2018JDE,HadzicRaphael2019JEMS,CollotMerleRaphael2020JAMS,CollotGhoulMasmoudiNguyen2022CPAM,
MerleRaphaelRodnianskiSzeftel2022-1Ann.Math.,MerleRaphaelRodnianskiSzeftel2022-2Ann.Math.,MerleRaphaelRodnianskiSzeftel2022Invention,Jeong2025Anal.PDE,JeongKim2024arXiv,ShaoWeiZhang2024arXiv,ShaoWeiZhang2025Pi,BuckmasterChen2024arXiv}.%\Red{Maybe one can also cite [Collot-Merle-Raphael on type-II anisotropic blow-up] at some point, and cite recent work on complex-valued defocusing NLW (if appropriate).}
%To our best knowledge, the present work provides the first classification result for quantized blow-up dynamics in a dispersive setting.
\end{rem}

\subsection{Gauge transform}

%\ \vspace{5bp}

To analyze the blow-up dynamics described in the previous sections, it is convenient to work in a gauge-transformed form of \eqref{CMdnls}, which reveals its self-dual structures.
We introduce a \emph{gauge transform} that converts \eqref{CMdnls} into a more convenient Hamiltonian form. 
Define
\begin{align*}
	v(t,x)=-\mathcal{G}(u)(t,x)\coloneqq-u(t,x)e^{-\frac{i}{2}\int_{-\infty}^x|u(y)|^2dy}.
\end{align*}
The map $\mathcal G$ is a diffeomorphism on $H^{s}(\mathbb{R})$ for every $s\geq0$, preserving the $L^2$ norm but \textit{not chirality}. Nonetheless, it transfers all continuous symmetries of \eqref{CMdnls} and will be used throughout the analysis. Under this change of variables, $v(t,x)$ satisfies
\begin{align}\label{CMdnls-gauged}
	\begin{cases}
		i\partial_tv+\partial_{xx}v+|D|(|v|^{2})v-\frac{1}{4}|v|^4v=0, \quad (t,x) \in \mathbb{R}\times\mathbb{R},
		\\
		v(0)=v_0,
	\end{cases}\tag{$\mathcal{G}$-CM}
\end{align}
where $|D|$ denotes the Fourier multiplier associated to symbol $|\xi|$. \eqref{CMdnls-gauged} is a Hamiltonian PDE endowed with the standard symplectic form $\omega(f,g)=\Im\int_{\bbR}f\ol{g}dx$, in contrast to \eqref{CMdnls} which is associated with a non-standard symplectic structure. This standard form allows us to fully exploit the \textit{self-dual structures} described below.

The gauge transform preserves all symmetries of \eqref{CMdnls}. Consequently, \eqref{CMdnls-gauged} also retains the conservation laws of energy, mass, and momentum:
\begin{align*}
    \begin{gathered}E(v)=\frac{1}{2}\int_{\mathbb{R}}\left|\partial_{x}v+\frac{1}{2}\mathcal{H}(|v|^{2})v\right|^{2}dx,\quad
    M(v)=\int_{\mathbb{R}}|v|^{2}dx,\quad P(v)=\int_{\mathbb{R}}\Im(\overline{v}\partial_{x}v)dx,
    \end{gathered}%\label{eq:intro energy gauge}
\end{align*}
where $\mathcal{H}$ is the Hilbert transform. The energy can be rewritten by a slightly different form:
\begin{align}
    E(v)=\frac{1}{2}\int_{\mathbb{R}}\big|\td{\bfD}_{v}v\big|^{2}dx,\qquad 
    \widetilde{\mathbf{D}}_{v}=\partial_{x} +\frac{1}{2}v\mathcal{H}(\overline{v}\,\cdot\,).\label{eq:intro Energy2 self dual}
\end{align}
Following the terminology of \cite{Bogomolnyi1976}, we call $v\mapsto\widetilde{\mathbf D}_v v$ the \emph{Bogomol'nyi operator}. The Hamiltonian structure of \eqref{CMdnls-gauged} and the complete square form of the energy \eqref{eq:intro Energy2 self dual} yield that \eqref{CMdnls-gauged} can be rewritten in the nonlinear factorization
\begin{align*}
    \partial_{t}v=-iL_{v}^{*}\td{\mathbf{D}}_{v}v,%\label{eq:CMdnls gauged SelfdualForm sec1}
\end{align*}
where $L_{v}^{*}$ is the formal $L^{2}$-adjoint of the linearized operator $L_{v}$ (at $v$) of the Bogomol'nyi operator $v\mapsto\td{\mathbf{D}}_{v}v$. This relation is referred to as the \emph{self-dual structure} of \eqref{CMdnls-gauged}.

The ground state soliton $\mathcal R$ of \eqref{CMdnls} is mapped to
\begin{align*}
    Q(x)\coloneqq-\mathcal{G}(\mathcal{R})(x)=\frac{\sqrt{2}}{\sqrt{1+x^{2}}}\in H^{1}(\mathbb{R}),\quad M(Q)=2\pi,\quad E(Q)=0.%\label{eq:solition gauged equation}
\end{align*}
Hence, the soliton $Q$ is the unique solution (up to symmetries) of the Bogomol'nyi equation $\widetilde{\mathbf D}_Q Q=0$. %Moreover, $Q$ is a real, positive, and even profile, exhibiting a $\frac{1}{x}$ decay. 

Moreover, for any sufficiently regular solution $v(t)$ of \eqref{CMdnls-gauged}, one has the Lax pair structure
\begin{align}
    \partial_t\widetilde{\mathbf{D}}_{v}=[-iH_{v},\widetilde{\mathbf{D}}_{v}], 
    \qquad 
    H_v=-\partial_{xx}+\frac{1}{4}|v|^{4}-v|D|\overline{v}. \label{eq:intro Lax gauge}
\end{align}
This Lax identity gives \eqref{CMdnls-gauged} the same integrable features as \eqref{CMdnls}, in particular, the hierarchy of conservation laws that we will use. Since $\mathcal{G}$ does not preserve the chirality, the explicit formula and its corollaries from the chiral setting are not directly available here. If the initial data for \eqref{CMdnls} satisfy $u_0\in L^2_+(\bbR)$, then the corresponding gauge-transformed solution $v=-\mathcal G(u)$ still admits the explicit formula, but we will not exploit this.

We shall primarily work with \eqref{CMdnls-gauged}, taking advantage of its standard symplectic structure and the resulting self-dual formulation, as well as the non-negativity of soliton \(Q\). 
All arguments will be carried out in the gauge-transformed variables and then transferred back to \eqref{CMdnls} through \(\mathcal G^{-1}\). Hence, establishing the classification for \eqref{CMdnls-gauged} (Theorem~\ref{thm:classfy GCM}), we transfer the result back through \(\mathcal G^{-1}\) to deduce Theorem~\ref{thm:classfy CMDNLS}.

\subsection{Strategy of the proof}

%\ \vspace{5bp}

We use the notation collected in Section~\ref{subsec:notation}. Since Corollary~\ref{cor:smooth} and~\ref{cor:two soliton under} are a direct consequence of Theorem~\ref{thm:classfy CMDNLS}, we only sketch the proof of Theorem~\ref{thm:classfy CMDNLS}. Our method, difficulties, and novelties at a conceptual level is mentioned in Remark~\ref{rem:1.7 method}.

Under even symmetry for \eqref{CMdnls-gauged}, the classification scheme of identifying quantized rates simplifies and becomes more transparent. Section~\ref{sec:radial} is devoted to conveying a heuristic picture of the classification scheme.

However, our main focus is the general case. While the discussion under even symmetry serves to illustrate the underlying mechanism, it is insufficient for the full classification problem.
In this non-symmetric setting, odd radiation modes contribute nontrivially to the modulation system; to capture them, we introduce radiation parameters $\{c_j\}$ and quadratic radiation interactions $\{\frakc_{i,j}\}$, which close the system at the desired order. 

%\Red{I don't like the phrasing in the following sentence because it may make this work feel less impressive. It seems better to rephrase and move it to another place. I think our radiation parameters and quadratic radiation interactions are more nontrivial than these works. Using the term ``near-soliton'' may also make this work feel less impressive.} In this sense, our analysis focuses on the general single-bubble dynamics in the spirit of \cite{MerleRaphael2005AnnMath,Raphael2005MathAnnalen,RS2011JAMSinhomo,MartelMerleRaphael2014Acta}.

\vspace{5bp}

\emph{1.~Gauge correspondence and single-bubble setting.}~
We first reduce the problem to the gauge-transformed formulation \eqref{CMdnls-gauged}.
Specifically, for a single-bubble blow-up solution $u(t)$ satisfying Assumption~\ref{assum:one soliton}, we apply the gauge transform $\calG$, which yields a corresponding solution $v(t)$ to \eqref{CMdnls-gauged} satisfying Assumption~\ref{assum:one soliton GCM}. The equivalence between the two formulations follows from Lemma~\ref{lem:gauge equiv}, whose proof is straightforward.
We then establish the classification result for \eqref{CMdnls-gauged} (stated later as Theorem \ref{thm:classfy GCM}) and transfer it back via the inverse transform to obtain Theorem \ref{thm:classfy CMDNLS}. We henceforth work with \eqref{CMdnls-gauged}.

Let $v$ be a solution to \eqref{CMdnls-gauged} satisfying Assumption~\ref{assum:one soliton GCM} and blowing up forwards in time at $T<+\infty$.  
The variational structure of the energy combined with the argument in the proof of soliton resolution~\cite{KimTKwon2024arxivSolResol} guarantees that, after rescaling, $v(t)$ remains close to the ground state $Q$ in the $\dot H^1$-topology as $t\to T$. 
Consequently, for $t$ sufficiently close to the blow-up time, $v(t)$ can be expressed as
\begin{equation*}
    v(t,x) =\frac{e^{i\gamma(t)}}{\sqrt{\lambda(t)}}
    \bigl(Q+\wt\eps(t)\bigr)\left(\frac{x-x(t)}{\lambda(t)}\right), \quad
    \lambda(t)\approx \frac{\|Q\|_{\dot H^1}}{\|v(t)\|_{\dot H^1}}.
\end{equation*}
Under the single-bubble blow-up assumption, we obtain the coercivity estimate
\begin{equation*}
    \|\wt\eps\|_{\dot H^1}^2\lesssim E(Q+\wt\eps)=\lambda^2E(v_0),
\end{equation*}
which provides the fundamental $\dot H^1$-level control that forms the starting point of the subsequent analysis.

\vspace{5bp}

\emph{2.~Nonlinear adapted variables and setup for the modulation analysis.}~Following the idea of \cite{KimKwonOh2020blowup,Kim2025JEMS} and of \cite{KimKimKwon2024arxiv,JeongKim2024arXiv} for \eqref{CMdnls}, we introduce the nonlinear variables $\td\bfD_v^j v = v_j$ by iteratively applying the Lax operator to $v$. 
In our setting, these nonlinear variables are designed to exploit the integrable structure of the equation. 
They preserve a nonlinear algebraic structure that allows higher-order energies to be expressed in terms of the first variable, making it possible to carry out the analysis without relying on coercivity for higher-order linearized operators.

The original equation \eqref{CMdnls-gauged} can be written in the compact form
\begin{align*}
    (\partial_t + iH_v)v = 0, %\label{eq:G-CM rewrite}
\end{align*}
which corresponds to the case $j=0$. 
Using the Lax structure \eqref{eq:intro Lax gauge}, one then verifies that each $v_j$ satisfies the same type of evolution equation,
\begin{align*}
    (\partial_t + iH_v)v_j = 0. %\label{eq:intro wk equ}
\end{align*}
To analyze the flow, we use the renormalized coordinates $(s,y)$ defined by 
\[
\frac{ds}{dt}=\frac{1}{\lambda^{2}(t)}\qquad\text{and}\qquad y=\frac{x-x(t)}{\lambda(t)},
\]
and the renormalized functions $w$ and $w_j$ defined by 
\begin{align*}
    w(s,y)&\coloneqq\lambda^{\frac{1}{2}}(t)e^{-i\gamma(t)}v(t,\lambda(t)y+x(t))|_{t=t(s)},
    \\
    w_{j}(s,y)&\coloneqq\td\bfD_{w}^jw=\lambda^{\frac{2j+1}{2}}(t) e^{-i\gamma(t)}(\td\bfD_{v}^jv)(t,\lambda(t)y+x(t))|_{t=t(s)}.
\end{align*}
One can find the evolution equations of $w$ and $w_{j}$ in \eqref{eq:wj equ}. The full family $w_j$ is introduced for two main reasons. First, these variables capture the hierarchy of conserved quantities arising from the Lax structure. Second, they explicitly capture the generalized kernel dynamics of the linearized operator around $Q$ in the $w$-variable. A more detailed discussion can be found in \cite{JeongKim2024arXiv}.

\vspace{5bp}

\emph{3.~Decompositions and coercivity estimates.}~We decompose the renormalized variables as
\begin{equation}
    \begin{split}
        w&=Q+\wt\eps
        \\
        &=Q+T+\eps,
    \end{split}
    \qquad
    \begin{split}
        w_j&=P_{j,e}+\wt\eps_j
        \\
        &=P_j+\eps_j,
    \end{split} \label{eq:intro decom}
\end{equation}
with profiles
\begin{align*}
    T=i\nu_0\frac{y}{2}Q-ib_1\frac{y^{2}}{4}Q-\eta_1\frac{1+y^{2}}{4}Q,
    \qquad
    P_j&=c_jQ+\bm{\beta}_jyQ,\qquad P_{j,e}=c_jQ,
\end{align*}
and
\begin{gather*}
    b_1\coloneqq -2\Im(\bm{\beta}_1),\quad \eta_1\coloneqq -2\Re(\bm{\beta}_1),\quad \bm{\beta}_1=-\tfrac12(ib_1+\eta_1),
    \\
    \nu_0\coloneqq 2\Im(c_1),\quad \mu_0\coloneqq 2\Re(c_1),\quad c_1=\tfrac12(i\nu_0+\mu_0).
\end{gather*}

As in \cite{KimKimKwon2024arxiv}, these profiles are defined without spatial cut-offs and we carry out the decompositions only in the topologies where each profile is well defined. The decompositions are fixed through orthogonality conditions on $\eps$ and $\eps_j$ transversal to the corresponding kernel elements (of relevant linearized operators), which allow coercivity estimates. This will allow us to transfer $L^2$-bounds on $w_j$ from the hierarchy of conservation laws to suitable Sobolev controls on $\eps$ and $\eps_j$.

As in \cite{JeongKim2024arXiv}, we begin with the controls of $w_j$ naturally obtained from the hierarchy of conservation laws. Conservation of the $L^2$-norm of each nonlinear variable $v_j$ yields the \emph{sharp scaling bounds} on $w_j$: 
\begin{equation*}
    \|w_j(s)\|_{L^2} = \lmb^j (s)\| v_j (s) \|_{L^2} \sim_{v_{j}(0)} \lambda^j(s).
\end{equation*}
Using the repulsivity identities \eqref{eq:BQBQstar equal I} and \eqref{eq:DQ BQ decomp}, one deduces
\begin{align*}
    \|\td{\bfD}_Q w_{j}\|_{L^2} \approx \|w_{j+1}\|_{L^2} \lesssim \lambda^{j+1},
    \qquad
    \|\td{\bfD}_Q^2 w_{j}\|_{L^2} \approx \|w_{j+2}\|_{L^2} \lesssim \lambda^{j+2}.
\end{align*}
After separating the (generalized) kernel directions of $\td{\bfD}_Q$, corresponding to the profile components $Q$ and $yQ$ encoded in $P_j$, one obtains the coercivity estimates for the radiation:
\begin{equation*}
    \|\eps_{j}\|_{\dot{\calH}^1} \lesssim \lambda^{j+1},
\end{equation*}
which provides quantitative higher-order control. This is a key advantage of the conservation laws, without requiring bootstrapping argument. 

At the $\dot H^2$ level, an additional degeneracy appears, since $(1+y^2)Q$ belongs to the kernel of $\td{\bfD}_Q^2$ in addition to the kernel directions of $\td{\bfD}_Q$.
Accordingly, to ensure coercivity at the $\dot H^2$ level, we refine the decomposition of $w_j$ by subtracting the profile
\begin{align*}
    (\bm{\beta}_{j+1}-\tfrac12\bm{\beta}_1c_j)(1+y^2)Q
\end{align*}
from the radiation part.
With this refinement, one can still work within the same decomposition framework and obtain the scale-sharp bound
\begin{equation*}
    \|\eps_j-(\bm{\beta}_{j+1}-\tfrac12\bm{\beta}_1c_j)(1+y^2)Q\|_{\dot{\calH}^2}\lesssim \lambda^{j+2}.
\end{equation*}

However, at the $\dot H^3$-level, the previous decompositions are insufficient to obtain the coercivity, motivating the introduction of new decompositions.
Although the hierarchy of conservation laws yields nonlinear higher-order energy bounds, these differ from the $\dot H^3$ coercivity bounds coming from the linear analysis by an error of order $\lambda^2$.
While this mismatch is harmless at the $\dot H^2$-level, it becomes a genuine obstruction for the part of the argument where we control the radiation dynamics without even symmetry. 

%More precisely, the coercivity requirements differ between the subsequent steps.
%For the modulation analysis in Step~4, it is sufficient to propagate $\dot H^2$-level control.
%By contrast, the analysis of radiation effects in Step~5 requires access to $\dot H^3$-level coercivity, though only for the even component.

To compensate for this $\lambda^2$-level mismatch, we introduce the \emph{check} $\check{\ }$-decomposition
\begin{equation*}
    w = Q+\check T+\check\eps,
\end{equation*}
where $\check T=T(\check\nu_0,\check b_1,\check\eta_1)$ has the same functional form as $T$, but the check parameters $(\check\nu_0,\check b_1,\check\eta_1)$ are defined directly at the level of the $w$-variable through suitable orthogonality conditions.
At the $\dot H^3$-level, an additional $Q$-component appears as a leading correction in the even part, and we therefore work with the compensated remainder
\begin{equation}
    \|\td{\bfD}_Q^3(\check{\eps}_e-\tfrac{\nu_0^2}{8}Q)\|_{L^2}\sim\|\check{\eps}_e -\tfrac{\nu_0^2}{8}Q\|_{\dot{\mathcal{H}}^{3}}\lesssim \lambda^3. \label{eq:H3 eps intro}
\end{equation}

For the higher-order variables $w_j$, we likewise refine the decomposition by extracting, in addition, a $(1+y^2)Q$-component at the next level:
\begin{equation}
    w_j = P_j+\check{\bm{\beta}}_{j+1}(1+y^2)Q+\check\eps_j \quad \text{with} \quad \|\check\eps_j\|_{\dot{\calH}^2}\lesssim \lambda^{j+2}, \quad \|\check\eps_{j,e}\|_{\dot{\calH}^3}\lesssim \lambda^{j+2+\frac12}, \label{eq:H3 epsj intro}
\end{equation}
and the refined coefficient satisfies $\check{\bm{\beta}}_{j+1}=\bm{\beta}_{j+1}-\tfrac12\bm{\beta}_1c_j+O(\lambda^{j+2})$.
This refinement restores $\dot H^3$-level coercivity beyond what can be achieved from the primary decomposition alone. %and is a key input in the non-symmetric radiation estimates.

\vspace{5bp}

\emph{4.~Modulation estimates.}
We extract the evolution laws of the modulation parameters
$\lambda,\gamma,x$ and $\{\bm{\beta}_j\}$ from the renormalized system for $w$ and $w_j$,
\begin{align}
    &(\partial_s-\tfrac{\lambda_s}{\lambda}\Lambda+\gamma_s i-\tfrac{x_s}{\lambda}\partial_y)w+iL_w^*w_1=0, \label{eq:w equ intro}
    \\
    &(\partial_s-\tfrac{\lambda_s}{\lambda}\Lambda_{-j}+\gamma_s i-\tfrac{x_s}{\lambda}\partial_y)w_j-iw_{j+2}
    +\tfrac{i}{2}[w_1\mathcal{H}(\overline{w}w_j)-w\mathcal{H}(\overline{w_1}w_j)]=0. \label{eq:wj equ intro}
\end{align}
We substitute the decompositions \eqref{eq:intro decom} and split the resulting identities into the explicit contributions produced by the profiles $(Q,T,P_{j,o}=\bm{\beta}_jyQ)$ and the contributions involving the radiation parts $(P_{j,e}=c_jQ,\eps,\eps_j)$.

A key point in the non-symmetric setting is that the nonlinearity produces additional leading-order terms coming from the nonlinear radiation terms, which are absent under even symmetry.
These terms deform the modulation laws already at the level of the main identities: they are quantified in \eqref{eq:main nonlinear eps j odd} and \eqref{eq:main nonlinear Q}.
After extracting these leading contributions, the remaining nonlinear remainders are controlled by the $\dot H^2$-level bounds from Step~3. Consequently,
\begin{equation}
    \begin{aligned}
        \eqref{eq:w equ intro}&\approx 
        -\left(\frac{\lambda_s}{\lambda}+b_1\right)\Lambda Q+\left(\gamma_s-\frac{\eta_1}{2}\right)iQ-\left(\frac{x_s}{\lambda}-\nu_0\right)Q_x\approx0,
        \\
        \eqref{eq:wj equ intro}_o&\approx \left(\left(\partial_s-\frac{\lambda_s}{\lambda}\left(j+\frac12\right)+i\gamma_s\right)\bm{\beta}_j
        -i\bm{\beta}_{j+2}-\nu_0\bm{\beta}_{j+1}+ic_j\bm{\beta}_2\right)yQ \approx 0.
    \end{aligned}\label{eq:mod esti intro}
\end{equation}
In particular, these relations yield the corresponding approximate modulation system,
whose quantitative form is recorded in Proposition~\ref{prop:modj}.

\vspace{5bp}

\emph{5.~Radiation estimates.}
Now we find the evolution equations for $c_j$ to close the modulation system \eqref{eq:mod esti intro}. Since $P_{j,e}=c_jQ$, taking the even part of \eqref{eq:wj equ intro} gives
\begin{align}
    \eqref{eq:wj equ intro}_e\approx 
    \left((c_j)_s-\tfrac{\lambda_s}{\lambda}jc_j-ic_{j+2}-a_j\right)Q\approx0, \label{eq:rad esti intro 1}
\end{align}
where $a_j$ consists of explicit combinations of the radiation parameters $c_j$ and \emph{quadratic radiation interactions} $\frakc_{i,j}$, which we explain below. The expression of $a_j$ is given in \eqref{eq:aj def}. The $\dot H^3$-level bounds \eqref{eq:H3 eps intro} and \eqref{eq:H3 epsj intro} obtained in Step~3 provide the smallness needed to justify \eqref{eq:rad esti intro 1}.

A key point is that $(c_j)_s$ has quadratic terms such as
\begin{equation}
    \int yQ^2\overline{\wt\eps_1}\wt\eps_j, \label{eq:intro quad interaction}
\end{equation}
whose size is comparable to $\lambda^{j+2}$ (or $c_{j+2}$).
These terms are not negligible and must be included in the leading-order dynamical system as new components (or variables). However, a direct treatment of these terms does not yield a closed system; differentiating them in time produces new quadratic terms that cannot be expressed solely in terms of $\{c_j\}$ and the original interactions \eqref{eq:intro quad interaction}.

To overcome this difficulty, we introduce modified quadratic quantities defined by 
\begin{equation*}
    \frakc_{i,j}\coloneqq \int \varphi\overline{w_i}w_j, \qquad \varphi=\frac12 yQ^2-3yQ^4+2yQ^6,
\end{equation*}
which we call \emph{quadratic radiation interactions}.
The weight $\varphi$ is chosen so that all profile contributions (and $c_j$ contributions) in $\overline{w_i}w_j$ cancel both in $\frakc_{i,j}$ and in $(\frakc_{i,j})_s$.
As a consequence, $\frakc_{i,j}$ still captures the dominant quadratic radiation,
\begin{equation*}
    \frakc_{i,j}\approx \frac12\int yQ^2 \overline{\wt\eps_i}\wt\eps_j,
\end{equation*}
but is designed so that its time derivative is (up to leading order) fully described by the previously introduced parameters, leading to a closed dynamical system.

In the derivation of $(\mathfrak c_{i,j})_s$, the nonlinear term gives rise to a 4-linear expression. After decomposing each factor into its profile and radiation parts, the resulting terms can be organized into three types: products of two quadratic radiation terms, single quadratic radiation terms, and perturbative errors.

More precisely, purely profile interactions are eliminated at the outset by the choice of the weight $\varphi$. Among the remaining contributions, the purely radiative case corresponds exactly to bilinear products of the variables $\mathfrak c_{i,j}$. This reduction follows from the symmetry of the underlying 4-linear expression together with a commutator identity for the Hilbert transform.
A specific mixed interaction produces weighted quadratic radiation terms. By an integration-by-parts argument with a suitably chosen weight, these terms can be rewritten in the same form as $\mathfrak c_{i,j}$.
All remaining contributions are perturbative in nature and can be treated as error terms.

As a result, up to perturbative errors, the evolution of $\mathfrak c_{i,j}$ closes within the system formed by the previously introduced radiation parameters $\{c_j\}$ and $\{\mathfrak c_{i,j}\}$, and no new quantities are generated. We finally arrive at the evolution laws for $\frakc_{i,j}$,
\begin{align}
    (\mathfrak{c}_{i,j})_s-\tfrac{\lambda_s}{\lambda}(i+j+1)\mathfrak{c}_{i,j}+i\frakc_{i,j+2} -i\frakc_{i+2,j}- \fraka_{i,j}\approx 0 \label{eq:rad esti intro 2}
\end{align}
where $\fraka_{i,j}$ is given explicitly in \eqref{eq:Aij def} and consists of linear combinations of the radiation parameters $\{c_j\}$ and $\{\frakc_{i,j}\}$.

Altogether, the evolution laws for $\{c_j\}$ and $\{\mathfrak c_{i,j}\}$ provide the missing radiation estimates that obstruct the closure of the modulation system.
In particular, combining the modulation equations for $(\lambda,\gamma,x,\{\bm\beta_j\})$ in \eqref{eq:mod esti intro} with the evolution laws for $\{c_j\}$ and $\{\frakc_{i,j}\}$ yields a finite-dimensional closed system for
\begin{equation*}
    (\lambda,\gamma,x,\{\bm\beta_j\},\{c_j\},\{\frakc_{i,j}\}).
\end{equation*}
This reduction allows us to pass from the original PDE dynamics to the analysis of the dynamics of this closed ODE system, which forms the basis of the classification argument.

\vspace{5bp}

\emph{6.~Identification of quantized blow-up rates.}
In the final step, we identify the blow-up rate by analyzing the finite-dimensional dynamical system obtained from the modulation estimates \eqref{eq:mod esti intro} and radiation estimates \eqref{eq:rad esti intro 1}, \eqref{eq:rad esti intro 2}.
The guiding idea is to reorganize the coupled system into variables that are uniformly bounded and whose evolution laws exhibit a simple recursive structure up to small errors.

We first renormalize the modulation and radiation parameters by factoring out their natural powers of $\lambda$. This leads, for instance, to the normalized variable
\begin{equation*}
    Z_0\coloneqq \lambda^{\frac12}e^{i\gamma},
\end{equation*}
as well as $Z_j$, $C_j$, $\frakC_{i,j}$, and $\frakA_{i,j}$.
Their precise definitions are given in \eqref{eq:def normalized parameters}.
In these variables, the approximate modulation equations \eqref{eq:mod esti intro} and the radiation evolution equations \eqref{eq:rad esti intro 1}--\eqref{eq:rad esti intro 2} take the following simplified form:
\begin{align*}
    (Z_0)_t+iZ_1&\approx 0,%\qquad x_t-2\Im(C_1)&\approx 0,
    &( Z_j)_t-i Z_{j+2}-2\Im (C_1) Z_{j+1}+iC_j Z_{2}&\approx 0,
    \\
    (C_j)_t-iC_{j+2}-A_{j} &\approx 0,
    &(\mathfrak{C}_{i,j})_t+i\frakC_{i,j+2} -i\frakC_{i+2,j}- \frakA_{i,j}&\approx 0.
\end{align*}
The structure of this system shows that the leading contribution to $(Z_0)_t$ is given by $-iZ_1$, and that the leading contribution to the time derivative of each quantity is expressed in terms of variables with higher indices.
This observation motivates the introduction of corrected modulation quantities $\{\Omega_k\}_{k\ge1}$, defined so as to extract these leading terms in a systematic way.

We set $\Omega_1\coloneqq -iZ_1$ and define $\Omega_{k+1}$ inductively so that $\Omega_{k+1}$ captures the leading contribution in $(\Omega_k)_t$ generated by the above normalized system.
With this choice, the evolution of $Z_0$ and of the corrected quantities $\{\Omega_k\}$ satisfies an approximated ODE system,
\begin{equation*}
    (Z_0)_t- \Omega_1\approx0,\qquad
    (\Omega_k)_t- \Omega_{k+1}\approx0. %\label{eq:intro Omega}
\end{equation*}
We then integrate the system as $t\to T$.
The limits $\Omega_k^*=\lim_{t\to T}\Omega_k(t)$ exist, and the smallest index $k$ for which $\Omega_k^*\neq 0$ determines the leading asymptotic behavior of $Z_0(t)$, and hence the blow-up rate.
This yields the quantized rates $\lambda(t)\sim (T-t)^{2k}$.
If all $\Omega_k^*$ vanish up to the available order, then $\lambda(t)$ decays faster than the last quantized rate, corresponding to the exotic regime.

%A key point is that \eqref{eq:intro Omega} are justified modulo an error of size $\lambda^{\frac12}$, as quantified in Lemma~\ref{lem:normalized estimate}. This control is essentially saturated for the present argument: the iterative integration of the system up to the blow-up time~$T$ can be closed with error terms of size at most $\lambda^{\frac12-}$.

\vspace{5pt}
\noindent\mbox{\textbf{Organization of the paper.~}}In Section~\ref{sec:notation}, we collect the
notation and preliminary facts for the analysis. In Section~\ref{sec:gague correspond}, we establish the gauge correspondence and set up the analytical framework for the single-bubble blow-up analysis. In Section~\ref{sec:radial}, we give a heuristic account of the classification scheme under even symmetry. In Section~\ref{sec:decom and coerc}, we introduce suitable decompositions and establish coercivity estimates. In Sections~\ref{sec:mod esti} and~\ref{sec:rad esti}, we derive the modulation and radiation estimates, respectively. Finally, we prove Theorem~\ref{thm:classfy CMDNLS} In Section~\ref{sec:classification}. For the remaining sections, we prove the coercivity estimates in Section~\ref{sec:coer pf} and the nonlinear estimates in Sections~\ref{sec:nonlin pf}, respectively.

\vspace{5pt}
\noindent\mbox{\textbf{Acknowledgements.~}}The authors would like to thank Charles Collot and Van Tien Nguyen for valuable discussions. K.~Kim's research was supported by the New Faculty Startup Fund from Seoul National University, the POSCO Science Fellowship of POSCO TJ Park Foundation, and the National Research Foundation of Korea, RS-2025-00523523. T.~Kim was partially supported by the National Research Foundation of Korea, RS-2019-NR040050 and RS-2022-NR069873 during the initial stage of this work and by a KIAS Individual Grant (MG105201) at Korea Institute for Advanced Study. U.~Jeong and S.~Kwon were partially supported by the National Research Foundation of Korea, RS-2019-NR040050 and RS-2022-NR069873.

\section{Notation and preliminaries}\label{sec:notation}
\subsection{Notations and definitions}\label{subsec:notation}
Here, we collect notations, function spaces, and operator definitions. Most of them were introduced in the authors' earlier works \cite{KimKimKwon2024arxiv,KimTKwon2024arxivSolResol,JeongKim2024arXiv} for \eqref{CMdnls}. 

\vspace{5pt}

\noindent\underline{\textit{General notations.}}
We denote by $\bbN$ the set of positive integers. We also denote $\bbN_0\coloneqq \bbN\cup\{0\}$ and $\frac{1}{2}\mathbb{N}\coloneqq \{\frac{n}{2}:n\in\mathbb{N}\}$. 
For quantities $A\in\bbC$ and $B\geq0$, we write $A\lesssim B$
if $|A|\leq CB$ holds for some implicit constant $C$. For $A,B\geq0$,
we write $A\sim B$ if $A\lesssim B$ and $B\lesssim A$. If $C$
is allowed to depend on some parameters, then we write them as subscripts of $\lesssim,\sim,\gtrsim$ to indicate the dependence. Unless otherwise specified, the dependence on the mass $M$ and the conserved energies $E_j$ will be left implicit, as making it explicit would unnecessarily complicate the notation.

We write $\langle x\rangle\coloneqq(1+x^{2})^{\frac{1}{2}}$. 
We fix a smooth even cut-off function $\chi$ such that $\chi(x)=1$ for $|x|\leq1$, and $\chi(x)=0$ for $|x|\geq2$.
For $R>0$, we define $\chi_{R}$ by $\chi_{R}(x)=\chi(R^{-1}x)$.
We also denote the sharp cut-off function on a set $A$ by $\mathbf{1}_{A}$.
We denote by $\delta_{jk}$ the Kronecker--Delta symbol, i.e., $\delta_{jk}=1$ if $j=k$ and $\delta_{jk}=0$ if $j\neq k$. 

Given a function $f$ on $\bbR$, let $f_o$ and $f_e$ denote its odd and even parts,
respectively.
For a function space $X$, we write $X_e$ and $X_o$ for its even and odd subspaces, respectively.

We denote by some functions $\calK_j$ for $1\leq j\leq 6$ and $\mathring \calK_i$ for $1\leq i\leq 3$
\begin{align}
    &\begin{split}
        \mathcal{K}_{1} & =\Lmb Q,\\
        \mathcal{K}_{4} & =-ix^{2}Q,
        \\\mathring \calK_1 &=Q,
    \end{split}
    &\begin{split}
        \mathcal{K}_{2} & =iQ,\\
        \mathcal{K}_{5} & =-(1+x^{2})Q,
        \\\mathring \calK_2 &=xQ,
    \end{split}
    &\begin{split}
        \mathcal{K}_{3} & =Q_{x},\\
        \mathcal{K}_{6} & =ixQ,
        \\\mathring \calK_3 &=(1+x^2)Q.
    \end{split} \label{eq:kernel element1}
\end{align}

Denote $L^{p}(\mathbb{R})$ and $H^{s}(\bbR)$ by the standard $L^{p}$ space and Sobolev spaces on $\bbR$, respectively. As we work on $\bbR$, we often omit $\bbR$ when there is no confusion. 

The Fourier transform (on $\bbR$) is denoted by 
\begin{align*}
	\mathcal{F}(f)(\xi)=\widehat{f}(\xi)\coloneqq\int_{\mathbb{R}}f(x)e^{-ix\xi}dx.
\end{align*}
The inverse Fourier transform is then given by $\mathcal{F}^{-1}(f)(x)\coloneqq\tfrac{1}{2\pi}\int_{\R}\wt f(\xi)e^{ix\xi}d\xi$.
We denote $|D|$ by the Fourier multiplier with symbol $|\xi|$, that
is, $|D|\coloneqq\mathcal{F}^{-1}|\xi|\mathcal{F}$. We denote by
$\mathcal{H}$ the Hilbert transform: 
\begin{align}
	\mathcal{H}f\coloneqq\left(\frac{1}{\pi}\text{p.v.}\frac{1}{x}\right)*f=\mathcal{F}^{-1}(-i\text{sgn}(\xi))\mathcal{F}f.\label{eq:HilbertFourierFormula}
\end{align}
In view of \eqref{eq:HilbertFourierFormula}, we have $\partial_{x}\mathcal{H}f=|D|f=\mathcal{H}\partial_{x}f$
for $f\in H^{1}$. 
%We denote by $\Pi_{+}$ the Cauchy--Szeg\H{o} projection from $L^{2}(\mathbb{R})$ onto the Hardy space $L_{+}^{2}(\mathbb{R})=\{f\in L^{2}(\bbR):\mathrm{supp}\wt f\subset[0,\infty)\}$:
%\begin{align*}
%	\Pi_{+}f=\mathcal{F}^{-1}\textbf{1}_{\xi>0}\mathcal{F}f.
%\end{align*}
%Then, we have $\Pi_{+}=\frac{1}{2}(1+i\mathcal{H}).$ We denote $D=-i\rd_{x}$, $D_{+}=D\Pi_{+}$, and $|D|=\calH\rd_{x}$. 
%We denote the Hardy-Sobolev space by $H_{+}^{s}\coloneqq H^{s}\cap L_{+}^{2}$, where $H^{s}$ is the usual Sobolev space based on $L^{2}$. 
%We also denote the Schwartz space $\mathcal{S}$.

We will use the \emph{real inner product} and \emph{complex inner product} defined by 
\begin{align*}
	(f,g)_{r}\coloneqq\Re\int_{\R}f\overline{g}dx,\qquad 
    (f,g)_{c}\coloneqq \int_{\R}f\overline{g}dx,
\end{align*}respectively.
We also use the modulated form of a function $f$ by parameters $\lmb$,
$\gmm$, and $x$: 
\begin{equation*}
	[f]_{\la,\ga,x}\coloneqq\frac{e^{i\ga}}{\la^{1/2}}f\lr{\frac{\cdot-x}{\la}}.
\end{equation*}
We denote by $\Lambda_{s}$ the $\dot{H}^{s}$-scaling generator in
$\R$ as 
\begin{align*}
	\Lambda_{s}f  \coloneqq\frac{d}{d\lambda}\bigg|_{\lambda=1}\lambda^{\frac{1}{2}-s}f(\lambda\cdot)=\left(\frac{1}{2}-s+x\partial_{x}\right)f,\quad 
	  \Lmb\coloneqq\Lmb_{0}.
\end{align*}

\vspace{5pt}

\noindent\underline{\textit{Function spaces.}}
We use the adapted Sobolev spaces $\dot{\mathcal{H}}^{k}$ for $k\in \bbN$ whose norm is defined by 
\begin{align*}
	\|f\|_{\dot{\mathcal{H}}^{k}}^{2}  \coloneqq \sum_{j=0}^{k}\left\Vert \langle x\rangle^{-j}\partial_{x}^{k-j}f\right\Vert_{L^{2}}^{2}.
\end{align*}
Note that $\dot{\mathcal{H}}^{k}\hookrightarrow\dot{H}^{k}$, $\dot{\mathcal{H}}^{k}\cap L^{2}=H^{k}$, and $\dot{\mathcal{H}}^{k}$ compactly embeds into $H_{\textnormal{loc}}^{k-1}$.
We also define 
\begin{align*}
	\|f\|_{\dot\calH^{k,\infty}}
	\coloneqq \sum_{j=0}^k\left\Vert\langle x\rangle^{-j}\partial_x^{k-j}f\right\Vert_{L^\infty}.
\end{align*}

\vspace{5pt}

\noindent\underline{\textit{Linear operators.}}
We denote some linear operator (at $v$ or at $Q$) by 
\begin{equation*}
    \mathbf{D}_{v}f  \coloneqq\partial_{x}f+\tfrac{1}{2}\mathcal{H}(|v|^{2})f,\qquad
        \widetilde{\mathbf{D}}_{v}f  \coloneqq\partial_{x}f+\tfrac{1}{2}v\mathcal{H}(\overline{v}f),
\end{equation*}
and
\begin{align*}
    L_{v}f & \coloneqq\partial_{x}f+\tfrac{1}{2}\mathcal{H}(|v|^{2})f+v\mathcal{H}(\text{Re}(\overline{v}f)), \\
    L_{v}^{*}f & \coloneqq-\partial_{x}f+\tfrac{1}{2}\mathcal{H}(|v|^{2})f-v\mathcal{H}(\text{Re}(\overline{v}f)), \\
    \widetilde{L}_{Q}f & \coloneqq \partial_x f+\tfrac{1}{2}\calH(Q^2)f+Q^{-1}\mathcal{H}\Re(Q^{3}f).
\end{align*}
Here, $L_{v}^{*}$ is the $L^{2}$-adjoint operator of $L_{v}$ with respect
to $(\cdot,\cdot)_{r}$,
and $\td{L}_Q$ is a modified operator of $L_Q$.

A second order operator $H_{v}$ is given by 
\begin{align*}
    H_{v}f\coloneqq\left(-\partial_{xx}+\tfrac{1}{4}|v|^{4}-v|D|\overline{v}\right)f.
\end{align*}
Linear operators $A_{Q}, B_Q$, and their adjoint operators $A_{Q}^{*},B_{Q}^{*}$ are defined by 
\begin{align*}
    A_{Q}f&\coloneqq\partial_{x}(x-\mathcal{H})\langle x\rangle^{-1}f,
    &A_{Q}^{*}f&\coloneqq-\langle x\rangle^{-1}(x+\mathcal{H})\partial_{x}f \\
    B_{Q}f&\coloneqq(x-\mathcal{H})\langle x\rangle^{-1}f,
    &B_{Q}^{*}f&\coloneqq\langle x\rangle^{-1}(x+\mathcal{H})f.
\end{align*}
Note that $A_{Q}=\partial_{x}B_{Q}$ and $A_{Q}^{*}=-B_{Q}^{*}\partial_{x}$.

We define
\begin{align*}
	\calB_j f\coloneqq& \partial_x^j (\tfrac{x}{\langle x\rangle}f)-\calH\partial_x^j(\tfrac{1}{\langle x\rangle}f),
\end{align*}
and
\begin{align*}
	\calL_{j}\coloneqq (\calL^{1}_j,\calL^{2}_j)=\calL^{1}_j\Re+\calL^{2}_j i\Im,
\end{align*}
where
\begin{align*}
	\begin{split}
		\calL^{1}_jf&\coloneqq 
		\partial_x^{j-1}[\tfrac{1}{\langle x\rangle}f+\partial_x(\tfrac{x}{\langle x\rangle}f)]
		-\calH\partial_x^{j}(\tfrac{1}{\langle x\rangle}f),
		\\
		\calL^{2}_jg&\coloneqq
		\partial_x^{j-1}[\tfrac{x}{\langle x\rangle^2}\partial_x(\langle x\rangle g)]-\calH\partial_x^{j-1}[\tfrac{1}{\langle x\rangle^2}\partial_x(\langle x\rangle g)].
	\end{split}
\end{align*}
One can check that, for sufficiently regular $\wt\eps$,
\begin{align*}
	\calB_j\wt\eps=\partial_x^jB_Q\wt\eps,
	\qquad
	\calL_{j}\wt\eps=\partial_x^{j-1}B_Q L_Q\wt\eps .
\end{align*}

Finally, we denote the nonlinear part by
\begin{align*}
    N_{v}(\eps) & \coloneqq\eps\mathcal{H}(\Re(\overline{v}\eps))+\tfrac{1}{2}(v+\eps)\mathcal{H}(|\eps|^{2}).
\end{align*}

\vspace{5pt} 

\noindent\underline{\textit{Decompositions and Profiles.}} 
We will decompose solutions into a leading profile described by several modulation parameters, and a radiation component. For the reader's convenience, we summarize below the key relations among these quantities. For renormalized variables $w$ and $w_j$, for $1\le j\le 2L$, we will use modulation decompositions.

For $w$, we use the following decompositions:
\begin{align*}
    w&=Q+\wt\eps& 
    \\
    &=Q+T+\eps& (\wt\eps=T+\eps) 
    \\
    &=Q+\check{T}+\check{\eps},& (\wt\eps=\check{T}+\check{\eps}) 
\end{align*}
where the profiles $T$ and $\check T$ are given by
\begin{align*}
    T &\coloneqq i\nu_0\tfrac{y}{2}Q-ib_1\tfrac{y^2}{4}Q-\eta_1\tfrac{1+y^2}{4}Q,
    \\
    \check{T}&\coloneqq i\check{\nu}_0\tfrac{y}{2}Q-i\check{b}_1\tfrac{y^2}{4}Q-\check{\eta}_1\tfrac{1+y^2}{4}Q.
\end{align*}
One can check its even and odd components:
\begin{align*}
    T_{o}\coloneqq i\nu_0\tfrac{y}{2}Q,\quad 
    T_{e}\coloneqq -ib_1\tfrac{y^2}{4}Q-\eta_1\tfrac{1+y^2}{4}Q.
\end{align*}

For $w_j$ with $1\le j \le 2L$, we use the following decompositions:
\begin{alignat*}{2}
    w_j &= P_{j,e} + \wt\eps_j 
    &\qquad& \\
    &= P_j + \eps_j 
    &\qquad& (\wt\eps_j = P_{j,o} + \eps_j) \\
    &= P_j + \check{\bm{\beta}}_{j+1}(1+y^2)Q + \check{\eps}_j, 
    &\qquad& (\eps_j = \check{\bm{\beta}}_{j+1}(1+y^2)Q + \check{\eps}_j)
\end{alignat*}
where
\begin{align*}
	P_{j}&\coloneqq c_{j}Q+ \bm{\beta}_{j} yQ,\quad   P_{j,e}\coloneqq c_{j}Q,\quad P_{j,o}\coloneqq  \bm{\beta}_{j} yQ
\end{align*}
for some $c_j,\bm{\beta}_{j} \in \bbC$.

We will have relations
\[
   \nu_0=2\Im(c_1),\quad \mu_0=2\Re(c_1),\quad 
   b_1=-2\Im(\bm{\beta}_{1}),\quad \eta_1=-2\Re(\bm{\beta}_{1}),
\]
so that
\[
   c_1=\tfrac{1}{2}(i\nu_0+\mu_0),\quad 
   \bm{\beta}_{1}=-\tfrac{1}{2}(ib_1+\eta_1).
\]
We also denote $\beta_j \coloneqq |\bm\beta_j|$.
We note that $c_1$ and $\bm{\beta}_{1}$ appear in both $w$ and $w_1$ decompositions.

In Section~\ref{subsec:frakc}, for notational convenience, we use some conventions for $j=0$ case
\begin{align*}
    P_{0}=Q+T_{o},\quad P_{0,e}=Q,\quad P_{0,o}=T_o, \quad \text{and} \quad c_0=1,\quad \bm{\beta}_0=\tfrac{i}{2}\nu_0.
\end{align*}

\subsection{Preliminaries} We collect some useful identities and formulas.

\vspace{5pt}

\noindent\underline{\textit{Hilbert transform.}}
We use some algebraic identities with respect to $\calH$ and $Q$: 
\begin{align}
    \mathcal{H}\big(\tfrac{1}{1+x^2}\big) & =\tfrac{x}{1+x^2}, \quad
    \textnormal{i.e., } \quad
    \calH(Q^2)=xQ^2, \label{eq:algebraic 1} \\
    \mathcal{H}\big(\tfrac{2x^{2}}{(1+x^2)^{2}}\big) & =\tfrac{x^{3}-x}{(1+x^2)^{2}},
    \quad
    \mathcal{H}\big(\tfrac{2x^{3}}{(1+x^2)^{2}}\big)=-\tfrac{3x^{2}+1}{(1+x^2)^{2}}. \label{eq:algebraic 2}
\end{align}
We also make frequent use of the following identities.
\begin{lem}[Formulas for Hilbert transform]\label{LemmaHilbertUsefulEquation} We have the following: 
    \begin{enumerate}
    \item For $f,g\in H^{\frac{1}{2}+}$, we have 
    \begin{align*}
        fg=\mathcal{H}f\cdot\mathcal{H}g-\mathcal{H}(f\cdot\mathcal{H}g+\mathcal{H}f\cdot g) %\label{eq:HilbertProductRule}
    \end{align*}
    in a pointwise sense. 
    \item For $f\in\langle x\rangle^{-1}L^{2}$, we have 
    \begin{align}
        [x,\mathcal{H}]f(x)=\frac{1}{\pi}\int_{\mathbb{R}}f(y)dy.\label{eq:CommuteHilbert}
    \end{align}
    \item If $f\in H^{1}$ and $\partial_{x}f\in\langle x\rangle^{-1}L^{2}$, then we have 
    \begin{align}
        \partial_{x}[x,\mathcal{H}]f=[x,\mathcal{H}]\partial_{x}f=0,\quad\text{i.e.,}\quad[\mathcal{H}\partial_{x},x]f=\mathcal{H}f.\label{eq:CommuteHilbertDerivative}
    \end{align}
    \end{enumerate}
\end{lem}

The proof of Lemma \ref{LemmaHilbertUsefulEquation} can be found in \cite[Appendix~C]{KimKimKwon2024arxiv}.

\vspace{5pt}

\noindent\underline{\textit{Lax pairs.}}
We present the Lax pair structure for \eqref{CMdnls-gauged}, along with the hierarchy of conservation laws, which serves as our key ingredient. This structure was rigorously proved for the first time in \cite{GerardLenzmann2024CPAM}.
\begin{prop}[Integrability for \eqref{CMdnls-gauged}]
\label{prop:Lax} For a solution $v\in C([0,T];H^{s})$ to \eqref{CMdnls-gauged} with $s\geq2$, the following hold true:
\begin{enumerate}
    \item(Lax equation) We have
    \begin{align}
        \partial_{t}\widetilde{\mathbf{D}}_{v}=[-iH_{v},\widetilde{\mathbf{D}}_{v}].\label{eq:LaxEqu Unconditional}
    \end{align}
    \item(Hierarchy of conservation laws) For $j\in \bbN_0$ with $0\leq j\leq 2s$, the quantities
    \begin{align}
        I_j(v)\coloneqq(\td \bfD_v^j v,v)_r \label{eq:hierarchy}
    \end{align}
    are conserved.

    \item Formally, for $j\in \bbN_{0}$, we have
    \begin{gather}
        \partial_{t}(\widetilde{\mathbf{D}}_{v}^jv)+iH_{v}\widetilde{\mathbf{D}}_{v}^jv=0, \label{eq:vk equ}
        \\
    	H_{v}h=-\widetilde{\mathbf{D}}_{v}^2h+\tfrac{1}{2}(\widetilde{\mathbf{D}}_{v}v)\mathcal{H}(\overline{v}h)-\tfrac{1}{2}v\mathcal{H}(\overline{\widetilde{\mathbf{D}}_{v}v}h).\label{eq:Dv square}
    \end{gather}
\end{enumerate}
\end{prop}
For the proof of \eqref{eq:LaxEqu Unconditional} and \eqref{eq:Dv square}, see Proposition~3.5 and the proof of Proposition~3.6 in \cite{KimKimKwon2024arxiv}, respectively. The equation \eqref{eq:vk equ} follows from a straightforward calculation using the equation for $v$ and \eqref{eq:LaxEqu Unconditional}. The proof of \eqref{eq:hierarchy} is standard; see \cite[Lemma~2.4]{GerardLenzmann2024CPAM} for example.

We define the initial higher-order energies associated with \eqref{eq:hierarchy} by
\begin{align*}
    E_j(v(t))=E_j(v_0)\coloneqq (-1)^jI_{2j}(v_0)=\|\td \bfD_{v_0}^j v_0\|_{L^2}^2.
\end{align*}
where $v_0$ denotes the initial data of \eqref{CMdnls-gauged}. In particular, we have $E_0=M(v_0)$ and $E_1=E(v_0)$.
By Proposition~\ref{prop:Lax}, these quantities are conserved in time. 

\vspace{5pt}

\noindent\underline{\textit{Conjugation identity.}} We record the conjugation identities for \eqref{CMdnls-gauged}. Similar identities first appeared in the study of the self-dual Chern–Simons–Schrödinger equation~\cite{KimKwon2020blowup,KimKwonOh2020blowup,Kim2025JEMS}, where they arose naturally from nonlinear conjugation. In contrast, for \eqref{CMdnls-gauged}, such a structure does not follow directly. Motivated by this analogy, \cite{KimKimKwon2024arxiv,JeongKim2024arXiv} introduced suitable operators $A_Q$ and $B_Q$ that realize an analogous factorization and capture the repulsive nature of the linearized flow. The precise identities are stated below.

\begin{prop}[Conjugation identity]\label{prop:AQ-BQ prop}
    The following hold true: 
    \begin{enumerate}
        \item For $f\in H^2$, we have
        \begin{align}
    		L_QiL_Q^*f=iH_Qf=iA_{Q}^*A_{Q}f=i\td\bfD_Q^*\td\bfD_Qf,\quad A_{Q}A_{Q}^{*}f=-\partial_{xx}f. \label{eq:conjugate identity}
    	\end{align}
        
        \item For $f\in L^2$,
        \begin{align}
        B_{Q}^{*}B_{Q}f=f-\frac{1}{2\pi}Q\int_{\R}Qfdx,\quad B_{Q}B_{Q}^{*}f=f.\label{eq:BQBQstar equal I}
        \end{align}
        
        \item For $f\in H^1$,
        \begin{align}\label{eq:DQ BQ decomp}
        	\widetilde{\D}_Qf=B_Q^*A_Qf=B_Q^*\partial_xB_Qf. 
        \end{align}
    \end{enumerate}
\end{prop}
The proof of Proposition~\ref{prop:AQ-BQ prop} can be found in \cite[Proposition~3.8]{KimKimKwon2024arxiv} and \cite[Proposition~2.3]{JeongKim2024arXiv}.

\section{Linear structure of gauge-transformed equation}\label{sec:profile}
We work with the gauge-transformed equation \eqref{CMdnls-gauged}, which is equivalent to the original equation \eqref{CMdnls} under the gauge transform.
In this section, we set up the analytical framework, partly recalling and partly refining previous results, to describe the single-bubble blow-up dynamics in this formulation.

We begin in Section~\ref{sec:gague correspond} by verifying the correspondence between \eqref{CMdnls} and \eqref{CMdnls-gauged} and restating the classification result for \eqref{CMdnls-gauged}. 
In Section~\ref{sec:linearized dynamics}, we discuss the linearized operators around the soliton $Q$ and their coercivity estimates.
Section~\ref{sec:one soliton configure} then introduces the single-bubble configuration for \eqref{CMdnls-gauged}, where a solution near $Q$ is decomposed into a modulated profile and a small radiation term satisfying a sharp scale bound.

\subsection{Gauge correspondence}\label{sec:gague correspond}

We clarify the correspondence between \eqref{CMdnls} and \eqref{CMdnls-gauged} in the single-bubble blow-up setting.
Specifically, we verify how the solution, modulation parameters, and asymptotic profiles transform under the gauge transform.
This correspondence allows us to translate statements between the two formulations and ensures that all results concerning single-bubble blow-up remain consistent under the gauge transform.

We begin by formulating the single-bubble blow-up assumption for the gauged equation, which corresponds to Assumption~\ref{assum:one soliton} for the original equation.
\begin{assum}[Single-bubble blow-up for \eqref{CMdnls-gauged}]\label{assum:one soliton GCM}
    We assume that the finite-time blow-up solution \( v(t) \) to \eqref{CMdnls-gauged} takes the following single-bubble form:
    \begin{equation}
        v(t) - \frac{e^{i\gamma(t)}}{\sqrt{\lambda(t)}} Q\left(\frac{\cdot - x(t)}{\lambda(t)}\right) \to z^* \quad \text{in } L^2 \text{ as } t \to T, \label{eq:one soliton cmdnls gauge}
    \end{equation}
    with the translation parameter $x(t)$ converging to $x^*$ as $t\to T$.
    In this case, we say that \(v(t)\) exhibits a \emph{single-bubble blow-up for \eqref{CMdnls-gauged}}.
\end{assum}
Next, we verify that the above notion is precisely equivalent to its counterpart in the original formulation. 
The following lemma shows that the gauge transform preserves the single-bubble structure and the associated parameters, thereby allowing us to pass freely between the two formulations.

\begin{lem}[Gauge correspondence of single-bubble blow-up]\label{lem:gauge equiv}
    Let $u(t)$ be a solution to \eqref{CMdnls} and $v(t)=-\mathcal{G}(u(t))$ denote its gauge-transformed counterpart, which solves \eqref{CMdnls-gauged}. Then $u(t)$ is a single-bubble blow-up solution if and only if so is $v(t)$. More precisely, the following two statements are equivalent:
    \begin{align}
        u(t)& - [\mathcal{R}]_{\lambda(t),\gamma_{\calG}(t),x(t)} \to u^* \quad \text{in } L^2 \text{ as } t \to T, \label{eq:u one soliton}
        \\
        v(t)& - [Q]_{\lambda(t),\gamma(t),x(t)} \to z^* \quad \text{in } L^2 \text{ as } t \to T, \label{eq:v one soliton}
    \end{align}
    where the phase parameters and the asymptotic profiles are related by
    \begin{align}
        \gamma_{\calG}(t)=\gamma(t)+\gamma^*_{\calG},\qquad \gamma^*_{\calG}=\frac{1}{2} \int_{-\infty}^{x^{*}}|z^{*}|^{2}dy =  \frac{1}{2} \int_{-\infty}^{x^{*}}|u^{*}|^{2}dy, \label{eq:gauge phase}
    \end{align}
    and 
    \begin{equation}
        \begin{aligned}
            u^{*}&= -\mathcal{G}^{-1}(z^{*})\mathbf{1}_{x<x^{*}} + \mathcal{G}^{-1}(z^{*})\mathbf{1}_{x>x^{*}},
            \\
            z^{*}&= -\mathcal{G}(u^{*})\mathbf{1}_{x<x^{*}} + \mathcal{G}(u^{*})\mathbf{1}_{x>x^{*}}.
        \end{aligned} \label{eq:gauge asymptotic profile}
    \end{equation}
    Moreover, the gauge transform $\calG$ is a diffeomorphism on $H^{2L+1}$. Consequently, $u_0\in H^{2L+1}$ if and only if $v_0=-\calG(u_0)\in H^{2L+1}$.
\end{lem}
\begin{proof}
    Throughout the proof of this lemma, we denote the spatial variable by $y$ instead of $x$ to avoid confusion with the modulation parameter $x(t)$.
    We note that
    \begin{align*}
        \calG(f)=f(y)\exp(-\tfrac{i}{2}{\textstyle\int}_{-\infty}^{y}|f(y^\prime)|^{2}dy^\prime),\quad
        \calG^{-1}(f)=f(y)\exp(\tfrac{i}{2}{\textstyle\int}_{-\infty}^{y}|f(y^\prime)|^{2}dy^\prime)
    \end{align*}
    and
    \begin{align*}
        \calG([f]_{\lambda,\gamma,x})=[\calG(f)]_{\lambda,\gamma,x}
        ,\qquad
        \calG^{-1}([f]_{\lambda,\gamma,x})=[\calG^{-1}(f)]_{\lambda,\gamma,x}.
    \end{align*}
    We assume \eqref{eq:v one soliton}. Then, we have
    \begin{align}
        u(&t,y)-[\mathcal{R}(y)]_{\lambda,\gamma_{\calG},x} \nonumber
        \\
        =&
        -\calG^{-1}(v-[Q]_{\lambda,\gamma,x})\exp(\tfrac{i}{2} {\textstyle\int}_{-\infty}^y |v|^2-|v-[Q]_{\lambda,\gamma,x}|^2dy^\prime ) \label{eq:gauge correpond 1}
        \\
        &-[\calG^{-1}(Q) \{\exp(\tfrac{i}{2} {\textstyle\int}_{-\infty}^y (|[v]_{\lambda^{-1},-\gamma,-\lambda^{-1}x}|^2-Q^2)dy^\prime )
        -e^{i\gamma_{\calG}^*}\}]_{\lambda,\gamma,x}. \label{eq:gauge correpond 2}
    \end{align}
    For \eqref{eq:gauge correpond 1}, we check $\calG^{-1}(v-[Q]_{\lambda,\gamma,x})\to \calG^{-1}(z^*)$ in $L^2$, and 
    \begin{align*}
        {\textstyle\int}_{-\infty}^y |v|^2-|v-[Q]_{\lambda,\gamma,x}|^2dy^\prime 
        =
        2\Re{\textstyle\int}_{-\infty}^y
        \ol{[Q]_{\lambda,\gamma,x}} (v-[Q]_{\lambda,\gamma,x})dy^\prime
        +{\textstyle\int}_{-\infty}^{\frac{y-x}{\lambda}}Q^{2}dy^\prime. 
    \end{align*}
    Thanks to \eqref{eq:v one soliton}, $\lambda\to 0$, and $x\to x^*$, we have the pointwise convergence
    \begin{align}
        {\textstyle\int}_{-\infty}^y
        \ol{[Q]_{\lambda,\gamma,x}} (v-[Q]_{\lambda,\gamma,x})dy^\prime\to 0, \label{eq:gauge cor interaction}
    \end{align}
    and
    \begin{align*}
        {\textstyle\int}_{-\infty}^{\frac{y-x}{\lambda}}Q^{2}dy^\prime \to 0 \quad\text{for} \quad y<x^* ,
        \quad \text{and} \quad
        {\textstyle\int}_{-\infty}^{\frac{y-x}{\lambda}}Q^{2}dy^\prime \to 2\pi \quad\text{for}\quad y>x^*.
    \end{align*}
    Hence we obtain
    \begin{align*}
        \eqref{eq:gauge correpond 1} \to 
        -\mathcal{G}^{-1}(z^{*})\mathbf{1}_{y<x^{*}} + \mathcal{G}^{-1}(z^{*})\mathbf{1}_{y>x^{*}},
    \end{align*}
    which yields \eqref{eq:gauge asymptotic profile}.
    
    We define $\gamma_{\calG}^*$ by the first equality in \eqref{eq:gauge phase}. With this choice, it suffices to prove the pointwise convergence
    \begin{equation}
        \begin{aligned}
            &\exp\{\tfrac{i}{2} {\textstyle\int}_{-\infty}^y |[v]_{\lambda^{-1},-\gamma,-\lambda^{-1}x}|^2-Q^2dy^\prime \} -e^{-i\gamma_{\calG}^*}
            \\
            &=
            \exp\{\tfrac{i}{2} {\textstyle\int}_{-\infty}^{\lambda y +x} |v|^2-|[Q]_{\lambda,\gamma,x}|^2 dy^\prime \} -e^{-i\gamma_{\calG}^*} \to 0,
        \end{aligned}\label{eq:gauge correpond 2-1}
    \end{equation}
    which leads to $\eqref{eq:gauge correpond 2} \to 0$ in $L^2$. By \eqref{eq:gauge cor interaction}, the convergence in \eqref{eq:gauge correpond 2-1} reduces to
    \begin{align*}
        \exp\{\tfrac{i}{2} {\textstyle\int}_{-\infty}^{\lambda y +x} |v-[Q]_{\lambda,\gamma,x}|^2 dy^\prime \} -e^{i\gamma_{\calG}^*} \to 0,
    \end{align*}
    which follows from \eqref{eq:v one soliton}, $\lambda\to 0$, $x\to x^*$.
    Thus, we arrive at \eqref{eq:u one soliton} after assuming \eqref{eq:v one soliton}. The converse follows by retracing the above argument. We also note that \eqref{eq:gauge asymptotic profile} implies $|u^*|=|z^*|$, which finishes the proof of \eqref{eq:gauge phase}.

    The regularity equivalence follows from the fact that the gauge transform $\calG$ is a diffeomorphism on $H^s$; see for example \cite[Lemma~3.2]{CollianderKeelStaffilaniTakaokaTao2001SIAM} and \cite[Lemma~7.2]{ErdoganGurelTzirakis2018AIHPC}.
\end{proof}

Having established the equivalence between the two formulations, we present the classification result for the gauge-transformed equation, which serves as the main analytic step toward Theorem~\ref{thm:classfy CMDNLS}.

\begin{thm}\label{thm:classfy GCM}
    Let $L \ge 1$ be an integer, and let $v(t)$ be a solution to \eqref{CMdnls-gauged} satisfying Assumption~\ref{assum:one soliton GCM}, with blow-up time $0<T<+\infty$ and initial data
    \begin{equation*}
        v(0)=v_0 \in H^{2L+1}(\bbR).
    \end{equation*}
    Then, the solution $v(t)$ falls into precisely one of the following scenarios:

    (1) \textbf{Quantized regime.} There exists an integer $1\leq k\leq L$ such that the modulation parameters $(\lambda(t),\gamma(t),x(t))$ in the decomposition \eqref{eq:one soliton cmdnls gauge} satisfy
    \begin{align*}
        \lambda(t)&=(\ell+o_{t\to T}(1))(T-t)^{2k}, 
        \\
        \gamma(t)&=\gamma^*+o_{t\to T}(1),
        \\
        x(t)&=x^*+2\nu^*(T-t)+O((T-t)^2),
    \end{align*}
    for some $\ell \in (0,\infty)$, $\gamma^*\in \bbR/2\pi\bbZ$, and $\nu^* \in \bbR$.

    (2) \textbf{Exotic regime.} The scale and translation parameters $(\lambda(t),x(t))$ in the decomposition \eqref{eq:one soliton cmdnls gauge} satisfy
    \begin{align*}
        \lambda(t)&\lesssim (T-t)^{2L+\frac32},\\
        x(t)&=x^*+2\nu^*(T-t)+O((T-t)^2),
    \end{align*}
    for some $\nu^* \in \bbR$.
    
\end{thm}
%\Red{- First, I deleted for some $x^*$ because this was already included in our assumption.}

%\Red{- I deleted the dependence $\ell=\ell (u)$ and $\gmm^*=\gmm^*(u)$. To mention the dependence, all $k, \ell, \gmm^*, \nu^*$ should depend on $u$.}
Lemma~\ref{lem:gauge equiv} bridges Theorem~\ref{thm:classfy GCM} and Theorem~\ref{thm:classfy CMDNLS}, thus completing the correspondence between the two formulations. The same lemma moreover shows that Assumption~\ref{assum:one soliton} and Assumption~\ref{assum:one soliton GCM} are equivalent through the gauge transform.
Henceforth, we work under Assumption~\ref{assum:one soliton GCM} and focus on establishing Theorem~\ref{thm:classfy GCM}.

\subsection{Linearized dynamics}\label{sec:linearized dynamics}
In this subsection, we briefly review the linearization of \eqref{CMdnls-gauged} around $Q$; see \cite[Section~3]{KimKimKwon2024arxiv} and \cite[Section~3]{JeongKim2024arXiv} for full details. 

As noted in \eqref{eq:intro Energy2 self dual}, we have
\begin{align*}
    E(v)=\frac{1}{2}\int_{\mathbb{R}}\big|\td{\mathbf{D}}_{v}v\big|^{2}dx,\quad  \td{\mathbf{D}}_{Q}Q=0,
\end{align*}
and linearizing $v\mapsto \td{\mathbf D}_v v$ yields
\begin{align*}
    \td{\mathbf{D}}_{v+\eps}(v+\eps)=\td{\mathbf{D}}_{v}v+L_{v}\eps+N_{v}(\eps).
\end{align*}
Since $i\nabla E(v)=iL_{v}^{*}\td\bfD_{v}v$ for the adjoint operator $L_{v}^{\ast}$ of $L_v$, we linearize \eqref{CMdnls-gauged} as
\begin{align*}
    iL_{w+\eps}^{*}\td\bfD_{w+\eps}(w+\eps)=iL_{w}^{*}\td\bfD_{w}w+i\calL_{w}\eps+R_{w}(\eps).
\end{align*}
Here, $i\calL_{w}\eps$ is the linear part, and $R_{w}(\eps)$ is the nonlinear part for $\eps$. Hence the linearized operator around $Q$ takes the self-dual form
\begin{align*}
i\calL_{Q}=iL_{Q}^{*}L_{Q}.
\end{align*}
We recall the kernel properties and the coercivity estimates from \cite{KimKimKwon2024arxiv,JeongKim2024arXiv}. The kernel of $L_Q$ is generated by the symmetries,
\begin{align*}
    \textnormal{ker}\,L_{Q}=\textnormal{ker}\,\mathcal{L}_{Q}
    =\textnormal{span}_{\mathbb{R}}\{iQ,\Lambda Q, Q_x\}, \quad \Lambda f\coloneqq\tfrac{1}{2}f+x\partial_{x}f.
\end{align*}
For higher-order variables $v_j=\widetilde{\mathbf D}_v^j v$, the corresponding linearized operator is $\widetilde{\mathbf D}_Q^{j-1}L_Q$, which is not directly defined on $\dot{H}^j$ because of the Hilbert transform $\calH$ in $L_Q$. Thus, following the observation in \cite{JeongKim2024arXiv}, we introduce higher-order linearized operators $\calL_j$ satisfying
\begin{align*}
	\calL_{j}\wt\eps
	=\partial_y^{j-1}B_Q L_Q\wt\eps \quad \text{for}\quad 1\leq j\leq 2L+1, \quad \wt\eps \in H^{2L+1}. %\label{eq:BQLQ and calL}
\end{align*}
More precisely, we define linear operators by
\begin{align*}
	\calL_{j}\coloneqq (\calL^{1}_j,\calL^{2}_j)=\calL^{1}_j\Re+\calL^{2}_j i\Im,
\end{align*}
where
\begin{align}
	\begin{split}
		\calL^{1}_jf&\coloneqq 
		\partial_x^{j-1}[\tfrac{1}{\langle x\rangle}f+\partial_x(\tfrac{x}{\langle x\rangle}f)]
		-\calH\partial_x^{j}(\tfrac{1}{\langle x\rangle}f),
		\\
		\calL^{2}_jg&\coloneqq
		\partial_x^{j-1}[\tfrac{x}{\langle x\rangle^2}\partial_x(\langle x\rangle g)]-\calH\partial_x^{j-1}[\tfrac{1}{\langle x\rangle^2}\partial_x(\langle x\rangle g)].
	\end{split} \label{eq:newdef BQLQ}
\end{align}
See \cite[Section~3]{JeongKim2024arXiv} for details. These operators are well defined on the weighted space $\dot{\calH}^j$ and inherit the spectral properties of $A_Q\widetilde L_Q$ in \cite{KimKimKwon2024arxiv}. In particular, the spectral properties of $A_Q\widetilde L_Q$ and $\calL_2$ coincide on $\dot{\calH}^2$, and thus
\begin{align*}
    \textnormal{ker}\,\calL_{2}
    =\textnormal{span}_{\mathbb{R}}\{iQ,\Lambda Q, Q_x, ixQ, ix^2Q,(1+x^2)Q\}.
\end{align*}
In addition, we expect that, for general $j\ge2$, the kernel of $\mathcal L_j$ on $\dot{\mathcal H}^j$ follows the same structural pattern,
\begin{align}
	\textnormal{ker}\calL_{j}
	=\textnormal{span}_{\mathbb{R}} \{iQ,\Lambda Q,Q_y,iyQ, iy^{\ell+1}Q, y^{\ell-1}(1+y^2)Q: 1\leq \ell\leq j-1\}. \label{eq:general kernel nonrad}
\end{align}
We recall that $\calK_j$ are given by \eqref{eq:kernel element1}. Now, we introduce coercivity lemmas for our analysis.

\begin{lem}[Linear coercivity estimates I]
\label{lem:coercivity for linear} The following hold true:
\begin{enumerate}
    \item (Coercivity for $L_{Q}$ on $\dot{\mathcal{H}}^{1}$) Let $\{\psi_{j}\}_{1\leq j\leq 3}$ be elements of the dual space $(\dot{\mathcal{H}}^{1})^{*}$ such that $(\psi_{i},\mathcal{K}_{j})_{r}=\delta_{i,j}$. Then, we have a coercivity estimate 
    \begin{align}
        \|v\|_{\dot{\mathcal{H}}^{1}}\sim\|L_{Q}v\|_{L^{2}},\quad\forall v\in\dot{\mathcal{H}}^{1}\cap\{\psi_{j}:1\leq j\leq 3\}^{\perp}. \label{eq:coercivity L1}
    \end{align}
    
    \item (Coercivity for $\calL_{2}$ on $\dot{\mathcal{H}}^{2}$) Let $\{\psi_{j}\}_{1\leq j\leq 6}$ be elements of the dual space $(\dot{\mathcal{H}}^{2})^{*}$ such that $(\psi_{i},\mathcal{K}_{j})_{r}=\delta_{i,j}$. Then, we have a coercivity estimate 
    \begin{align}
        \|v\|_{\dot{\mathcal{H}}^{2}}\sim\|\calL_{2}v\|_{L^{2}},
        \quad\forall v\in\dot{\mathcal{H}}^{2}\cap\{ \psi_{j}:1\leq j\leq 6\}^{\perp}. \label{eq:coercivity L2}
    \end{align}

    \item (Coercivity for $\calL_{3}$ on $\dot{\mathcal{H}}^{3}_e$) Let $\psi_{1},\psi_{2},\psi_{4},\psi_{5}$ be elements of the dual space $(\dot{\mathcal{H}}^{3}_e)^{*}$ such that $(\psi_{i},\mathcal{K}_{j})_{r}=\delta_{i,j}$. Then, we have a coercivity estimate 
    \begin{align}
        \|v\|_{\dot{\mathcal{H}}^{3}_e}\sim\|\calL_{3}v\|_{L^{2}},
        \quad\forall v\in\dot{\mathcal{H}}^{3}_e\cap\{ \psi_{1},\psi_{2},\psi_{4},\psi_{5}\}^{\perp}. \label{eq:coercivity L3}
    \end{align}

\end{enumerate}
\end{lem}
We note that, due to the degeneracy of the linearized operators, the adapted Sobolev spaces $\dot{\mathcal{H}}^{j}$ are used instead of the usual ones.
\begin{proof}[Proof of Lemma~\ref{lem:coercivity for linear}]
    The proofs of \eqref{eq:coercivity L1} and \eqref{eq:coercivity L2} are given in \cite{KimKimKwon2024arxiv}. For \eqref{eq:coercivity L3}, we need to find the kernel of $\calL_3$ on $\dot{\mathcal{H}}^{3}_e$. Formally, one can observe that $\calL_3=\partial_y \calL_2$. 
    Thus, we may expect that additional kernel elements $f_1,f_2$ satisfy
    \begin{align*}
       \calL_2 f_1 = 1, \quad \calL_2 f_2 = i.
    \end{align*}
    Following the argument of \cite[Lemma~3.2]{JeongKim2024arXiv}, one finds
    \begin{align*}
       f_1 = \tfrac{1}{3\sqrt{2}}y(1+y^2)Q,
       \quad 
       f_2 = \tfrac{i}{3\sqrt{2}}y^3 Q.
    \end{align*}
    However, both are odd, and hence irrelevant in the even sector. 
    Therefore, on $\dot{\mathcal{H}}^3_e$, we expect
    \begin{align}
       \ker \calL_3 
       = \mathrm{span}_{\mathbb{R}}\{\, iQ, \Lambda Q,\, i x^2 Q,\,(1+x^2)Q \,\}. \label{eq:kernel calL3 even}
    \end{align}
    A rigorous justification of this kernel description can be obtained by setting $h=\frac{1}{(1+x^2)^5}f$ and applying the Fourier transform. 
    This reduces $\calL_3 f=0$ to a sixth-order ODE with a discontinuous coefficient at $\xi=0$. The resulting analysis parallels the method used in Propositions~3.2 and~3.12 of \cite{KimKimKwon2024arxiv}, which establish the kernels of $\calL_1$ and $\calL_2$. We omit the details, as the argument is purely technical, but note that this reduction rigorously confirms the expected kernel structure. In particular, one concludes that the kernel of $\calL_3$ in $\dot\calH^3$ coincides with \eqref{eq:general kernel nonrad} for $j=3$, and upon restriction to the even subspace $\dot\calH^3_e$, it reduces precisely to \eqref{eq:kernel calL3 even}. 
    Therefore, combining the subcoercivity of $\calL_3$ on $\dot{\mathcal{H}}^{3}_e$, Lemma~\ref{lem:calL3 subcoer}, and the kernel characterization \eqref{eq:kernel calL3 even}, we obtain \eqref{eq:coercivity L3} by a standard argument; see \cite[Lemma~A.15]{KimKwon2020blowup}.
\end{proof}

Finally, we recall that the operator $A_Q$, introduced through the conjugation identity \eqref{eq:conjugate identity}, plays an essential role in establishing propagation of smallness of higher-order energies. In the analysis of \cite{JeongKim2024arXiv}, only even perturbations were considered, so the coercivity was formulated for $A_Q$ acting on odd functions, particularly on variables of the form $A_Q v_{2k-1}$. In the present setting, however, we treat the general case without parity restrictions and hence require coercivity for all variables $v_j$.

To this end, we introduce another operators $\calB_j$ that extend the action of $A_Q$ to higher orders, defined for $f\in H^j$ by
\begin{align}
	\partial_{x}^jB_Q f=\partial_{x}^j(\tfrac{x}{\langle x \rangle}f)-\calH\partial_{x}^j(\tfrac{1}{\langle x \rangle}f)\eqqcolon \calB_j f. \label{eq:def Bj}
\end{align}
We note that $\calB_1=A_Q$. As before, $\mathring \calK_i$ are given by \eqref{eq:kernel element1}. The following lemma provides coercivity estimates for these operators, which generalize those of $A_Q$ to higher regularity levels.

\begin{lem}[Linear coercivity estimates II]
\label{lem:coercivity for linear wj} The following hold true:
\begin{enumerate}
    \item (Coercivity for $\calB_1 $ on $\dot{\mathcal{H}}^{1}$) Let $\psi_{1}, \psi_2$ be elements of the dual space $(\dot{\mathcal{H}}^{1})^{*}$ such that $(\psi_{i},\mathring \calK_j )_{r}=\delta_{i,j}$. Then, we have a coercivity estimate 
    \begin{align}
        \|v\|_{\dot{\mathcal{H}}^{1}}\sim \|\calB_1 v\|_{L^{2}},\quad\forall v\in\dot{\mathcal{H}}^{1}\cap\{\psi_{j},i\psi_{j}: 1\leq j\leq 2\}^{\perp}. \label{eq:coercivity A1}
    \end{align}

    \item (Coercivity for $\calB_2 $ on $\dot{\mathcal{H}}^{2}$) Let $\{\psi_{j}\}_{1\leq j\leq 3}$ be elements of the dual space $(\dot{\mathcal{H}}^{2})^{*}$ such that $(\psi_{i},\mathring \calK_j )_{r}=\delta_{i,j}$. Then, we have a coercivity estimate 
    \begin{align}
        \|v\|_{\dot{\mathcal{H}}^{2}}\sim \|\calB_2 v\|_{L^{2}},\quad\forall v\in\dot{\mathcal{H}}^{2}\cap\{\psi_{j},i\psi_{j}: 1\leq j\leq 3\}^{\perp}. \label{eq:coercivity A2}
    \end{align}

    \item (Coercivity for $\calB_3 $ on $\dot{\mathcal{H}}^{3}_e$) Let $\psi_{1}, \psi_3$ be elements of the dual space $(\dot{\mathcal{H}}^{3}_e)^{*}$ such that $(\psi_{i},\mathring \calK_j )_{r}=\delta_{i,j}$. Then, we have a coercivity estimate 
    \begin{align}
        \|v\|_{\dot{\mathcal{H}}^{3}}\sim \|\calB_3 v\|_{L^{2}},\quad\forall v\in\dot{\mathcal{H}}^{3}_e\cap\{\psi_{j},i\psi_{j}: j=1, 3\}^{\perp}. \label{eq:coercivity A3}
    \end{align}
\end{enumerate}
\end{lem}
\begin{proof}%[Proof of Lemma~\ref{lem:coercivity for linear}]
    For $(1)$, the proof can be found in \cite[Proposition~3.8]{KimKimKwon2024arxiv}. As noted in the proof of Lemma~\ref{lem:coercivity for linear}, a Fourier transform computation shows that
    \begin{equation*}
        \begin{alignedat}{2}
            \ker \calB_2 
            &= \mathrm{span}_{\mathbb{C}}\{  Q, xQ, (1+x^2)Q  \} &\quad& \text{on} \quad \dot \calH^2,
            \\
            \ker \calB_3 
            &=\mathrm{span}_{\mathbb{C}}\{ Q, (1+x^2)Q \}  &\quad& \text{on} \quad \dot \calH^3_e.
        \end{alignedat}
    \end{equation*}
    Together with Lemmas~\ref{lem:subcoer calB2} and~\ref{lem:subcoer calB3}, this leads by a standard argument to the conclusion.
    %Together with Lemmas~\ref{lem:subcoer calB2} and~\ref{lem:subcoer calB3}, this gives the proof by a standard argument.
\end{proof}

\subsection{Single-bubble configuration}\label{sec:one soliton configure}
In this subsection we reformulate the single-bubble blow-up assumption in a form suited to \eqref{CMdnls-gauged}. Our objective is to describe how a solution near the soliton $Q$ can be expressed as a modulated profile plus a small radiation term that is orthogonal to the symmetry directions.
To achieve this, we introduce a family of test functions $\calZ_k$ satisfying suitable transversality relations with the kernel functions $\calK_j$.

In \cite{KimKimKwon2024arxiv}, the functions $\calZ_k\in C^\infty_c$ were constructed to satisfy
\begin{align*}
    (\mathcal{K}_{j},\mathcal{Z}_{k})_r=d_{j}\delta_{jk},
\end{align*}
where $\delta_{jk}$ denotes the Kronecker delta and each $d_j$ is a nonzero constant. This condition was sufficient for the decomposition of the solution $w$ in that setting.
However, in the present single-bubble setting, we further impose
\begin{align*}
    (Q,\mathcal{Z}_{k})_r=0,
\end{align*}
a condition that will play a crucial role in the coercivity analysis carried out in the proof of Lemma~\ref{lem:coercivity H3 epsj} (see Section~\ref{sec:coer pf}).

To construct such functions, we introduce a first-order operator
\begin{align*}
    \bfD\coloneqq x\partial_x -1+\frac{5x^2}{1+x^2}.
\end{align*}
A direct computation reveals
\begin{align*}
    \bfD (\bfD_Q \{\calK_1\})=\bfD(Q^2Q_x)=0, \quad
    \bfD [\bfD_Q \{\calK_5\}]=-\bfD(2xQ)=-4x^3Q^3,
\end{align*}
which guides the design of the new family $\{\calZ_k\}$.
Explicitly, we set
\begin{align*}
	\mathcal{Z}_{1} & =\bfD_Q^*\left[Q_x \chi-\tfrac{(Q_x\chi,xQ)_r}{(\bfD^*(xQ\chi), xQ)_r}\bfD^*(xQ\chi)\right],\\
	\mathcal{Z}_{2} & =i(1+x^{2})Q\chi-\tfrac{(x^{2}Q,(1+x^{2})Q\chi)_{r}}{2(xQ,xQ\chi)_{r}}\mathbf{D}_{Q}^{*}(ixQ\chi),\\
	\mathcal{Z}_{3} & =xQ\chi,\\
	\mathcal{Z}_{4} & =\mathbf{D}_{Q}^{*}(ixQ\chi),\qquad
	\mathcal{Z}_{5}  =\bfD_Q^*\bfD^*(xQ\chi),\qquad
	\mathcal{Z}_{6}  =ixQ\chi.
\end{align*}
By construction, this family satisfies the transversality relations
\begin{align}
	(\mathcal{K}_{j},\mathcal{Z}_{k})_r=d_{j}\delta_{jk}, \quad (Q,\mathcal{Z}_{k})_r=0, \quad 1\leq j,k\leq 6,\label{eq:transversality}
\end{align}
where each $d_j$ is a nonzero constant.
The above orthogonality relations ensure that the modulation parameters can be uniquely determined for solutions sufficiently close to $Q$. We summarize this decomposition in the following proposition. 

\begin{prop}[Decomposition]\label{prop:Decomposition}
	For any $M>0$, there exist $0<\alpha^*\ll \eta \ll R^{-1}\ll M^{-1}$ such that the following hold for all $v\in H^1$ with $\|v\|_{L^2}\leq M$ and small energy condition $\sqrt{E(v)}\leq \alpha^*\|v\|_{\dot H^1}$:
	\begin{enumerate}
		\item (Decomposition) There exists unique $(\wt\lambda,\wt\gamma,\wt x,\wt\eps)\in \bbR_+\times\bbR/2\pi\bbZ\times\bbR\times H^1$ satisfying
		\begin{align}
			v=[Q+\wt\eps]_{\wt\lambda,\wt\gamma,\wt x},\quad (\wt\eps,\mathcal{Z}_k)_r=0 \quad \text{for}\quad k=1,2,3, \label{eq:Decomposition and orthogonal}
		\end{align}
		and smallness $\|\wt\eps\|_{\dot\calH^1}<\eta$.

        \item ($C^{1}$-regularity) The map $v\mapsto(\wt\lambda,\wt\gamma,\wt x)$
		for each decomposition is $C^{1}$. 
        
		\item (Estimate for $\lambda$) We have 
		\begin{align}
			\left|\frac{\|v\|_{\dot H^1}}{\|Q\|_{\dot H^1}}\wt\lambda -1 \right|\lesssim \|\wt\eps\|_{\dot H^1}. \label{eq:Decompose lambda estimate}
		\end{align}
		
		\item (Refined estimate for $\wt\eps$) We have
		\begin{align*}
			\|\wt\eps\|_{\dot\calH^1_R}^2+E(\varphi_R\wt\eps)
			\lesssim_M \wt\lambda^2E(v). %\label{eq:Nolinear estimate truncated}
		\end{align*}
	\end{enumerate}
\end{prop}
For the proof of Proposition~\ref{prop:Decomposition}, refer to \cite[Proposition~3.1]{KimTKwon2024arxivSolResol}. Based on this proposition, we claim the following proposition, which asserts that the single-bubble configuration yields a quantitative coercivity for the radiation term.
\begin{prop}\label{prop:one soliton} 
    Under Assumption~\ref{assum:one soliton GCM}, let $(\lambda,\gamma,x)$ be the modulation parameters appearing in the assumption, and let $v$ be decomposed as in Proposition~\ref{prop:Decomposition}, with corresponding parameters $(\wt\lambda,\wt\gamma,\wt x)$. Then, we have
    \begin{align}
    	\wt\lambda^2 E_1=E(Q+\wt\eps) \gtrsim \|\wt\eps\|_{\dot\calH^1}^2. \label{eq:coercivity H1 from one sol}
    \end{align}
    Moreover, we have
    \begin{align}
        \bigg|\frac{\wt\lambda}{\lambda} -1 \bigg|\lesssim \wt\lambda,\quad |\wt x-x|\lesssim \wt\lambda, \quad \lim_{t\to T}|\wt\ga-\ga|=0, \label{eq:mod para change}
    \end{align}
    and
    \begin{align}
        \wt\lambda+\lambda\lesssim T-t. \label{eq:lambda bdd one sol}
    \end{align}
\end{prop}
\begin{proof}
    We show \eqref{eq:coercivity H1 from one sol} by contradiction. If \eqref{eq:coercivity H1 from one sol} fails, then by \cite[Lemma~4.1, Lemma~4.2, Lemma~4.3, Lemma~4.6]{KimTKwon2024arxivSolResol}, there exist $N\geq2$, $(\lambda_{j},\gamma_{j},x_{j})\in \bbR_+\times\bbR/2\pi\bbZ\times\bbR$, and $\td\eps_N\in H^1$ such that
    \begin{align*}
        v(t)=[Q]_{\wt\lambda,\wt\gamma,\wt x}+\sum_{j=2}^{N}[Q]_{\lambda_{j},\gamma_{j},x_{j}}+\td \eps_{N}
    \end{align*}
    with $\|\eps_{N}\|_{H^1}\lesssim 1$, and
    \begin{align*}
        \int_{\bbR}\psi|v|^{2}dx
        =\int_{\bbR}\psi|[Q]_{\wt\lambda,\wt\gamma,\wt x}|^2
        +\sum_{j=1}^{N}\int_{\bbR}\psi|[Q]_{\lambda_{j},\gamma_{j},x_{j}}|^{2}
        +\psi|\td\eps_{N}|^{2}dx+o_{t\to T}(1)\cdot\|\psi\|_{L^{\infty}}%\label{eq:soliton decoupling}
    \end{align*}
    for any $\psi \in L^\infty$. Moreover, $\lim_{t\to T}\wt x(t)=\wt x^*$ and $\lim_{t\to T}x_{j}(t)=x_{j}^*$ exist for all $j$. Therefore, taking $\psi=\chi_{r}(x-\wt x^*)$ with sufficiently small $r>0$, we reach a contradiction unless $N=1$ and $x^*=\wt x^*$. 
    
    By \eqref{eq:coercivity H1 from one sol} and \eqref{eq:Decompose lambda estimate}, we directly deduce $|\frac{\wt\lambda}{\lambda} -1 |\lesssim \wt\lambda$. From \eqref{eq:one soliton cmdnls gauge} and \eqref{eq:Decomposition and orthogonal}, we have
    \begin{align*}
        v(t)=[Q]_{\la,\ga,x}+\eps_{\text{aux}}=[Q]_{\wt\la,\wt\ga,\wt x}+[\wt\eps]_{\wt\la,\wt\ga,\wt x}.
    \end{align*}
    Moreover, from \textbf{Step 1} of the proof of \cite[Theorem~1.2]{KimTKwon2024arxivSolResol}, we obtain $[\wt\eps]_{\wt\la,\wt\ga,\wt x}\to \wt z^*$ in $L^2$. Since the assumption \eqref{eq:one soliton cmdnls gauge} implies $\eps_{\text{aux}}\to z^*$, we derive
    \begin{align*}
        [Q]_{\wt\la,\wt\ga,\wt x}-[Q]_{\la,\ga,x}=\eps_{\text{aux}}-[\wt\eps]_{\wt\la,\wt\ga,\wt x}\to z^*-\wt z^* \quad \text{in} \quad L^2.
    \end{align*}
    Since $\lambda,\wt\lambda\to 0$, their difference cannot converge strongly in $L^2$ to a nonzero limit. Hence, we obtain $z^*-\wt z^*=0$, and thus we conclude \eqref{eq:mod para change}.

    Finally, the bound \eqref{eq:lambda bdd} in Theorem~\ref{thm:Soliton resolution} directly yields \eqref{eq:lambda bdd one sol}.
\end{proof}
Thanks to \eqref{eq:mod para change}, the difference between the two parametrizations is negligible. From this point on, we work with $(\wt\lambda,\wt\gamma,\wt x)$ in place of $(\lambda, \gamma, x)$, but keep the original notation $(\lambda, \gamma, x)$ for notational convenience.
\begin{rem}
    In view of soliton resolution (Theorem~\ref{thm:Soliton resolution}), the convergence of $x(t)$ as $t\to T$  holds even without being assumed in Assumption~\ref{assum:one soliton GCM}. Here, this can be justified directly by \eqref{eq:mod para change}.
\end{rem}

\section{Heuristics: the radial case}\label{sec:radial}

%\vspace{10pt}
%\textit{In the remainder of this section, we assume $w$ is radial and present a heuristic description of the blow-up classification in this symmetric setting.}
%\vspace{10pt}

This section provides a brief structural overview that prepares the reader for the detailed analysis to follow. Under the radial assumption, we outline the classification of blow-up rates in a simplified setting, highlighting the main ideas and presenting the core argument in a concise and direct manner.

To keep the presentation clear and focused, we introduce notation tailored to radial solutions. These notational choices and symmetry assumptions are used only for this heuristic discussion. When we return to the full, non-symmetric analysis, we switch back to the general notation introduced in Section~\ref{sec:notation}.

For completeness, we introduce the renormalized coordinates and variables used in this section. 
We define the rescaled time and space variables $(s,y)$ and the renormalized solution $w$ by
\begin{align*}
	\frac{ds}{dt}=\frac{1}{\lambda(t)^{2}},\qquad 
    y=\frac{x}{\lambda(t)},\qquad 
    w(s,y)=\lambda^{\frac{1}{2}}e^{-i\gamma}v(t,\lambda(t)y)\big|_{t=t(s)}. 
\end{align*}
For notational convenience, we write $\lambda(s)=\lambda(t(s))$ and $\gamma(s)=\gamma(t(s))$, 
and we adopt the same convention for the other parameters. 
The corresponding (renormalized) nonlinear adapted derivatives are defined by
\begin{align*}
	w_{j}\coloneqq\td\bfD_{w}^{j}w
    =\lambda^{\frac{2j+1}{2}} e^{-i\gamma}(v_{j})(t,\lambda(t)y)\big|_{t=t(s)}.
\end{align*}
Under radial symmetry, the translation parameter $x(t)$ in Proposition~\ref{prop:Decomposition} is absent. 
Rewriting \eqref{CMdnls-gauged} and \eqref{eq:vk equ} together with \eqref{eq:Dv square} in the renormalized variables, we obtain the renormalized equations
\begin{equation}
    \begin{aligned}
        &(\partial_s-\tfrac{\lambda_s}{\lambda}\Lambda +\gamma_s i)w+iL_w^*w_1=0,
        \\
        &(\partial_s-\tfrac{\lambda_s}{\lambda}\Lambda_{-j} +\gamma_s i)w_{j}-iw_{j+2}
        +\tfrac{i}{2}[w_1\mathcal{H}(\overline{w}w_{j})-w\mathcal{H}(\overline{w_1}w_{j})]=0. 
    \end{aligned} \label{eq:w equ rad}
\end{equation}
In particular, only the odd-index $(j=2k-1)$ equations will be used:
\begin{equation}
    (\partial_s-\tfrac{\lambda_s}{\lambda}\Lambda_{-(2k-1)} +\gamma_s i)w_{2k-1}-iw_{2k+1}
    +\tfrac{i}{2}[w_1\mathcal{H}(\overline{w}w_{2k-1})-w\mathcal{H}(\overline{w_1}w_{2k-1})]=0. \label{eq:wj equ rad}
\end{equation}
Following \cite[Section~3.2]{JeongKim2024arXiv}, we decompose
\begin{equation}
    \begin{aligned}
        w&=Q+\wt\eps=Q+T + \eps ,
        \\
        w_{2k-1}&=P_{2k-1}+\eps_{2k-1},
    \end{aligned}\label{eq:decom rad}
\end{equation}
where
\begin{align*}
    T=-ib_1\tfrac{y^2}{4}Q-\eta_1\tfrac{1+y^2}{4}Q,\quad P_{2k-1}=-\tfrac{1}{2}(-i)^{k-1}(ib_k+\eta_k)yQ.
\end{align*}

As observed in \cite{JeongKim2024arXiv}, inserting only the profiles into \eqref{eq:w equ rad} and \eqref{eq:wj equ rad} and matching the coefficients along the kernel directions yield the formal ODE system of modulation parameters. 
Equation \eqref{eq:w equ rad} gives 
\begin{align*}
    \text{LHS}\eqref{eq:w equ rad}\approx
    -\left(\frac{\lambda_{s}}{\lambda}+b_1\right)\Lambda Q+\left(\gamma_{s}-\frac{\eta_1}{2}\right)iQ+\text{h.o.t.},
\end{align*}
which suggests 
\begin{align*}
    \frac{\lambda_{s}}{\lambda}+b_1=0,\qquad \gamma_{s}-\frac{\eta_1}{2}=0.%\label{eq:first mod equ}
\end{align*}
The odd-index equation \eqref{eq:wj equ rad} for $1\le k\le L$ (with the convention $Q\mathcal H(1)=0$) gives 
\begin{align*}
	\text{LHS}\eqref{eq:wj equ rad} & \approx
	((b_k)_s-b_{k+1}+(2k-\tfrac{1}{2})b_1b_{k}+\tfrac{1}{2}\eta_1\eta_k)(-i)^{k}\tfrac{y}{2}Q
	\\
	& \quad -
	((\eta_k)_s-\eta_{k+1}+(2k-\tfrac{1}{2})b_1\eta_{k}-\tfrac{1}{2}\eta_1b_k)(-i)^{k-1}\tfrac{y}{2}Q+\text{h.o.t.},
\end{align*}
which suggests that, for each $k$,
\begin{align}
    \begin{split}
        (b_k)_s &= b_{k+1}-(2k-\tfrac12)b_1b_k-\tfrac12\eta_1\eta_k,\\
        (\eta_k)_s&=\eta_{k+1}-(2k-\tfrac12)b_1\eta_k+\tfrac12\eta_1 b_k.
    \end{split}\label{eq:mod eqn}
\end{align}

The modulation system \eqref{eq:mod eqn} is obtained by inserting only the profile parts in the decomposition \eqref{eq:decom rad} into \eqref{eq:wj equ rad} and matching the coefficients along the kernel directions.
To justify \eqref{eq:mod eqn} rigorously, it therefore remains to control the contribution of the radiation parts $\wt\eps$, $\eps$ and $\eps_{2k-1}$ appearing in \eqref{eq:decom rad}.

More precisely, after subtracting the profile contributions, the remaining terms in \eqref{eq:wj equ rad} involving the radiation parts must be shown to give rise only to error terms that are sufficiently small compared with the leading-order profile dynamics.
This requires quantitative estimates on $\wt\eps$, $\eps$ and $\eps_{2k-1}$, ensuring that their influence on the system of modulation parameters can be regarded as perturbative errors.

In what follows, we derive such bounds by combining the hierarchy of conservation laws with linear coercivity estimates for the operators $\calB_1$ and $\calB_2$.
These estimates allow us to control the radiation part and thereby justify that the formal modulation system \eqref{eq:mod eqn} holds up to negligible errors.

We begin by recalling the consequences of the hierarchy of conservation laws \eqref{eq:hierarchy}. Assuming $v_0\in H^{2L+1}$, we have for $1\le j\le 2L$,
\begin{align*}
    \|\td{\bfD}_w w_j\|_{L^2}=\|w_{j+1}\|_{L^2}=\lambda^{j+1}\sqrt{E_{j+1}},\quad 
    \|\td{\bfD}_w^2 w_j\|_{L^2}=\|w_{j+2}\|_{L^2}=\lambda^{j+2}\sqrt{E_{j+2}}.
\end{align*}
Moreover, combining \eqref{eq:BQBQstar equal I}, \eqref{eq:DQ BQ decomp}, and \eqref{eq:coercivity H1 from one sol}, we deduce\footnote{The non-radial version is stated in Lemma~\ref{lem:coercivity estimate 2}, and its rigorous proof is given in Section~\ref{sec:coer pf}.}
\begin{align*}
    \|\calB_1 w_j\|_{L^2}\lesssim \lambda^{j+1},\qquad \|\calB_2 w_j\|_{L^2}\lesssim \lambda^{j+2}.
    %\label{eq:boots energy AQ}
\end{align*}
Since $w_{2k-1}$ is odd in the radial setting, it suffices to analyze the kernels of $\calB_1$ and $\calB_2$ restricted to the odd subspace.
One can show that $\ker \calB_1=\mathrm{span}_{\bbC}\{yQ\}$ on $\dot{\calH}^1_{o}$ and $\ker \calB_2=\mathrm{span}_{\bbC}\{yQ\}$ on $\dot{\calH}^2_{o}$.
Accordingly, it is enough to extract the $yQ$-mode coefficients, namely $b_k$ and $\eta_k$, which correspond to the profile $P_{2k-1}$. As a consequence of the linear coercivity estimates
(Lemma~\ref{lem:coercivity for linear wj}), one obtains that, for $v_0\in H^{2L+1}$ and $1\le k\le L$,
\begin{align}
    \|\eps_{2k-1}\|_{\dot{\calH}^{1}}\lesssim \lambda^{2k}
    \quad\text{and} 
    \quad
    \|\eps_{2k-1}\|_{\dot{\calH}^{2}}\lesssim \lambda^{2k+1}.
    \label{eq:e2k-1 coercivity rad}
\end{align}
At this point, we record these bounds without proofs.\footnote{A rigorous justification of these bounds in the non-radial setting is given in Lemma~\ref{lem:coercivity estimate 2}.} These bounds are sufficient to control the error terms arising in the modulation estimates. 
We introduce
\begin{align*}
    \Mod_0 \coloneqq 
    \begin{pmatrix}
        \frac{\lambda_s}{\lambda}+b_1 \\
        \gamma_s-\frac{\eta_1}{2}
    \end{pmatrix},
    \quad
    \Mod_k \coloneqq 
    \begin{pmatrix}
        (b_k)_s - b_{k+1} + (2k-\frac12)b_1 b_k + \frac12 \eta_1 \eta_k \\
        (\eta_k)_s - \eta_{k+1} + (2k-\frac12)b_1 \eta_k - \frac12 \eta_1 b_k
    \end{pmatrix}.
\end{align*}
We also define the complexified modulation parameters by
\begin{align*}
    \bm{\beta}_{2k-1}\coloneqq -\tfrac{1}{2}(-i)^{k-1}(i b_k+\eta_k),\quad \beta_{2k-1}\coloneqq |\bm{\beta}_{2k-1}|.
\end{align*}
Then, by interpolating the hierarchy of conservation laws with \eqref{eq:e2k-1 coercivity rad}, and using the decomposition $w_{2k-1}=P_{2k-1}+\eps_{2k-1}$, we obtain
\begin{align*}
    \beta_{2k-1} R^{1/2}\sim\|\chi_R P_{2k-1}\|_{L^2}
    &\lesssim \|\chi_R w_{2k-1}\|_{L^2}+\|\chi_R \eps_{2k-1}\|_{L^2}
    \lesssim \lambda^{2k-1}+\lambda^{2k}R.
\end{align*}
Choosing $R=\lambda^{-1}$ then yields the scale bound
\begin{align}
    \beta_{2k-1}\lesssim \lambda^{2k-\frac{1}{2}}. \label{eq:beta scale bdd rad}
\end{align}
In addition, one can observe that
\begin{align}
    |\Mod_k|=|(\bm{\beta}_{2k-1})_s +
    (2k-\tfrac{1}{2})b_1\bm{\beta}_{2k-1} + \tfrac{i}{2}\eta_1\bm{\beta}_{2k-1} -i\bm{\beta}_{2k+1}|. \label{eq:mod equ rad 3}
\end{align}
Then, as a result of the coercivity estimates \eqref{eq:e2k-1 coercivity rad}, for $1\le k\le L-1$, one can prove\footnote{The non-radial analogue of \eqref{eq:mod equ rad 2} is given in Proposition~\ref{prop:modj}.}
\begin{equation}
    \begin{aligned}
        |\Mod_0|\lesssim \lambda^2,\quad
        |\Mod_k|
        \lesssim \beta_1\lambda^{2k}+\lambda^{2k+2}.
    \end{aligned} \label{eq:mod equ rad 2}
\end{equation}
We reformulate \eqref{eq:mod equ rad 2} in a scale-invariant form. Denote normalized parameters
\begin{align*}
    Z_0\coloneqq \lambda^{\frac12}e^{i\gamma},\quad
    Z_{2k-1}\coloneqq \frac{\bm{\beta}_{2k-1}}{\lambda^{2k-\frac12}} e^{i\gamma}
    =-\frac{1}{2}(-i)^{k-1}\cdot\frac{i b_k+\eta_k}{\lambda^{2k-\frac12}} e^{i\gamma},
\end{align*}
so that \eqref{eq:beta scale bdd rad} and \eqref{eq:mod equ rad 2} become
\begin{align}
    |Z_{2k-1}|\lesssim 1,\quad
    |(Z_0)_t+iZ_1|\lesssim \lambda^{\frac12},\quad
    |( Z_{2k-1})_t-i Z_{2k+1}|\lesssim |Z_1|+\lambda^{\frac12}. \label{eq:mod equ rad normalize}
\end{align}
Neglecting the error terms, the system reads
\begin{align*}
    (Z_0)_t &\approx - i Z_1,\\
    (Z_{2k-1})_t &\approx i Z_{2k+1} \quad \text{for } k\ge 1 .
\end{align*}
From \eqref{eq:mod equ rad normalize}, one checks that the trajectories $Z_{2k-1}$ are \textit{uniformly} Lipschitz on $[t_0,T)$ for $1\le k\le L-1$. 
Consequently, each $Z_{2k-1}$ with $1\leq k\le L-1$ admits a finite limit as $t\to T$, which we denote by
\[
   z_{2k-1}^*\coloneqq \lim_{t\to T} Z_{2k-1}(t).
\]
With these in hand, the asymptotics of $\lambda$ are determined by the first index $2k-1$ ($k\le L-1$) for which $z_{2k-1}^*\neq 0$.

Roughly speaking, the relations can be viewed as generating a Taylor expansion of $Z_0$ around $t=T$: 
each identity expresses a higher derivative of $Z_0$ in terms of the next coefficient. 
The first index $2k-1$ for which $z_{2k-1}^*\neq0$ determines the leading non-vanishing term in this expansion, so that $Z_0$ admits a polynomial profile in $(T-t)$ of order $k$. 
In turn, this directly fixes the quantized blow-up rate for $\lambda$. 

Within this viewpoint, the error in \eqref{eq:mod equ rad normalize} is precisely what is needed to carry out the expansion to arbitrary order.
For the convergence of each $Z_{2k-1}$, uniform Lipschitz continuity, namely $|(Z_{2k-1})_t-iZ_{2k+1}|\lesssim 1$, already suffices.
However, to proceed inductively and extract higher-order coefficients in the Taylor-type expansion, the error terms must be sufficiently small so that each time integration yields an additional polynomial gain in $(T-t)$.

From the Taylor-expansion viewpoint described above, it becomes clear why the $\lambda^{\frac12}$-level error in \eqref{eq:mod equ rad normalize} is sufficient for our argument.
At this stage, we already know that $\lambda(t)\to0$ as $t\to T$ and that the normalized variables $Z_{2k-1}$ remain uniformly bounded.

If $Z_1$ converges to a nonzero limit, then the first blow-up rate follows immediately from the relation $|(Z_0)_t+iZ_1|\lesssim\lambda^{\frac12}$.
If instead $Z_1(t)\to0$, then the same estimate shows that $|Z_1(t)|\sim \lambda^{\frac12}(t)/(T-t)$; this bound is also sufficient for our purpose.

More generally, each step of the Taylor-type expansion corresponds to additional time integration.
As a consequence, proceeding to the $k$-th step requires the error terms to gain a factor of $(T-t)$ at each iteration.
Since $Z_0$ has size $\lambda^{\frac12}$ by definition, an error of order $\lambda^{\frac12}$ precisely ensures this polynomial gain at every step, and thus allows the expansion to be continued to arbitrarily high order.
Errors of smaller, for instance of order $\lambda^{\frac12-\epsilon}$ for some $\epsilon>0$, would fail to produce sufficient decay after finitely many iterations, preventing further refinement of the expansion.

More explicitly, this Taylor-type expansion shows that $Z_0$ behaves like a Taylor polynomial in $(T-t)$ whose coefficients are given by the limits $z_{2k-1}^*$.
If $z_{2k-1}^*$ is the first non-vanishing coefficient, then\footnote{A rigorous proof for non-radial case can be found in Section~\ref{sec:classification}.}
\begin{equation*}
    \begin{aligned}
        Z_0(t)&=-\tfrac{(-i)^{k}}{k!}\,(z_{2k-1}^*+o_{t\to T}(1))(T-t)^{k},
        \\
        Z_{2j-1}(t)&=-\tfrac{(-i)^{k-j}}{(k-j)!}(z_{2k-1}^*+o_{t\to T}(1))(T-t)^{k-j},\qquad 1\leq i \leq k,
        \\
        \lambda(t)&=\tfrac{1}{(k!)^2}(|z_{2k-1}^*|^2+o_{t\to T}(1))(T-t)^{2k}.
    \end{aligned}
\end{equation*}
For example,
\begin{align*}
  z_1^*\neq0 &\Longrightarrow \lambda(t)\sim (T-t)^2,\\
  z_1^*=0,\; z_3^*\neq0 &\Longrightarrow \lambda(t)\sim (T-t)^4, \\
  & \quad \vdots 
\end{align*}
If $z_{2k-1}^*=0$ for all $1\le k\le L-1$, then $\lambda(t)\lesssim (T-t)^{2L-2}$. To capture the precise asymptotics at the final level $k=L$, one needs refined modulation estimates for the last mode, $Z_{2L-1}$. 
This is achieved by introducing a refined parameter; see Sections~\ref{sec:refined mod esti} and \ref{sec:classification} for details. 
In particular,
\begin{align*}
  \lambda(t)=\tfrac{1}{(L!)^2}(|\td z_{2L-1}^*|^2+o_{t\to T}(1))(T-t)^{2L},
\end{align*}
and if moreover $\td z_{2L-1}^*=0$, then
\begin{align*}
  \lambda(t)\lesssim (T-t)^{2L+\frac32}.
\end{align*}

To summarize, the radial reduction provides a simplified picture of the classification scheme. In this setting, the modulation system closes at the level of the parameters $Z_{2k-1}$, and the blow-up rate is determined solely by the first non-vanishing limit $z_{2k-1}^*$.
The Taylor-type structure of the normalized system already captures the essential mechanism underlying the quantized blow-up rates.

In the general non-symmetric setting, however, this simplification is no longer available, and several genuinely new difficulties arise.
First, one must simultaneously handle all components $w_j$ for both odd and even indices, rather than restricting to the odd ones only.
Moreover, slow odd radiations in the $w$-variable produce nontrivial contributions that are completely invisible in the radial case, which leads to the introduction of additional dynamical parameters.

The analytic origin of this phenomenon lies in the structure of the kernel of the operator $A_Q=\calB_1$.
Since $\ker A_Q=\mathrm{span}_{\bbC}\{Q,yQ\}$, the radial setting eliminates the $Q$-direction for odd indices, and coercivity can be recovered by tracking only the $yQ$-mode coefficients.
In contrast, in the non-symmetric case the $Q$-component is no longer absent.
To recover coercivity, one must explicitly extract this component in the decomposition of $w_j$, which naturally leads to the introduction of additional radiation parameters $c_j$ through the projection onto the $Q$-direction.

As a consequence, the classification problem is no longer governed solely by the parameters $\{Z_{2k-1}\}$.
To close the system, one must derive precise evolution equations for the radiation parameters $c_j$.
At this level, quadratic radiation interaction terms arise as leading-order contributions rather than mere error terms, and their dynamics must be taken into account.

Since this procedure cannot be iterated indefinitely, these quadratic interactions are reorganized into the quantities $\frakc_{i,j}$.
A key point is to choose $\frakc_{i,j}$ so that their evolution forms a closed system, allowing the radiation effects to be incorporated without introducing new parameters at each step.
This substantially complicates the analysis and constitutes one of the main technical challenges beyond the radial setting.
A detailed treatment of this mechanism is carried out in Section~\ref{sec:rad esti}.

\section{Decompositions and coercivity estimates}\label{sec:decom and coerc}

In this section, we introduce several types of decompositions and coercivity estimates that will be used throughout the paper. The analysis here is \emph{time–independent} in the sense that we do not study the time evolution of the renormalized variables.
The main objective is twofold: (i) to decompose the solution into a profile part and a radiation part and (ii) to establish coercivity estimates that yield quantitative controls of higher-order energies associated with the radiation. 

To set up the framework, we introduce the renormalized coordinates and variables that will be used throughout the remainder of the paper.
We define
\begin{align*}
	\frac{ds}{dt}=\frac{1}{\lambda(t)^{2}},\quad 
    y=\frac{x-x(t)}{\lambda(t)},\quad 
    w(s,y)=\lambda^{\frac{1}{2}}e^{-i\gamma}v(t,\lambda(t)y+x(t))\big|_{t=t(s)}. %\label{eq:renormal}
\end{align*}
For notational convenience, we write $\lambda(s)=\lambda(t(s))$, $\gamma(s)=\gamma(t(s))$, and $x(s)=x(t(s))$, and adopt the same convention for other parameters. 
The (renormalized) nonlinear adapted derivatives $w_j$ are defined by
\begin{align*}
	w_j \coloneqq \td{\bfD}_w^j w
    = \lambda^{\frac{2j+1}{2}} e^{-i\gamma}(v_j)(t,\lambda(t)y+x(t))\big|_{t=t(s)}.
    %\label{eq:wk def}
\end{align*}
The evolution equations satisfied by $w$ and $w_j$ will be introduced later in Section~\ref{sec:mod esti}.
In the present section, we only focus on the structural properties of these renormalized variables, which will be used to establish the decompositions and coercivity estimates.

With these renormalized variables in place, we now focus on the decompositions and coercivity estimates.
To achieve this, we employ two complementary strategies. First, we use different decompositions of solutions depending on the norms under consideration. This strategy, introduced in \cite{KimKimKwon2024arxiv}, is crucial when dealing with slowly  decaying tails of corrector profiles. Second, we employ nonlinear higher-order variables, a method first introduced in \cite{KimKwonOh2020blowup} and subsequently refined in \cite{KimKimKwon2024arxiv, JeongKim2024arXiv} for systems with a hierarchy of conservation laws.

However, the hierarchy of conservation laws alone is not sufficient for the present classification problem. 
Although it yields nonlinear higher-order energy estimates, these differ from the higher-order energy bounds derived from the linear coercivity by an error of order $\lambda^2$. 
Such an error is negligible at the $\dot H^2$-level, which suffices for the construction of a particular blow-up solution, but becomes significant at the $\dot H^3$-level required for classification.
We therefore introduce an additional decomposition to compensate for the $\lambda^2$-level mismatch between the hierarchy of conservation laws and the linear higher-order energy estimates.
This auxiliary decomposition, denoted by the \textit{check} $\check{}$-notation, is designed to recover $\dot H^3$-level coercivity that is not guaranteed by the primary decomposition alone.
The \textit{check} decomposition isolates the non-coercive directions and establishes higher-order energy estimates.

For the reader's convenience, we collect here the decomposition formulas and the corresponding scaling bounds that will be used throughout the paper. 
Let $w$ and $w_j$ denote the renormalized solutions in $(s,y)$-variables to \eqref{CMdnls-gauged}, satisfying \eqref{eq:wj equ} for each $j=0,1,\dots$, with the convention $w=w_0$. 

More precisely, we use the following decompositions for $w$:
\begin{equation}
    \begin{aligned}
        w&=Q+\wt\eps, 
        \\
        &=Q+T+\eps, 
        \\
        &=Q+\check{T}+\check{\eps},
    \end{aligned} \label{eq:w decomposition}
\end{equation}
that satisfies the scaling bounds
\begin{equation*}
    \|\wt\eps\|_{\dot \calH^1}\lesssim \lambda,\quad
    \|\eps\|_{\dot \calH^2}\lesssim \lambda^2, \quad 
     \|\check{\eps}_e\|_{\dot \calH^3}\lesssim \lambda^3,
\end{equation*}
as stated in \eqref{eq:check eps decom} and \eqref{eq:w decom h2 new} with Lemma~\ref{lem:proximity}, where
\begin{align*}
    T&=T(\nu_0,b_1,\eta_1)
    =i\nu_0\tfrac{y}{2}Q-ib_1\tfrac{y^2}{4}Q-\eta_1\tfrac{1+y^2}{4}Q,
    \\
    \check{T}&=T(\check{\nu}_0,\check{b}_1,\check{\eta}_1)
    =i\check{\nu}_0\tfrac{y}{2}Q-i\check{b}_1\tfrac{y^2}{4}Q-\check{\eta}_1\tfrac{1+y^2}{4}Q.
\end{align*}
For higher-order variables $w_j$, we will decompose them into 
\begin{equation}
    \begin{aligned}
        w_j &= P_{j,e}+\wt\eps_j \\
        &= P_j+\eps_j,\\
        &= P_j+\check{\bm{\beta}}_{j+1}(1+y^2)Q+\check\eps_j,
    \end{aligned} \label{eq:wj decomposition}
\end{equation}
that satisfy bounds
\begin{gather*}
    \|\eps_j\|_{\dot\calH^1}\lesssim \lambda^{j+1},\quad
    \|\check\eps_{j}\|_{\dot\calH^2}\lesssim \lambda^{j+2},\quad \|\check\eps_{j,e}\|_{\dot\calH^3}\lesssim \lambda^{j+2+\frac12},
    \\
    \check{\bm{\beta}}_{j+1}=\bm{\beta}_{j+1}-\tfrac{1}{2}\bm{\beta}_1c_j+O(\lambda^{j+2}),
\end{gather*}
as in \eqref{eq:wj decom 1} and \eqref{eq:eps je further decom} with Lemma~\ref{lem:coercivity estimate 2}, where
\begin{align*}
    P_{j}=c_jQ+\bm{\beta}_{j}yQ,\quad P_{j,e}=c_jQ,\quad P_{j,o}=\bm{\beta}_{j}yQ.
\end{align*}

We begin in Section~\ref{subsec:modified profile P1} with a refinement of the basic decomposition $w=Q+\widehat\eps$, introducing an auxiliary \emph{check} decomposition that yields $\dot H^3$-level coercivity. 
In Section~\ref{sec:wj decomposition} we extend this framework to the nonlinear higher-order variables $w_j$, decomposing each $w_j$ into its profile part $P_j$, modulation remainder $\eps_j$, and radiation component parameterized by $c_j$. 
The $\dot H^3$ coercivity estimates for the even part of $\check\eps$ and $\check\eps_j$ are established at the end of Section~\ref{sec:wj decomposition}; these estimates will later play a key role in the radiation analysis of Section~\ref{sec:rad esti}.

\subsection{Decomposition for \texorpdfstring{$w$}{w} and modified profile \texorpdfstring{$P_1$}{P1}}\label{subsec:modified profile P1}
In this subsection, we take the decomposition for $w$ and obtain its coercivity estimates. This part was originally introduced in \cite{KimKimKwon2024arxiv}. Starting from the basic decomposition $w = Q + \widehat{\eps}$ in Proposition~\ref{prop:Decomposition}, the authors in \cite{KimKimKwon2024arxiv} introduced a further decomposition $\widehat{\eps} = \check{T} + \check{\eps}$ to establish coercivity at the $\dot{H}^2$-level. 
In the present work, we recall this construction because it provides the appropriate framework for deriving $\dot{H}^3$-level estimates directly for the $w$ variable, which will be essential for the later analysis of radiation.
If one parametrizes $w$ only by the quantities defined from $w_1$ (to be introduced in Section~\ref{sec:wj decomposition}), namely $(b_1,\eta_1,\nu_0)$, then the proximity relations \eqref{eq:proximity} enforce an optimal gap of size $O(\lambda^2)$, and consequently one only obtains $\|\eps\|_{\dot{\calH}^3}\lesssim \lambda^2$ as in \eqref{eq:coercivity H2 and H1 inf}. 
For this reason, following \cite{KimKimKwon2024arxiv}, we impose orthogonality directly on $w$ and use the resulting decomposition $\widehat{\eps} = \check{T} + \check{\eps}$, under which we establish quantitative coercivity estimates for both $\widehat{\eps}$ and $\check{\eps}$.

We now recall the precise definition of this decomposition and the associated modulation parameters.
We define new modulation parameters,
\begin{align}
    \check{b}_1\coloneqq \frac{(\widehat{\eps},\mathcal{Z}_{4})_{r}}{(\calK_4,\mathcal{Z}_{4})_{r}}
    ,\quad
    \check{\eta}_1 \coloneqq \frac{(\widehat{\eps},\mathcal{Z}_{5})_{r}}{(\calK_5,\mathcal{Z}_{5})_{r}}, \quad
    \check{\nu}_0\coloneqq \frac{(\widehat{\eps},\mathcal{Z}_{6})_{r}}{(\calK_6,\mathcal{Z}_{6})_{r}},\label{eq:check definition j=1}
\end{align} 
and then, 
\begin{align*}
    \check{T}=T(\check{\nu}_0,\check{b}_1,\check{\eta}_1)
    =i\check{\nu}_0\tfrac{y}{2}Q-i\check{b}_1\tfrac{y^2}{4}Q-\check{\eta}_1\tfrac{1+y^2}{4}Q. %\label{eq:check T}
\end{align*}
Combining Proposition~\ref{prop:Decomposition}, we have a decomposition given by the orthogonality conditions,
\begin{align}
    \wt\eps=\check{T}+\check{\eps},\quad (\check{\eps},\mathcal{Z}_{k})_{r}=0\quad\text{for} \quad 1\leq k\leq 6. \label{eq:check eps decom}
\end{align}
Note that the map $v\mapsto(\check{b}_1, \check{\eta}_1, \check{\nu}_0)$ is $C^{1}$. 
Here, main point is to define $\check{b}_1, \check{\eta}_1,$ and $\check{\nu}_0$ by imposing orthogonality conditions in $w$. 

For temporary use, following \cite{KimKimKwon2024arxiv}, we decompose $w_1$ as follows;
\begin{align*}
	w_1=\check{P}_1+\td{\eps}_1, \qquad \check{P}_{1}=P_1(\check{\nu}_0,\check{\mu}_0,\check{b}_1,\check{\eta}_1)
    \coloneqq
	(i\check{\nu}_0+\check{\mu}_0)\tfrac{1}{2}Q
    -(i\check{b}_1+\check{\eta}_1)\tfrac{y}{2}Q %\label{eq:w1 decompose}
\end{align*}
with
\begin{equation}
	\check{\mu}_0\coloneqq -\tfrac{1}{\pi}{\textstyle\int_{\bbR}}\Re(yQ^{3}\check{\eps}_{o})dy
	-\tfrac{1}{2\pi} {\textstyle\int_{\bbR}} yQ^2|\wt\eps|^2dy.\label{eq:def mu}
\end{equation}
Moreover,
\begin{align}
	\td{\eps}_1 =\widetilde{L}_{Q}\check{\eps}+\td N_{Q}(\widehat{\eps})-yQ\cdot\tfrac{1}{2\pi}{\textstyle \int_{\R}\Re(Q^{3}\check{\eps}_e )dy},\label{eq:e prime relation}
\end{align}
with
\begin{equation}
    \begin{aligned}
    	\td L_Q&=\partial_y+\tfrac{1}{2}\calH(Q^2)+Q^{-1}\calH\Re(Q^3\ \cdot \ ),
    	\\
    	\td N_Q(\wt\eps)
    	&=\wt\eps\mathcal{H}(\Re(Q\wt\eps))+\tfrac{1}{2}\wt\eps\mathcal{H}(|\wt\eps|^{2})+\tfrac{1}{4}yQ\mathcal{H}(yQ^2|\wt\eps|^{2})
    	+\tfrac{1}{4}Q\mathcal{H}(Q^2|\wt\eps|^{2}).
    \end{aligned}\label{eq:L N tilde}
\end{equation}
We emphasize that the $Q$-direction in \eqref{eq:def mu} and \eqref{eq:e prime relation} does not come from the linearized kernel alone; see \cite[Section~5.3]{KimKimKwon2024arxiv} for the motivation of \eqref{eq:def mu}. It arises from the strong interaction between the radiation and the main profile, and, technically, from the action of the Hilbert transform.

\begin{lem}[Coercivity estimates~I]
	\label{lem:coercivity w1} Let $v(t)$ be a solution on $[0,T)$ that blows up at time $T$. Also, let $M$ and $E_1$ be the initial mass and energy of $v(t)$, respectively. Then, there exists a time $t^*>0$, depending on $M$ and $E_1$, such that the following estimates hold on the interval $[T-t^*, T)$: 
	\begin{enumerate}
		\item We have 
		\begin{align}
			\|\widehat{\eps}\|_{\dot{\mathcal{H}}^{1}}+\|w_{1}\|_{L^{2}}\lesssim \lambda,\label{eq:coercivity H1}
		\end{align}
		from which we deduce 
		\begin{align}
			\|\widehat{\eps}\|_{L^{\infty}}\lesssim \lambda^{\frac{1}{2}}.\label{eq:coercivity L inf}
		\end{align}
		\item We have 
		\begin{align}
			\|\check{\eps}\|_{\dot{\mathcal{H}}^{2}} & \lesssim \lambda^{2}.\label{eq:coercivity H2 check}
		\end{align}
		In addition, we have interpolation estimates, 
		\begin{align}
			\|\check{\eps}\|_{\dot\calH^{1,\infty}}
			\lesssim \lambda^{\frac32}.\label{eq:coercivity H1 inf check}
		\end{align}
		\item We have estimates for modulation parameters: 
		\begin{equation}
            |\check{b}_1|+|\check{\eta}_1|\lesssim \lambda^{\frac32}, \quad 
			|\check{\nu}_0|+|\check{\mu}_0|\lesssim \lambda.\label{eq:mu check estimate}
		\end{equation}
		
	\end{enumerate}
\end{lem}
The proof of Lemma \ref{lem:coercivity w1} is given in Section~\ref{sec:coer pf}.

\subsection{Nonlinear higher-order variables and their decompositions}\label{sec:wj decomposition}
Now, we introduce the main modulation decomposition for higher-order nonlinear variables $w_j$. Motivated by the Lax pair structure \eqref{eq:hierarchy} and \eqref{eq:vk equ}, we understand \eqref{CMdnls-gauged} as a system of $w$ together with $w_j$ for $j=1,2,\cdots, 2L+1$. Each $w_j$ is decomposed into a profile part $P_j$ and an error part $\eps_j$, for which the necessary coercivity is ensured. The argument of decomposition is similar to earlier work \cite{JeongKim2024arXiv}, which was carried out under a radial (evenness) assumption. In that setting, the dynamics are simplified due to the absence of odd components. However, for general solutions, there arise issues coming from \textit{slow odd radiations} in the $w$-variable. To address this, we introduce additional parameters $c_j$, referred to as \emph{radiation parameters}, whose evolutions are to be tracked. See Section~\ref{sec:rad esti} for details. 

Let us recall the linear operators $\calL_{\ell}$ introduced in Section~\ref{sec:linearized dynamics} (see \eqref{eq:newdef BQLQ}) and their kernel (see \eqref{eq:general kernel nonrad}):
\begin{align*}
	\textnormal{ker}\calL_{\ell}
	=\textnormal{span}_{\mathbb{R}} \{iQ,\Lambda Q,Q_y,iyQ, iy^{j+1}Q, y^{j-1}(1+y^2)Q: 1\leq j\leq \ell-1\}.
\end{align*}
The directions $iy^{j+1}Q$ and $y^{j-1}(1+y^2)Q$ for $w$ correspond, respectively, to the profiles $iyQ$ and $yQ$ for $w_j$. We therefore decompose $w_j$ into the $iyQ$ and $yQ$ parts and the rest. Thus, for $1\le j\le 2L$, we define complex-valued parameters $\bm{\beta}_j$ by
%We therefore decompose $w_j$ into the profiles $iyQ$ and $yQ$ parts and the rest. Thus, for $1\le j\le 2L$, we define complex-valued parameters $\bm{\beta}_j$ by
\begin{align*}
    \bm{\beta}_j\coloneqq \frac{(w_{j},\calZ_h)_c}{(yQ,\calZ_h)_r},\qquad \calZ_h\coloneqq yQ\chi,\qquad  \beta_j\coloneqq |\bm{\beta}_j|.
\end{align*}
In addition to this $yQ$-direction, since the kernel of $\td\bfD_Q$ contains $yQ$ and $Q$, one also needs to control the $Q$-direction in $w_j$ in order to exploit coercivity for $w_j$. This naturally suggests the introduction of additional parameters corresponding to the $Q$-direction in $w_j$. This can be viewed as a generalization of the introduction of $\check{\mu}_0$ for $w_1$ in the previous section (or in \cite{KimKimKwon2024arxiv}). We also emphasize that, just as $\check{\mu}_0$ did not originate from the kernel direction of $w$ but from nonlinear interactions, the present $Q$-direction likewise does not arise from the linearized dynamics alone, but emerges from the interaction between the odd radiation and the main profile through the Hilbert transform.

Fix $\calZ_c\in C_c^\infty$ satisfying 
\begin{equation*}
    (Q,\calZ_c)_r\neq 0 \quad \text{and} \quad ((1+y^2)Q,\calZ_c)_r=0.
\end{equation*}
We define, for $1\le j\le 2L+1$, the (complex-valued) radiation parameter $c_j$ by
\begin{align}
    c_j\coloneqq \frac{(w_{j},\calZ_{c} )_c}{(Q,\calZ_{c})_r}. \label{eq:def cj}
\end{align}
These choices are equivalent to imposing the orthogonality conditions
\begin{align}
    (\wt\eps_j,\calZ_c)_c=0, \qquad (\eps_j,\calZ_c)_c=(\eps_j,\calZ_h)_c=0, \label{eq:epsj orthogonal}
\end{align}
and obtain the two-step decomposition
\begin{equation}
    \begin{aligned}
        P_j&\coloneqq c_jQ+\bm{\beta}_j yQ,\qquad P_{j,e}=c_jQ,\qquad P_{j,o}=\bm{\beta}_jyQ, \\
        w_{j}&=P_{j,e}+\wt\eps_{j}=P_{j}+ \eps_{j},
        \qquad
    	\wt\eps_j=P_{j,o}+\eps_j.
    \end{aligned} \label{eq:wj decom 1}
\end{equation}

When $j=1$, The decomposition \eqref{eq:wj decom 1} for $w_1$ has a proximity to earlier decompositions \eqref{eq:check definition j=1},\eqref{eq:check eps decom}, and \eqref{eq:def mu} for $w=Q+\check{T}+\check{\eps}$. More precisely, the linearization of $w_1$ is given by
\begin{align*}
    w_1=L_Q\wt\eps + O(|\wt\eps|^2),
\end{align*}
and a formal computation shows
\begin{align*}
    L_Q\check{T}=\tfrac{1}{2}i\check{\nu}_0 Q-\tfrac12(i\check{b}_1+\check{\eta}_1)yQ.
\end{align*}
Thus, we can expect the proximity between  $(\check{\nu}_0,\check{b}_1,\check{\eta}_1)$ and $(\nu_0,b_1,\eta_1)$, which is defining 
\begin{align}
    b_1\coloneqq -2\Im(\bm{\beta}_1),&\quad \eta_1\coloneqq -2\Re(\bm{\beta}_1),\quad \bm{\beta}_1=-\tfrac12(ib_1+\eta_1), \label{eq:b para def} \\
    	\nu_0\coloneqq 2\Im(c_1),&\quad \mu_0\coloneqq 2\Re(c_1),\quad \quad \, c_1=\tfrac12(i\nu_0+\mu_0). \label{eq: nu mu new def}
\end{align} 
In the case of $\check{\mu}_0$, this parameter was designed to remove non-trivial interactions arising from the Hilbert transform, which ensures that $\td{\eps}_1$ has sufficiently small values, $\|\td{\eps}_1\|_{\dot\calH^1}\lesssim \lambda^2$. See also \cite{KimKimKwon2024arxiv}. Since the definition of $\mu_0$ leads to the coercivity of $\eps_1$, and hence we have  $\|\eps_1\|_{\dot\calH^1}\lesssim \lambda^2$, we can also expect the proximity between $\mu_0$ and $\check{\mu}_0$.

%\begin{rem}\label{rem:dif quantize construct}
%    In \cite{JeongKim2024arXiv}, 
%    the authors used the modulation parameters by
%    \begin{equation*}
%    	\begin{aligned}
%    		b_k&\coloneqq \frac{(w_{2k-1},(-i)^{k}\calZ_h)_r}{(yQ,\calZ_h)_r},\quad 
%    		\eta_k\coloneqq \frac{(w_{2k-1},-(-i)^{k-1}\calZ_h)_r}{(yQ,\calZ_h)_r},
%    	\end{aligned}\quad k\geq 1.
%    \end{equation*}
%    Accordingly, the authors used the following profiles,
%    \begin{align*}
%	   P_{2k-1,o}&=-(-i)^{k-1}(ib_k+\eta_k)\tfrac{y}{2}Q.
%    \end{align*}
%    In our context, this coincides with
%    \begin{align*}
%    	\bm{\beta}_{2k-1}&= -(-i)^{k-1}(ib_k+\eta_k)\tfrac{1}{2}.
%    \end{align*}
%    However, instead of $b_k$ and $\eta_k$, we mainly use the complex-valued form, $\bm{\beta}_{j}$ since the separation of $b_k $ and $\eta_k$ is unnecessary.
%\end{rem}

We also use a reparametrization for the base variable $w$: define a new decomposition for $\wt\eps$ by
\begin{align}
    \eps\coloneqq \wt\eps - T(\nu_0, b_1, \eta_1),\qquad w=Q+T+\eps  \label{eq:w decom h2 new}
\end{align}
where $(\nu_0,b_1,\eta_1)$ is given by \eqref{eq: nu mu new def} and \eqref{eq:b para def}. 
The check parameters $(\check b_1,\check\eta_1,\check\nu_0,\check\mu_0)$ introduced in Section~\ref{subsec:modified profile P1} will be used when we reach $H^3$-level coercivity; for $H^2$-level precision, they can be replaced by $(b_1,\eta_1,\nu_0)$ through the following proximity estimates.

\begin{lem}[Proximity between decompositions]\label{lem:proximity}
	\begin{align}
		|\check{\nu}_0-\nu_0|+|\check{\mu}_0-\mu_0|
		+|\check{b}_1-b_1|
		+|\check{\eta}_1-\eta_1|
		\lesssim \lambda^2. \label{eq:proximity}
	\end{align}
	We also have
	\begin{align}
		\|\eps\|_{\dot \calH^2}\lesssim\lambda^2, \quad
		\|\eps\|_{\dot \calH^{1,\infty}}\lesssim   \lambda^{\frac32}. \label{eq:coercivity H2 and H1 inf}
	\end{align}
    In particular, we also have
    \begin{align}
        |b_1|+|\eta_1|\lesssim \lambda^{\frac32},\quad |\nu_0|+|\mu_0|\lesssim \lambda. \label{eq: b1 estimate}
    \end{align}
\end{lem}
\begin{proof}
	We first estimate the $\nu_0$- and $\mu_0$-terms of \eqref{eq:proximity} from the even part of \eqref{eq:e prime relation}:
	\begin{align*}
		c_1Q+\eps_{1,e}=w_{1,e}=\check{P}_{1,e} +\td{\eps}_{1,e} 
		=\check{P}_{1,e} +\td L_Q\check{\eps}_o +[\td N_Q(\wt\eps)]_e,
	\end{align*}
	where $\td L_Q$ and $\td N_Q(\wt\eps)$ are given by \eqref{eq:L N tilde}.
	Collecting all $Q$-terms, we have
	\begin{align}
		c_1 Q -\check{P}_{1,e} =\tfrac{1}{2}[i(\nu_0-\check{\nu}_0)+(\mu_0-\check{\mu}_0)]Q
		=-\eps_{1,e}+\td L_Q\check{\eps}_o +[\td N_Q(\wt\eps)]_e. \label{eq:nnu mu proximity}
	\end{align}
	Taking the $L^2$ norm after multiplying \eqref{eq:nnu mu proximity} by $Q$, we obtain
	\begin{align*}
		|\check{\nu}_0-\nu_0|+|\check{\mu}_0-\mu_0|
		&\lesssim \|\eps_1\|_{\dot\calH^1}+\|\check{\eps}_o \|_{\dot\calH^2}+\|Q[\td N_Q(\wt\eps)]_e\|_{L^2}
		\lesssim \lambda^2+\|Q[\td N_Q(\wt\eps)]_e\|_{L^2}.
	\end{align*}
	Here, we used \eqref{eq:coercivity H2 check} and \eqref{eq:coercivity H1 epsj}. For the nonlinear part, we can check that
	\begin{align*}
		\|Q[\td N_Q(\wt\eps)]_e\|_{L^2}
		&\lesssim \|Q\wt\eps\|_{L^\infty}\|Q\wt\eps\|_{L^2}+\|Q\wt\eps\|_{L^2}\|\calH(|\wt\eps|^2)\|_{L^\infty}+\|[yQ^2\calH(yQ^2|\wt\eps|^2)]_e\|_{L^2}
		\\
		&\lesssim \lambda^2 +\lambda \|\calH(|\wt\eps|^2)\|_{L^\infty}+\|yQ^2\calH(yQ^2[|\wt\eps|^2]_o)\|_{L^2}.
	\end{align*}
    Using Gagliardo--Nirenberg inequality, we have
	\begin{align*}
		\|\calH(|\wt\eps|^2)\|_{L^\infty}\lesssim \||\wt\eps|^2\|_{L^2}^{\frac{1}{2}} \|\partial_y|\wt\eps|^2\|_{L^2}^{\frac{1}{2}}
		\lesssim \|\wt\eps\|_{L^2}\|\partial_y\wt\eps\|_{L^2}\lesssim \lambda.
	\end{align*}
	For the last term, thanks to \eqref{eq:CommuteHilbert}, we observe that
	\begin{align*}
		\|yQ^2\calH(yQ^2[|\wt\eps|^2]_o)\|_{L^2}=\|y^2Q^2\calH(Q^2[|\wt\eps|^2]_o)\|_{L^2}
		\lesssim 
		\|Q\wt\eps\|_{L^\infty}\|Q\wt\eps\|_{L^2}\lesssim  \lambda^2
	\end{align*}
	%Therefore, we have conclude 
	%\begin{align*}
	%	\|Q[\td N_Q(\wt\eps)]_e\|_{L^2}\lesssim \lambda^2,
	%\end{align*}
	and we conclude $|\check{\nu}_0-\nu_0|+|\check{\mu}_0-\mu_0| \lesssim \lambda^2$.
	%\begin{align*}
	%	|\check{\nu}_0-\nu_0|+|\check{\mu}_0-\mu_0|
	%	\lesssim \lambda^2.
	%\end{align*}
    
    Similarly, we extract the $b_1$- and $\eta_1$-terms of \eqref{eq:proximity} from the odd part of \eqref{eq:e prime relation}:
	\begin{align*}
		-\tfrac{1}{2}(i(b_1-\check{b}_1)+(\eta_1-\check{\eta}_1))yQ=-\eps_{1,o}+\widetilde{L}_{Q}\check{\eps}_e +[\td N_{Q}(\widehat{\eps})]_{o}-yQ\cdot\tfrac{1}{2\pi}{\textstyle \int_{\R}\Re(Q^{3}\check{\eps}_e )dy}.
	\end{align*}
	We again take the $L^2$ norm after multiplying by $Q$; we have
	\begin{align*}
		|\check{b}_1-b_1|
		+|\check{\eta}_1-\eta_1|
		&\lesssim \|\eps_{1}\|_{\dot\calH^1}+\|\check{\eps}\|_{\dot\calH^2}+\|Q[\td N_{Q}(\widehat{\eps})]_{o}\|_{L^2}
		\lesssim \lambda^2 +\|Q[\td N_{Q}(\widehat{\eps})]_{o}\|_{L^2}.
	\end{align*}
	We note that
	\begin{align*}
		{\textstyle \int_{\R}\Re(Q^{3}\check{\eps}_e )dy}\leq \|Q\|_{L^2}\|Q^2\check{\eps}_e \|_{L^2}\lesssim \|\check{\eps}\|_{\dot\calH^2}.
	\end{align*}
	For the nonlinear term, we proceed as above and obtain
	\begin{align*}
		\|Q[\td N_{Q}(\widehat{\eps})]_{o}\|_{L^2}\lesssim \lambda^2+\|yQ^2\calH(yQ^2[|\wt\eps|^2]_e)\|_{L^2}.
	\end{align*}
	For the last term, using \eqref{eq:CommuteHilbert}, we derive
	\begin{align*}
		\|yQ^2\calH(yQ^2[|\wt\eps|^2]_e)\|_{L^2}
		&=\|y^2Q^2\calH(Q^2[|\wt\eps|^2]_e)- \tfrac{1}{\pi}yQ^2\cdot{\textstyle\int_\bbR} Q^2[|\wt\eps|^2]_e dy\|_{L^2}\\
        &\lesssim \|Q^2|\wt\eps|^2\|_{L^2}+\|Q\wt\eps\|_{L^2}^2\lesssim \lambda^2.
	\end{align*}
	Therefore, we conclude \eqref{eq:proximity} by obtaining $|\check{b}_1-b_1|
		+|\check{\eta}_1-\eta_1|
		\lesssim \lambda^2$.
	%\begin{align*}
	%	|\check{b}_1-b_1|
	%	+|\check{\eta}_1-\eta_1|
	%	\lesssim \lambda^2,
	%\end{align*}
    
	Now, we show \eqref{eq:coercivity H2 and H1 inf}. By the definition of $\eps$, we have
	\begin{align*}
		\eps=\wt\eps-T (\nu_0,b_1,\eta_1)=T (\check{\nu}_0-\nu_0,\check{b}_1-b_1, \check{\eta}_1-\eta_1)+\check{\eps}.
	\end{align*}
	Thus, \eqref{eq:coercivity H2 check}, \eqref{eq:coercivity H1 inf check}, and \eqref{eq:proximity} imply \eqref{eq:coercivity H2 and H1 inf}.

    Also, \eqref{eq:proximity} and \eqref{eq:mu check estimate} yield \eqref{eq: b1 estimate}. Thus, we finish the proof of Lemma~\ref{lem:proximity}.
\end{proof}
Thanks to Lemma~\ref{lem:proximity}, in all statements at the $H^2$-level we freely replace the parameters $(\check{\nu}_0,\check{b}_1,\check{\eta}_1)$ coming from the $w$-decomposition in Section~\ref{subsec:modified profile P1} by $(\nu_0,b_1,\eta_1)$; this removes unnecessary \emph{check} notation in what follows.

As an analogy, we use an auxiliary \emph{check} decomposition for $w_j$. To obtain $H^2$-level coercivity, we further decompose $\eps_j$ by extracting the $(1+y^2)Q$ profile. Choose $\calZ_{c,2}\in C_c^\infty$ satisfying $(Q,\calZ_{c,2})_r=0$ and $((1+y^2)Q,\calZ_{c,2})_r\neq 0$, and define
\[
	\check{\bm{\beta}}_{j+1}\coloneqq \frac{(w_{j},\calZ_{c,2})_c}{((1+y^2)Q,\calZ_{c,2})_r}.
\]
We then decompose
\begin{equation}
	\eps_{j}=\check{\eps}_{j}+\check{\bm{\beta}}_{j+1}(1+y^2)Q,
	\qquad 
	\check\eps_j\in \{\calZ_c,\calZ_h,\calZ_{c,2}\}^\perp. \label{eq:eps je further decom}
\end{equation}
Since $(1+y^2)Q$ is even, this is equivalent to $\eps_{j,e}=\check{\eps}_{j,e}+\check{\bm{\beta}}_{j+1}(1+y^2)Q$ and $\eps_{j,o}=\check{\eps}_{j,o}$, and hence we have
\[
	w_j=P_{j,e}+P_{j,o}+\check{\bm{\beta}}_{j+1}(1+y^2)Q+\check\eps_j.
\]
The parameter $\check{\bm{\beta}}_{j+1}$ is tied to $\bm{\beta}_{j+1}$ (not to $\bm{\beta}_j$): it is the linear extraction of the $(1+y^2)Q$-component at the $w_j$-stage, whereas $\bm{\beta}_{j+1}$ is the next-iterate parameter arising at the $w_{j+1}$ stage. This choice is motivated by a formal observation. 
\begin{align*}
    w_{j+1}=\td{\bfD}_w w_j\approx \td{\bfD}_Q w_j
    \approx \check{\bm{\beta}}_{j+1}\td{\bfD}_Q[(1+y^2)Q]=\check{\bm{\beta}}_{j+1} yQ,
\end{align*}
which suggests $\bm{\beta}_{j+1}\approx \check{\bm{\beta}}_{j+1}$. But this is purely formal, more precisely, there is a correction. Lemma~\ref{lem:coercivity estimate 2} provides the quantitative relation
\[
    \bigl|\check{\bm{\beta}}_{j+1}-\bm{\beta}_{j+1}+\tfrac{1}{2}\bm{\beta}_1 c_j\bigr| 
    \lesssim \lambda^{j+2}.
\]
Here, $\tfrac{1}{2}\bm{\beta}_1 c_j$ comes from the interaction between the profile part and the radiation.
Consequently, in every estimate involving $H^2$-level, $\check{\bm{\beta}}_{j+1}$ can be substituted with $\bm{\beta}_{j+1}-\tfrac{1}{2}\bm{\beta}_1c_j$, allowing the residual to be absorbed. This approach follows the method of using $(\nu_0,b_1,\eta_1)$ instead of $(\check\nu_0,\check b_1,\check\eta_1)$ as described in Lemma~\ref{lem:proximity}.

With the decompositions \eqref{eq:check eps decom}, \eqref{eq:wj decom 1}, \eqref{eq:eps je further decom} and the replacements guaranteed by Lemma~\ref{lem:proximity}%in place
, we now ready to prove the $\dot H^1$- and $\dot H^2$-level coercivity for $w_j$.

\begin{lem}[Coercivity estimates~II]\label{lem:coercivity estimate 2} Assume $v_0 \in H^{2L+1} $. That is, $w_j \in L^2 $ for $1\leq j\leq 2L+1$. Then, the following estimates hold: 
    \begin{enumerate}
        \item (Consequence of the hierarchy of conservation law) For $1\leq j \leq 2L+1$, we have the conservation laws,
    	\begin{align}
    		\|B_Qw_j\|_{L^2}\lesssim \|w_j\|_{L^2}\lesssim \lambda^{j}. \label{eq:coercivity H0 epsj}
    	\end{align}

        \item ($\dot H^1$-level coercivity estimates) For $1\leq j \leq 2L$, the coercivity estimates for nonlinear variables hold true:
        \begin{align}
    		\|\eps_{j}\|_{\dot\calH^{1}}\sim \|A_Qw_{j}\|_{L^2}\lesssim \lambda^{j+1}. \label{eq:coercivity H1 epsj}
    	\end{align}
        We also have the interpolation estimate,
    	\begin{align}
    		\|\eps_j\|_{L^\infty}\lesssim \lambda^{j+\frac12}. \label{eq:coercivity Linf epsj}
    	\end{align}

        \item ($\dot H^2$-level coercivity estimates) For $1\leq j \leq 2L-1$, we have
        \begin{equation}
            \begin{aligned}
                \|\td\bfD_{Q}^{2}w_{j,o}\|_{L^2}\sim \|\eps_{j,o}\|_{\dot\calH^2}
                \lesssim \lambda^{j+2}, 
                \\
                \|\eps_{j,e}\|_{\dot\calH^2}+\|\eps_{j}\|_{\dot\calH^{1,\infty}}
                \lesssim \lambda^{j+\frac32},
            \end{aligned} \label{eq:coercivity H2 epsj}
        \end{equation}
    	In addition, we have 
    	\begin{align}
    		\|\check{\eps}_{j,e}\|_{\dot\calH^2}\lesssim \lambda^{j+2},\quad \|\check{\eps}_{j,e}\|_{\dot\calH^{1,\infty}}\lesssim \lambda^{j+\frac32}. \label{eq:coercivity H2 epsj check}
    	\end{align}

        \item (Proximity for $\check{\bm{\beta}}_{j+1}$) For $1\leq j \leq 2L-1$, we also have
        \begin{align}
            |\check{\bm{\beta}}_{j+1}-\bm{\beta}_{j+1}+\tfrac{1}{2}\bm{\beta}_1c_j|
            \lesssim \lambda^{j+2}. \label{eq:betaj and prime relation}
        \end{align}

        \item (Scale bounds for modulation parameters) For $1\leq j\leq 2L$, we have
    	\begin{align}
    		|\bm{\beta}_j|+|\check{\bm{\beta}}_{j} |\lesssim \lambda^{j+\frac{1}{2}},\label{eq:beta j estimate}
    	\end{align}
        with a convention $\check{\bm{\beta}}_{1}=-\frac{1}{2}(i\check{b}_1+\check{\eta}_1)$.

        \item (Scale bounds for radiation parameters) For $1\leq j\leq 2L+1$, we have
    	\begin{align}
    		|c_j|\lesssim  \lambda^{j}. \label{eq:c j estimate}
    	\end{align}
    \end{enumerate}
\end{lem}
The proof of Lemma~\ref{lem:coercivity estimate 2} is given in Section~\ref{sec:coer pf}.

%\subsection{Coercivity at \texorpdfstring{$\dot{H}^3$}{H3}-level}\label{sec:coercivity H3}
We next extend the coercivity estimate to the $\dot H^3$-level. 
As explained in Section~\ref{subsec:modified profile P1}, the auxiliary \textit{check} decomposition is used to obtain $\dot H^3$-coercivity, which is required for the radiation analysis in Section~\ref{sec:rad esti}. 
Recall the decompositions \eqref{eq:check eps decom} and \eqref{eq:eps je further decom}, and note that the transversality condition \eqref{eq:transversality} for the functions $\calZ_i$ ($i=1,\dots,6$) will be used in Section~\ref{sec:coer pf}.

\begin{lem}[$\dot H^3$-level coercivity estimates]\label{lem:coercivity H3 epsj} We have
	\begin{align}
		\|\check{\eps}_e -\tfrac{\nu_0^2}{8}Q\|_{\dot{\mathcal{H}}^{3}} & \lesssim\lambda^{3}.\label{eq:coercivity H3}
	\end{align}
    For $L\geq 2$ and $1\leq j \leq 2L-2$, we have
    \begin{align}
        \|\check \eps_{j,e}\|_{\dot\calH^3}\lesssim \lambda^{j+2+\frac12}. \label{eq:coercivity H3 epsj check}
    \end{align}
\end{lem}
\begin{rem}
    While \eqref{eq:coercivity H3 epsj check} is slightly weaker than the scaling bound by a factor of $\lambda^{\frac12}$, it is sufficient for our purpose. 
    At the technical level, this loss originates from the term \eqref{eq:coer2 even 1} in the proof given in Section~\ref{sec:coer pf}. 
    If one could further establish $\dot{H}^3$-level coercivity for the odd part of $\check{\eps}$, the full scaling bound would be recovered.
\end{rem}
The proof of Lemma~\ref{lem:coercivity H3 epsj} is also deferred to Section~\ref{sec:coer pf}.

\section{Modulation estimates}\label{sec:mod esti}
In our analysis, we track the dynamics of modulation parameters $\lambda,\gamma,x,\bm{\beta}_j$ and the radiation parameters $c_j$. They form an approximate finite-dimensional ODE system that governs the evolution of the solution. Our ultimate goal is to identify this system precisely and to show that it admits only quantized blow-up rates, thereby establishing the rigidity of the blow-up dynamics.

In the present section, we establish the evolution equations for the modulation parameters. 
These modulation laws are obtained by differentiating in time the orthogonality conditions \eqref{eq:Decomposition and orthogonal} and \eqref{eq:epsj orthogonal}. 
The first one, \eqref{eq:Decomposition and orthogonal}, yields equations for $\lambda_s$, $\gamma_s$, and $x_s$. The second, \eqref{eq:epsj orthogonal}, produces evolution equations for $(\bm{\beta}_j)_s$. In these relations, each $(\bm{\beta}_j)_s$ interacts only with nearby indices (through $\bm{\beta}_{j+1}$ and $\bm{\beta}_{j+2}$, and with the radiation parameter $c_j$, giving an approximate ODE system for $(\lambda,\gamma,x)$ coupled to the $\{\bm{\beta}_j\}$.

Although the PDE flow may, a priori, exhibit complicated blow-up behavior, the modulation procedure reduces a decisive part of the dynamics to a finite-dimensional approximate ODE system for $(\lambda,\gamma,x)$ coupled to $\{\bm{\beta}_j\}$. Since an ODE system admits only restricted types of trajectories compared to the full PDE, this reduction narrows the admissible blow-up scenarios and provides the lever for the classification result. In the previous work of the first and second authors, \cite{JeongKim2024arXiv}, one can select a particular trajectory of modulation parameters governed by the ODE, which in turn yields a quantized blow-up rate. The present paper, instead, aims at a classification of behaviors of modulation parameters; we cannot preselect a specific ODE solution, so we must derive sharper information in modulation laws that are uniform across all admissible trajectories. 

On the other hand, in the general setting without the symmetry assumption, it turns out that there is a nontrivial radiation part involved. In contrast to the symmetric case treated in \cite{JeongKim2024arXiv}, this radiation cannot be eliminated from the dynamics, and it contributes at leading order. Consequently, the radiation parameters $c_j$ must be tracked, and their evolution laws have to be incorporated together with those of the modulation parameters. The evolution of radiation is deferred to Section~\ref{sec:rad esti}.

In Section~\ref{subsec:modulation estimate}, we prove the modulation estimates. Section~\ref{sec:refined mod esti} develops refined modulation estimates, providing the additional control needed at the index $j=2L-1$. The details of this refinement are presented there.

\subsection{Modulation estimate}\label{subsec:modulation estimate}
To derive the modulation laws, we start from the renormalized formulation of the equation.  
Rewriting \eqref{CMdnls-gauged} and \eqref{eq:vk equ} together with \eqref{eq:Dv square} in the renormalized variables $(s,y)$, we obtain
\begin{equation}
    \begin{aligned}
        &(\partial_s-\tfrac{\lambda_s}{\lambda}\Lambda+\gamma_s i-\tfrac{x_s}{\lambda}\partial_y)w+iL_w^*w_1=0,\\
        &(\partial_s-\tfrac{\lambda_s}{\lambda}\Lambda_{-j}+\gamma_s i-\tfrac{x_s}{\lambda}\partial_y)w_j-iw_{j+2}
        =-\tfrac{i}{2}[w_1\mathcal{H}(\overline{w}w_j)-w\mathcal{H}(\overline{w_1}w_j)].
    \end{aligned}
    \label{eq:wj equ}
\end{equation}

Differentiating in time the orthogonality conditions \eqref{eq:Decomposition and orthogonal} and \eqref{eq:epsj orthogonal} 
and testing the resulting identities against the equations above, 
we derive the quantitative modulation laws for the parameters $(\lambda,\gamma,x)$ and $\{\bm{\beta}_j\}$.  
For convenience, we rewrite \eqref{eq:wj equ} in the following equivalent form:
\begin{equation}
    \begin{aligned}
        &(\partial_s-\tfrac{\lambda_s}{\lambda}\Lambda_{-j} +\gamma_s i)w_j-\tfrac{x_s}{\lambda}w_{j+1}-iw_{j+2}
        +\tfrac{1}{2}\big(\tfrac{x_s}{\lambda}-\nu_0\big)w\calH(\ol{w}w_j)
        \\
        &+\tfrac{i}{2}\big[\big(w_1-\tfrac{i\nu_0}{2}w\big)\mathcal{H}(\overline{w}w_{j})-w\mathcal{H}\big\{\overline{\big(w_1-\tfrac{i\nu_0}{2}w\big)}w_{j}\big\}\big]
        =0.
    \end{aligned} \label{eq:wj equ in section mod esti}
\end{equation}
Combining this equation with the time-differentiated orthogonality relations,
we obtain the evolution laws for $(\lambda, \gamma, x)$ and $\{\bm{\beta}_j\}$. The decompositions of $w$ and $w_j$, as well as the parameters $(\lambda, \gamma, x)$, $\{\bm{\beta}_j\}$, and $\{c_j\}$ are recalled in Section~\ref{sec:one soliton configure} and~\ref{sec:wj decomposition}. We remark that $c_1=\frac{1}{2}(i\nu_0+\mu_0)$ and $\bm{\beta}_1=\frac{1}{2}(ib_1+\eta_1)$. In differentiating the orthogonality conditions, the modified profiles for $w$ and $w_j$ are computed in the nonlinear term, and their dominant components influence the modulation dynamics. This leads to the modulation laws.
To determine the dynamical behavior, these laws must be sufficiently sharp to justify the reduction to an approximate ODE system.
The resulting estimates are summarized below.

\begin{prop}[Modulation estimates]\label{prop:modj}
    We have the following estimates:
    \begin{enumerate}
        \item (First modulation estimate) We have
    	\begin{align}
    		\left|\frac{\lambda_s}{\lambda}+b_1\right|
            +\left|\gamma_s-\frac{\eta_1}{2}\right|
            +\left|\frac{x_s}{\lambda}-\nu_0\right|\lesssim \lambda^2, \label{eq:mod0 esti}
    	\end{align}
        Moreover, we have stronger bounds:
        \begin{equation}
            \begin{aligned}
                \left|\frac{\lambda_s}{\lambda}+b_1+\frac{\mu_0\nu_0}{2}\right|
                +\left|\gamma_s-\frac{\eta_1}{2}-\frac{\nu_0^2}{2}+\frac{\mu_0^2}{4}-2\Re(c_2)\right|
                &\lesssim \lambda^{2+\frac{1}{2}},
                \\
                \left|\frac{x_s}{\lambda}-\nu_0\right|
                &\lesssim \beta_1\lambda+\lambda^{3}.
            \end{aligned}\label{eq:mod0 esti refine}
        \end{equation}

        \item (Second modulation estimate) For $L\geq 2$ and $1\leq j\leq 2L-2$, we have
        \begin{align}
            \left|\left(\partial_s - \frac{\lambda_s}{\lambda} 
            \left(j+\frac{1}{2}\right) + \gamma_si\right)\bm{\beta}_j -i\bm{\beta}_{j+2}
            -\nu_0\bm{\beta}_{j+1} +ic_j\bm{\beta}_2 \right|
            \lesssim  \beta_1 \lambda^{j+1}+\lambda^{j+3}. \label{eq:modj esti}
        \end{align}
    \end{enumerate}
\end{prop}
Estimates \eqref{eq:mod0 esti} and \eqref{eq:mod0 esti refine} determine modulation laws for $\lambda,\gamma, $ and $x$. The refined bound \eqref{eq:mod0 esti refine} is sharper and therefore allows the modulation laws to incorporate additional correction terms. Together with $\beta_1\lesssim \lambda^{\frac{3}{2}}$ from \eqref{eq:beta j estimate}, this leads to an overall improvement of the error bound in \eqref{eq:mod0 esti refine} to $\lambda^{2+\frac{1}{2}}$. Estimates \eqref{eq:modj esti} control the higher-order time derivatives of the modulation parameters $(\lambda,\gamma,x)$. Indeed, \eqref{eq:modj esti} shows that each equation for $(\bm{\beta}_j)_s$ is coupled to its neighboring parameters $\bm{\beta}_{j+1}$, $\bm{\beta}_{j+2}$, and to the radiation parameter $c_j$. Through the relations \eqref{eq:mod0 esti} and $\bm{\beta}_1=\frac{1}{2}(ib_1+\eta_1)$, these couplings connect the sequence $\{\bm{\beta}_j\}$ directly to the higher-order behavior of $(\lambda,\gamma,x)$. This coupling structure reflects the complexity of the non-radial dynamics. For comparison, in the radial setting, the modulation laws are much simpler as shown in \eqref{eq:mod equ rad 2} and \eqref{eq:mod equ rad 3}.

In this subsection, we prove Proposition~\ref{prop:modj} under the assumption of the required nonlinear estimates, whose proofs are deferred to Section~\ref{sec:nonlin pf}. In deriving the modulation laws, we separate the leading profiles arising in the nonlinear interactions from the higher-order remainders. We begin by introducing the reorganized nonlinear terms.
\begin{equation}
    \begin{aligned}
        \textnormal{NL}_0
    	:=\,&(w_1-\tfrac{i\nu_0}{2}w)\mathcal{H}(\overline{w}\wt\eps)-w\mathcal{H}\{\overline{(w_1-\tfrac{i\nu_0}{2}w)}\wt\eps\},
        \\
        \textnormal{NL}_j
        :=\,& (w_1-\tfrac{i\nu_0}{2}w)\mathcal{H}(\overline{w}\wt\eps_{j})-w\mathcal{H}\{\overline{(w_1-\tfrac{i\nu_0}{2}w)}\wt\eps_{j}\}
        \\
        &-P_{1,o}\mathcal{H}(QP_{j,o})
        +yQ\mathcal{H}\{\tfrac{y}{1+y^2}\overline{P_{1,o}}P_{j,o}\}
        +Q\mathcal{H}\{\tfrac{1}{1+y^2}\overline{P_{1,o}}P_{j,o}\},
        \\
        \textnormal{NL}_{Q}
        :=\,&
        (w_1-\tfrac{i\nu_0}{2}w)\mathcal{H}(\overline{w}Q)-w\mathcal{H}\{\overline{(w_1-\tfrac{i\nu_0}{2}w)}Q\}
        \\
        &-(P_1-\tfrac{i\nu_0}{2}Q)\mathcal{H}(Q^2)
        +Q\mathcal{H}\{\overline{(P_1-\tfrac{i\nu_0}{2}Q)}Q\}
        \\
        &+c_1T_{o}\calH(Q^2)-Q\calH(\ol{c_1 T_{o}} Q).
    \end{aligned} \label{eq:def NL}
\end{equation}
We now state the first nonlinear estimates. In their analysis, the parity decomposition should be in consideration. 

\begin{lem}[Nonlinear estimates I]\label{lem:nonlinear esti I}
	Let $f\in C^\infty$ with $Q^{-1}f\in L^2$. Then the following estimates hold:

    (1) For the terms of the form $w\calH(\ol{w}w_j)$,
        \begin{equation}
    		\begin{aligned}
    		      |(w\calH(|w|^2)-Q\calH(Q^2),f)_c|+|(w\calH(\ol{w} Q)-Q\calH(Q^2),f)_c|\lesssim \|Q^{-1}f\|_{L^2}\lambda,
    		\end{aligned} \label{eq:wH w square}
    	\end{equation}
    	and for $1\leq j\leq 2L$,
    	\begin{equation}
    		\begin{aligned}
                |(w\calH(\ol{w} \wt\eps_{j}),f_o)_c|\lesssim \|Q^{-1}f_o\|_{L^2}\lambda^{j+1},
                \\
    			|(w\calH(\ol{w} \wt\eps_{j}),f_e)_c|\lesssim \|Q^{-1}f_e\|_{L^2}\lambda^{j+\frac{1}{2}}.
    		\end{aligned} \label{eq:wH w wj}
    	\end{equation}

    (2) For the nonlinear terms, for $1\leq j\leq 2L-1$
        \begin{align}
    		|(\textnormal{NL}_j-\tfrac{i}{2}\mu_0\nu_0 \bm{\beta}_j yQ^3,f_o)_c|\lesssim \|Q^{-1}f_o\|_{L^2}(\beta_1\lambda^{j+1}+\lambda^{j+3}). \label{eq:main nonlinear eps j odd}
    	\end{align}
        For $j=0$,
        \begin{equation}
    		\begin{aligned}
    			|(\textnormal{NL}_0, f_o)_c|&\lesssim  \|Q^{-1}f_o\|_{L^2}(\beta_1\lambda+\lambda^3),
                \\
                |(\textnormal{NL}_0
    		+(2\ol{c_{2}}+\tfrac{1}{2}(\nu_0^2-\mu_0^2)+i\mu_0\nu_0)Q, f_e)_c|&\lesssim \|Q^{-1}f_e\|_{L^2}\lambda^{2+\frac12}.
    		\end{aligned} \label{eq:main nonlinear wt eps}
    	\end{equation}
        Finally,
        \begin{equation}
            \begin{aligned}
                |(\textnormal{NL}_{Q}-2\bm{\beta}_2 yQ, f_o)_c|&\lesssim \|Q^{-1}f_o\|_{L^2}(\beta_1\lambda+\lambda^{3}),
                \\
                |(\textnormal{NL}_{Q}, f_e)_c|&\lesssim \|Q^{-1}f_e\|_{L^2}\lambda^{2+\frac12}
                .
            \end{aligned}\label{eq:main nonlinear Q}
        \end{equation}
    
\end{lem}

With these estimates at hand,  we are now ready to prove Proposition~\ref{prop:modj}. 

\begin{proof}[Proof of Proposition~\ref{prop:modj} assuming Lemma~\ref{lem:nonlinear esti I}]
    We first control $(\lambda_s, \gamma_s, x_s)$ and then derive the equations for $(\bm{\beta}_j)_s$. In \underline{\textbf{Step 1}}, we prove \eqref{eq:mod0 esti} and \eqref{eq:mod0 esti refine}; the proof of the rough bound \eqref{eq:mod0 esti} is similar to \cite[Lemma 5.15]{KimKimKwon2024arxiv}. In \underline{\textbf{Step 2}}, we establish \eqref{eq:modj esti}. 

    In the proof, we derive the modulation laws \eqref{eq:mod0 esti refine} and \eqref{eq:modj esti}. To obtain these laws, we must isolate the nontrivial profile contributions hidden inside the nonlinear terms. This extraction is already carried out in the nonlinear estimate lemma, and in what follows, we simply use its consequences to identify the profile terms and control the remaining nonlinear contributions.

    We emphasize that this separation is essential: certain quadratic terms that at first appear to be part of the radiation actually contain nontrivial contributions from the profile. These contributions must be isolated and inserted into the modulation laws. The relation between these quadratic radiation terms and the parameters $c_j$ is recorded in Section~\ref{sec:nonlin pf}.

    \underline{\textbf{Step 1.}}
	For simplicity, we denote 
    \begin{align*}
    	\Mod_0\coloneqq \left(\tfrac{\lambda_s}{\lambda}+b_1,\gamma_s-\tfrac{\eta_1}{2}, \tfrac{x_s}{\lambda}-\nu_0 \right).
    \end{align*}
	For \eqref{eq:mod0 esti refine}, we recall \eqref{eq:wj equ in section mod esti} with $j=0$,
	\begin{equation}
		\begin{aligned}
			&(\partial_s-\tfrac{\lambda_s}{\lambda}\Lambda +\gamma_s i)w-\tfrac{x_s}{\lambda}w_{1}-iw_{2} +\tfrac{1}{2}\left(\tfrac{x_s}{\lambda}-\nu_0\right)w\calH(\ol{w}w)
			\\
			&
			+\tfrac{i}{2}((w_1-\tfrac{i\nu_0}{2}w)\mathcal{H}\overline{w}-w\mathcal{H}\overline{(w_1-\tfrac{i\nu_0}{2}w)})w=0.
		\end{aligned} \label{eq:w equ in mod esti}
	\end{equation}
	We use the decompositions, $w=Q+\wt{\eps}$ and $w_j=P_j+\eps_j$. We first observe the profile terms in \eqref{eq:w equ in mod esti}. The profiles in the first line of \eqref{eq:w equ in mod esti} are
	\begin{align}
		-\tfrac{\lambda_s}{\lambda}\Lambda Q+\gamma_s iQ-\tfrac{x_s}{\lambda}P_{1,e}-iP_{2,e}
		+\tfrac{1}{2}(\tfrac{x_s}{\lambda}-\nu_0)Q\calH(Q^2). \label{eq:w equ mod time profile}
	\end{align}
	The second line of \eqref{eq:w equ in mod esti} is decomposed into
	\begin{align}
		&\tfrac{i}{2}[(w_1-\tfrac{i\nu_0}{2}w)\mathcal{H}(\overline{w}\wt\eps)-w\mathcal{H}\{\overline{(w_1-\tfrac{i\nu_0}{2}w)}\wt\eps\}]  \label{eq:w equ in mod nonlin 1}
		\\
		&+
		\tfrac{i}{2}[(w_1-\tfrac{i\nu_0}{2}w)\mathcal{H}(\overline{w}Q)-w\mathcal{H}\{\overline{(w_1-\tfrac{i\nu_0}{2}w)}Q\}], \label{eq:w equ in mod nonlin 2}
	\end{align}
	In view of \eqref{eq:main nonlinear wt eps}, the profiles in \eqref{eq:w equ in mod nonlin 1} are
	\begin{align}
		-\tfrac{i}{2}(2\ol{c_{2}}+\tfrac{1}{2}(\nu_0^2-\mu_0^2)+i\mu_0\nu_0)Q
        =
        [-i\ol{c_{2}}-\tfrac{i}{4}(\nu_0^2-\mu_0^2)+\tfrac{1}{2}\mu_0\nu_0]Q. \label{eq:w equ in mod nonlin 1 prof goal}
	\end{align}
	Moreover, by \eqref{eq:main nonlinear Q}, the profiles in \eqref{eq:w equ in mod nonlin 2} are
	\begin{align}
		\tfrac{i}{2}\big[&(P_1-\tfrac{i\nu_0}{2}Q)\mathcal{H}(Q^2)-Q\mathcal{H}\{\overline{(P_1-\tfrac{i\nu_0}{2}Q)}Q\}
		-c_1T_{o}\calH(Q^2)+Q\calH(\ol{c_1 T_{o}} Q)\big]. \label{eq:w equ in mod nonlin 2 prof}
	\end{align}
    Using \eqref{eq:algebraic 1}, we derive
    \begin{align*}
        (P_1-\tfrac{i\nu_0}{2}Q)\mathcal{H}(Q^2)-Q\mathcal{H}\{\overline{(P_1-\tfrac{i\nu_0}{2}Q)}Q\}
        &=-\tfrac{1}{2}(ib_1+\eta_1)y^2Q^3+\tfrac{1}{2}(ib_1-\eta_1) Q^3\\
        &=2ib_1\Lambda Q-\eta_1 Q.
    \end{align*}
    Also, we can calculate the other terms:
    \begin{align*}
        -c_1T_{o}\calH(Q^2)+Q\calH(\ol{c_1 T_{o}} Q)=
        -\tfrac{i}{2}\nu_0c_1 y^2Q^3 +\tfrac{i}{2}\nu_0\ol{c_1} Q^3
        =i\mu_0\nu_0 \Lambda Q+\tfrac{1}{2}\nu_0^2 Q.
    \end{align*}
    Thus, we have
    \begin{align}
        \eqref{eq:w equ in mod nonlin 2 prof}
        =-(b_1+\tfrac{1}{2}\mu_0\nu_0)\Lambda Q -(\tfrac{1}{2}\eta_1-\tfrac{1}{4} \nu_0^2) iQ. \label{eq:w equ in mod nonlin 2 prof goal}
    \end{align}
	Therefore, from \eqref{eq:w equ mod time profile}, \eqref{eq:w equ in mod nonlin 1 prof goal}, and \eqref{eq:w equ in mod nonlin 2 prof goal}, the collection of profiles in \eqref{eq:w equ in mod esti} is written by
	\begin{align*}
		&\eqref{eq:w equ mod time profile}+\eqref{eq:w equ in mod nonlin 1 prof goal}+\eqref{eq:w equ in mod nonlin 2 prof}
		\\
        &=-(\tfrac{\lambda_s}{\lambda}+b_1+\tfrac{1}{2}\mu_0\nu_0)\Lambda Q
		+(\gamma_s-\tfrac{1}{2}\eta_1-\tfrac{1}{2}\nu_0^2+\tfrac{1}{4}\mu_0^2-2\Re(c_2)) iQ
        \\
        &\quad \,+(\tfrac{x_s}{\lambda}-\nu_0)(\tfrac{1}{2}yQ^3-c_1Q).
	\end{align*}
	Thus, we write \eqref{eq:w equ in mod esti} as
    \begin{equation}
        \begin{aligned}
            &-(\tfrac{\lambda_s}{\lambda}+b_1+\tfrac{1}{2}\mu_0\nu_0)\Lambda Q
    		+(\gamma_s-\tfrac{1}{2}\eta_1-\tfrac{1}{2}\nu_0^2+\tfrac{1}{4}\mu_0^2-2\Re(c_2)) iQ
    		\\
    		&+(\tfrac{x_s}{\lambda}-\nu_0)(\tfrac{1}{2}yQ^3-c_1Q)
            =-\calM_0,
        \end{aligned} \label{eq:mod0 equ rewrite}
    \end{equation}
	where $\calM_0=\calM_{0,s}+\calM_{0,1}+\calM_{0,2}+\calM_{0,3}$, and
	\begin{align*}
		\calM_{0,s}
		&=(\partial_s-\tfrac{\lambda_s}{\lambda}\Lambda +\gamma_s i)\wt\eps-\tfrac{x_s}{\lambda}\wt\eps_1-i\wt\eps_2,
		\\
		\calM_{0,1}&=\tfrac{1}{2}\left(\tfrac{x_s}{\lambda}-\nu_0\right)(w\calH(|w|^2)-Q\calH(Q^2)),
		\\
		\calM_{0,2}&=
        \tfrac{i}{2}\textnormal{NL}_{0}
        +\tfrac{i}{2}(2\ol{c_{2}}+\tfrac{1}{2}(\nu_0^2-\mu_0^2)+i\mu_0\nu_0)Q,
		\\
		\calM_{0,3}&=
		\tfrac{i}{2}\textnormal{NL}_{Q}.
	\end{align*}
	We recall $w_{j}=P_{j,e}+\wt\eps_{j}$.
    Since $(\wt\eps,\calZ_\ell)_r=0$ for $\ell=1,2,3$, we have
    \begin{align*}
        (\calM_{0,s},\calZ_\ell )_r
        \lesssim\,& (|\tfrac{\lambda_s}{\lambda} +b_1|+|\gamma_s-\tfrac{1}{2}\eta_1|
        +\beta_1)\|\wt\eps\|_{\dot\calH^1}
        +|\tfrac{x_s}{\lambda}|\|\wt\eps_1\|_{\dot\calH^1}+\|\wt\eps_2\|_{\dot\calH^1}
        \\
        \lesssim\,& |\Mod_0|\lambda +\lambda^{2+\frac12}.
    \end{align*}
	By \eqref{eq:wH w square}, \eqref{eq:main nonlinear wt eps}, and \eqref{eq:main nonlinear Q}, we deduce
    \begin{align*}
        (\calM_{0,1},\calZ_\ell )_r\lesssim |\Mod_0|\lambda, \quad  (\calM_{0,2},\calZ_\ell )_r\lesssim \lambda^{2+\frac12},\quad (\calM_{0,3},\calZ_\ell )_r\lesssim \lambda^{2+\frac12},
    \end{align*}
    %Since $(Q,f_o)_r=0$, \eqref{eq:main nonlinear wt eps} implies
    %\begin{align*}   (\calM_{0,2},\calZ_\ell )_r\lesssim \lambda^{2+\frac12}.\end{align*}
    %For $\calM_{0,3}$, using \eqref{eq:main nonlinear Q}, we arrive at
    %\begin{align*}    (\calM_{0,3},\calZ_\ell )_r\lesssim \lambda^{2+\frac12}.\end{align*}
    respectively. We remark \eqref{eq:kernel element1}, \eqref{eq:transversality}, and $Q_y=-yQ^3$. Therefore, taking the inner product of \eqref{eq:mod0 equ rewrite} with $\calZ_\ell$, $\ell=1,2,3$, we conclude
    \begin{align*}
        |\tfrac{\lambda_s}{\lambda}+b_1+\tfrac{1}{2}\mu_0\nu_0|
		\!+\!|\gamma_s-\tfrac{1}{2}\eta_1-\tfrac{1}{2}\nu_0^2\!+\!\tfrac{1}{4}\mu_0^2-2\Re(c_2)|
        \!+\!|\tfrac{x_s}{\lambda}-\nu_0|\lesssim |\Mod_0|\lambda\!+\!\lambda^{2+\frac12}\!.
    \end{align*}
    Since $\mu_0,\nu_0=O(\lambda)$ and $c_2=O(\lambda^2)$, we have
    %\begin{align*}
     %   |\mu_0\nu_0|\lesssim \lambda^2, \quad |-\tfrac{1}{2}\nu_0^2+\tfrac{1}{4}\mu_0^2-2\Re(c_2)|\lesssim \lambda^2,
    %\end{align*}
    %we obtain
    \begin{align*}
        |\tfrac{\lambda_s}{\lambda}+b_1|
		+|\gamma_s-\tfrac{1}{2}\eta_1|
        +|\tfrac{x_s}{\lambda}-\nu_0|
        \lesssim |\Mod_0|\lambda+\lambda^{2},
    \end{align*}
    which yields $|\Mod_0|\lesssim \lambda^2$, \eqref{eq:mod0 esti}.
    Using this, we derive
    \begin{align*}
        |\tfrac{\lambda_s}{\lambda}+b_1+\tfrac{1}{2}\mu_0\nu_0|
		+|\gamma_s-\tfrac{1}{2}\eta_1-\tfrac{1}{2}\nu_0^2+\tfrac{1}{4}\mu_0^2-2\Re(c_2)|
        +|\tfrac{x_s}{\lambda}-\nu_0|
        \lesssim \lambda^{2+\frac12}.
    \end{align*}
    Thus, we conclude the first line of \eqref{eq:mod0 esti refine}. Now, we further upgrade the estimate for $\frac{x_s}{\lambda}$. 
    Since $\calZ_3$ is odd and $(yQ^3,\calZ_3)_r\neq 0$, we have
    \begin{align*}
        (\text{LHS of } \eqref{eq:mod0 equ rewrite},\calZ_3)_r
        =\tfrac{1}{2}((\tfrac{x_s}{\lambda}-\nu_0)yQ^3,\calZ_3)_r.
    \end{align*}
    For the remainder term $\calM_0$, by \eqref{eq:wH w square} and \eqref{eq:main nonlinear wt eps}, we obtain
    \begin{align*}
        (\calM_{0,1},\calZ_3 )_r
        \lesssim
        |\tfrac{x_s}{\lambda}-\nu_0|\lambda\lesssim \lambda^3,\qquad 
        (\calM_{0,2},\calZ_3 )_r
        \lesssim \beta_1\lambda+\lambda^3.
    \end{align*}
    For $\calM_{0,s}+\calM_{0,3}$, we can check that
    \begin{align*}
        (\calM_{0,s}+\calM_{0,3},\calZ_3 )_r
        \lesssim \beta_1\lambda+\lambda^3+|(-i\wt\eps_2 + \tfrac{i}{2}\textnormal{NL}_{Q},\calZ_3 )_r|.
    \end{align*}
    From the decomposition $\wt\eps_2=\beta_2yQ+\eps_2$ with \eqref{eq:coercivity H1 epsj} and \eqref{eq:main nonlinear Q}, we have
    \begin{align*}
        |(-i\wt\eps_2 + \tfrac{i}{2}\textnormal{NL}_{Q},\calZ_3 )_r|
        \leq
        |(i\eps_{2},\calZ_3 )_r|
        +
        |(\tfrac{i}{2}[\textnormal{NL}_{Q}-2\beta_2yQ],\calZ_3 )_r|
        \lesssim \beta_1\lambda+\lambda^3.
    \end{align*}
    Hence, we finish the proof of the second line of \eqref{eq:mod0 esti refine}.

    \underline{\textbf{Step 2.}} Next, we show \eqref{eq:modj esti}.
    We again recall \eqref{eq:wj equ in section mod esti},
    \begin{equation}
    	\begin{aligned}
    		&(\partial_s-\tfrac{\lambda_s}{\lambda}\Lambda_{-j} +\gamma_s i)w_j-\tfrac{x_s}{\lambda}w_{j+1}-iw_{j+2}
            +\tfrac{1}{2}\left(\tfrac{x_s}{\lambda}-\nu_0\right)w\calH(\ol{w}w_j)
    		\\
    		&+\tfrac{i}{2}[(w_1-\tfrac{i\nu_0}{2}w)\mathcal{H}(\overline{w}w_{j})-w\mathcal{H}\{\overline{(w_1-\tfrac{i\nu_0}{2}w)}w_{j}\}]
            =0.
    	\end{aligned} \label{eq:wj equ in mod proof}
    \end{equation}
    We take the inner product of \eqref{eq:wj equ in mod proof} and $\calZ_h$. Since $\calZ_h$ is odd, it suffices to consider the odd part of \eqref{eq:wj equ in mod proof}. The odd profiles in the first line of \eqref{eq:wj equ in mod proof} are
    \begin{align}
        (\partial_s-\tfrac{\lambda_s}{\lambda}\Lambda_{-j} +\gamma_s i)\bm{\beta}_j yQ-\tfrac{x_s}{\lambda}\bm{\beta}_{j+1} yQ-i\bm{\beta}_{j+2} yQ
        +\tfrac{1}{2}(\tfrac{x_s}{\lambda}-\nu_0)c_jQ\calH(Q^2). \label{eq:mod j pf profile s}
    \end{align}
    For the second line of \eqref{eq:wj equ in mod proof}, since $w_j=c_jQ+\wt\eps_j$, we decompose this into
    \begin{align}
        &\tfrac{i}{2}\big[(w_1-\tfrac{i\nu_0}{2}w)\mathcal{H}(\overline{w}\wt\eps_{j})-w\mathcal{H}\{\overline{(w_1-\tfrac{i\nu_0}{2}w)}\wt\eps_{j}\}\big] \label{eq:mod j proof nonlin 1}
        \\
        & +
        \tfrac{i}{2}c_j\big[(w_1-\tfrac{i\nu_0}{2}w)\mathcal{H}(\overline{w}Q)-w\mathcal{H}\{\overline{(w_1-\tfrac{i\nu_0}{2}w)}Q\}\big]. \label{eq:mod j proof nonlin 2}
    \end{align}
    In view of \eqref{eq:main nonlinear eps j odd} and \eqref{eq:main nonlinear Q}, the odd profiles corresponding to \eqref{eq:mod j proof nonlin 1} and \eqref{eq:mod j proof nonlin 2} are, respectively,
    \begin{align}
        &\begin{aligned}
            \eqref{eq:mod j proof nonlin 1}:\, &\tfrac{i}{2}[P_{1,o}\mathcal{H}(QP_{j,o})
            -yQ\mathcal{H}\{\tfrac{y}{1+y^2}\overline{P_{1,o}}P_{j,o}\}
            -Q\mathcal{H}\{\tfrac{1}{1+y^2}\overline{P_{1,o}}P_{j,o}\}
            \\
            &\quad +\tfrac{i}{2}\mu_0\nu_0\bm{\beta}_j yQ^3],
        \end{aligned}  \label{eq:mod j pf profile 1}
        \\
        &\eqref{eq:mod j proof nonlin 2}:\, \tfrac{i}{2}c_j[
        (P_{1,e}-\tfrac{i\nu_0}{2}Q)\mathcal{H}(Q^2)-Q\mathcal{H}\{\overline{(P_{1,e}-\tfrac{i\nu_0}{2}Q)}Q\}
        +2\bm{\beta}_2 yQ  ].  \label{eq:mod j pf profile 2}
    \end{align}
    We note that $\Lambda (yQ)=\frac{1}{2}yQ+\frac{1}{2}yQ^3$ and 
    \begin{align*}
        P_{1,o}\mathcal{H}(QP_{j,o})
        -yQ\mathcal{H}\{\tfrac{y}{1+y^2}\overline{P_{1,o}}P_{j,o}\}-Q\mathcal{H}\{\tfrac{1}{1+y^2}\overline{P_{1,o}}P_{j,o}\}
        &=ib_1\bm{\beta}_jyQ^3,
        \\
        (P_{1,e}-\tfrac{i\nu_0}{2}Q)\mathcal{H}(Q^2)-Q\mathcal{H}\{\overline{(P_{1,e}-\tfrac{i\nu_0}{2}Q)}Q\}&=0.
    \end{align*}
    Thus, the collection of profiles is
    \begin{align*}
        &\eqref{eq:mod j pf profile s}+\eqref{eq:mod j pf profile 1}+\eqref{eq:mod j pf profile 2}
        \\
        &=[(\bm{\beta}_j)_s - \tfrac{\lambda_s}{\lambda} (j+\tfrac{1}{2})\bm{\beta}_j + i\gamma_s\bm{\beta}_j 
        -\tfrac{x_s}{\lambda}\bm{\beta}_{j+1} -i\bm{\beta}_{j+2}
        +ic_j\bm{\beta}_2 ] yQ
        \\
        &\quad 
        -[\tfrac{\lambda_s}{\lambda}+b_1 +\tfrac{1}{2}\mu_0\nu_0]\tfrac{1}{2}\bm{\beta}_j yQ^3 +\tfrac{1}{2}(\tfrac{x_s}{\lambda}-\nu_0)c_j yQ^3.
    \end{align*}
    Therefore, the odd part of \eqref{eq:wj equ in mod proof} can be written as
    \begin{equation}
        \begin{aligned}
            &[(\bm{\beta}_j)_s - \tfrac{\lambda_s}{\lambda} (j+\tfrac{1}{2})\bm{\beta}_j + i\gamma_s\bm{\beta}_j 
            -\nu_0\bm{\beta}_{j+1} -i\bm{\beta}_{j+2}
            +ic_j\bm{\beta}_2 ] yQ
            \\
            &=[\tfrac{\lambda_s}{\lambda}+b_1 +\tfrac{1}{2}\mu_0\nu_0]\tfrac{1}{2}\bm{\beta}_j yQ^3 
            -(\tfrac{x_s}{\lambda}-\nu_0)[\tfrac{1}{2}c_j yQ^3 -\bm{\beta}_{j+1} yQ]
            -\calM_j^{\text{odd}},
        \end{aligned} \label{eq:modj equ rewrite}
    \end{equation}
    where $\calM_j^{\text{odd}}=
        \calM_{j,s}^{\text{odd}}+\calM_{j,1}^{\text{odd}}
        +\calM_{j,2}^{\text{odd}}+\calM_{j,3}^{\text{odd}}$, and
    \begin{equation}
        \begin{aligned}
            \calM_{j,s}^{\text{odd}}
    		&=(\partial_s-\tfrac{\lambda_s}{\lambda}\Lambda +\gamma_s i)\eps_{j,o}-\tfrac{x_s}{\lambda}\eps_{j+1,o}-i\eps_{j+2,o},
    		\\
    		\calM_{j,1}^{\text{odd}}&=\tfrac{1}{2}\left(\tfrac{x_s}{\lambda}-\nu_0\right)[w\calH(\ol{w}w_j)-c_jQ\calH(Q^2)]_{o},
    		\\
    		\calM_{j,2}^{\text{odd}}&=
            \tfrac{i}{2}[\textnormal{NL}_{j}-\tfrac{i}{2}\mu_0\nu_0\bm{\beta}_j yQ^3]_{o},
    		\\
    		\calM_{j,3}^{\text{odd}}&=
    		\tfrac{i}{2}c_j[\textnormal{NL}_{Q}-2\bm{\beta}_2 yQ]_{o}.
        \end{aligned} \label{eq:calM odd def}
    \end{equation}
    Now, we estimate the right-hand side of \eqref{eq:modj equ rewrite}.
    The first two terms can be regarded as a perturbative term; by \eqref{eq:mod0 esti refine},
    \begin{align*}
        ((\tfrac{x_s}{\lambda}-\nu_0)[\tfrac{1}{2}c_j yQ^3 -\bm{\beta}_{j+1} yQ], \calZ_h)_c
        &\lesssim \beta_1\lambda^{j+1}+\lambda^{j+3},
        \\
        ([\tfrac{\lambda_s}{\lambda}+b_1 +\tfrac{1}{2}\mu_0\nu_0]\tfrac{1}{2}\bm{\beta}_j yQ^3), \calZ_h)_c
        &\lesssim \lambda^{j+3}.
    \end{align*}
    For $\calM_{j,s}^{\text{odd}}$, since $(\eps_{j,o}, \calZ_h)_c$, we derive
    \begin{align*}
        (\calM_{j,s}^{\text{odd}}, \calZ_h )_c=
        -(\tfrac{\lambda_s}{\lambda}\Lambda \eps_{j,o},  \calZ_h)_c
        \lesssim \beta_1 \lambda^{j+1}.
    \end{align*}
    To estimate $\calM_{j,1}^{\text{odd}}$, we first remark the following decomposition,
    \begin{align*}
        w\calH(\ol{w}w_j)
        =w\calH(\ol{w}\wt\eps_j)+c_j w\calH(\ol{w}Q).
    \end{align*}
    Thus, using \eqref{eq:wH w square} and \eqref{eq:wH w wj} with \eqref{eq:mod0 esti refine}, we obtain
    \begin{align*}
        (\calM_{j,1}^{\text{odd}}, \calZ_h )_c
        \lesssim (\beta_1 \lambda+\lambda^{3})\lambda^{j+1}
        \lesssim \lambda^{j+3}.
    \end{align*}
    Similarly, thanks to \eqref{eq:main nonlinear eps j odd} and \eqref{eq:main nonlinear Q}, we derive
    \begin{align*}
        |(\calM_{j,2}^{\text{odd}}, \calZ_h )_c|+|(\calM_{j,3}^{\text{odd}}, \calZ_h )_c|
        \lesssim \beta_1 \lambda^{j+1}+\lambda^{j+3}.
    \end{align*}
    Therefore, we conclude \eqref{eq:modj esti}.
\end{proof}	

\subsection{Refined modulation estimates}\label{sec:refined mod esti}

In this subsection, we establish refined modulation estimates. The estimate \eqref{eq:modj esti} holds for $1\le j\le 2L-2$, but it does not provide a control on $(\bm{\beta}_{2L-1})_s$. To handle this final parameter, we replace $\bm\beta_{2L-1}$ with a refined version. Specifically, we introduce a new parameter $\td{\bm{\beta}}_{2L-1}$, chosen to approximate $\bm{\beta}_{2L-1}$ closely enough, and we will derive its modulation law, namely the equation governing $(\td{\bm{\beta}}_{2L-1})_s$. Our construction follows the approach used in earlier works on \eqref{CMdnls} \cite{KimKimKwon2024arxiv,JeongKim2024arXiv}, originally introduced in \cite{JendrejLawrieRodriguez2022ASENS,KimKwonOh2024arXiv}.

The refinement is implemented by testing the $w_{2L-1}$-equation against the almost-null space elements $B_Q^*\chi_R\approx \frac{1}{\sqrt{2}}yQ$, with $R$ chosen near the self-similar scale. As pointed out in \cite[Section~3.1]{Kim2025JEMS}, the condition $b_1>0$ is not available in the blow-up classification problem, and therefore, unlike \cite{KimKimKwon2024arxiv,JeongKim2024arXiv}, one cannot take $R=b_1^{-1/2}$. Instead, using the relation $(\lambda^2)_t\sim |b_1|$, we are naturally led to take $R=\lambda^{-1}\sqrt{T-t}$.

\begin{lem}[Refined modulation estimates]\label{lem:modj refine}
	Define the refined modulation parameters as
    \begin{align}
        \td{\bm{\beta}}_{2L-1}\coloneqq \tfrac{1}{AR}(w_{2L-1},B_Q^*\chi_R)_c, \label{eq:def refine mod}
    \end{align}
	where $R=\lambda^{-1}\sqrt{T-t}$, and $A$ is given by
	\begin{align*}
		A\coloneqq \tfrac{1}{R}(\tfrac{\sqrt{2}}{2},\chi_R)_r=(\tfrac{\sqrt{2}}{2},\chi)_r.
	\end{align*}
	Then, we have a proximity,
	\begin{align}
        |\bm{\beta}_{2L-1}-\td{\bm{\beta}}_{2L-1}|
        \lesssim \lambda^{2L-\frac12}(T-t)^{\frac34}, \label{eq:refined mod gap}
	\end{align}
    and a modulation estimate, 
    \begin{equation}
        |(\partial_s-\tfrac{\lambda_s}{\lambda}(2L-\tfrac{1}{2})+\gamma_si)\td{\bm{\beta}}_{2L-1}|\lesssim \lambda^{2L+\frac32}[(T-t)^{-\frac14}+\lambda(T-t)^{-\frac32}]. \label{eq:refined mod esti}
    \end{equation}
\end{lem}
We note that the extra terms in \eqref{eq:refined mod esti} compared to \eqref{eq:modj esti} are negligible since %absorbed into the error since
\begin{align*}
   |\nu_0\bm{\beta}_{2L}|+|c_{2L-1}\bm{\beta}_{2}|\lesssim \lambda^{2L+\frac32}\ll \lambda^{2L+\frac32}(T-t)^{-\frac14}.
\end{align*}
%which shows that this term is negligible.

\begin{proof}[Proof of Lemma~\ref{lem:modj refine} assuming Lemma~\ref{lem:nonlinear esti I}]
	We first prove \eqref{eq:refined mod gap}.
	Using the decomposition $w_{2L-1}=P_{2L-1}+\eps_{2L-1}$ with $B_QP_{2L-1}=\tfrac{\sqrt{2}}{2}\bm{\beta}_{2L-1}$,
	we have
	\begin{align*}
		AR \td{\bm{\beta}}_{2L-1}
		=\bm{\beta}_{2L-1}(\tfrac{\sqrt{2}}{2},\chi_R)_c
		+(\eps_{2L-1,o},B_Q^*\chi_R)_c.
	\end{align*}
    We observe that, from \eqref{eq:CommuteHilbert},
    \begin{align*}
        (1+y^2)B_Q^*\chi_R=
        y\langle y \rangle\chi_R+\langle y \rangle\calH\chi_R.
    \end{align*}
	Thanks to this, $(\tfrac{\sqrt{2}}{2},\chi_R)_c=AR$, and $\lambda\lesssim T-t$, we have
    \begin{equation}
        \begin{aligned}
            |\bm{\beta}_{2L-1}-\td{\bm{\beta}}_{2L-1}|
            &\lesssim \tfrac{1}{R} \big[\|\eps_{2L-1,o}\|_{\dot\calH^2}\|y^2\chi_R\|_{L^2}
            +\|\calH(Q\eps_{2L-1,o})\|_{L^{\infty}}\|\chi_R\|_{L^1} \big]
            \\
            &\lesssim \lambda^{2L+1}R^{\frac{3}{2}}+\lambda^{2L+\frac12}
            \lesssim \lambda^{2L-\frac12}(T-t)^{\frac34}.
        \end{aligned} \label{eq:eps BQ techniq}
    \end{equation}
	Now, we show \eqref{eq:refined mod esti}. 
	From the definition of $\td b_{L}$, we compute
	\begin{align*}
		&AR(\partial_s-\tfrac{\lambda_s}{\lambda}(2L-\tfrac{1}{2})+\gamma_si)\td{\bm{\beta}}_{2L-1} 
        \\
        &=
		((\partial_s-\tfrac{\lambda_s}{\lambda}(2L-\tfrac{1}{2})+\gamma_si)w_{2L-1},B_Q^*\chi_R)_c
        -\tfrac{R_s}{R}[(w_{2L-1}, B_Q^*(y\partial_y\chi_R+\chi_R))_c]. 
	\end{align*}
    Similarly to \eqref{eq:modj equ rewrite}, from \eqref{eq:wj equ in section mod esti}, we have
    \begin{align*}
        &[(\partial_s-\tfrac{\lambda_s}{\lambda}(2L-\tfrac{1}{2})+\gamma_si)w_{2L-1}]_o
        \\
        &=[\tfrac{\lambda_s}{\lambda}y\partial_yw_{2L-1}+\tfrac{x_s}{\lambda}w_{2L}+iw_{2L+1}]_o-\tfrac{1}{2}(\tfrac{x_s}{\lambda}-\nu_0)c_{2L-1}yQ^3
        \\
        &\quad +\tfrac{1}{4}\mu_0\nu_0\bm{\beta}_{2L-1}yQ^3-i\bm{\beta}_{2}c_{2L-1}yQ
        -\td{\calM}_{2L-1}^{\text{odd}},
    \end{align*}
    where $  \td{\calM}_{2L-1}^{\text{odd}}=\calM_{2L-1,1}^{\text{odd}}+\calM_{2L-1,2}^{\text{odd}}+\calM_{2L-1,3}^{\text{odd}}$, and the terms $\calM_{2L-1,j}^{\text{odd}}$ are defined in \eqref{eq:calM odd def}.
    Hence, we arrive at
    \begin{align}
        &AR(\partial_s-\tfrac{\lambda_s}{\lambda}(2L-\tfrac{1}{2})+\gamma_si)\td{\bm{\beta}}_{2L-1} \nonumber
        \\
        &=\tfrac{\lambda_s}{\lambda}(y\partial_yw_{2L-1},B_Q^*\chi_R)_c
        -\tfrac{R_s}{R}(w_{2L-1}, B_Q^*(y\partial_y\chi_R+\chi_R))_c \label{eq:refined mod equ b-1}
        \\
        &\quad +(\tfrac{x_s}{\lambda}w_{2L}+iw_{2L+1}, B_Q^*\chi_R)_c
        - (\tfrac{1}{2}(\tfrac{x_s}{\lambda}-\nu_0)c_{2L-1}yQ^3, B_Q^*\chi_R)_c \label{eq:refined mod equ b-2}
        \\
        &\quad +(\tfrac{1}{4}\mu_0\nu_0\bm{\beta}_{2L-1}yQ^3-i\bm{\beta}_{2}c_{2L-1}yQ, B_Q^*\chi_R)_c- (\td{\calM}_{2L-1}^{\text{odd}}, B_Q^*\chi_R)_c. \label{eq:refined mod equ b-3}
    \end{align}
    A direct computation yields
    \begin{align}
        \eqref{eq:refined mod equ b-1}
        =\tfrac{(\lambda R)_s}{\lambda R}(y\partial_yw_{2L-1},B_Q^*\chi_R)_c
        -\tfrac{R_s}{R}(w_{2L-1}, [B_Q^*,\partial_yy]\chi_R)_c. \label{eq:refined mod equ b-1-1}
    \end{align}
    Here, from $R= \lambda^{-1}\sqrt{T-t}$, we obtain
	\begin{align*}
		\tfrac{(\lambda R)_s}{\lambda R}=\tfrac{\lambda^2}{2(T-t)},\quad
        \left|\tfrac{R_s}{R}\right|\lesssim |b_1|+\tfrac{\lambda^2}{T-t}.
	\end{align*}
    From $w_{2L-1}=P_{2L-1}+\eps_{2L-1}$, we estimate the first term of \eqref{eq:refined mod equ b-1-1} as
	\begin{align*}
		&\tfrac{(\lambda R)_s}{\lambda R}(y\partial_yw_{2L-1},B_Q^*\chi_R)_c
        \lesssim 
        \tfrac{\lambda^2}{T-t}[|\bm{\beta}_{2L-1}(y\partial_y(yQ),B_Q^*\chi_R)_c|
        \!+\!|(y\partial_y\eps_{2L-1},B_Q^*\chi_R)_c|].
	\end{align*}
	Using $\partial_y(yQ)=\tfrac{1}{2}Q^3$, \eqref{eq:algebraic 1}, and \eqref{eq:algebraic 2}, we have
	\begin{align}
		(y\partial_y(yQ),B_Q^*\chi_R)_c
		=(\tfrac{1}{2}yQ^3,B_Q^*\chi_R)_c=\tfrac{\sqrt{2}}{2}(\langle y\rangle^{-2},\chi_R)_c\lesssim 1. \label{eq:BQ yQ3}
	\end{align}
	As in \eqref{eq:eps BQ techniq}, and using $\lambda\lesssim T-t$, we also deduce
	\begin{align*}
		|(y\partial_y\eps_{2L-1},B_Q^*\chi_R)_c|
		\lesssim\lambda^{2L+1}R^{\frac{5}{2}}+\lambda^{2L+\frac12}R
        \lesssim \lambda^{2L-\frac32}(T-t)^{\frac54}.
	\end{align*}
    We hence arrive at
    \begin{align}
        |\tfrac{(\lambda R)_s}{\lambda R}(y\partial_yw_{2L-1},B_Q^*\chi_R)_c|
        \lesssim \lambda^{2L+\frac12}[\tfrac{\lambda}{T-t}+(T-t)^{\frac14}]. \label{eq:refined mod equ b goal 1}
    \end{align}
    To estimate the second term of \eqref{eq:refined mod equ b-1-1}, by \eqref{eq:CommuteHilbertDerivative}, we note that
    \begin{align*}
    	[B_Q^*,\partial_yy]=-y\langle y\rangle^{-3}+y^2\langle y\rangle^{-3}\calH.
    \end{align*}
    Thus, using $|\calH(y^3Q^4)|\lesssim Q^2$, we obtain
    \begin{equation}
        \begin{aligned}
            \tfrac{R_s}{R}(w_{2L-1}, [B_Q^*,\partial_yy]\chi_R)_c
            &\lesssim (|b_1|+\tfrac{\lambda^2}{T-t})(|\bm{\beta}_{2L-1}|+\lambda^{2L}R^{\frac{1}{2}})\\
            &\lesssim \lambda^{2L-\frac12}(\lambda^{\frac32}+\tfrac{\lambda^2}{T-t}).
        \end{aligned} \label{eq:refined mod equ b goal 2}
    \end{equation}
    By \eqref{eq:refined mod equ b goal 1} and \eqref{eq:refined mod equ b goal 2}, we have
    \begin{align}
        |\eqref{eq:refined mod equ b-1}|
        \lesssim 
        \lambda^{2L+\frac12}[\tfrac{\lambda}{T-t}+(T-t)^{\frac14}+\lambda^{\frac12}]. \label{eq:refined mod equ b-1 goal}
    \end{align}
    
    We now control \eqref{eq:refined mod equ b-2}. The first term of \eqref{eq:refined mod equ b-2} becomes
    \begin{align*}
        |(\tfrac{x_s}{\lambda}w_{2L}+iw_{2L+1}, B_Q^*\chi_R)_c|
        \lesssim \lambda^{2L+1} R^{\frac12}.
    \end{align*}
	Thanks to \eqref{eq:BQ yQ3} and \eqref{eq:mod0 esti refine}, we estimate the second part as
    \begin{align*}
        |(\tfrac{1}{2}(\tfrac{x_s}{\lambda}-\nu_0)c_{2L-1}yQ^3, B_Q^*\chi_R)_c|
        \lesssim \lambda^{2L+\frac32}.
    \end{align*}
    Thus, we have
    \begin{align}
        |\eqref{eq:refined mod equ b-2}|
        \lesssim \lambda^{2L+\frac12}(T-t)^{\frac14}+ \lambda^{2L+\frac32}. \label{eq:refined mod equ b-2 goal}
    \end{align}

    For \eqref{eq:refined mod equ b-3}, again using \eqref{eq:BQ yQ3}, the first term becomes
    \begin{align*}
        |(\tfrac{1}{4}\mu_0\nu_0\bm{\beta}_{2L-1}yQ^3-i\bm{\beta}_{2}c_{2L-1}yQ, B_Q^*\chi_R)_c|
        \lesssim \lambda^{2L+1}+ \lambda^{2L+\frac32}R^{\frac12} \lesssim \lambda^{2L+1}.
    \end{align*}
    As in the proof of Proposition~\ref{prop:modj}, by using \eqref{eq:wH w square}, \eqref{eq:wH w wj}, \eqref{eq:main nonlinear Q}, and \eqref{eq:main nonlinear eps j odd}, the nonlinear part of \eqref{eq:refined mod equ b-3} becomes
    \begin{align*}
        |(\td{\calM}_{2L-1}^{\text{odd}}, B_Q^*\chi_R)_c|
        \lesssim \lambda^{2L+\frac32}R^{\frac12}\lesssim \lambda^{2L+1}.
    \end{align*}
	Therefore, we have $\eqref{eq:refined mod equ b-3} =O(\lambda^{2L+1})$. From this, \eqref{eq:refined mod equ b-1 goal}, and \eqref{eq:refined mod equ b-2 goal} 
	%\begin{align}
		%|\eqref{eq:refined mod equ b-3}|	\lesssim\lambda^{2L+1}. \label{eq:refined mod equ b-3 goal}
	%\end{align}
	%By \eqref{eq:refined mod equ b-1 goal}, \eqref{eq:refined mod equ b-2 goal}, and \eqref{eq:refined mod equ b-3 goal} 
    with $R= \lambda^{-1}\sqrt{T-t}$ and $\lambda\lesssim T-t$, we conclude
	\begin{align*}
		&|(\partial_s-\tfrac{\lambda_s}{\lambda}(2L-\tfrac{1}{2})+\gamma_si)\td{\bm{\beta}}_{2L-1}|
		\\
        &\lesssim 
        \tfrac{1}{R}\big[\lambda^{2L+\frac12}[\tfrac{\lambda}{T-t}+(T-t)^{\frac14}+\lambda^{\frac12}]
        +\lambda^{2L+\frac12}(T-t)^{\frac14}+ \lambda^{2L+\frac32}
        \big]
        \\
        &\lesssim \lambda^{2L+\frac32}[(T-t)^{-\frac14}+\lambda(T-t)^{-\frac32}]. \qedhere
	\end{align*}
\end{proof}

\section{Radiation estimates}\label{sec:rad esti}
In this section we establish quantitative analysis for the radiation parameters $\{c_j\}$ defined by \eqref{eq:def cj} in Section~\ref{sec:wj decomposition}.
As discussed therein, these parameters capture the slow odd radiation modes of the $w$-variable that are not removed by the orthogonality conditions.
While the modulation parameters $(\lambda,\gamma,x)$ and $\{\bm{\beta}_j\}$ describe the core dynamics of the blow-up profile, the radiation parameters $\{c_j\}$ account for the residual radiation and influence on the modulation system.

Recall from the modulation estimate \eqref{eq:modj esti} in Section~\ref{sec:mod esti} that each $(\bm{\beta}_j)_s$ depends on the neighboring parameters $\bm{\beta}_{j+1}$ and $\bm{\beta}_{j+2}$, as well as on the radiation parameter $c_j$.
To close that system, it is therefore necessary to control the evolution of $\{c_j\}$ themselves.
We show that ${c_j}$ satisfy a coupled system of the form
\[
    (c_j)_s - \tfrac{\lambda_s}{\lambda}j c_j - i c_{j+2} \approx a_j,
\]
which is consistent with the size information $|c_j|\lesssim \lambda^j$ from \eqref{eq:c j estimate}. The source term $a_j$ is expressed in terms of the radiation parameters $\{c_j\}$ and a new family of quantities $\{\frakc_{i,j}\}$, which arise from quadratic interactions among the radiation modes.
These $\mathfrak{c}_{i,j}$ are introduced to record certain quadratic radiation interactions that appear in the derivation of $(c_j)_s$ and are comparable in size to the leading radiation parameters. 
Since they cannot be absorbed into the error, isolating them allows the system for $\{c_j\}$ to close at the desired order. 

We begin in Section~\ref{subsec:cj} by deriving the evolution equations for $\{c_j\}$. In this part, the quadratic radiation interactions $\{\mathfrak{c}_{i,j}\}$ naturally appear in the computation of $(c_j)_s$. Their detailed construction and quantitative properties will be established in Section~\ref{subsec:frakc}.

\subsection{Evolution of \texorpdfstring{$c_j$}{cj}}\label{subsec:cj}
In this subsection, we derive the evolution of the radiation parameters $c_j$.
During the computation of $(c_j)_s$, one encounters quadratic radiation contributions such as 
\begin{equation}
    \frac{1}{2}\int yQ^2 \ol{\wt{\eps}_1} \wt{\eps}_j, \label{eq:quadratic raw form}
\end{equation}
whose size is comparable to the leading radiation parameters. 
These terms cannot be discarded, and in order to obtain a closed ODE system for $\{c_j\}$, one must also track the evolution of the corresponding quadratic interactions.

A direct attempt to use the quantity \eqref{eq:quadratic raw form} is problematic: differentiating it in time $s$ produces a complicated expression containing not only $\{c_j\}$ and quadratic interactions of the same type, but also undesirable profile contributions coming from $\overline{w_i} w_j$ that do not fit into the structure needed to close the modulation system.
If kept in this raw form, these additional terms prevent us from closing the inductive system for $(c_j)_s$.

To overcome this, we introduce the modified quadratic quantities
\begin{align} 
	\frakc_{i,j}\coloneqq \int\varphi\ol{w_i} w_j,\qquad \varphi\coloneqq \frac{1}{2}yQ^2-3yQ^4+2yQ^6. \label{eq: frakc ij def}
\end{align}
The weight $\varphi$ is designed so that the profile contributions present in $\overline{w_i}w_j$ cancel out. 
As a consequence, the time derivative $(\frakc_{i,j})_s$ can be expressed solely in terms of the radiation parameters $\{c_j\}$ and the quantities $\{\frakc_{i,j}\}$. Although $\mathfrak{c}_{i,j}$ is a modified version of the raw quadratic term, it still approximates the non-negligible quadratic interactions 
\begin{equation*}
    \mathfrak{c}_{i,j}
    \approx \frac{1}{2}\int yQ^2 \ol{\wt{\eps}_i} \wt{\eps}_j, \quad
    \text{up to} \quad O(\lambda^{i+j+1+\frac12}).
\end{equation*}
A detailed explanation of $\varphi$, further properties of $\frakc_{i,j}$, and the evolution equations of $\frakc_{i,j}$ will be given in Section~\ref{subsec:frakc}.

With these preparations, we now state the main estimates for the radiation parameters.

\begin{prop}[Radiation estimates I]\label{prop:cj mod esti} 
    Let $L\geq 1$ and $v_0\in H^{2L+1}$. For $1\leq j \leq 2L-1$, we have
    \begin{align}
        |(c_j)_s-\tfrac{\lambda_s}{\lambda} j c_j|
        \lesssim \lambda^{j+2}. \label{eq:cj mod aux}
    \end{align}
    If $L\geq 2$ and $1\leq j \leq 2L-2$, we have
	\begin{align}
		|(c_j)_s-\tfrac{\lambda_s}{\lambda} j c_j-ic_{j+2}-a_{j}|
        \lesssim \lambda^{j+2+\frac{1}{2}}, \label{eq:cj mod}
	\end{align}
    where
    \begin{align}
        a_{j}=\nu_0 c_{j+1} 
        -i(\tfrac{3}{4} \nu_0^2-\tfrac{1}{4}\mu_0^2+2\Re(c_2))c_j -\tfrac{i}{2\pi}(\frakc_{1,j}- \ol{c_1}\frakc_{0,j}). \label{eq:aj def}
    \end{align}
\end{prop}
To prove Proposition~\ref{prop:cj mod esti}, we need a refined control of the nonlinear term $\textnormal{NL}_j$, defined in \eqref{eq:def NL}, appearing in the $w_j$ equation, in addition to $\textnormal{NL}_Q$ and $\textnormal{NL}_0$.
\begin{lem}[Nonlinear estimates II]\label{lem:nonlinear esti II}
    Let $f\in C^\infty$ be such that $Q^{-1}f\in L^2$. Then, for $1\leq j\leq 2L-1$, we have
    \begin{align}
		|(\textnormal{NL}_j-\tfrac{1}{\pi}(\frakc_{1,j}- \ol{c_1}\frakc_{0,j})Q,f_e)_c|\lesssim \lambda^{j+2+\frac12}. \label{eq:main nonlinear eps j even}
	\end{align}
\end{lem}
Indeed, using \eqref{eq:nonli rad cj relation}, the term $\frakc_{0,j}$ can be rewritten by $\pi(-2c_{j+1}+i\nu_0c_j)$; see Remark~\ref{rem:convert rad to mod} for further details. Since the proofs of these nonlinear estimates involve only technical computations, we postpone them to Section~\ref{sec:nonlin pf}.

\begin{proof}[Proof of Proposition~\ref{prop:cj mod esti} assuming Lemmas~\ref{lem:nonlinear esti I} and \ref{lem:nonlinear esti II}] From \eqref{eq:wj equ in section mod esti}, we have \begin{subequations}\label{eq:wj equ rewrite-1}
		\begin{align}
			&(\partial_s-\tfrac{\lambda_s}{\lambda}\Lambda_{-j} +\gamma_s i)w_j-\tfrac{x_s}{\lambda}w_{j+1}-iw_{j+2}+\tfrac{1}{2}\left(\tfrac{x_s}{\lambda}-\nu_0\right)w\calH(\ol{w}w_j)
            \label{eq:wj equ rewrite-1-1}\\
			&+\tfrac{i}{2}\big[(w_1-\tfrac{i\nu_0}{2}w)\mathcal{H}(\overline{w}w_{j})-w\mathcal{H}\{\overline{(w_1-\tfrac{i\nu_0}{2}w)}w_{j}\}\big]
            =0.\label{eq:wj equ rewrite-1-2}
		\end{align}
	\end{subequations}
	We only consider the even part of \eqref{eq:wj equ rewrite-1} to track $c_j$ due to 
	\begin{align*}
		w_j=P_j+\eps_j, \quad \text{where} \quad P_{j,e}=c_jQ, \quad P_{j,o}=\bm{\beta}_jyQ.
	\end{align*}
	Thus, we find out the profile part as follows: The even profiles of $\eqref{eq:wj equ rewrite-1-1}$ are
	\begin{align*}
		&(c_j)_sQ-\tfrac{\lambda_s}{\lambda}c_j\Lambda_{-j}Q+\gamma_s c_j iQ
		-\tfrac{x_s}{\lambda}c_{j+1}Q
		\\
		&=[(c_j)_s-\tfrac{\lambda_s}{\lambda} j c_j-\tfrac{x_s}{\lambda}c_{j+1}]Q-\tfrac{\lambda_s}{\lambda}c_j\Lambda Q+\gamma_s c_j iQ.
	\end{align*}
    For \eqref{eq:wj equ rewrite-1-2}, we extract certain profile terms from the nonlinear parts in order to apply \eqref{eq:main nonlinear Q} and \eqref{eq:main nonlinear eps j even}. This shows that the even profiles appearing in \eqref{eq:wj equ rewrite-1-2} are
    \begin{align*}
        &\tfrac{i}{2}c_j\big[P_{1,o}\mathcal{H}(Q^2)
        -Q\mathcal{H}\{\overline{P_{1,o}}Q\}
        -c_1T_{o}\calH(Q^2)+Q\calH(\ol{c_1 T_{o}} Q)\big]
        +\tfrac{i}{2\pi}[\frakc_{1,j}- \ol{c_1}\frakc_{0,j}]Q
        \\
        &=\tfrac{i}{2}c_j[2ib_1 \Lambda Q-\eta_1 Q+\tfrac{1}{2}\nu_0^2 Q+i\mu_0\nu_0 \Lambda Q]
        +\tfrac{i}{2\pi}[\frakc_{1,j}- \ol{c_1}\frakc_{0,j}]Q
        \\
        &=-(b_1+\tfrac{1}{2} \mu_0\nu_0)c_j \Lambda Q -\tfrac{1}{2}\eta_1 c_j iQ
        +[\tfrac{i}{4} \nu_0^2c_j +\tfrac{i}{2\pi}(\frakc_{1,j}- \ol{c_1}\frakc_{0,j})]Q.
    \end{align*}
	Here, we used $\frac{1}{4}(1-y^2)Q^3=\Lambda Q$. Thus, the profile parts become
    \begin{align*}
        &[(c_j)_s-\tfrac{\lambda_s}{\lambda} j c_j-\nu_0 c_{j+1}
        +\tfrac{i}{4} \nu_0^2c_j +\tfrac{i}{2\pi}(\frakc_{1,j}- \ol{c_1}\frakc_{0,j})]Q
        \\
        &+(\tfrac{1}{2}\nu_0^2-\tfrac{1}{4}\mu_0^2+2\Re(c_2))c_jiQ -(\tfrac{x_s}{\lambda}-\nu_0)c_{j+1}Q
        \\
        &-(\tfrac{\lambda_s}{\lambda}+b_1+\tfrac{1}{2} \mu_0\nu_0)c_j\Lambda Q
        +(\gamma_s -\tfrac{1}{2}\eta_1-\tfrac{1}{2}\nu_0^2+\tfrac{1}{4}\mu_0^2-2\Re(c_2) )c_j iQ.
    \end{align*}
	Therefore, we rewrite the even part of \eqref{eq:wj equ rewrite-1} as
    \begin{equation}
        \begin{aligned}
            &[(c_j)_s-\tfrac{\lambda_s}{\lambda} j c_j-\nu_0 c_{j+1}
            +i(\tfrac{3}{4} \nu_0^2-\tfrac{1}{4}\mu_0^2+2\Re(c_2))c_j +\tfrac{i}{2\pi}(\frakc_{1,j}- \ol{c_1}\frakc_{0,j})]Q
    		\\
            &=(\tfrac{\lambda_s}{\lambda}+b_1+\tfrac{\mu_0\nu_0}{2})c_j \Lambda Q-(\gamma_s -\tfrac{1}{2}\eta_1-\tfrac{1}{2}\nu_0^2+\tfrac{1}{4}\mu_0^2-2\Re(c_2) )c_j iQ
            \\
            &\quad +(\tfrac{x_s}{\lambda}-\nu_0)c_{j+1}Q
            -\calM_{j}^{\text{even}},
        \end{aligned} \label{eq:wj equ mod even}
    \end{equation}
	where $\calM_j^{\text{even}}= \calM_{j,s}^{\text{even}}+\calM_{j,1}^{\text{even}} +\calM_{j,2}^{\text{even}}+\calM_{j,3}^{\text{even}}$, and
	\begin{align*}
		\calM_{j,s}^{\text{even}}&=
		(\partial_s-\tfrac{\lambda_s}{\lambda}\Lambda_{-j} +\gamma_s i)\eps_{j,e}-\tfrac{x_s}{\lambda}\eps_{j+1,e}-iw_{j+2,e},
		\\
		\calM_{j,1}^{\text{even}}&=\tfrac{1}{2}\left(\tfrac{x_s}{\lambda}-\nu_0\right)[w\calH(\ol{w}w_j)]_e,
		\\
		\calM_{j,2}^{\text{even}}&=\tfrac{i}{2}[\textnormal{NL}_j-\tfrac{1}{\pi}(\frakc_{1,j}- \ol{c_1}\frakc_{0,j})Q]_e,
        \\
		\calM_{j,3}^{\text{even}}&=\tfrac{i}{2}c_j[\textnormal{NL}_{Q}]_{e}.
	\end{align*}
	The modulation parts, the first three terms of the right side of \eqref{eq:wj equ mod even} are controlled by \eqref{eq:mod0 esti refine}. For $\calM_{j,s}^{\text{even}}$, using $(\eps_{j,e}, \calZ_c)_c=0$ and Lemma~\ref{lem:coercivity estimate 2}, we first derive
    \begin{align}
        (\calM_{j,s}^{\text{even}},\calZ_c)_c
        =((-\tfrac{\lambda_s}{\lambda}\Lambda_{-j} +\gamma_s i)\eps_{j,e} -\tfrac{x_s}{\lambda}\eps_{j+1,e}-iw_{j+2,e},\calZ_c)_c
        \lesssim \lambda^{j+2}. \label{eq:wj equ mod even s 1}
    \end{align}
    If $1\leq j\leq 2L-2$ with $L\geq2$, from \eqref{eq:coercivity H1 epsj}, we have
    \begin{align}
        (\calM_{j,s}^{\text{even}}+ic_{j+2}Q,\calZ_c)_c
        &=((-\tfrac{\lambda_s}{\lambda}\Lambda_{-j} +\gamma_s i)\eps_{j,e} -\tfrac{x_s}{\lambda}\eps_{j+1,e}-i\eps_{j+2,e},\calZ_c)_c \nonumber\\
        &\lesssim \lambda^{j+2+\frac12}. \label{eq:wj equ mod even s 2}
    \end{align}
    For the rest, we use nonlinear estimates. From \eqref{eq:wH w square}, \eqref{eq:wH w wj}, and \eqref{eq:mod0 esti refine}, we obtain
    \begin{align*}
        (\calM_{j,1}^{\text{even}},\calZ_c)_c
        =\tfrac{1}{2}([\tfrac{x_s}{\lambda}-\nu_0][w\calH(\ol{w}w_j)-c_jQ\calH(Q^2)]_e,\calZ_c)_c\lesssim \lambda^{j+2+\frac12}.
    \end{align*}
    From \eqref{eq:main nonlinear eps j even} and \eqref{eq:main nonlinear Q}, we have
    \begin{equation}
        \begin{aligned}
            (\calM_{j,2}^{\text{even}},\calZ_c)_c
            \lesssim   \lambda^{j+2+\frac12},
            \quad
            (\calM_{j,3}^{\text{even}},\calZ_c)_c
            \lesssim   \lambda^{j+2+\frac12},
        \end{aligned}\label{eq:wj equ mod even 12}
    \end{equation}
    respectively. Therefore, collecting \eqref{eq:wj equ mod even}, \eqref{eq:wj equ mod even s 2}, and \eqref{eq:wj equ mod even 12}, we conclude \eqref{eq:cj mod}. For \eqref{eq:cj mod aux}, we combine \eqref{eq:wj equ mod even}, \eqref{eq:wj equ mod even s 1}, \eqref{eq:wj equ mod even 12}, \eqref{eq:c j estimate}, together with the forthcoming estimate \eqref{eq:frakc esti}, which yields the desired bound.
\end{proof}

\subsection{Evolution of \texorpdfstring{$\frakc_{i,j}$}{cij}}\label{subsec:frakc}

In this subsection, we establish the quantitative properties of the modified quadratic
radiation parameters. Recall $\frakc_{i,j}$ and $\varphi$ defined in \eqref{eq: frakc ij def}.
The quantity $\frakc_{i,j}$ is designed to capture the leading-order quadratic radiation as the raw expression
\begin{align*}
    \int yQ^2 \ol{\wt\eps_i} \wt\eps_j,
\end{align*}
which appears in the derivation of $(c_j)_s$. 
However, differentiating it in time $s$ necessarily brings in
\[
    (\wt\eps_j)_s = (w_j)_s - (c_j)_s Q,
\]
and therefore requires precise control of both $(w_j)_s$ and $(c_j)_s$.  
Unlike the evolution of the parameters $\bm\beta_j$ and $c_j$, where the terms involving $(\wt\eps_j)_s$ can be eliminated by differentiating the orthogonality conditions, no such mechanism is available for the quadratic quantity considered here.
Consequently, differentiating the raw quadratic quantity produces additional quadratic radiation terms of size $O(\lambda^{i+j+3})$, so they cannot be treated as errors and must be included at the same level as the leading contributions in $(\frakc_{i,j})_s$.

Controlling these contributions would require incorporating them into the leading part of $(c_j)_s$, which in turn amounts to improving the estimate \eqref{eq:cj mod} from $O(\lambda^{j+2+\frac12})$ to $O(\lambda^{j+3})$. 
At this level, however, the evolution of $c_j$ generates new quadratic radiation terms of size $O(\lambda^{j+2+\frac12})$ that are no longer of the form $\tfrac12\int yQ^2\,\overline{\wt\eps_i}\wt\eps_j$.
Analyzing the time evolution of these new terms would then require a sharper equation for $\bm\beta_j$, and this refinement generates further quadratic radiation contributions at the corresponding error scales of $(\bm\beta_j)_s$. 
These new contributions are not of the same form as the original quadratic interaction, and each gain in accuracy introduces additional terms at the next scale, preventing the computation from closing. 
In this way, the procedure does not stabilize and fails to produce a closed system. For this reason, the raw quantity $\tfrac12\int yQ^2\,\overline{\wt\eps_i}\wt\eps_j$ is not suitable for deriving a \textit{closed} ODE system.

The modified quantity $\frakc_{i,j}$ avoids this issue.  
It is comparable in size to the raw quadratic term, while its time derivative does not require computing $(\wt\eps_j)_s$, since it involves only $w_i$ and $w_j$.
A key feature of this construction is the choice of weight
\[
    \varphi = \tfrac12 yQ^2 - 3yQ^4 + 2yQ^6,
\]
which cancels all profile contributions both in $\frakc_{i,j}$ and in its derivative $(\frakc_{i,j})_s$.  
This cancellation enables us to obtain a closed ODE system for the family $\{\frakc_{i,j}\}$.
%This cancellation is what makes it possible to obtain a closed ODE system for the family $\{\frakc_{i,j}\}$.

As discussed so far, our first goal is to derive the quantitative relation
\begin{align}
    \frakc_{i,j}=
    \frac{1}{2}\int yQ^2 \ol{\wt\eps_i} \wt\eps_j
    +O(\lambda^{i+j+1+\frac12}), \label{eq:frakc radiatin relation}
\end{align}
which verifies that $\frakc_{i,j}$ indeed provides the correct modified form of the dominant quadratic radiation.
To see why the definition of $\frakc_{i,j}$ produces the desired cancellation,
we begin by recalling the explicit identities
\begin{align}
	{\textstyle\int} y^2Q^4=2\pi,\quad 
	{\textstyle \int} y^2Q^6={\textstyle \int} y^2Q^8=\pi,\quad
	{\textstyle \int} y^2Q^{10}=\tfrac{5}{4}\pi,\quad
	{\textstyle \int} y^2Q^{12}=\tfrac{7}{4}\pi.
    \label{eq:yQ2n integration}
\end{align}	
These yield the cancellation relations
\begin{align}
    \tfrac{1}{2}{\textstyle \int} y^2Q^4-3{\textstyle \int} y^2Q^6+2{\textstyle \int} y^2Q^8=0, \quad
    \tfrac{1}{2}{\textstyle \int} y^2Q^6-3{\textstyle \int} y^2Q^8+2{\textstyle \int} y^2Q^{10}=0. \label{eq:profile cancel 1}
\end{align}
Moreover, using $2y^2Q^{2n}=y^2(1+y^2)Q^{2n+2}$, we also obtain
\begin{align}
    \tfrac{1}{2}{\textstyle \int} y^4Q^6-3{\textstyle \int} y^4Q^8+2{\textstyle \int} y^4Q^{10}=0. \label{eq:profile cancel 3}
\end{align}

The motivation behind the weight $\varphi$ is as follows.
Among the weights $yQ^{2m}$, the term $yQ^2$ carries the dominant contribution of the
interaction $\overline{w_i}w_j$, while $yQ^4$ and $yQ^6$ decay faster and may be used
as perturbative corrections.   
The profile part of $\overline{w_i}w_j$ always takes the form
\[
    \ol{P_{i,o}}P_{j,e}+\ol{P_{i,e}}P_{j,o}=(\overline{\bm\beta_i}c_j+\overline{c_i}\bm\beta_j) yQ^2,
\]
and thus, every profile contribution in the quadratic interaction $\int yQ^{2m}\overline{w_i} w_j$ is equal to
\[
    (\overline{\bm\beta_i}c_j+\overline{c_i}\bm\beta_j){\textstyle \int} y^2 Q^{2(m+1)}.
\]
The identities \eqref{eq:profile cancel 1} show that the coefficients $\frac12,-3,2$ are chosen so that the corresponding linear combination of these profile integrals vanishes.  
Consequently, the weight $\varphi$ annihilates the profile part of $\overline{w_i}w_j$ while retaining the dominant radiation contribution.

More generally, \eqref{eq:yQ2n integration} implies that any coefficients $(p,q)$ satisfying $p-q=1$ produce a linear combination
\[
    \tfrac12 yQ^2 - p yQ^4 + q yQ^6
\]
that cancels the profile contribution in $\overline{w_i} w_j$ itself.  
However, for our purposes, we must control not only $\frakc_{i,j}$ but also its time derivative $(\frakc_{i,j})_s$, since both are needed to obtain a closed evolution system for the radiation parameters.
Most choices of $(p,q)$ eliminate the profile term in $\frakc_{i,j}$, but leave residual profile contributions in $(\frakc_{i,j})_s$, which obstructs obtaining a precise evolution equation for $(\frakc_{i,j})_s$.
A direct computation shows that the choice $(p,q)=(3,2)$ achieves cancellation in both $\frakc_{i,j}$ and $(\frakc_{i,j})_s$, which is why the weight~$\varphi$ is defined in this particular form.

To analyze $\frakc_{i,j}$, we begin by isolating the genuine radiation component
in the product $\overline{w_i} w_j$. For all $i,j\ge0$, define
\begin{align*}
    g_{i,j}=\ol{w_i}w_j-\ol{P_{i,e}}P_{j,e}-\ol{P_{i,o}}P_{j,e}-\ol{P_{i,e}}P_{j,o},
\end{align*}
with a convention
\begin{align*}
    P_{0,e}\coloneqq Q,\quad P_{0,o}\coloneqq T_{o}.
\end{align*}
Note that, for $i=j=0$, we have
\begin{align*}
    g_{0,0}=|w|^2-Q^2-\ol{T_{o}}Q-QT_{o}=|w|^2-Q^2.
\end{align*}

\begin{rem}
    Among the quadratic profiles, the term $\overline{P_{i,o}}P_{j,o}$ is not subtracted in the definition of $g_{i,j}$. Since the quadratic profile term
    \[
        yQ^2 \overline{P_{i,o}}P_{j,o}
        = \overline{\bm{\beta}_i}\bm{\beta}_j y^3 Q^4
    \]
    is not in $L^1$, it must be regarded as part of the radiation.  
    For this reason, $\overline{P_{i,o}}P_{j,o}$ is left inside $g_{i,j}$ rather than 
    being subtracted as a profile term.     
\end{rem}

The following lemma records a basic estimate that will be used repeatedly.
\begin{lem}\label{lem:gij estimate}
    For any odd function $\psi_o$ satisfying $|\psi_o|\lesssim Q^3$, we have
    \begin{align}
        \int \psi_o g_{i,j}=O(\lambda^{i+j+2}). \label{eq:ij radiation g odd}
    \end{align}
    
\end{lem}

With this lemma in hand, we now derive \eqref{eq:frakc radiatin relation}.
Using the oddness of $\varphi$ together with \eqref{eq:profile cancel 1}, one checks that
\begin{align}
    \frakc_{i,j}=
    {\textstyle \int} \varphi g_{i,j}=\tfrac{1}{2} {\textstyle \int} yQ^2 g_{i,j}
    + {\textstyle \int} (-3yQ^4+2yQ^6)g_{i,j}. \label{eq:frakc rad aux}
\end{align}
For the last term on the right-hand side of \eqref{eq:frakc rad aux}, Lemma~\ref{lem:gij estimate} gives
\begin{align*}
    \frakc_{i,j}= \tfrac{1}{2}{\textstyle \int} yQ^2 g_{i,j}+O(\lambda^{i+j+2}).
\end{align*}
Recalling the definition of $g_{i,j}$, we write
\begin{align*}
    g_{i,j}= \overline{\wt\eps_i}\wt\eps_j
    + \overline{P_{i,e}}\eps_j
    + \overline{\eps_i}P_{j,e}.
\end{align*}
By inserting this expression for $g_{i,j}$ and applying \eqref{eq:coercivity H2 epsj}, we obtain \eqref{eq:frakc radiatin relation}.
Furthermore, using \eqref{eq:coercivity H1 epsj}, \eqref{eq:beta j estimate}, and $\wt\eps_j = P_{j,o} + \eps_j$, we obtain the bound
\begin{align}
    |\frakc_{i,j}|+\big|\tfrac{1}{2}\textstyle\int yQ^2 \ol{\wt\eps_i} \wt\eps_j\big|+ 
    \left|\tfrac{1}{2}\textstyle\int yQ^2 g_{i,j}\right|
    \lesssim \lambda^{i+j+1}. \label{eq:frakc esti}
\end{align}

\begin{proof}[Proof of Lemma~\ref{lem:gij estimate}]
    We note that, for $i,j\geq 1$,
    \begin{equation}
        \begin{aligned}
            g_{i,j}&=\ol{\eps_i}P_{j,e}+\ol{\wt{\eps}_i}P_{j,o}
            +\ol{P_{i,e}}\eps_j+\ol{\wt\eps_i}\wt\eps_j,
            \\
            g_{i,0}&=\ol{\eps_i}Q+\ol{\wt{\eps}_i}T_{o}+\ol{P_{i,e}}(\wt\eps-T_{1,o})+\ol{\wt\eps_i}\wt\eps,
            \\
            g_{0,j}&=\ol{(\wt\eps-T_{1,o})}P_{j,e}+\ol{\wt\eps} P_{j,o}+Q\eps_j+\ol{\wt\eps}\wt\eps_j,
            \\
            g_{0,0}&=2\Re[Q(\wt\eps-T_{1,o})]+|\wt\eps|^2.
        \end{aligned} \label{eq: gij rewrite}
    \end{equation}
    By \eqref{eq:coercivity Linf epsj}, \eqref{eq:coercivity H2 epsj}, \eqref{eq:beta j estimate}, \eqref{eq:c j estimate}, and the decomposition $\wt\eps_j=P_{j,o}+\eps_j$, we derive
    \begin{align*}
        {\textstyle\int} \psi_o g_{i,j} 
        &\!\lesssim\!
        \|\eps_{i,o}\|_{\dot \calH^2}\!|c_j|\!+\!\|\eps_{i,e}\|_{\dot \calH^2}\beta_j
        \!+\!|c_i|\|\eps_{j,o}\|_{\dot \calH^2}
        \!+\!\|\wt\eps_{i,o}\|_{\dot\calH^1}\!\|\eps_{j,e}\|_{\dot \calH^2}\!+\!\|\eps_{i,e}\|_{\dot\calH^2}\!\|\wt\eps_{j,o}\|_{\dot \calH^1}
        \\
        &\!\lesssim \lambda^{i+j+2}.
    \end{align*}
   For $g_{i,0}$, further using \eqref{eq:coercivity H2 and H1 inf}, we derive
    \begin{align*}
        {\textstyle\int} \psi_o g_{i,0}
        &\lesssim  
        \|\eps_{i,o}\|_{\dot \calH^2}+\|\eps_{i,e}\|_{\dot \calH^2}|\nu_0|
        +|c_i|\|\eps_{o}\|_{\dot \calH^2}
        +\|\wt\eps_{i,o}\|_{\dot\calH^1}\|\wt\eps_{e}\|_{\dot \calH^2}+\|\eps_{i,e}\|_{\dot\calH^2}\|\wt\eps_{o}\|_{\dot \calH^1}
        \\
        &\lesssim  \lambda^{i+2}.
    \end{align*}
    By symmetry, we also have ${\int} \psi_o g_{0,j} \lesssim \lambda^{j+2}$.
    For $g_{0,0}$, we have
    \begin{align*}
        {\textstyle\int} \psi_o g_{0,0}\lesssim \|\eps_o\|_{\dot\calH^2}+\|Q\wt\eps\|_{L^2}^2\lesssim \lambda^2,
    \end{align*}
    which conclude \eqref{eq:ij radiation g odd} for all $i,j\geq 0$.
\end{proof}

Having established the structural properties of $\frakc_{i,j}$ and its relation to the raw quadratic radiation, we now turn to the evolution of $\frakc_{i,j}$.
The next proposition provides the precise radiation estimate for $(\frakc_{i,j})_s$.

\begin{prop}[Radiation estimates II]\label{prop:radiation estimate 2}
    Let $L\geq 2$ with $v_0\in H^{2L+1}$. For each $0\leq i,j \leq 2L-2$, we have
	\begin{align}
		|(\mathfrak{c}_{i,j})_s-\tfrac{\lambda_s}{\lambda}(i+j+1)\mathfrak{c}_{i,j}+i\frakc_{i,j+2} -i\frakc_{i+2,j}- \fraka_{i,j}|\lesssim \lambda^{i+j+3+\frac12}, \label{eq:frakcij mod}
	\end{align}
    where $\fraka_{i,j}$ is given explicitly in \eqref{eq:Aij def} and satisfies $|\mathfrak{a}_{i,j}| \lesssim \lambda^{i+j+3}$.
\end{prop}
\begin{rem}
    The explicit expression of $\fraka_{i,j}$ in \eqref{eq:Aij def} is quite lengthy, but its precise formula is not needed for our purposes. What matters here is only that each $\fraka_{i,j}$ can be written in terms of the parameters $c_j$ and $\mathfrak{c}_{i,j}$ (together with their complex conjugates). Moreover, while $\fraka_{i,j}$ has the same size as $\mathfrak{c}_{i,j+2}$ and $\mathfrak{c}_{i+2,j}$, it involves only indices that increase by at most one.  More precisely, $\fraka_{i,j}$ consists of combinations of $\mathfrak{c}_{\ell_1,\ell_2}$ and $c_{\ell_3}$ with $\ell_1,\ell_2 \in \{0,1,2,i,i+1,j,j+1\}$ and $\ell_3 \in \{1,i,j\}$.
\end{rem}

In order to derive a closed system for $(\frakc_{i,j})_s$, it is necessary to separate the remaining profile contributions that persist inside $g_{i,j}$.  
For this purpose, we introduce a modified quantity $\check g_{i,j}$, obtained by subtracting the additional profile components as follows. For $i,j\neq 0$, we set
\begin{align*}
    \check g_{i,j}
    \coloneqq\,& g_{i,j}-\ol{P_{i,e}}\,\check{\bm{\beta}}_{j+1}(1+y^2)Q
     - \ol{\check{\bm{\beta}}_{i+1}}(1+y^2)Q \cdot P_{j,e}
    \\
    =\,&\ol{w_i}w_j
    -\ol{P_{i,e}}P_{j,e}-\ol{P_{i,o}}P_{j,e}-\ol{P_{i,e}}P_{j,o}
    - \ol{P_{i,o}}P_{j,o}  \qquad (i,j\neq 0).
    \\
    &-\ol{P_{i,e}}\,\check{\bm{\beta}}_{j+1}(1+y^2)Q
     - \ol{\check{\bm{\beta}}_{i+1}}(1+y^2)Q \cdot P_{j,e},
\end{align*}
When one of the indices is zero, we set
\begin{align*}
    \check g_{i,0}&\coloneqq g_{i,0}-\tfrac{\nu_0^2}{8}\ol{P_{i,e}}Q
    -\ol{P_{i,o}}T_{o}-\ol{P_{i,e}}\check{T}_e 
    -\ol{\check{\bm{\beta}}_{i+1}}(1+y^2)Q\cdot Q,
    \\
    \check g_{0,j}&\coloneqq g_{0,j}-\tfrac{\nu_0^2}{8}Q P_{j,e}-\ol{T_{o}}P_{j,o}
    -Q\cdot \check{\bm{\beta}}_{j+1}(1+y^2)Q
    -\ol{\check{T}_e }\cdot P_{j,e},
    \\
    \check g_{0,0}&\coloneqq g_{0,0} -2\cdot \tfrac{\nu_0^2}{8}Q^2-|\check T_{o}|^2-2\Re(Q\cdot \check{T}_e ).
\end{align*}
We now state the key lemma for the second radiation estimates.

\begin{lem}\label{lem:gij induction}
    For $0\le \ell_1,\ell_2,\ell_3,\ell_4 \le 2L-2$, whenever at least one among $\ell_1,\ell_2,\ell_3,\ell_4$ is zero, we have
    \begin{align}
        \begin{aligned}
            {\textstyle\int}\varphi \ol{w_{\ell_1}}w_{\ell_2}\calH( \ol{w_{\ell_3}}w_{\ell_4})=&{\textstyle\int}\tfrac{y}{1+y^2}g_{\ell_1,\ell_2}\calH(g_{\ell_3,\ell_4})
            +\tfrac{1}{2}\ol{c_{\ell_3}}c_{\ell_4} {\textstyle\int}y^2Q^4 \check g_{\ell_1,\ell_2}
            \\
            &+O(\lambda^{\ell_1+\ell_2+\ell_3+\ell_4+2+\frac12}).
        \end{aligned} \label{eq:ijs quartic nonlin goal varphi}
    \end{align}
    For arbitrary indices $0\le \ell_1,\ell_2,\ell_3,\ell_4 \le 2L-2$, we have
    \begin{equation}
        \begin{aligned}
            {\textstyle\int}\tfrac{y}{1+y^2}g_{\ell_1,\ell_2}\!\calH(g_{\ell_3,\ell_4}\!)
            +
            {\textstyle\int}\tfrac{y}{1+y^2}g_{\ell_3,\ell_4}\!\calH(g_{\ell_1,\ell_2}\!)
            =-\tfrac{1}{\pi}\mathfrak{c}_{\ell_1,\ell_2}\mathfrak{c}_{\ell_3,\ell_4}
            +O(\lambda^{\ell_1+\ell_2+\ell_3+\ell_4+2+\frac12}).
        \end{aligned}
        \label{eq:ijs quartic commute}
    \end{equation}
    Moreover, for $0\leq i,j\leq 2L-2$, we have
    \begin{equation}
        \begin{aligned}
            \tfrac{1}{2}{\textstyle\int} y^2Q^4\check g_{i,j} = &
            -\tfrac{1}{2}\ol{c_i}{\textstyle \int}  y^2Q^4 \check g_{0,j}
            -\tfrac{1}{2}{c_j}{\textstyle \int}  y^2Q^4 \check g_{i,0}
            \\
            &+2\mathfrak{c}_{i+1,j}
            +2\mathfrak{c}_{i,j+1}
            +\tfrac{1}{\pi}\mathfrak{c}_{0,j}\mathfrak{c}_{i,0}
            +O(\lambda^{i+j+2+\frac12}).
        \end{aligned}\label{eq:ijs F nonlin goal}
    \end{equation}
\end{lem}

The goal of the proposition is to derive a closed ODE for the coefficients $\{ \mathfrak{c}_{i,j} \}$. In the computation of $(\mathfrak{c}_{i,j})_s$, the nonlinear contribution produces a 4-linear expression of the form 
\begin{equation*}
   {\textstyle\int}\varphi \ol{w_{\ell_1}}w_{\ell_2}\calH( \ol{w_{\ell_3}}w_{\ell_4})
   +{\textstyle\int}\varphi \ol{w_{\ell_3}}w_{\ell_4}\calH( \ol{w_{\ell_1}}w_{\ell_2}).
\end{equation*} 
A first observation is that the expression is symmetric in $(\ell_1,\ell_2)$ and $(\ell_3,\ell_4)$.  
Using $(1+y^2)\langle y\rangle^{-2}=1$ together with \eqref{eq:CommuteHilbert}, this symmetry yields a structural cancellation at the level of the 4-linear form, reducing it to a product of two quadratic radiation terms; see \eqref{eq:ijs quartic commute}.  

Crucially, this symmetric cancellation can be exploited only after the leading profile pieces inside $\overline{w_i}w_j$ have been separated out.  
This motivates the introduction of the modified quantities $g_{i,j}$ and $\check g_{i,j}$.  
The first identity \eqref{eq:ijs quartic nonlin goal varphi} implements this separation and rewrites the 4-linear expression entirely in terms of $g_{\ell_1,\ell_2}$, up to an additional quadratic term weighted by $y^2Q^4$.

To understand the structure of \eqref{eq:ijs quartic nonlin goal varphi}, 
we expand each $w_{\ell_k}$ into its profile part and its radiation part inside 
the 4-linear expression.  
By the choice of $\varphi$, the pure profile contribution (when all four 
factors are profiles) cancels and does not appear in the identity.  
Among the remaining terms, two kinds of structures emerge, while the rest are absorbed into the error.

On one hand, if all four factors are radiative, we obtain the genuine 4-linear radiation 
interaction
\[
    {\textstyle\int} \tfrac{y}{1+y^2} g_{\ell_1,\ell_2} 
        \calH(g_{\ell_3,\ell_4}),
\]
which will be the first term on the right-hand side of 
\eqref{eq:ijs quartic nonlin goal varphi}. 
This term is then handled through the symmetric cancellation described earlier, by applying \eqref{eq:ijs quartic commute}, which reduces it to a product of $\frakc_{\ell,\ell^\prime}$'s.

On the other hand, when one side of $\mathcal H$ carries a quadratic profile 
and the other side is purely radiative (schematically, ``$(1,2)$ profile, $(3,4)$ 
radiation''), the resulting mixed interaction generates the weighted quadratic 
term
\[
    \tfrac{1}{2}\overline{c_{\ell_3}}c_{\ell_4}
    {\textstyle\int} y^2 Q^4\check g_{\ell_1,\ell_2}.
\]
In this case, the decay cannot be freely exchanged through the Hilbert transform.  
When viewed with the stronger weight $y^2Q^4$, some profile parts of $g_{\ell_1,\ell_2}$ no longer behave as radiation and must be subtracted.  
The remaining radiative part is denoted by $\check g_{\ell_1,\ell_2}$.

Aside from these two main cases, all other mixed configurations already contain enough decay. For instance, this happens when each side of $\mathcal H$ has one profile factor and one radiation factor. These contributions can be estimated directly and are absorbed into the $O(\lambda^{\ell_1+\ell_2+\ell_3+\ell_4+2+\frac12})$ remainder in \eqref{eq:ijs quartic nonlin goal varphi}.

To close the ODE system, it remains to express
\[
{\textstyle\int} y^2Q^4\check g_{i,j}
\]
in terms of known quantities such as $\frakc_{i,j}$. 
This is achieved in the last identity \eqref{eq:ijs F nonlin goal}. 
Here, another technical trick comes in. We replace the weight $\frac{1}{2}y^2Q^4$ by a more convenient form $\Phi$ satisfying
\[
\Phi = \tfrac{1}{2} y^2Q^4 + O(Q^4), 
\qquad 
\Phi = \partial_y\phi,
\]
so that an integration by parts converts the term into one involving $\phi$. We note that $\phi$ plays a role similar to $\varphi$, but with a \textit{different} weight.
Because the correction $O(Q^4)$ decays faster, we may replace the weight $y^2 Q^4$ by $\Phi$, and we obtain %the estimate
\begin{align*}
    {\textstyle\int} \tfrac{1}{2}y^2Q^4\check g_{i,j} = {\textstyle\int} \Phi \check g_{i,j}+ O(\lambda^{i+j+2+\frac12})
\end{align*}
Lastly, integrating by parts with $\Phi = \partial_y\phi$ gives an expression in terms of $\frakc_{i,j}$, \eqref{eq:ijs F nonlin goal}.%applying integration by parts with $\Phi = \partial_y\phi$ expresses the result in terms of $\frakc_{i,j}$, yielding \eqref{eq:ijs F nonlin goal}.

Before embarking on the proof of Lemma \ref{lem:gij induction}, we explain the motivation for the specific choice of the auxiliary weights $\Phi$ and $\phi$, and we fix them to fulfill the required cancellation identities.
We set
\begin{align*}
	\Phi=&\tfrac{1}{2}y^2Q^4+a_1 y^2Q^6+a_2 y^2Q^8+a_3 y^2Q^{10}+a_4 y^2Q^{12}.
\end{align*}
We seek to choose $\Phi$ as the derivative of an odd rational function $\phi$,
\begin{align}
    \partial_y \phi=\Phi, \quad \phi\in L^2, \quad \phi \text{ is an odd rational function,} \label{eq:odd rational function}
\end{align}
and to impose three cancellation conditions:
\begin{equation}
    \begin{aligned}
        {\textstyle \int} \Phi=0,
        \qquad
        {\textstyle \int} \Phi Q^2=0,
        \qquad
        {\textstyle \int} \phi yQ^2=0.
    \end{aligned} \label{eq:prof cancel choice}
\end{equation}
By integration by parts, one can check
\begin{align}
    {\textstyle \int} \phi yQ^4={\textstyle \int} \Phi Q^2=0. \label{eq:prof cancel int by part}
\end{align}
Then, based on \eqref{eq:yQ2n integration} and the identity ${\textstyle \int} y^2Q^{12} = \frac{21}{8}\pi$, we can choose the coefficients $(a_1,a_2,a_3,a_4)$ to satisfy \eqref{eq:odd rational function} and \eqref{eq:prof cancel choice}, using the four available degrees of freedom.

In particular, if we take
\begin{align*}
    (a_1,a_2,a_3,a_4)=(-13, \tfrac{86}{3}, -\tfrac{40}{3},0),
\end{align*}
then, for this explicit choice, $\Phi$ and $\phi$ satisfy
\begin{align*}
    \partial_y \phi=\Phi, \quad 
    \phi= -\tfrac{1}{24}y^3(3y^4-42y^2+35)Q^{8} = - yQ^2+y\cdot O(Q^4).
\end{align*}
This choice also satisfies 
\eqref{eq:prof cancel choice} and \eqref{eq:prof cancel int by part}:
\begin{align}
    {\textstyle \int} \Phi
    = {\textstyle \int} \Phi Q^2
    = {\textstyle \int} \phi yQ^2
    = {\textstyle \int} \phi yQ^4
    = {\textstyle \int} \phi y^3Q^4 = 0.
    \label{eq:profile cancel phi}
\end{align}
We now turn to the proof of Lemma~\ref{lem:gij induction}.

\begin{proof}[Proof of Lemma~\ref{lem:gij induction}]
    We first claim that, for any even function $\psi_e$ with $|\psi_e|\lesssim Q^4$,
    \begin{align}
        {\textstyle\int} \psi_e \check g_{i,j}=O(\lambda^{i+j+2+\frac12}). \label{eq:ij radiation}
    \end{align}
    The proof of \eqref{eq:ij radiation} is similar to Lemma~\ref{lem:gij estimate}.
    We first observe, for $i,j\geq 1$,
    \begin{equation}
        \begin{aligned}
            (\check g_{i,j})_e
            &=
            \ol{P_{i,e}} \check\eps_{j,e} 
            +
            \ol{\check\eps_i} P_{j,e}
            +
            \ol{P_{i,o}}\eps_{j,o}
            +
            \ol{\eps_{i,o}}P_{j,o}
            +
            (\ol{\eps_i} \eps_{j})_e,
            \\
            (\check g_{i,0})_e
            &=
            \ol{P_{i,e}} (\check{\eps}_e -\tfrac{\nu_0^2}{8}) 
            +
            \ol{\check\eps_i} Q
            +
            \ol{P_{i,o}}\eps_{o}
            +
            (\ol{\eps_i} \wt\eps)_e,
            \\
            (\check g_{0,j})_e
            &=
            Q \check\eps_{j,e} 
            +
            \ol{(\check{\eps}_e -\tfrac{\nu_0^2}{8})} P_{j,e}
            +
            \ol{\eps_{o}}P_{j,o}
            +
            (\ol{\wt\eps} \eps_{j})_e,
            \\
            (\check g_{0,0})_e
            &=2\Re [Q(\check{\eps}_e -\tfrac{\nu_0^2}{8})]
            +2\Re(\ol{\eps_{o}}T_{o})
            +
            (|\wt\eps-T_{o}|^2)_e.
        \end{aligned} \label{eq:gij check rewrite}
    \end{equation}
    Using $|\psi_e|\lesssim Q^4$, \eqref{eq:coercivity H3}, and \eqref{eq:coercivity H3 epsj check}, together with term-by-term estimates for \eqref{eq:gij check rewrite} as in Lemma~\ref{lem:gij estimate}, we deduce \eqref{eq:ij radiation}.

    Moreover, we claim
    \begin{align}
        \|Q (g_{i,j}-\ol{P_{i,o}}P_{j,o})\|_{L^2}\lesssim \lambda^{i+j+1+\frac12}. \label{eq:ij rad g H1}
    \end{align}
    Recall the convention,
    \begin{align*}
        P_{0}=Q+T_o,\quad P_{0,e}=Q,\quad P_{0,o}=T_o,\quad c_0=1 ,\quad \bm{\beta}_0=\tfrac{i}{2}\nu_0.
    \end{align*}
    From \eqref{eq: gij rewrite}, we have
    \begin{equation*}
        \begin{aligned}
            g_{i,j}-\ol{P_{i,o}}P_{j,o}&=\ol{\eps_i}P_{j,e}+\ol{{\eps}_i}P_{j,o}
            +\ol{P_{i,e}}\eps_j+\ol{\wt\eps_i}\wt\eps_j,
            \\
            g_{i,0}-\ol{P_{i,o}}T_{o}&=\ol{\eps_i}Q+\ol{{\eps}_i}T_{o}+\ol{P_{i,e}}(\wt\eps-T_{1,o})+\ol{\wt\eps_i}\wt\eps,
            \\
            g_{0,j}-\ol{T_{o}}P_{j,o}&=\ol{(\wt\eps-T_{1,o})}P_{j,e}+\ol{(\wt\eps-T_{1,o})} P_{j,o}+Q\eps_j+\ol{\wt\eps}\wt\eps_j,
            \\
            g_{0,0}-|T_o|^2&=2\Re[Q(\wt\eps-T_{1,o})]+|\wt\eps|^2-|T_o|^2.
        \end{aligned} 
    \end{equation*}
    By \eqref{eq:coercivity Linf epsj}, \eqref{eq:coercivity H2 epsj}, \eqref{eq:beta j estimate}, and \eqref{eq:c j estimate}, we conclude \eqref{eq:ij rad g H1}.
    
%\vspace{5pt} 
   \underline{\textbf{Step 1.}} We prove \eqref{eq:ijs quartic nonlin goal varphi}.
    In addition to \eqref{eq:ijs quartic nonlin goal varphi}, we will employ the following analogue for $\phi$, which will play a role in the proof of \eqref{eq:ijs F nonlin goal}: if at least one of $\ell_1,\ell_2,\ell_3,\ell_4$ is zero, then
    \begin{equation}
        \begin{aligned}
            {\textstyle\int}(-\tfrac12\phi) \ol{w_{\ell_1}}w_{\ell_2}\calH( \ol{w_{\ell_3}}w_{\ell_4})=&
            {\textstyle\int}\tfrac{y}{1+y^2}g_{\ell_1,\ell_2}\calH(g_{\ell_3,\ell_4})
            +\tfrac12\ol{c_{\ell_3}}c_{\ell_4} {\textstyle\int}y^2Q^4 \check g_{\ell_1,\ell_2}
            \\
            &+O(\lambda^{\ell_1+\ell_2+\ell_3+\ell_4+2+\frac12}).
        \end{aligned} \label{eq:ijs quartic nonlin goal phi}
    \end{equation}
    For notational convenience, we use $\psi$ to cover both cases; either $\psi=\varphi$ or $\psi=-\frac12 \phi$. Although Lemma~\ref{lem:gij induction} is stated only for $\varphi$, the same argument applies to both cases. Hence, in what follows, we establish the claim for $\psi$, which simultaneously yields \eqref{eq:ijs quartic nonlin goal varphi} and \eqref{eq:ijs quartic nonlin goal phi}.
    Using the decomposition
    \begin{align*}
        \ol{w_{\ell_3}}w_{\ell_4}
        =g_{\ell_3,\ell_4}
        +\ol{c_{\ell_3}}c_{\ell_4}Q^2
        +(\ol{c_{\ell_3}}\bm{\beta}_{\ell_4}+\ol{\bm{\beta}_{\ell_3}}c_{\ell_4})yQ^2,
    \end{align*}
    we rewrite the $4$-linear term ${\textstyle\int}\psi \ol{w_{\ell_1}}w_{\ell_2}\calH( \ol{w_{\ell_3}}w_{\ell_4})$ as
    \begin{align}
		{\textstyle\int}\psi \ol{w_{\ell_1}}w_{\ell_2}\calH( \ol{w_{\ell_3}}w_{\ell_4})
        =&{\textstyle \int} \psi \ol{w_{\ell_1}}w_{\ell_2}\calH(g_{\ell_3,\ell_4}) \label{eq:ijs 3-1 without prof}
		\\
		&+\ol{c_{\ell_3}}c_{\ell_4}{\textstyle \int} \psi yQ^2 \ol{w_{\ell_1}}w_{\ell_2} \label{eq:ijs 3-1-1 with prof}
		\\
		&-(\ol{c_{\ell_3}}\bm{\beta}_{\ell_4}+\ol{\bm{\beta}_{\ell_3}}c_{\ell_4}){\textstyle \int} \psi Q^2 \ol{w_{\ell_1}}w_{\ell_2}, \label{eq:ijs 3-1-2 with prof}
	\end{align}
    with the convention $c_0=1$ and $\bm{\beta}_{0}=\tfrac{i}{2}\nu_0$.
    
	We consider \eqref{eq:ijs 3-1 without prof}. We further decompose this into
	\begin{align}
		\eqref{eq:ijs 3-1 without prof}
		=&{\textstyle \int} \psi g_{\ell_1,\ell_2}\calH(g_{\ell_3,\ell_4}) \nonumber
		\\
		&+\ol{c_{\ell_1}}c_{\ell_2}{\textstyle \int} \psi Q^2 \calH(g_{\ell_3,\ell_4}) 
        +(\ol{c_{\ell_1}}\bm{\beta}_{\ell_2}+\ol{\bm{\beta}_{\ell_1}}c_{\ell_2}){\textstyle \int} \psi yQ^2\calH(g_{\ell_3,\ell_4}). \label{eq:ijs 3-1 rem-2}
	\end{align}
	We claim that
	\begin{align}
		\eqref{eq:ijs 3-1 rem-2}=O(
		\lambda^{\ell_1+\ell_2+\ell_3+\ell_4+2+\frac12}). \label{eq:ijs 3-1 rem claim}
	\end{align}
	Before showing this, we first claim the following properties for weight functions:
    \begin{equation}
		\begin{aligned}
			{\textstyle \int} \calH[\psi Q^2]  =0,\qquad
            {\textstyle \int} \calH[\psi Q^2] \cdot y^2Q^2=0, \qquad
			{\textstyle \int} \calH[\psi Q^2] \cdot Q^2 =0.
		\end{aligned} \label{eq:profile cancel additional}
	\end{equation}
   For \eqref{eq:profile cancel additional}, using \eqref{eq:profile cancel 1} and \eqref{eq:profile cancel 3} for $\varphi$, and \eqref{eq:profile cancel phi} for $\phi$, we have $\psi$ is odd and
    \begin{align}
        {\textstyle \int} \psi yQ^2
        = {\textstyle \int} \psi yQ^4
        = {\textstyle \int} \psi y^3Q^4=0. \label{eq:profile cancel psi}
    \end{align}
    Thus, we have $\calH(\psi Q^2)$ is even and
    \begin{align*}
        \calH[\psi Q^2] \sim y^{-4} \quad \text{at}\quad y\to \pm \infty.
    \end{align*}
	Using \eqref{eq:CommuteHilbert} first, then \eqref{eq:profile cancel 1} for $\varphi$ and \eqref{eq:profile cancel phi} for $\phi$, we obtain
	\begin{equation*}
		\begin{aligned}
			{\textstyle \int} \calH[\psi Q^2] \cdot y^2Q^2
			&=
			{\textstyle \int} \calH[\psi yQ^2] \cdot yQ^2=-{\textstyle \int} [\psi yQ^2] \cdot \calH(yQ^2)
			= {\textstyle \int} \psi yQ^4=0, 
		\end{aligned} 
	\end{equation*}
	Therefore, we conclude \eqref{eq:profile cancel additional} by obtaining
	\begin{align*}
		{\textstyle \int} \calH[\psi Q^2] 
		=
		\tfrac{1}{2}{\textstyle \int} \calH[\psi Q^2] \cdot (1+y^2)Q^2
        =
		\tfrac{1}{2}{\textstyle \int} \calH[\psi Q^2] \cdot Q^2
		=
		-\tfrac{1}{2}{\textstyle \int} \psi yQ^4=0.
	\end{align*}

	To show \eqref{eq:ijs 3-1 rem claim}, we control the first part of \eqref{eq:ijs 3-1 rem-2}. %, $\ol{c_{\ell_1}}c_{\ell_2}{\textstyle \int} \psi Q^2 \calH(g_{\ell_3,\ell_4})$. 
    The terms $g_{i,j}-\check g_{i,j}$ contain only the profile part, and for convenience, we denote $g_{\ell,\ell^\prime}-\check g_{\ell,\ell^\prime}\eqqcolon P_{\ell,\ell^\prime}^g$ by
    \begin{equation*}
        \begin{aligned}
            P_{i,j}^g&\coloneqq \ol{P_{i,o}}P_{j,o} +\ol{P_{i,e}}\,\check{\bm{\beta}}_{j+1}(1+y^2)Q
        + \ol{\check{\bm{\beta}}_{i+1}}(1+y^2)Q \cdot P_{j,e},
        \\
        P_{i,0}^g&\coloneqq \tfrac{\nu_0^2}{8}\ol{P_{i,e}}Q
        +\ol{P_{i,o}}T_{o}-\ol{P_{i,e}}\check{T}_e 
        +\ol{\check{\bm{\beta}}_{i+1}}(1+y^2)Q\cdot Q,
        \\
        P_{0,j}^g&\coloneqq \tfrac{\nu_0^2}{8}Q P_{j,e}-\ol{T_{o}}P_{j,o}
        +Q\cdot \check{\bm{\beta}}_{j+1}(1+y^2)Q
        +\ol{\check{T}_e }\cdot P_{j,e},
        \\
        P_{0,j}^g&\coloneqq 2\cdot \tfrac{\nu_0^2}{8}Q^2+|\check T_{o}|^2+2\Re(Q\cdot \check{T}_e ),
        \end{aligned} 
    \end{equation*}
    for $i,j\geq 1$.
    Then, we rewrite the first part of \eqref{eq:ijs 3-1 rem-2}, $\ol{c_{\ell_1}}c_{\ell_2}{\textstyle \int} \psi Q^2 \calH(g_{\ell_3,\ell_4})$ as follows:
	\begin{align*}
		\ol{c_{\ell_1}}c_{\ell_2}{\textstyle \int} \psi Q^2 \calH(g_{\ell_3,\ell_4})
		=&
		-\ol{c_{\ell_1}}c_{\ell_2}{\textstyle \int} \calH[\psi Q^2]  \cdot P_{\ell_3,\ell_4}^g
        -\ol{c_{\ell_1}}c_{\ell_2}{\textstyle \int} \calH[\psi Q^2]  \cdot  (\check g_{\ell_3,\ell_4})_e.
	\end{align*}
	Thanks to \eqref{eq:profile cancel additional}, we deduce ${\textstyle \int} \calH[\psi Q^2]  \cdot P_{\ell_3,\ell_4}^g=0$. 
	%\begin{align*}
	%	-\ol{c_{\ell_1}}c_{\ell_2}{\textstyle \int} \calH[\psi Q^2]  \cdot P_{\ell_3,\ell_4}^g=0.
	%\end{align*}
	For $(\check g_{\ell_3,\ell_4})_e$, we remark \eqref{eq:gij check rewrite}.
	By \eqref{eq:c j estimate} and \eqref{eq:coercivity H3}, we arrive at
	\begin{align*}
		{\textstyle \int} \calH[\psi Q^2] \cdot [(\ol{\check{\eps}_e -\tfrac{\nu_0^2}{8}Q}) P_{\ell,e}]
		\lesssim \lambda^{\ell}\|\check{\eps}_e  -\tfrac{\nu_0^2}{8}Q\|_{\dot\calH^3}
		\lesssim 
		\lambda^{\ell+3}.
	\end{align*}
	For $\ol{\check\eps_{\ell_3,e}} P_{\ell_4,e} $ type terms, \eqref{eq:c j estimate} and \eqref{eq:coercivity H3 epsj check} imply
	\begin{align*}
		{\textstyle \int} \calH[\psi Q^2] \cdot [\ol{\check\eps_{\ell_3,e}} P_{\ell_4,e} ]
		\lesssim
		\lambda^{\ell_4}\|\check\eps_{\ell_3,e}\|_{\dot\calH^3}
		\lesssim 
		\lambda^{\ell_3+\ell_4+2+\frac12}.
	\end{align*}
    For $\ol{\eps_o} P_{\ell_4,o}$ and $\ol{\eps_{\ell_3,o}} P_{\ell_4,o}$ type terms, thanks to \eqref{eq:beta j estimate}, \eqref{eq:coercivity H2 and H1 inf}, and \eqref{eq:coercivity H2 epsj}, we have
    \begin{align*}
        {\textstyle \int} \calH[\psi Q^2] \cdot [\ol{\eps_o} P_{\ell_4,o}]
		&\lesssim \lambda^{\ell_4+\frac12}\|\eps_o\|_{\dot\calH^2}
		\lesssim 
		\lambda^{\ell_4+2+\frac12},
        \\
        {\textstyle \int} \calH[\psi Q^2] \cdot [\ol{\eps_{\ell_3,o}} P_{\ell_4,o}]
		&\lesssim \lambda^{\ell_4+\frac12}\|\eps_{\ell_3,o}\|_{\dot\calH^2}
		\lesssim 
		\lambda^{\ell_3+\ell_4+2+\frac12}.
    \end{align*}
	By \eqref{eq:coercivity H1} and \eqref{eq:coercivity H2 epsj}, the $(\ol{\wt\eps}\eps_{\ell_4})_e$ and $(\ol{\eps_{\ell_3}}\eps_{\ell_4})_e$ type terms are estimated as follows:%nonlinear terms $(\ol{\wt\eps}\eps_{\ell_4})_e$ and $(\ol{\eps_{\ell_3}}\eps_{\ell_4})_e$ are estimated as follows:
	\begin{align*}
		|{\textstyle \int} \calH[\psi Q^2] \cdot (\ol{\wt\eps}\eps_{\ell_4})_e|
		&\lesssim \|\wt\eps\|_{\dot\calH^1}\|\eps_{\ell_4}\|_{\dot\calH^2}
		\lesssim 
		\lambda^{\ell_4+2+\frac12},
        \\
        |{\textstyle \int} \calH[\psi Q^2] \cdot (\ol{\eps_{\ell_3}}\eps_{\ell_4})_e|
		&\lesssim \|\eps_{\ell_3}\|_{\dot\calH^2}\|\eps_{\ell_4}\|_{\dot\calH^2}
		\lesssim 
		\lambda^{\ell_3+\ell_4+2+\frac12}.
	\end{align*}
	Therefore, we conclude $\ol{c_{\ell_1}}c_{\ell_2}{\textstyle \int} \psi Q^2 \calH(g_{\ell_3,\ell_4}) =O(
		\lambda^{\ell_1+\ell_2+\ell_3+\ell_4+2+\frac12})$.
	%\begin{align*}
	%	\ol{c_{\ell_1}}c_{\ell_2}{\textstyle \int} \psi Q^2 \calH(g_{\ell_3,\ell_4}) =O(
	%	\lambda^{\ell_1+\ell_2+\ell_3+\ell_4+2+\frac12}).
	%\end{align*}
	
	Next, we consider the second term of \eqref{eq:ijs 3-1 rem-2}.
	We observe that
	\begin{align*}
		\calH[\psi yQ^2]
		=y\calH[\psi Q^2]
		\sim y^{-3} \quad \text{at} \quad y\to \pm\infty.
	\end{align*}
	We also note that
	\begin{align*}
		(g_{\ell_3,\ell_4})_o=
		\ol{\eps_{\ell_3,o}} P_{\ell_4,e} +\ol{P_{\ell_3,e}}\eps_{\ell_4,o}+\ol{\wt\eps_{\ell_3,e}} P_{\ell_4,o}
        +\ol{P_{\ell_3,o}} \eps_{\ell_4,e}
        +(\ol{\eps_{\ell_3}}\eps_{\ell_4})_o,
	\end{align*}
    with the convention $\eps_{0}=\wt\eps-T_{o}$.
    Thus, from \eqref{eq:ij radiation g odd}, we obtain 
	\begin{align*}
		{\textstyle \int} \calH[\psi yQ^2]\cdot (g_{\ell_3,\ell_4})_o
		\lesssim \lambda^{\ell_3+\ell_4+2},
	\end{align*}
	which yields $(\ol{c_{\ell_1}}\bm{\beta}_{\ell_2}+\ol{\bm{\beta}_{\ell_1}}c_{\ell_2}){\textstyle \int} \psi yQ^2\calH(g_{\ell_3,\ell_4})=O(\lambda^{\ell_1+\ell_2+\ell_3+\ell_4+2+\frac12})$.
	%\begin{align*}
	%	(\ol{c_{\ell_1}}\bm{\beta}_{\ell_2}+\ol{\bm{\beta}_{\ell_1}}c_{\ell_2}){\textstyle \int} \psi yQ^2\calH(g_{\ell_3,\ell_4})=O(\lambda^{\ell_1+\ell_2+\ell_3+\ell_4+2+\frac12}).
	%\end{align*}
    
	Therefore, we have finished proving the claim \eqref{eq:ijs 3-1 rem claim}, and we conclude
	\begin{align*}
		\eqref{eq:ijs 3-1 without prof}
		={\textstyle \int} \psi g_{\ell_1,\ell_2}\calH(g_{\ell_3,\ell_4})
		+O(\lambda^{\ell_1+\ell_2+\ell_3+\ell_4+2+\frac12}).
	\end{align*}
    Next, we claim that
    \begin{equation}
        \begin{aligned}
            {\textstyle \int} \psi g_{\ell_1,\ell_2}\calH(g_{\ell_3,\ell_4})
            =&{\textstyle \int} \tfrac{1}{2}yQ^2 g_{\ell_1,\ell_2}\calH(g_{\ell_3,\ell_4})
            +\ol{\bm{\beta}_{\ell_1}}\bm{\beta}_{\ell_2}\ol{\bm{\beta}_{\ell_3}}\bm{\beta}_{\ell_4}{\textstyle \int} \tfrac{1}{2}y^4Q^6
            \\
            &+O(\lambda^{\ell_1+\ell_2+\ell_3+\ell_4+2+\frac12}),
        \end{aligned} \label{eq:ijs psi to yQ square claim}
    \end{equation}
    with the convention $\bm{\beta}_{0}=\tfrac{i}{2}\nu_0$. Thus, if one of $\ell_1,\ell_2,\ell_3,\ell_4$ is zero, then we have
    \begin{align*}
        \eqref{eq:ijs 3-1 without prof}
        ={\textstyle \int} \tfrac{1}{2}yQ^2 g_{\ell_1,\ell_2}\calH(g_{\ell_3,\ell_4})
        +O(\lambda^{\ell_1+\ell_2+\ell_3+\ell_4+2+\frac12}).
    \end{align*}
    Let $\td \psi \coloneqq \psi -\tfrac{1}{2}yQ^2=y\cdot O(Q^4)$. %$\td \psi$ be
    %\begin{align*}
     %   \td \psi \coloneqq \psi -\tfrac{1}{2}yQ^2=y\cdot O(Q^4).
    %\end{align*}
    For our claim \eqref{eq:ijs psi to yQ square claim}, it suffices to show
    \begin{align*}
        {\textstyle \int} \td \psi g_{\ell_1,\ell_2}\calH(g_{\ell_3,\ell_4})
		=\ol{\bm{\beta}_{\ell_1}}\bm{\beta}_{\ell_2}\ol{\bm{\beta}_{\ell_3}}\bm{\beta}_{\ell_4}{\textstyle \int} \tfrac{1}{2}y^4Q^6
        +O(\lambda^{\ell_1+\ell_2+\ell_3+\ell_4+2+\frac12}).
    \end{align*}
    We extract the profile term in ${\textstyle \int} \td \psi g_{\ell_1,\ell_2}\calH(g_{\ell_3,\ell_4})$ as follows:
    \begin{align}
        {\textstyle \int} \td \psi g_{\ell_1,\ell_2}\calH(g_{\ell_3,\ell_4})
        =
        &-{\textstyle \int} \calH[\td \psi (g_{\ell_1,\ell_2}-\ol{P_{\ell_1,o}}P_{\ell_2,o})]  \cdot g_{\ell_3,\ell_4} \label{eq:ijs td psi quartic 1}
        \\
        &-{\textstyle \int} \calH(\td \psi \ol{P_{\ell_1,o}}P_{\ell_2,o}) \cdot (g_{\ell_3,\ell_4}-\ol{P_{\ell_3,o}}P_{\ell_4,o}) \label{eq:ijs td psi quartic 2}
        \\
        &-{\textstyle \int} \calH(\td \psi \ol{P_{\ell_1,o}}P_{\ell_2,o}) \cdot  \ol{P_{\ell_3,o}}P_{\ell_4,o}. \label{eq:ijs td psi quartic 3}
    \end{align}
	Using \eqref{eq:CommuteHilbert}, we have
	\begin{align*}
		\eqref{eq:ijs td psi quartic 1}
		=&
		{\textstyle \int} y\td \psi (g_{\ell_1,\ell_2}\!-\!\ol{P_{\ell_1,o}}P_{\ell_2,o}) \calH(\tfrac{y}{1+y^2}g_{\ell_3,\ell_4})
		+{\textstyle \int} \td \psi (g_{\ell_1,\ell_2}\!-\!\ol{P_{\ell_1,o}}P_{\ell_2,o}) \calH(\tfrac{1}{1+y^2}g_{\ell_3,\ell_4})
		\\
		&-
		\tfrac{1}{\pi}{\textstyle \int} \td \psi (g_{\ell_1,\ell_2}-\ol{P_{\ell_1,o}}P_{\ell_2,o}) {\textstyle \int} \tfrac{y}{1+y^2}g_{\ell_3,\ell_4}.
	\end{align*}
	By \eqref{eq:ij radiation g odd} and \eqref{eq:ij rad g H1}, we have
	\begin{align}
		{\textstyle \int} \tfrac{y}{1+y^2}g_{i,j}
		=O(\lambda^{i+j+1}),& \quad 
        {\textstyle \int} \td \psi(g_{i,j}-\ol{P_{i,o}}P_{j,o}) ={\textstyle \int} \td \psi g_{i,j}
		=O(\lambda^{i+j+2}),\nonumber\\
		\|Qg_{i,j}\|_{L^2}\lesssim \lambda^{i+j+1},&\quad 
		\|Q(g_{i,j}-\ol{P_{i,o}}P_{j,o})\|_{L^2}\lesssim \lambda^{i+j+1+\frac{1}{2}}.
        \label{eq:ij g estimate}
	\end{align}
	Hence, we derive
    \begin{align}
        \eqref{eq:ijs td psi quartic 1}=O(\lambda^{\ell_1+\ell_2+\ell_3+\ell_4+2+\frac12}). \label{eq:ijs psi to yQ square claim goal 1}
    \end{align}

    For \eqref{eq:ijs td psi quartic 2}, since $\td \psi y^2Q^2$ is odd, we get $\calH(\td \psi \ol{P_{\ell_1,o}}P_{\ell_2,o}) = \ol{\bm{\beta}_{\ell_1}}\bm{\beta}_{\ell_2} O(Q^2)$. By this and \eqref{eq:ij g estimate}, we have
    %\begin{align*}
        %\calH(\td \psi \ol{P_{\ell_1,o}}P_{\ell_2,o}) = \ol{\bm{\beta}_{\ell_1}}\bm{\beta}_{\ell_2}\cdot O(Q^2). 
    %\end{align*}
    %Thus, using this and \eqref{eq:ij g estimate},
    \begin{align}
        \eqref{eq:ijs td psi quartic 2}
        \lesssim \beta_{\ell_1}\beta_{\ell_2} \lambda^{\ell_3+\ell_4+1+\frac12}
        \lesssim\lambda^{\ell_1+\ell_2+\ell_3+\ell_4+2+\frac12}. \label{eq:ijs psi to yQ square claim goal 2}
    \end{align}

    For \eqref{eq:ijs td psi quartic 3}, thanks to \eqref{eq:CommuteHilbert} and the oddness of $\td \psi$, we compute
	\begin{align}
		\eqref{eq:ijs td psi quartic 3}&=
        -\ol{\bm{\beta}_{\ell_1}}\bm{\beta}_{\ell_2}\ol{\bm{\beta}_{\ell_3}}\bm{\beta}_{\ell_4}
        {\textstyle \int} \calH(\td \psi y^3 Q^2) \cdot yQ^2 \nonumber\\
        &=
        -\ol{\bm{\beta}_{\ell_1}}\bm{\beta}_{\ell_2}\ol{\bm{\beta}_{\ell_3}}\bm{\beta}_{\ell_4}
        {\textstyle \int} \td \psi y^3 Q^4= \ol{\bm{\beta}_{\ell_1}}\bm{\beta}_{\ell_2}\ol{\bm{\beta}_{\ell_3}}\bm{\beta}_{\ell_4}
        {\textstyle \int} \tfrac{1}{2}y^4 Q^6. \label{eq:ijs psi to yQ square claim goal 3}
	\end{align}
    Here, we used the following identity from \eqref{eq:profile cancel 1} and \eqref{eq:profile cancel 3} for $\varphi$, and \eqref{eq:profile cancel phi} for $\phi$, %we have
    \begin{align*}
        {\textstyle \int} \td \psi y^3 Q^4=
        {\textstyle \int} (\psi -\tfrac{1}{2}yQ^2)y^3 Q^4=
        -{\textstyle \int} \tfrac{1}{2}y^4 Q^6.%,
    \end{align*}
   % which yields
    %\begin{align}
     %   \eqref{eq:ijs td psi quartic 3}
      %  =
       % \ol{\bm{\beta}_{\ell_1}}\bm{\beta}_{\ell_2}\ol{\bm{\beta}_{\ell_3}}\bm{\beta}_{\ell_4}
        %{\textstyle \int} \tfrac{1}{2}y^4 Q^6. \label{eq:ijs psi to yQ square claim goal 3}
    %\end{align}
    
	Therefore, by \eqref{eq:ijs psi to yQ square claim goal 1}, \eqref{eq:ijs psi to yQ square claim goal 2}, and \eqref{eq:ijs psi to yQ square claim goal 3}, we conclude \eqref{eq:ijs psi to yQ square claim}. In particular, if at least one among $\ell_1,\ell_2,\ell_3,\ell_4$ is zero, then $|\ol{\bm{\beta}_{\ell_1}}\bm{\beta}_{\ell_2}\ol{\bm{\beta}_{\ell_3}}\bm{\beta}_{\ell_4}|\lesssim \lambda^{\ell_1+\ell_2+\ell_3+\ell_4+2+\frac12}$ since we have $|\bm{\beta}_0|=|\frac12 \nu_0|\lesssim \lambda$ instead of $|\bm\beta_j|\lesssim \lambda^{j+\frac12}$. Hence, we derive
	\begin{align}
		\eqref{eq:ijs 3-1 without prof}
		=\tfrac{1}{2}{\textstyle \int} yQ^2g_{\ell_1,\ell_2}\calH(g_{\ell_3,\ell_4})
		+O(\lambda^{\ell_1+\ell_2+\ell_3+\ell_4+2+\frac12}). \label{eq:ijs 3-1 goal}
	\end{align}

    Next, we estimate \eqref{eq:ijs 3-1-1 with prof}. By \eqref{eq:profile cancel psi}, we have ${\textstyle \int} \psi yQ^2  P_{\ell_1,\ell_2}^g=0$, 
    %\begin{align*}
     %   {\textstyle \int} \psi yQ^2  P_{\ell_1,\ell_2}^g=0,
    %\end{align*}
    which implies
    \begin{align*}
        \eqref{eq:ijs 3-1-1 with prof}=\ol{c_{\ell_3}}c_{\ell_4}{\textstyle \int} \psi yQ^2 (\check g_{\ell_1,\ell_2}+ P_{\ell_1,\ell_2}^g)=\ol{c_{\ell_3}}c_{\ell_4}{\textstyle \int} \psi yQ^2 \check g_{\ell_1,\ell_2}.
    \end{align*}
    Applying \eqref{eq:ij radiation} with $|\td{\psi}yQ^2|\lesssim Q^4$, we arrive at
    \begin{align*}
        {\textstyle \int} \psi yQ^2 \check g_{\ell_1,\ell_2}=
        \tfrac{1}{2}{\textstyle \int} y^2Q^4 \check g_{\ell_1,\ell_2}+
        {\textstyle \int} \td{\psi} yQ^2 \check g_{\ell_1,\ell_2}
        =\tfrac{1}{2}{\textstyle \int} y^2Q^4 \check g_{\ell,\ell^\prime}+O(\lambda^{\ell_1+\ell_2+2+\frac12}).
    \end{align*}
    From \eqref{eq:c j estimate}, we deduce
    \begin{align}
        \eqref{eq:ijs 3-1-1 with prof}
        =\tfrac{1}{2}\ol{c_{\ell_3}}c_{\ell_4}{\textstyle \int} y^2Q^4 \check g_{\ell,\ell^\prime}+O(\lambda^{\ell_1+\ell_2+\ell_3+\ell_4+2+\frac12}). \label{eq:ijs 3-1-1 goal}
    \end{align}

    Finally, we consider \eqref{eq:ijs 3-1-2 with prof}. Thanks to \eqref{eq:profile cancel psi}, \eqref{eq:ij radiation g odd} with $|\psi Q^2|\lesssim Q^3$, \eqref{eq:c j estimate}, and \eqref{eq:beta j estimate}, we obtain
    \begin{align}
        \eqref{eq:ijs 3-1-2 with prof}
        =-(\ol{c_{\ell_3}}\bm{\beta}_{\ell_4}+\ol{\bm{\beta}_{\ell_3}}c_{\ell_4}){\textstyle \int} \psi Q^2 g_{\ell_1,\ell_2}
        =O(\lambda^{\ell_1+\ell_2+\ell_3+\ell_4+2+\frac12}). \label{eq:ijs 3-1-2 goal}
    \end{align}

    Collecting \eqref{eq:ijs 3-1 goal}, \eqref{eq:ijs 3-1-1 goal}, \eqref{eq:ijs 3-1-2 goal}, and $\psi \in \{\varphi,-\frac{1}{2}\phi\}$, we obtain \eqref{eq:ijs quartic nonlin goal varphi} and \eqref{eq:ijs quartic nonlin goal phi}. %the fact that $\psi$ is either $\varphi$ or $-\frac{1}{2}\phi$, we conclude \eqref{eq:ijs quartic nonlin goal varphi} and \eqref{eq:ijs quartic nonlin goal phi}.
    
    %\vspace{5pt}
     \underline{\textbf{Step 2.}} We prove \eqref{eq:ijs quartic commute}.
    By \eqref{eq:CommuteHilbert} with $(1+y^2)\langle y\rangle^{-2}=1$ and \eqref{eq:ij g estimate}, we have
    \begin{align*}
        {\textstyle \int} \tfrac{y}{1+y^2}  g_{\ell_1,\ell_2} \calH(g_{\ell_3,\ell_4})
        =&
        {\textstyle \int}  g_{\ell_1,\ell_2} \calH(\tfrac{y}{1+y^2}g_{\ell_3,\ell_4})
        -\tfrac{1}{\pi}{\textstyle \int} \tfrac{y}{1+y^2}  g_{\ell_1,\ell_2} {\textstyle \int} \tfrac{y}{1+y^2}g_{\ell_3,\ell_4}
        \\
        &-
        {\textstyle \int} \tfrac{1}{1+y^2}  g_{\ell_1,\ell_2} \calH(\tfrac{y}{1+y^2}g_{\ell_3,\ell_4})
        +
        {\textstyle \int} \tfrac{y}{1+y^2}  g_{\ell_1,\ell_2} \calH(\tfrac{1}{1+y^2}g_{\ell_3,\ell_4})
        \\
        =&
        {\textstyle \int}  g_{\ell_1,\ell_2} \calH(\tfrac{y}{1+y^2}g_{\ell_3,\ell_4})
        -\tfrac{1}{\pi}\mathfrak{c}_{\ell_1,\ell_2}\mathfrak{c}_{\ell_3,\ell_4}
        +O(\lambda^{\ell_1+\ell_2+\ell_3+\ell_4+2+\frac12}).
    \end{align*}
    That is, we conclude \eqref{eq:ijs quartic commute}.%,
    %\begin{align*}
     %   {\textstyle \int} \tfrac{y}{1+y^2}   g_{\ell_1,\ell_2}\! \calH(g_{\ell_3,\ell_4}\!)
      %  +{\textstyle \int} \tfrac{y}{1+y^2} g_{\ell_3,\ell_4} \!\calH(g_{\ell_1,\ell_2}\!)
       % =-\tfrac{1}{\pi}\mathfrak{c}_{\ell_1,\ell_2}\mathfrak{c}_{\ell_3,\ell_4}
        %+O(\lambda^{\ell_1+\ell_2+\ell_3+\ell_4+2+\frac12}).
    %\end{align*}
    
     %\vspace{5pt}
     \underline{\textbf{Step 3.}} We prove \eqref{eq:ijs F nonlin goal}.
    First, by \eqref{eq:profile cancel phi} and \eqref{eq:ij radiation}, we obtain
    \begin{align}
        {\textstyle \int} \Phi \ol{w_i} w_j
        ={\textstyle \int} \Phi \check g_{i,j}
        =\tfrac{1}{2}{\textstyle \int} y^2Q^4 \check g_{i,j}
        +O(\lambda^{i+j+2+\frac12}). \label{eq:g check ij convert}
    \end{align}
    Thus, it suffices to consider ${\textstyle \int} \Phi \ol{w_i} w_j$ instead of $\frac{1}{2}{\textstyle \int} y^2Q^4 \check g_{i,j}$. We have
    \begin{align}
        {\textstyle \int} \Phi \ol{w_i} w_j
        ={\textstyle \int} (\partial_y \phi) \ol{w_i} w_j
        =&
        -{\textstyle \int} \phi \ol{ w_{i+1}} w_j
        -
        {\textstyle \int} \phi \ol{ w_i} w_{j+1} \label{eq:ijs F}
        \\
        &+\tfrac{1}{2}{\textstyle \int} \phi \ol{w} w_j \calH(\ol{w_{i}}w) 
        +\tfrac{1}{2}{\textstyle \int} \phi \ol{w_{i}}w \calH(\ol{w} w_j ). \label{eq:ijs F nonlin}
    \end{align}
    We remark that $ g_{i,j}=\ol{w_i}w_j-\ol{P_{i,e}}P_{j,e}-\ol{P_{i,o}}P_{j,e}-\ol{P_{i,e}}P_{j,o}$. 
    %\begin{align*}
     %   g_{i,j}=\ol{w_i}w_j-\ol{P_{i,e}}P_{j,e}-\ol{P_{i,o}}P_{j,e}-\ol{P_{i,e}}P_{j,o},
    %\end{align*}
    By \eqref{eq:profile cancel phi}, we have
    \begin{align*}
        \eqref{eq:ijs F}
        =
        -{\textstyle \int} \phi g_{i+1,j}
        -
        {\textstyle \int} \phi  g_{i,j+1}
        =
        {\textstyle \int} yQ^2 g_{i+1,j}
        +
        {\textstyle \int} yQ^2  g_{i,j+1}
        +O(\lambda^{i+j+2+\frac12}).
    \end{align*}
    Here, we used \eqref{eq:ij radiation g odd} for the second equality.
    Also, this implies that
    \begin{align*}
        {\textstyle \int} [-6yQ^4  +4 yQ^6 ]g_{i+1,j}&=O(\lambda^{i+j+2+\frac12}).
    \end{align*}
    Therefore, we have
    \begin{align*}
        {\textstyle \int} yQ^2  g_{i+1,j}
        &=
        {\textstyle \int} [yQ^2 -6yQ^4  +4 yQ^6 ] g_{i+1,j}
        +O(\lambda^{i+j+2+\frac12})
        \\
        &= 2\mathfrak{c}_{i+1,j}+O(\lambda^{i+j+2+\frac12}).
    \end{align*}
    By symmetry, we conclude
    \begin{align}
        \eqref{eq:ijs F}
        =
        2\mathfrak{c}_{i+1,j}+2\mathfrak{c}_{i,j+1}
        +O(\lambda^{i+j+2+\frac12}). \label{eq:ijs F goal}
    \end{align}
    For \eqref{eq:ijs F nonlin}, we apply \eqref{eq:ijs quartic nonlin goal phi}. We have 
    \begin{equation}
        \begin{aligned}
            \tfrac{1}{2}{\textstyle \int} \phi \ol{w} w_j \calH(\ol{w_{i}}w) 
            =&
            -{\textstyle \int} \tfrac{y}{1+y^2}  g_{0,j} \calH(g_{i,0})
            -\tfrac{1}{2}\ol{c_i}{\textstyle \int}  y^2Q^4 \check g_{0,j}+O(\lambda^{i+j+2+\frac12}),
            \\
            \tfrac{1}{2}{\textstyle \int} \phi \ol{w_{i}}w \calH(\ol{w} w_j) 
            =&
            -{\textstyle \int} \tfrac{y}{1+y^2}  g_{i,0} \calH(g_{0,j})
            -\tfrac{1}{2}{c_j}{\textstyle \int}  y^2Q^4 \check g_{i,0}+O(\lambda^{i+j+2+\frac12}).
        \end{aligned}\label{eq:ijs F nonlin aux}
    \end{equation} 
    By \eqref{eq:ijs F nonlin aux}, we obtain
    \begin{align*}
        \eqref{eq:ijs F nonlin}
        =&
        -\tfrac{1}{2}\ol{c_i}{\textstyle \int}  y^2Q^4 \check g_{0,j}
        -\tfrac{1}{2}{c_j}{\textstyle \int}  y^2Q^4 \check g_{i,0}
        \\
        &-{\textstyle \int} \tfrac{y}{1+y^2}  g_{0,j} \calH(g_{i,0})
        -{\textstyle \int} \tfrac{y}{1+y^2} g_{i,0} \calH(g_{0,j})
        +O( \lambda^{i+j+2+\frac12}).
    \end{align*}
    Thus, \eqref{eq:ijs quartic commute} implies
    \begin{equation}
        \begin{aligned}
            \eqref{eq:ijs F nonlin}
            =-\tfrac{1}{2}\ol{c_i}{\textstyle \int}  y^2Q^4 \check g_{0,j}
            -\tfrac{1}{2}{c_j}{\textstyle \int}  y^2Q^4 \check g_{i,0}
            +\tfrac{1}{\pi}\mathfrak{c}_{0,j}\mathfrak{c}_{i,0}
            +O(\lambda^{i+j+2+\frac12}).
        \end{aligned} \label{eq:ijs F nonlin goal pf}
    \end{equation}
    Therefore, combining \eqref{eq:g check ij convert}, \eqref{eq:ijs F goal}, and \eqref{eq:ijs F nonlin goal pf}, we conclude \eqref{eq:ijs F nonlin goal}.    
\end{proof}

With Lemma~\ref{lem:gij induction} in place, we can now analyze the nonlinear term in $(\frakc_{i,j})_s$ appearing in Proposition~\ref{prop:radiation estimate 2}.  
Combining the three identities in the lemma, we first decompose the 4-linear nonlinear term into its radiative components (the 4-linear radiation part and the weighted quadratic part).  
We then use the symmetry identity \eqref{eq:ijs quartic commute} to collapse the 4-linear radiation term into a product of quadratic radiation quantities, and finally apply \eqref{eq:ijs F nonlin goal} to convert the weighted quadratic contribution into $\frakc_{i,j}$.  
These reductions provide exactly the structure required to obtain the closed ODE for $(\frakc_{i,j})_s$ in Proposition~\ref{prop:radiation estimate 2}.

\begin{proof}[Proof of Proposition~\ref{prop:radiation estimate 2}]
	From \eqref{eq:wj equ}, we have
    \begin{equation}
        \begin{aligned}
            (\partial_s&-\tfrac{\lambda_s}{\lambda}\Lambda_{-j} +\gamma_s i)w_j-\tfrac{x_s}{\lambda}w_{j+1}-iw_{j+2}
    		\\
    		=&-\tfrac{1}{2}\tfrac{x_s}{\lambda}w\calH(\ol{w}w_j)
    		-\tfrac{i}{2}(w_1\mathcal{H}\overline{w}-w\mathcal{H}\overline{w_1})w_{j}.
        \end{aligned} \label{eq:sec7prop wj equ}
    \end{equation}
	A straightforward computation based on the definition of $\frakc_{i,j}$ and \eqref{eq:sec7prop wj equ} yields the following identity:
	\begin{align}
		(\frakc_{i,j})_s
		=\,&\tfrac{\lambda_s}{\lambda}(i+j)\frakc_{i,j}
		+\tfrac{\lambda_s}{\lambda}{\textstyle \int} \varphi\ol{w_i} \Lambda w_j
		+\tfrac{\lambda_s}{\lambda}{\textstyle \int} \varphi \ol{\Lambda w_i}  w_j \label{eq: cijs 1}
		\\
		&+\tfrac{x_s}{\lambda} (\frakc_{i,j+1}+\frakc_{i+1,j}) +i\frakc_{i,j+2} -i\frakc_{i+2,j} \label{eq: cijs 2}
		\\
		&-\tfrac{1}{2}\tfrac{x_s}{\lambda}[{\textstyle \int} \varphi\ol{w_i} w\calH(\ol{w}w_j)
		+{\textstyle \int} \varphi\ol{w}w_j \calH(\ol{w_i} w) ] \label{eq: cijs 3}
		\\
		&-\tfrac{i}{2}[{\textstyle \int} \varphi\ol{w_i} w_1\calH(\ol{w}w_j)
		+{\textstyle \int} \varphi\ol{w}w_j \calH(\ol{w_i}w_1) ] \label{eq: cijs 4}
		\\
		&+\tfrac{i}{2}[{\textstyle \int}  \varphi \ol{w_i} w\calH(\ol{w_1}w_j)
		+{\textstyle \int} \varphi\ol{w_1}w_j \calH(\ol{w_i}w) ]. \label{eq: cijs 5}
	\end{align}
	For \eqref{eq: cijs 1}, we can check that
	\begin{align*}
		{\textstyle \int} \varphi\ol{w_i} \Lambda w_j
		+{\textstyle \int} \varphi\ol{\Lambda w_i}  w_j
		&=
		-{\textstyle \int}  \Lambda(\varphi\ol{w_i}) w_j
		+{\textstyle \int} \varphi\ol{\Lambda w_i}  w_j
		\\
		&=c_{i,j}
		-
		{\textstyle \int} (\tfrac{13}{2}yQ^4-14yQ^6+6yQ^8)\ol{w_i} w_j.
	\end{align*}
	Thus, by \eqref{eq:ij radiation g odd}, \eqref{eq:yQ2n integration}, and \eqref{eq:mod0 esti}, the first line \eqref{eq: cijs 1} becomes
	\begin{align}
		\eqref{eq: cijs 1}=\,&\tfrac{\lambda_s}{\lambda}(i+j+1)\frakc_{i,j} -
		\tfrac{1}{2}\tfrac{\lambda_s}{\lambda}(\ol{c_i}\bm{\beta}_j+\ol{\bm{\beta}_i}c_j){\textstyle \int}(\tfrac{13}{2}y^2Q^6-14y^2Q^8+6y^2Q^{10}) \nonumber
        \\
        &+O(\lambda^{i+j+3+\frac12}) %\nonumber
		=\tfrac{\lambda_s}{\lambda}(i+j+1)\frakc_{i,j}+O(\lambda^{i+j+3+\frac12}). \label{eq: cijs 1 goal}
	\end{align}
	Now, we consider $\eqref{eq: cijs 3}$. Thanks to \eqref{eq:ijs quartic nonlin goal varphi}, we derive
    \begin{align*}
        {\textstyle \int} \varphi\ol{w_i} w\calH(\ol{w}w_j)&=
        {\textstyle\int}\tfrac{y}{1+y^2}g_{i,0}\calH(g_{0,j})
        +\tfrac{1}{2}c_{j} {\textstyle\int}y^2Q^4 \check g_{i,0}
        +O(\lambda^{i+j+2+\frac12}),
        \\
        {\textstyle \int} \varphi\ol{w}w_j \calH(\ol{w_i} w)&=
        {\textstyle\int}\tfrac{y}{1+y^2}g_{0,j}\calH(g_{i,0})
        +\tfrac{1}{2}\ol{c_{i}} {\textstyle\int}y^2Q^4 \check g_{0,j}
        +O(\lambda^{i+j+2+\frac12}).
    \end{align*}
    Thus, applying \eqref{eq:ijs quartic commute} and \eqref{eq:mod0 esti}, we have
    \begin{align}
        \eqref{eq: cijs 3}=
        -\tfrac{1}{2}\nu_0[
        \tfrac{1}{2}c_{j} {\textstyle\int}y^2Q^4 \check g_{i,0}
        +\tfrac{1}{2}\ol{c_{i}} {\textstyle\int}y^2Q^4 \check g_{0,j}
        -\tfrac{1}{\pi}\mathfrak{c}_{i,0}\mathfrak{c}_{0,j}
        ]+O(\lambda^{i+j+3+\frac12}). \label{eq: cijs 3-1}
    \end{align}
    Before applying \eqref{eq:ijs F nonlin goal}, for notational convenience, we define
    \begin{equation}
        \mathfrak{g}_{i,j}\coloneqq \mathfrak{c}_{i+1,j}+\mathfrak{c}_{i,j+1}
                +\tfrac{1}{2\pi}\mathfrak{c}_{0,j}\mathfrak{c}_{i,0}. \label{eq:frakg def}
    \end{equation}
	By \eqref{eq:ijs F nonlin goal} with either $j=0$ or $i=0$ and $c_0=1$, we obtain
    \begin{align}
        \tfrac{1}{2}{\textstyle\int} y^2Q^4\check g_{i,0}&=
        -\tfrac{1}{2}\ol{c_i}{\textstyle \int}  y^2Q^4 \check g_{0,0}
        +2\mathfrak{g}_{i,0}
        +O(\lambda^{i+2+\frac12}), \label{eq: cijs 3-2} 
        \\
        \tfrac{1}{2}{\textstyle\int} y^2Q^4\check g_{0,j}&=
        -\tfrac{1}{2}{c_j}{\textstyle \int}  y^2Q^4 \check g_{0,0}
        +2\mathfrak{g}_{0,j}
        +O(\lambda^{j+2+\frac12}). \label{eq: cijs 3-3}
    \end{align}
    For ${\textstyle \int}  y^2Q^4 \check g_{0,0} $, taking $i=j=0$ in \eqref{eq:ijs F nonlin goal}, we derive
    \begin{align}
        {\textstyle \int}  y^2Q^4 \check g_{0,0}
        =\tfrac{4}{3}\mathfrak{g}_{0,0}. \label{eq:00s F}
    \end{align}
    Inserting \eqref{eq: cijs 3-2}, \eqref{eq: cijs 3-3}, and \eqref{eq:00s F} into \eqref{eq: cijs 3-1}, we deduce
    \begin{equation}
        \begin{aligned}
            \eqref{eq: cijs 3}=
            -\tfrac{1}{2}\nu_0\big[
            -\tfrac{4}{3}\ol{c_i}c_j\mathfrak{g}_{0,0}
            +2c_j\mathfrak{g}_{i,0}
            +2\ol{c_i}\mathfrak{g}_{0,j}
            -\tfrac{1}{\pi}\mathfrak{c}_{i,0}\mathfrak{c}_{0,j}
            \big]+O(\lambda^{i+j+3+\frac12}).
        \end{aligned} \label{eq: cijs 3 goal}
    \end{equation}
	Similarly to the first step of \eqref{eq: cijs 3}, by \eqref{eq:ijs quartic nonlin goal varphi}, each part of \eqref{eq: cijs 4} becomes
    \begin{equation}
        \begin{aligned}
            {\textstyle \int} \varphi\ol{w_i} w_1\calH(\ol{w}w_j)&=
            {\textstyle\int}\tfrac{y}{1+y^2}g_{i,1}\calH(g_{0,j})
            +\tfrac{1}{2}c_{j} {\textstyle\int}y^2Q^4 \check g_{i,1}
            +O(\lambda^{i+j+2+\frac12}),
            \\
            {\textstyle \int} \varphi\ol{w}w_j \calH(\ol{w_i} w_1)&=
            {\textstyle\int}\tfrac{y}{1+y^2}g_{0,j}\calH(g_{i,1})
            +\tfrac{1}{2}\ol{c_{i}} c_1{\textstyle\int}y^2Q^4 \check g_{0,j}
            +O(\lambda^{i+j+2+\frac12}).
        \end{aligned} \label{eq: cijs 4-1}
    \end{equation}
    In this case, using \eqref{eq:ijs F nonlin goal}, \eqref{eq: cijs 3-2}, \eqref{eq: cijs 3-3} with $j=1$, and \eqref{eq:00s F}, we have
    \begin{equation}
        \begin{aligned}
            \tfrac{1}{2}{\textstyle\int} y^2Q^4\check g_{i,1} = &
            -\tfrac{1}{2}\ol{c_i}{\textstyle \int}  y^2Q^4 \check g_{0,1}
            -\tfrac{1}{2}{c_1}{\textstyle \int}  y^2Q^4 \check g_{i,0} 
            +2\mathfrak{g}_{i,1}
            +O(\lambda^{i+3+\frac12}) \\
            = &
            \,\ol{c_i}[\tfrac{2}{3}{c_1}\mathfrak{g}_{0,0}
            -2\mathfrak{g}_{0,1}
            ]
            +c_1[\tfrac{2}{3}\ol{c_i}\mathfrak{g}_{0,0}
            -2\mathfrak{g}_{i,0}
            ]
            +2\mathfrak{g}_{i,1}
            +O(\lambda^{i+3+\frac12}).
        \end{aligned} \label{eq: cijs 4-2}
    \end{equation}
	Thus, by \eqref{eq:ijs quartic commute}, \eqref{eq: cijs 4-1}, \eqref{eq: cijs 4-2}, and \eqref{eq: cijs 3-3} with \eqref{eq:00s F}, we have
    \begin{equation}
        \begin{aligned}
            \eqref{eq: cijs 4}=-\tfrac{i}{2}
            \big[&-\tfrac{1}{\pi}\frakc_{i,1}\frakc_{0,j}
            +\tfrac{1}{\pi}\ol{c_{i}} c_1\mathfrak{c}_{0,j}\mathfrak{c}_{0,0}
            +\tfrac{2}{3}\ol{c_i}c_jc_1\mathfrak{g}_{0,0}
            -2\ol{c_i}c_j\mathfrak{g}_{0,1}
            \\
            &-2c_j{c_1}\mathfrak{g}_{i,0}
            +2c_j\mathfrak{g}_{i,1}
            +2\ol{c_{i}} c_1\mathfrak{g}_{0,j}
            \big]+O(\lambda^{i+j+3+\frac12}).
        \end{aligned} \label{eq: cijs 4 goal}
    \end{equation}
    By symmetry, exchanging $i$ and $j$ in \eqref{eq: cijs 4 goal} and taking the complex conjugate yields
    \begin{equation}
        \begin{aligned}
            \eqref{eq: cijs 5}=\tfrac{i}{2}
            \big[&-\tfrac{1}{\pi}\frakc_{i,0}\frakc_{1,j}
            +\tfrac{1}{\pi}\ol{c_{1}}c_{j}\mathfrak{c}_{i,0}\mathfrak{c}_{0,0}
            +\tfrac{2}{3}\ol{c_i}c_j\ol{c_1}\mathfrak{g}_{0,0}
            -2\ol{c_i}c_j\mathfrak{g}_{1,0}
            \\
            &-2\ol{c_i c_1}\mathfrak{g}_{0,j}
            +2\ol{c_i}\mathfrak{g}_{1,j}
            +2\ol{c_{1}} c_j\mathfrak{g}_{i,0}
            \big]+O(\lambda^{i+j+3+\frac12}).
        \end{aligned} \label{eq: cijs 5 goal}
    \end{equation}
    Now, collecting \eqref{eq: cijs 2}, \eqref{eq: cijs 1 goal}, \eqref{eq: cijs 3 goal}, \eqref{eq: cijs 4 goal}, and \eqref{eq: cijs 5 goal}, we arrive at
    \begin{align*}
        (\frakc_{i,j})_s=\tfrac{\lambda_s}{\lambda}(i+j+1)\frakc_{i,j}+i\frakc_{i,j+2} -i\frakc_{i+2,j} +\fraka_{i,j}+O(\lambda^{i+j+3+\frac12}),
    \end{align*}
    where
    \begin{equation}
        \begin{aligned}
            \fraka_{i,j}=\,
            &\nu_0(\frakc_{i,j+1}+\frakc_{i+1,j}) 
            -\tfrac{1}{2}\nu_0\big[
                -\tfrac{2}{3}\ol{c_i}c_j\mathfrak{g}_{0,0}
                +2c_j\mathfrak{g}_{i,0}
                +2\ol{c_i}\mathfrak{g}_{0,j}
                -\tfrac{1}{\pi}\mathfrak{c}_{i,0}\mathfrak{c}_{0,j}
                \big]
            \\
            &-\tfrac{i}{2}
                \big[-\tfrac{1}{\pi}\frakc_{i,1}\frakc_{0,j}
                +\tfrac{1}{\pi}\ol{c_{i}} c_1\mathfrak{c}_{0,j}\mathfrak{c}_{0,0}
                +\tfrac{2}{3}\ol{c_i}c_jc_1\mathfrak{g}_{0,0}
                -2\ol{c_i}c_j\mathfrak{g}_{0,1}
                \\
                &\qquad
                -2c_j{c_1}\mathfrak{g}_{i,0}
                +2c_j\mathfrak{g}_{i,1}
                +2\ol{c_{i}} c_1\mathfrak{g}_{0,j}
                \big]
            \\
            &+\tfrac{i}{2}
                \big[-\tfrac{1}{\pi}\frakc_{i,0}\frakc_{1,j}
                +\tfrac{1}{\pi}\ol{c_{1}}c_{j}\mathfrak{c}_{i,0}\mathfrak{c}_{0,0}
                +\tfrac{2}{3}\ol{c_i}c_j\ol{c_1}\mathfrak{g}_{0,0}
                -2\ol{c_i}c_j\mathfrak{g}_{1,0}
                \\
                &\qquad
                -2\ol{c_i c_1}\mathfrak{g}_{0,j}
                +2\ol{c_i}\mathfrak{g}_{1,j}
                +2\ol{c_{1}} c_j\mathfrak{g}_{i,0}
                \big].
        \end{aligned}\label{eq:Aij def}
    \end{equation}
    We recall the definition of $\mathfrak{g}_{i,j}$ \eqref{eq:frakg def}. 
    We can observe $|\mathfrak{g}_{i,j}|\lesssim \lambda^{i+j+2}$ and $|\fraka_{i,j}|\lesssim \lambda^{i+j+3}$ from \eqref{eq:c j estimate} and \eqref{eq:frakc esti}.
\end{proof}

\section{Classification of single-bubble blow-up}\label{sec:classification}

In this section, we complete the proof of Theorem~\ref{thm:classfy CMDNLS} and Theorem~\ref{thm:classfy GCM} by providing a rigid classification of the blow-up rates.
Building on the estimates established in Sections~\ref{sec:mod esti} and \ref{sec:rad esti}, we analyze the resulting dynamical system to determine the precise asymptotic behavior of the solution.

Our argument proceeds in three steps.
First, we reformulate the dynamics in terms of \emph{normalized parameters} defined in \eqref{eq:def normalized parameters} below. By doing this, we obtain a dynamical system where all the variables become of size $O(1)$, while their evolution equations satisfy small error estimates controlled by $\lambda^{\frac12}$ (hence being uniform in $j$); see Lemma~\ref{lem:normalized estimate}.

Next, to capture the precise dynamics of the scaling parameter $\lambda(t)$, we refine the modulation system by correcting the influence of radiation. Unlike the radial case (cf. the heuristic discussion in Section~\ref{sec:radial}) where the dynamics are governed solely by the modulation parameters $\bm{\beta}_{2k-1}$, the general setting involves nontrivial radiation effects mediated by $c_j$ and the quadratic radiation interactions $\frakc_{i,j}$.
To incorporate these effects coherently into the modulation dynamics, we construct a sequence $\{\Omega_k\}_{k\ge1}$ that successively corrects the modulation variables by the relevant radiation components.
Each $\Omega_k$ thus represents a renormalized modulation quantity that satisfies a recursive evolution law; see Lemma~\ref{lem:Omega ode system}.

Finally, in Lemma~\ref{lem:blowup rate classify}, we analyze the limiting behavior of $\{\Omega_k\}$ as $t \to T$. We show that the limit $\Omega_k^*$ exists and the first $k$ with $\Omega_k^* \neq 0$ determines the rate $\lambda(t) \sim (T-t)^{2k}$, concluding (essentially) the proof of the main theorems.

We define the normalized parameters
\begin{equation}
    \begin{alignedat}{4}
        Z_0&\coloneqq \lambda^{\frac12}e^{i\gamma},
        \\
        Z_j&\coloneqq \frac{\bm{\beta}_j}{\lambda^{j+\frac12}}e^{i\gamma},  &\quad&(1\leq j\leq 2L),&
        \td{Z}_{2L-1}&\coloneqq  \frac{\td{\bm{\beta}}_{2L-1}}{\lambda^{2L-\frac12}}e^{i\gamma},
        \\
        C_j&\coloneqq \frac{c_j}{\lambda^{j}},   &&(1\leq j \leq 2L+1), &\mathfrak{C}_{i,j}&\coloneqq  \frac{\frakc_{i,j}}{\lambda^{i+j+1}}, &&(0\leq i, j \leq 2L),
        \\
        A_j&\coloneqq \frac{a_j}{\lambda^{j+2}},   &&(1\leq j \leq 2L), &\frakA_{i,j}&\coloneqq \frac{\fraka_{i,j}}{\lambda^{i+j+3}}, &&(0\leq i, j \leq 2L-1).
    \end{alignedat}  \label{eq:def normalized parameters}
\end{equation}
Here, $\td{\bm{\beta}}_{2L-1}$ is given by \eqref{eq:def refine mod}.
From the definitions of $a_j$ \eqref{eq:aj def} and $\fraka_{i,j}$ \eqref{eq:Aij def}, one sees that $A_j$ and $\frakA_{i,j}$ can be expressed as finite linear combinations of $C_m$, $\mathfrak{C}_{n,m}$, and their complex conjugates.

We rewrite previously obtained modulation and radiation estimates in these new variables:
\begin{lem}\label{lem:normalized estimate} The following normalized estimates hold true:
    \begin{enumerate}
        \item (Bounds for normalized parameters) For the indices specified in the definitions \eqref{eq:def normalized parameters}, we have
        \begin{align}
            |Z_j|+ |C_j|+|\mathfrak{C}_{i,j}|+|A_j|+|\frakA_{i,j}|\lesssim 1. \label{eq:A frakA esti}
        \end{align}
        \item (Modulation estimates) We have
        \begin{align}
            |(Z_0)_t+iZ_1|\lesssim \lambda^{\frac12}, \label{eq:mod0 normalize}
        \end{align}
        and, for $L\geq 2$ and $1\leq j\leq 2L-2$,
        \begin{align}
            |( Z_j)_t-i Z_{j+2}-2\Im (C_1) Z_{j+1}+iC_j Z_{2}|\lesssim| Z_1|+\lambda^{\frac12} . \label{eq:modj normalize}
        \end{align}
        Moreover, the translation parameter satisfies
        \begin{align}
            |x_t-2\Im(C_1)|\lesssim \lambda^{\frac32}. \label{eq:trans normalize}
        \end{align}
        
        \item (Radiation estimates) For $1\leq j\leq 2L-1$, we have
        \begin{align}
            |(C_j)_t| \lesssim 1. \label{eq:cj mod aux normalize}
        \end{align}
        Moreover, when $L\ge2$ and $1\le i,j\le 2L-2$, we have
        \begin{align}
            |(C_j)_t-iC_{j+2}-A_{j}|
            & \lesssim \lambda^{\frac{1}{2}}, \label{eq:cj mod normalize} \\
            |(\mathfrak{C}_{i,j})_t+i\frakC_{i,j+2} -i\frakC_{i+2,j}- \frakA_{i,j}| & \lesssim \lambda^{\frac12}. \label{eq:frakcij mod normalize}
        \end{align}

        \item (Refined modulation estimates) We have
        \begin{equation}
            \begin{aligned}
                |Z_{2L-1}-\td{Z}_{2L-1}|&\lesssim (T-t)^{\frac34},\quad |\td{Z}_{2L-1}|\lesssim 1,
                \\
                |(\td{Z}_{2L-1})_t|&\lesssim  (T-t)^{-\frac14}+ \lambda(T-t)^{-\frac32}.
            \end{aligned}\label{eq:refine Z mod normalize}
        \end{equation}
        
    \end{enumerate}
\end{lem}
\begin{proof}
    For \eqref{eq:A frakA esti}, we immediately obtain this by \eqref{eq:beta j estimate}, \eqref{eq:c j estimate}, \eqref{eq:frakc esti} with \eqref{eq:coercivity H1 epsj}, and the definitions of $a_j$ \eqref{eq:aj def} together with $\fraka_{i,j}$ \eqref{eq:Aij def}.
    
    We show \eqref{eq:mod0 normalize}. We observe that, from $\bm{\beta}_{1}=-\tfrac{1}{2}(ib_1+\eta_1)$,
    \begin{align*}
        ( Z_0)_t+i  Z_1
        =[((\lambda^{\frac12})_t+\tfrac{b_1}{2\lambda^{3/2}})
        +i(\lambda^{\frac12}\gamma_t-\tfrac{\eta_1}{2\lambda^{3/2}})]e^{i\gamma}.
    \end{align*}
    By \eqref{eq:mod0 esti} and $\partial_s=\lambda^2\partial_t$, we conclude \eqref{eq:mod0 normalize} by obtaining
    \begin{align*}
        |2(\lambda^{\frac12})_t+\tfrac{b_1}{\lambda^{3/2}}|
        +|2\lambda^{\frac12}\gamma_t-\tfrac{\eta_1}{\lambda^{3/2}}|
        \lesssim \lambda^{\frac12}.
    \end{align*}
    For \eqref{eq:modj normalize}, we observe that
    \begin{align*}
        \partial_s( Z_j)
        =((\bm{\beta}_j)_s-\tfrac{\lambda_s}{\lambda}(j+\tfrac12)\bm{\beta}_j+i\gamma_s\bm{\beta}_j )e^{i\gamma}\lambda^{-(j+\frac12)}.
    \end{align*}
    Thus, we rewrite \eqref{eq:modj esti} as
    \begin{align*}
        | \lambda^{j+\frac12}e^{-i\gamma} \partial_s( Z_j) -i\bm{\beta}_{j+2}
        -\nu_0\bm{\beta}_{j+1} +ic_j\bm{\beta}_2 |
        \lesssim  \beta_1 \lambda^{j+1}+\lambda^{j+3}.
    \end{align*}
    Dividing both sides by $\lambda^{j+2+\frac12}$ and using $\partial_s=\lambda^2\partial_t$, we derive \eqref{eq:modj normalize},
    \begin{align*}
        |( Z_j)_t-iZ_{j+2}-2\Im (C_1) Z_{j+1}+iC_j Z_{2}|\lesssim |Z_1|+\lambda^{\frac12}.%|( Z_j)_t-iZ_{j+2}-\tfrac{\nu_0}{\lambda} Z_{j+1}+iC_j Z_{2}|\lesssim |Z_1|+\lambda^{\frac12}.
    \end{align*}
    Here, we note that $\beta_1 \lambda^{-\frac{3}{2}} =| Z_1|$ and $\tfrac{\nu_0}{\lambda}=\tfrac{2\Im (c_1)}{\lambda}=2\Im (C_1)$. %Since $c_1=\tfrac{1}{2}(i\nu_0+\mu_0)$, we have
    %\begin{align}
     %   \tfrac{\nu_0}{\lambda}=\tfrac{2\Im (c_1)}{\lambda}=2\Im (C_1), \label{eq:nu0 normal}
    %\end{align}
    %which yields \eqref{eq:modj normalize}. 
    In addition, thanks to \eqref{eq:mod0 esti refine} and $\tfrac{\nu_0}{\lambda}=2\Im (C_1)$, we obtain \eqref{eq:trans normalize}. %\eqref{eq:nu0 normal}, we obtain \eqref{eq:trans normalize}.
    For \eqref{eq:cj mod aux normalize}, \eqref{eq:cj mod normalize}, and \eqref{eq:frakcij mod normalize}, since
    \begin{align*}
        (C_j)_s=\tfrac{1}{\lambda^j}[(c_j)_s-\tfrac{\lambda_s}{\lambda} j c_j], \quad
        (\mathfrak{C}_{i,j})_s
        =\tfrac{1}{\lambda^{i+j+1}}[(\mathfrak{c}_{i,j})_s-\tfrac{\lambda_s}{\lambda}(i+j+1)\mathfrak{c}_{i,j}],
    \end{align*}
    \eqref{eq:cj mod aux}, \eqref{eq:cj mod}, and \eqref{eq:frakcij mod} directly imply \eqref{eq:cj mod aux normalize}, \eqref{eq:cj mod normalize}, and \eqref{eq:frakcij mod normalize} respectively.

    In addition, from \eqref{eq:refined mod gap} and \eqref{eq:refined mod esti}, we conclude \eqref{eq:refine Z mod normalize}.
\end{proof}

To simplify the presentation of parameters, we introduce the notion of \emph{admissible combinations} and assign an \emph{order} to each parameter.

\begin{definition}[Admissible combinations]\label{def:admissible}
    For $m \in \frac{1}{2}\bbN$, let $\mathscr P_m$ denote the set of all finite complex linear
    combinations of products of the parameters
    \[
    \{Z_j,\,\overline{Z_j}\}_{j\ge1},\quad
    \{C_{j},\,\overline{C_{j}}\}_{j\ge0},\quad
    \{\frakC_{i,j},\,\overline{\frakC_{i,j}}\}_{i,j\ge0},
    \]
    such that every product (monomial) has order $m$.
    Here the order is defined by
    \begin{equation}
        \begin{gathered}
            \mathrm{ord}(Z_{j})\coloneqq j+\tfrac12,\quad
            \mathrm{ord}(C_{j})\coloneqq j,\quad
            \mathrm{ord}(\frakC_{i,j})\coloneqq i+j+1, 
            \\
            \mathrm{ord}(\overline{z})\coloneqq \mathrm{ord}(z)
            ,\quad
            \mathrm{ord}(pq)\coloneqq \mathrm{ord}(p)+\mathrm{ord}(q).
        \end{gathered}\label{eq:order def}
    \end{equation}
    Hence, $\bm{\beta}_je^{i\gmm}$, $c_j$, and $\frakc_{i,j}$ can be written as $P\lambda^{\mathrm{ord}(P)}$, where $P \in \{Z_j,C_j,\frakC_{i,j}\}$ is the corresponding normalized parameter in \eqref{eq:def normalized parameters}. 
    
    We call $\mathscr P_m$ the class of \emph{admissible combinations of order $m$},
    and write $\mathcal P \in \mathscr P_m$ for a generic admissible element. Let $\mathscr P=\bigcup_{m\ge 1}\mathscr P_m$ denote the collection of all admissible combinations.
\end{definition}

\begin{rem}
Equivalently, $\mathscr P_m$ consists of complex-coefficient polynomials 
of order $m$ in the variables $\{\Re(Z_{j+1}),\Im(Z_{j+1}) ,\Re (C_j),\Im (C_j),\Re(\frakC_{i,j}),\Im(\frakC_{i,j})\}_{i,j\ge0}$.
%$\{\Re(Z_j),\Im(Z_j)\}_{j\ge1}$, $\{\Re (C_j),\Im (C_j)\}_{j\ge0}$, and $\{\Re(\frakC_{i,j}),\Im(\frakC_{i,j})\}_{i,j\ge0}$.
For example,
\begin{align*}
    \frakC_{1,1}+\Re(\frakC_{2,0})+C_1^{3}+\Im(Z_1)\ol{Z_1}  \in \mathscr P_3.
\end{align*}
Moreover, by the definition of $a_j$ \eqref{eq:aj def} and $\fraka_{i,j}$ \eqref{eq:Aij def}, we have
\begin{align}
    A_j\in \mathscr P_{j+2},\quad \frakA_{i,j}\in  \mathscr P_{i+j+3}. \label{eq:order a}
\end{align}
\end{rem}

To classify the blow-up dynamics, we introduce the sequence $\{\Omega_k\}_{k\ge1}$ of admissible combinations.
Each $\Omega_k$ represents a dynamical quantity that coherently incorporates all modulation parameters \emph{and} additional parameters $c_j$ and $\frakc_{i,j}$ from the radiation effects.
By introducing these \textit{corrected parameters} $\Omega_k$, we are able to reduce the coupled ODE system for the dynamical variables to a simpler recursive form, which ultimately enables a complete classification of the possible dynamics.

The construction of $\Omega_k$s is based on the normalized modulation and radiation estimates \eqref{eq:mod0 normalize}--\eqref{eq:frakcij mod normalize} 
and the order structure from Definition~\ref{def:admissible}.

\begin{definition}[Derivation $\mathbb D$ on admissible combinations]
    Define a linear derivation $\mathbb D:\mathscr P\to\mathscr P$ by setting, for the generators,
    \begin{equation*}
        \begin{aligned}
            &\mathbb D(Z_j)\coloneqq iZ_{j+2}+2\Im(C_1)Z_{j+1}-iC_jZ_2,\\
            &\begin{aligned}[t]
                \mathbb D(C_j)&\coloneqq iC_{j+2}+A_j,\\
                \mathbb D(\overline{u})&\coloneqq \overline{\mathbb D(u)},
            \end{aligned}
            \qquad
            \begin{aligned}[t]
                \mathbb D(\frakC_{i,j})&\coloneqq -i\frakC_{i,j+2}+i\frakC_{i+2,j}+\frakA_{i,j},\\
                \mathbb D(uv)&\coloneqq \mathbb D(u)v+u\mathbb D(v).
            \end{aligned}
        \end{aligned}
    \end{equation*}
    and extending to linear combinations by linearity. We note that $\mathbb D|_{\mathscr P_m} : \mathscr P_m \to \mathscr P_{m+2}$ for $1\le m\in \frac{1}{2}\bbN$ by \eqref{eq:order a}.
\end{definition}

Define $\{\Omega_k\}_{1\leq k\leq L}$ by
\begin{align*}
    \Omega_1 \coloneqq -iZ_1,\qquad 
    \Omega_{k} \coloneqq \mathbb D(\Omega_{k-1}).
\end{align*}
From \eqref{eq:order def}, one can check that
\begin{align*}
    \Omega_1 \in \mathscr P_{\frac32},\qquad \Omega_k \in  \mathscr P_{2k-\frac12}.
\end{align*}

\begin{lem}\label{lem:well def omega}
    The indices of the parameters appearing in $\Omega_k$ satisfy the following:
    \begin{enumerate}
        \item \label{enumi:Omega 1} The unique factor with the largest index is $-i^{k}Z_{2k-1}$, and all other monomials in $\Omega_k$ involve only parameters with indices strictly less than $2k-1$.
        
        \item \label{enumi:Omega 3} If any factor of the form $C_j$ or $\frakC_{i,j}$ appears in a monomial of $\Omega_k$, its indices are at most $2k-3$.
    \end{enumerate}
\end{lem}
\begin{proof}
    To show \eqref{enumi:Omega 1} and \eqref{enumi:Omega 3}, we use induction.
    For the base case $k=1$, these follow from the definition $\Omega_1=-iZ_1$. For $k=2$ and $k=3$, from the definition of $\mathbb D$, we can check \eqref{enumi:Omega 1} and \eqref{enumi:Omega 3} directly from the definition.
    
    Assume the claim holds for some $k$. By the definition of $\mathbb D$, the unique highest-index term $-i^{k}Z_{2k-1}$ in $\Omega_k$ produces $-i^{k+1}Z_{2k+1}$ together with additional factors including $C_{2k-1}$. This $C_{2k-1}$ has index $2k-1$, which is consistent with the bounds in \eqref{enumi:Omega 3} for $\Omega_{k+1}$.
    Similarly, for the factor $Z_{2k-2}$, the action of $\mathbb D$ does not generate any term of the form $Z_{2k+1}$, and the resulting factors remain consistent with the index bounds in \eqref{enumi:Omega 3} for $\Omega_{k+1}$.
    
    Since the indices can increase by at most two under the action of $\mathbb D$, the remaining factors $Z_j$, $C_j$, and $\frakC_{i,j}$ with indices bounded by $2k-3$ generate only terms with indices bounded by $2k-1$. 
    Therefore, in $\Omega_{k+1}=\mathbb D(\Omega_k)$, the highest-index term is uniquely $-i^{k+1}Z_{2k+1}$, and the properties \eqref{enumi:Omega 1} and \eqref{enumi:Omega 3} remain valid for $k+1$.
\end{proof}
In addition, we need to consider a refined version of $\Omega_L$, which contains $Z_{2L-1}$.
As discussed in Section~\ref{sec:refined mod esti}, the modulation estimate for $\bm{\beta}_{2L-1}$ itself is weak, and hence so is its normalized parameter $Z_{2L-1}$. 
Thus in $\Omega_L$ we would like to replace $Z_{2L-1}$ by its refined version to improve the modulation control; we define
\begin{align}
    \td{\Omega}_{L}\coloneqq \Omega_L+i^{L}Z_{2L-1}-i^{L}\td{Z}_{2L-1}. \label{eq:def Omega tilde}
\end{align} 
In view of the scale of $\td{Z}_{2L-1}$, the order of $\td{\Omega}_L$ is the same as that of $\Omega_L$.
We now record the time evolution properties of the sequence $\{\Omega_k\}_{1\le k\le L}$ and $\td{\Omega}_L$.

\begin{lem}\label{lem:Omega ode system}
    Let $L\geq 1$ and $1\leq k\leq L$. Consider $\{\Omega_k\}_{1\le k\le L}$ and $\td{\Omega}_L$. The following hold true:
    \begin{enumerate}
        \item (Boundedness of parameters) We have
        \begin{align}
            |\Omega_k| + |\td{\Omega}_{L}|\lesssim 1. \label{eq:Omega esti}
        \end{align}
            
        \item (Corrected modulation estimates) We have
        \begin{align}
            |(Z_{0})_t-\Omega_{1}|\lesssim \lambda^{\frac12}. \label{eq:mod0 omega}
        \end{align}
        Moreover, for $L\geq 2$ and $1\leq k\leq L-1$, we have
        \begin{align}
            |(\Omega_{k})_t-\Omega_{k+1}|\lesssim |\Omega_1|+\lambda^{\frac12}. \label{eq:modk omega}
        \end{align}

        \item (Proximity of refined parameter) We have
        \begin{align}
            |{\Omega}_{L}-\td{\Omega}_{L}|\lesssim (T-t)^{\frac34}, \label{eq:Omega refine esti}
        \end{align}

        \item (Refined parameter estimates) We have
        \begin{align}
            |(\td{\Omega}_{L})_t|\lesssim (T-t)^{-\frac14}+|Z_0|^2(T-t)^{-\frac32}. \label{eq:modk omega refine}
        \end{align}
        
    \end{enumerate}
\end{lem}
\begin{proof} 
    The definition of $\Omega_k$ together with \eqref{eq:A frakA esti} yields $|\Omega_k|\lesssim 1$. The boundedness of $\td{\Omega}_L$ follows from the refined estimate \eqref{eq:Omega refine esti}, whose proof will be given below. By $\Omega_1=-iZ_1$ and \eqref{eq:mod0 normalize}, we obtain \eqref{eq:mod0 omega}. Similarly, \eqref{eq:modj normalize}, \eqref{eq:cj mod normalize}, and \eqref{eq:frakcij mod normalize}, together with the definition of $\Omega_k$, imply \eqref{eq:modk omega}.
    
    Next, from \eqref{eq:refine Z mod normalize} together with the definition of $\td{\Omega}_L$, we derive \eqref{eq:Omega refine esti} and \eqref{eq:modk omega refine}, respectively. Thanks to \eqref{eq:Omega refine esti}, we also have $|\td{\Omega}_L|\lesssim 1$, which implies \eqref{eq:Omega esti}.
\end{proof}

Finally, the rigidity of $\Omega_k(t)$ follows by integrating the estimates for $\Omega_k$ established in the previous lemma.
\begin{lem}[Identification of quantized blow-up rates]\label{lem:blowup rate classify} 
    The limits $\Omega_k^*\coloneqq \lim_{t\to T}\Omega_k(t)$ exist for all $1\leq k\leq L$, and the following hold true:
    \begin{enumerate}
        \item (Quantized regime) Assume $\Omega_{1}^*=\Omega_{2}^*=\cdots=\Omega_{k-1}^*=0$ for some $1\leq k \leq L$.
        Then, we have 
        \begin{equation}
            \begin{aligned}
                Z_0(t)&=\frac{(-1)^k}{k!}(\Omega_{k}^*+o_{t\to T}(1))(T-t)^{k},
        		\\
        		\Omega_i(t)&=\frac{(-1)^{k-i}}{(k-i)!}(\Omega_{k}^*+o_{t\to T}(1))(T-t)^{k-i},\qquad 1\leq i \leq k.
            \end{aligned} \label{eq:classify z0}
        \end{equation}
        In particular, we have
        \begin{align}
            \lambda(t) = \frac{1}{(k!)^2}(|\Omega_k^*|^2+o_{t\to T}(1))(T-t)^{2k}. \label{eq:classify lambda}
        \end{align}

        \item (Exotic regime) If $\Omega_{k}^*=0$ for all $1\leq k\leq L$, then 
        \begin{align}
            \lambda(t)\lesssim (T-t)^{2L+\frac32}. \label{eq:classify exotic}
        \end{align}
    \end{enumerate}
\end{lem}
\begin{proof}
    Note that \eqref{eq:classify lambda} is a direct consequence of $\lmb = |Z_0|^2$ and \eqref{eq:classify z0}.

    We first establish the convergence of $\Omega_k$. By \eqref{eq:Omega esti} and \eqref{eq:modk omega}, for $k=1$ or $1\leq k \leq L-1$, we have $|(\Omega_k)_t|\lesssim 1$. Hence, $\Omega_k$ converges as $t\to T$. Similarly, $\td{\Omega}_{L}$ converges as $t\to T$ by \eqref{eq:modk omega refine} and \eqref{eq:lambda bdd one sol}. By \eqref{eq:Omega refine esti} and the convergence of $\td{\Omega}_{L}$, $\Omega_{L}$ also converges to the same limit.

    We show \eqref{eq:classify z0} for $k=1$. By \eqref{eq:mod0 omega} and $\lmb(t) \to 0$, we get $(Z_0)_t = \Omega_1^* + o_{t\to T}(1)$. Integrating this, we get $Z_0 (t) = -(\Omega_1^* + o_{t\to T}(1))(T-t)$ as desired.

    We show \eqref{eq:classify z0} for $2\leq k \leq L$. For a function $z(t)$ such that $\lim_{t\to T} z(t)$ exists, we introduce a notation 
    \begin{align*}
        \mathcal{I}z=\mathcal{I}z(t)\coloneqq -{\textstyle\int}_t^T z(\tau)d\tau,\qquad \|\mathcal{I}^j\|_{L^\infty\to L^\infty} \leq \tfrac{1}{j!}(T-t)^j.
    \end{align*}
    Using $\Omega_1^*=\cdots=\Omega_{k-1}^*=0$, we can rewrite \eqref{eq:mod0 omega} and \eqref{eq:modk omega} as 
    \begin{align}
        Z_0 & = \mathcal{I}\Omega_1 + \mathcal{I}r_0, & |r_0| & \aleq |Z_0|, \label{eq:Z_0 induction} \\
        \Omega_i & = \mathcal{I}\Omega_{i+1}+\mathcal{I}r_i, & |r_i| & \aleq|\Omega_1| + |Z_0|,\qquad 1\leq i \leq k-1. \label{eq:Omega induction}
    \end{align}
    Integrating \eqref{eq:Z_0 induction}, we have 
    \begin{align*}
        \sup_{[t,T)}|Z_0|\lesssim (T-t)\sup_{[t,T)}|\Omega_1|+(T-t)\sup_{[t,T)}|Z_0|,
    \end{align*}
    so for all $t$ close to $T$ we obtain  
    \begin{equation}
        \sup_{[t,T)}|Z_0|\lesssim (T-t)\sup_{[t,T)}|\Omega_1|. \label{eq:Z_0 bound}
    \end{equation}
    Iterating \eqref{eq:Omega induction}, we get 
    \begin{align}
        \Omega_i=\mathcal{I}^{k-i} \Omega_{k}+{\textstyle\sum_{\ell=i}^{k-1}} \mathcal{I}^{\ell} r_{\ell}\qquad \quad 1\leq i\leq k. \label{eq:Omega i j}
    \end{align}
    Substituting $i=1$ into the above and using \eqref{eq:Z_0 bound}, we get 
    \begin{equation*}
        \sup_{[t,T)}|\Omega_1| \aleq (T-t)^{k-1} \sup_{[t,T)}|\Omega_k| + (T-t) \sup_{[t,T)}|\Omega_1|,
    \end{equation*}
    so for all $t$ close to $T$ we obtain 
    \begin{equation}
        \sup_{[t,T)}|\Omega_1|\lesssim (T-t)^{k-1}\sup_{[t,T)}|\Omega_k|. \label{eq:Omega1 bound}
    \end{equation}
    Now, substituting the bounds \eqref{eq:Z_0 bound} and \eqref{eq:Omega1 bound} into \eqref{eq:Omega i j}, we get for $1\leq i\leq k$
    \begin{equation}
        \Omega_i = \mathcal{I}^{k-i} \Omega_{k} + {\textstyle\sum_{\ell=i}^{k-1}} O((T-t)^{\ell+k-1}) = \tfrac{(-1)^{k-i}}{(k-i)!}(\Omega_k^* + o_{t\to T}(1))(T-t)^{k-i}. \label{eq:8.25}
    \end{equation}
    Furthermore, substituting this (when $i=1$) and \eqref{eq:Z_0 bound} into  \eqref{eq:Z_0 induction}, we get 
    \begin{equation}
        Z_0 = \mathcal{I} \Omega_{1}+O((T-t)^{k+1}) = \tfrac{(-1)^k}{k!}(\Omega_k^* + o_{t\to T}(1))(T-t)^k. \label{eq:8.26}
    \end{equation}
    This completes the proof of \eqref{eq:classify z0}.
    
    Finally, we prove \eqref{eq:classify exotic}. By \eqref{eq:modk omega refine} and $|Z_0|\lesssim (T-t)^L$ (from \eqref{eq:classify z0}), $|(\td{\Omega}_L)_t| \aleq (T-t)^{-\frac14}$.
    Since $\lim_{t\to T}\td{\Omega}_L={\Omega}_{L}^*=0$, we get $|\td{\Omega}_L| \aleq (T-t)^{\frac34}$. By \eqref{eq:Omega refine esti}, 
    \begin{equation*}
        |\Omega_L| \aleq (T-t)^{\frac34}.
    \end{equation*}
    This combined with \eqref{eq:8.25} for $i=1$ and $k=L$ gives $|\Omega_1| \aleq (T-t)^{L-\frac14}$. Then, by \eqref{eq:8.26} for $k=L$, we get $|Z_0| \aleq (T-t)^{L+\frac34}$. Recalling $\lmb = |Z_0|^2$, we get \eqref{eq:classify exotic}.
\end{proof}

We are ready to finalize the proof of our main theorem.
\begin{proof}[Proof of Theorem~\ref{thm:classfy GCM} and Theorem~\ref{thm:classfy CMDNLS}]
    As mentioned in Section~\ref{sec:gague correspond}, Theorem~\ref{thm:classfy CMDNLS} is equivalent to Theorem~\ref{thm:classfy GCM} by Lemma~\ref{lem:gauge equiv}. Thus it suffices to prove Theorem~\ref{thm:classfy GCM}.

    Under Assumption~\ref{assum:one soliton GCM}, Proposition~\ref{prop:one soliton} allows us to identify the modulation parameters $(\lambda,\gamma,x)$ with the parameters $(\wt\lambda,\wt\gamma,\wt x)$ appearing in Proposition~\ref{prop:Decomposition}, without affecting any of the conclusions.
    Accordingly, we work with a solution $v$ admitting a decomposition of the form \eqref{eq:Decomposition and orthogonal}.
    For notational convenience, we drop the hat and denote the parameters again by $(\lambda,\gamma,x)$.
    
    By Lemma~\ref{lem:well def omega} and \eqref{eq:def Omega tilde}, for any single-bubble blow-up solution $v(t)$ to \eqref{CMdnls-gauged} with $v_0\in H^{2L+1}$, we can define $\{\Omega_k\}_{1\leq j\leq L}$ and $\td \Omega_L$ on $[T-t^*,T)$ where $T$ is the blow-up time and $t^*>0$ is a universal constant depending on $E_j$. 
    
    By Lemma~\ref{lem:Omega ode system} and Lemma~\ref{lem:blowup rate classify}, one of the following cases hold:
    \begin{enumerate}
        \item\label{eq:pf main quantize} $Z_0(t)=\frac{(-1)^k}{k!}({\Omega}_{k}^*+o_{t\to T}(1))(T-t)^{k}$ for some $1\leq k\leq L$, where ${\Omega}_{k}^*\neq0$.

        \item\label{eq:pf main exotic} $\lambda(t)\lesssim (T-t)^{2L+\frac32}$.
    \end{enumerate}
    In the quantized case \eqref{eq:pf main quantize}, since $Z_0(t)=\lambda(t)^{\frac12}e^{i\gamma(t)}$, we obtain
    \begin{align*}
        \lambda(t)=(\ell +o_{t\to T}(1))(T-t)^{2k},\quad \gamma(t)=\gamma^*+o_{t\to T}(1),
    \end{align*}
    where
    \begin{equation*}
        (\ell, \gamma^*) =
        \big(\tfrac{1}{(k!)^2}|\Omega_k^*|^2, \arg\big((-1)^k \Omega_k^*\big) \big).
    \end{equation*}

    It remains to find the asymptotics of $x(t)$. 
    By \eqref{eq:cj mod aux normalize} and \eqref{eq:A frakA esti}, $2\Im(C_1)$ converges to some $2\nu^*\in \bbR$.
    Integrating \eqref{eq:trans normalize}, there exists $x^*\in \bbR$ such that 
    \begin{equation*}
        x(t)=x^*+2\nu^*(T-t)+O((T-t)^2).
        \qedhere
    \end{equation*}
\end{proof}

\section{Proof of coercivity estimates}\label{sec:coer pf}
In this section, we prove the coercivity estimates Lemmas~\ref{lem:coercivity w1}, \ref{lem:coercivity estimate 2}, and \ref{lem:coercivity H3 epsj}.

\begin{rem}\label{rem:no circularity}
    In the proof of \eqref{eq:coercivity H2 check}, we use the coercivity of $w_1$ (the case $j=1$ of \eqref{eq:coercivity H1 epsj} in Lemma~\ref{lem:coercivity estimate 2}), namely
    \begin{align*}
        \|A_Q w_1\|_{L^2}\lesssim \lambda^2.
    \end{align*}
    The proof of that lemma is postponed to a later point in this section and is independent of the estimates established here, so there is no circularity.
\end{rem}

\begin{proof}[Proof of Lemma~\ref{lem:coercivity w1}]
    We basically follow the proof of \cite[Lemma~5.14]{KimKimKwon2024arxiv}. We use the coercivity $\|L_{Q}\widehat{\eps}\|_{L^{2}}\sim\|\widehat{\eps}\|_{\dot{\mathcal{H}}^{1}}$ and $\|\check{\eps}\|_{\dot{\mathcal{H}}^{2}}\sim\|\calL_2\check{\eps}\|_{L^{2}}$ (Lemma~\ref{lem:coercivity for linear}) with $\calL_2\check{\eps}=A_Q\td L_Q\check{\eps}$ without mentioning.

    \underline{\textbf{Step 1.}} We show \eqref{eq:coercivity H1} and \eqref{eq:coercivity L inf}. We note that $\|\widehat{\eps}\|_{\dot{\mathcal{H}}^{1}}$ directly follows from \eqref{eq:coercivity H1 from one sol}.
    Moreover, the energy conservation law implies $\|w_{1}\|_{L^{2}}\lesssim \lambda$ on $[0, T)$. Plus, \eqref{eq:coercivity L inf} is obtained by interpolating between the mass and energy conservation laws.

	\underline{\textbf{Step 2.}} We show \eqref{eq:coercivity H2 check}. Due to $A_{Q}(yQ)=0$, we have
	\begin{align}
		A_{Q} w_{1}=A_{Q}\widetilde{L}_{Q}\check{\eps}+A_{Q}\td N_{Q}(\widehat{\eps}). \label{eq:AQw1 w level decom}
	\end{align}
	We decompose $\td N_Q (\widehat{\eps})=\td N_1+ \td N_2$, where
    \begin{equation}
        \begin{aligned}
            \td N_1=\wt\eps\mathcal{H}(\Re(Q\wt\eps))+\tfrac{1}{2}\wt\eps\mathcal{H}(|\wt\eps|^{2}),
    		\quad \td N_2=\tfrac{1}{4}yQ\mathcal{H}(yQ^2|\wt\eps|^{2})
    		+\tfrac{1}{4}Q\mathcal{H}(Q^2|\wt\eps|^{2}).
        \end{aligned} \label{eq:N12 tilde def}
    \end{equation}
	For $\td N_1$, since $ \|A_{Q}\td N_1\|_{L^{2}}\sim\|\td N_1\|_{\dot{\mathcal{H}}^{1}}$, we have 
        \begin{align}\label{eq:td N1 goal}
            \|A_{Q}\td N_1\|_{L^{2}}
            \lesssim   \|\wt\eps\mathcal{H}(\Re(Q\wt\eps))\|_{\dot{\mathcal{H}}^{1}}  + \|\wt\eps\mathcal{H}(|\wt\eps|^{2})\|_{\dot{\mathcal{H}}^{1}} \lesssim\lambda^{2}+ \lambda^{\frac{1}{2}}\|A_{Q}\widetilde{L}_{Q}\eps\|_{L^{2}},
        \end{align}
    from the bounds obtained by Gagliardo--Nirenberg inequality and $\partial_{y}\Re(QP)=0$:
    \begin{align*}
      \|\wt\eps\mathcal{H}(\Re(Q\wt\eps))\|_{\dot{\mathcal{H}}^{1}}
        & \lesssim \|\widehat{\eps}\|_{\dot{\calH}^{1}}\|\mathcal{H}(\Re(Q\widehat{\eps}))\|_{L^{\infty}}+\|\widehat{\eps}\|_{L^{\infty}}\|\partial_{y}\Re(Q\widehat{\eps})\|_{L^{2}}
        \\
         & \lesssim\|\widehat{\eps}\|_{\dot{\mathcal{H}}^{1}}^{2}+\|\widehat{\eps}\|_{L^{\infty}}
         \|\check{\eps}\|_{\dot{\mathcal{H}}^{2}}
          \lesssim \lambda^{2}+ \lambda^{\frac{1}{2}}\|A_{Q}\widetilde{L}_{Q}\check{\eps}\|_{L^{2}},\\
          \|\wt\eps\mathcal{H}(|\wt\eps|^{2})\|_{\dot{\mathcal{H}}^{1}} & \lesssim \|\widehat{\eps}\|_{\dot{\mathcal{H}}^{1}}\|\mathcal{H}(|\widehat{\eps}|^{2})\|_{L^{\infty}}
        +\|\widehat{\eps}\|_{L^{\infty}}\|\partial_{y}(|\widehat{\eps}|^{2})\|_{L^{2}}
        \lesssim \lambda^{2}.
    \end{align*}
    
    For $\td N_2$, by \eqref{eq:CommuteHilbert} with $A_Q(yQ)=0$ and the definition of $A_Q$ with \eqref{eq:CommuteHilbertDerivative}, we derive
    \begin{align*}
    	2\sqrt{2}A_Q\td N_2=\sqrt{2}A_Q[Q^{-1}\mathcal{H}(Q^2|\wt\eps|^{2})]
        =\mathcal{H}(Q^2|\wt\eps|^{2})+\mathcal{H}[y\partial_y(Q^2|\wt\eps|^{2})]+\partial_y(Q^2|\wt\eps|^{2}).
    \end{align*}
    Hence, Gagliardo--Nirenberg inequality implies
    \begin{align}
        \|A_Q\td N_2\|_{L^2}&\lesssim \|Q^2|\wt\eps|^{2} \|_{L^2} \!+\! \|y\partial_y(Q^2|\wt\eps|^{2}) \|_{L^2}\! +\! \|\partial_y(Q^2|\wt\eps|^{2}) \|_{L^2} 	\lesssim \|\wt\eps\|_{\dot \calH^1}^2
            \lesssim  \lambda^2.\label{eq:td N2 goal}
    \end{align}
    Collecting \eqref{eq:td N1 goal} and \eqref{eq:td N2 goal}, we deduce 
    \begin{align}
        \|A_{Q}w_{1}-A_{Q}\widetilde{L}_{Q}\check{\eps}\|_{L^{2}} =  \|A_{Q}\td N_{Q}(\widehat{\eps})\|_{L^{2}}
        \lesssim
        \lambda^{2}+\lambda^{\frac{1}{2}}\|A_{Q}\widetilde{L}_{Q}\check{\eps}\|_{L^{2}}.\label{eq:5.3 proof coercivity 1}
    \end{align}
    By \eqref{eq:coercivity H1 epsj}\footnote{As noted in Remark~\ref{rem:no circularity}, this use of coercivity for $w_1$ is independent of the present estimates and hence involves no circularity.} and \eqref{eq:5.3 proof coercivity 1}, we arrive at 
    \begin{align*}
        (1-C_{M,E_1}\lambda^{\frac{1}{2}})\|A_{Q}\widetilde{L}_{Q}\check{\eps}\|_{L^{2}}
        \lesssim \lambda^{2}+\|A_{Q}w_{1}\|_{L^{2}}
        \lesssim \lambda^{2}.
    \end{align*}
    Since $\lambda \lesssim T-t$, after taking small $t^*$ depending on $M$ and $E_1$ if necessary, we have
    \begin{align*}
        \|A_{Q}\widetilde{L}_{Q}\check{\eps}\|_{L^{2}}\lesssim \lambda^{2} \quad \text{on} \quad [T-t^*, T).
    \end{align*}
    Therefore, we conclude \eqref{eq:coercivity H2 check} since $\|\check{\eps}\|_{\dot{\mathcal{H}}^{2}}\sim\|A_{Q}\widetilde{L}_{Q}\check{\eps}\|_{L^{2}}$.
    
    \underline{\textbf{Step 3.}} We show \eqref{eq:mu check estimate}.
    We first prove for $\check{\nu}_0$. Since $\wt\eps_o=\tfrac{i}{2}\check{\nu}_0 yQ+\check{\eps}_o$, we have
	\begin{align*}
		|\check{\nu}_0|\sim \|\chi_{10} \check{\nu}_0 yQ\|_{L^2}\lesssim \|\chi_{10}\wt\eps_o\|_{L^2}+\|\chi_{10}\check{\eps}_o \|_{L^2}
		\lesssim \lambda+\lambda^2\lesssim \lambda.
	\end{align*}
    For $\check{b}_1$ and $\check{\eta}_1$, using $\wt\eps_e=-\tfrac{i}{4}\check{b}_1 iy^2Q-\tfrac{1}{4}\check{\eta}_1(1+y^2)Q+\check{\eps}_e$, we obtain
    \begin{align*}
        (\check{b}_1^2+\check{\eta}_1^2)^{\frac12}\lambda^{-\frac32}\!
         \sim \!\|\chi_{\lambda^{-1}}[\check{b}_1 iy^2Q+\check{\eta}_1(1+y^2)Q]\|_{L^2}\lesssim 
        \|\chi_{\lambda^{-1}} \wt\eps_e\|_{L^2} +\|\chi_{\lambda^{-1}}\check{\eps}_e\|_{L^2}
        \lesssim 1.
    \end{align*}
    For $\check{\mu}_0$, by the definition of $\check{\mu}_0$ \eqref{eq:def mu}, we have
    \begin{align*}
        |\check{\mu}_0|
        \lesssim |{\textstyle \int} \Re (yQ^3\check{\eps})|+|{\textstyle \int} yQ^2|\wt\eps|^2|
        \lesssim \lambda +\lambda^{\frac32} \lesssim \lambda.
    \end{align*}

    \underline{\textbf{Step 4.}} We show \eqref{eq:coercivity H1 inf check}.
    %Using $\wt\eps=\check{T}+\check{\eps}$, 
    We define $ \mathring{\eps}\coloneqq \check{T}(1-\chi_{\lambda^{-1}})+\check{\eps}$. Using Gagliardo--Nirenberg inequality with \eqref{eq:mu check estimate} and \eqref{eq:coercivity H2 check}, we deduce \eqref{eq:coercivity H1 inf check} as follows:
    \begin{align*}
          \|\check{\eps}\|_{\dot{\calH}^{1,\infty}} 
            &\leq\|\check{T}(1-\chi_{\lambda^{-1}})\|_{\dot{\calH}^{1,\infty}}
            \!+\!\|\mathring{\eps}\|_{\dot{\calH}^{1,\infty}} \\
            &\lesssim |\check{\nu}_0|\lambda  +|\check{b}_1|+|\check{\eta}_1| + \|\mathring{\eps}\|_{\dot{\calH}^{1}}^{\frac{1}{2}}\|\mathring{\eps}\|_{\dot{\calH}^2}^{\frac{1}{2}} \\
            &\lesssim \lambda^{\frac{3}{2}} + (\|\widehat{\eps}\|_{\dot{\calH}^{1}}+\|\check{T}\chi_{\lambda^{-1}}\|_{\dot{\calH}^{1}})^{\frac{1}{2}}(\|\check{\eps}\|_{\dot{\calH}^{2}}+\|\check{T}(1-\chi_{\lambda^{-1}})\|_{\dot{\calH}^{2}})^{\frac{1}{2}}\lesssim \lambda^{\frac{3}{2}}. \qedhere
    \end{align*}
   %we conclude the proof of Lemma~\ref{lem:coercivity w1}
\end{proof}

Next, we prove Lemma~\ref{lem:coercivity estimate 2}.

\begin{proof}[Proof of Lemma~\ref{lem:coercivity estimate 2}]
    \underline{\textbf{Step 1.}} We prove \eqref{eq:coercivity H0 epsj}, \eqref{eq:c j estimate}, \eqref{eq:coercivity H1 epsj}, the $\bm{\beta}_{j}$ part of \eqref{eq:beta j estimate}, and \eqref{eq:coercivity Linf epsj}. The estimate \eqref{eq:coercivity H0 epsj} follows directly from Proposition~\ref{prop:Lax} and the definition of $B_Q$. We also directly obtain \eqref{eq:c j estimate} by the definition of $c_j$, 
    \begin{align*}
        |c_j|\lesssim \|w_j\|_{L^2}\lesssim \lambda^{j}.
    \end{align*}
    Next, we show \eqref{eq:coercivity H1 epsj}.
    Thanks to \eqref{eq:BQBQstar equal I} and \eqref{eq:DQ BQ decomp}, we obtain
    \begin{align*}
    	\|A_Qw_j\|_{L^2}
        =\|B_Q\td \bfD_Qw_j\|_{L^2}
        \lesssim \|w_{j+1}\|_{L^2}+\|B_Q(w_{j+1}-\td\bfD_Qw_j)\|_{L^2}.
    \end{align*}
    By \eqref{eq:coercivity H0 epsj}, the proof of \eqref{eq:coercivity H1 epsj} is reduced to showing
    \begin{align}
        \|B_Q(w_{j+1}-\td\bfD_Qw_j)\|_{L^2}\lesssim \lambda^{j+1}. \label{eq:coerciv step 2 goal}
    \end{align}
    From $w_{j+1}=\td\bfD_ww_{j}$ and $w=Q+\wt\eps$, we decompose $B_Q(w_{j+1}-\td\bfD_Qw_j)$ into
    \begin{align}\label{eq:coerciv step 2 goal decomp}
    	2B_Q(w_{j+1}-\td\bfD_Qw_j)=B_Q[Q\calH (\ol{\wt\eps}w_j)+\wt\eps\calH (Qw_j)+\wt \eps\calH (\ol{\wt\eps}w_j)].
    \end{align}
    Using \eqref{eq:CommuteHilbert}, $B_QQ=0$, \eqref{eq:coercivity L inf} and \eqref{eq:coercivity H0 epsj}, we have
    \begin{align*}
        \|\eqref{eq:coerciv step 2 goal decomp}\|_{L^2}\lesssim \,& \|B_Q[Q\calH(\tfrac{1}{1+y^2}\ol{\wt\eps} w_j)]\|_{L^2}+\|B_Q[yQ\calH(\tfrac{y}{1+y^2}\ol{\wt\eps} w_j)]\|_{L^2} \\
        &+ \|\tfrac{1}{1+y^2}\wt\eps\calH (Qw_j)\|_{L^2}+\|\tfrac{y}{1+y^2}\wt\eps\calH(yQw_j)\|_{L^2}+\|\tfrac{y}{1+y^2}\wt\eps{\textstyle\int} Qw_j\|_{L^2} \\
        &+  \|\wt\eps\calH (\ol{\wt\eps} w_j)\|_{L^2} \lesssim \|\wt\eps\|_{\dot\calH^1}\|w_j\|_{L^2} +\|\wt\eps\|_{L^\infty}^2\|w_j\|_{L^2}\lesssim \lambda^{j+1}.
    \end{align*}
    Thus, we have \eqref{eq:coerciv step 2 goal}, which implies \eqref{eq:coercivity H1 epsj}. Next, the $\bm{\beta}_{j}$ part of \eqref{eq:beta j estimate} comes from
	\begin{align*}
		\beta_j \lambda^{-1/2}\sim\|\chi_{\lambda^{-1}} P_{j,o}\|_{L^2}
		&\lesssim \|\chi_{\lambda^{-1}} w_{j}\|_{L^2}+\|\chi_{\lambda^{-1}} \eps_{j}\|_{L^2}
		\lesssim \lambda^{j}.
	\end{align*}
	Gagliardo--Nirenberg inequality and $w_j-P_j\chi_{\lambda^{-1}} = P_j(1-\chi_{\lambda^{-1}}) + \eps_j$ yield \eqref{eq:coercivity Linf epsj}:
	\begin{align*}
		\|\eps_j\|_{L^\infty}
		\leq \|\eps_j+P_j(1-\chi_{\lambda^{-1}})\|_{L^\infty}+\|P_j(1-\chi_{\lambda^{-1}})\|_{L^\infty}
        \lesssim\lambda^{j+\frac12}.
	\end{align*}

    \underline{\textbf{Step 2.}} In this step, before showing \eqref{eq:coercivity H2 epsj}, \eqref{eq:coercivity H2 epsj check}, \eqref{eq:betaj and prime relation} and the $\check{\bm{\beta}}_{j}$ part of \eqref{eq:beta j estimate}, we prove the following estimates:
    \begin{subequations}\label{eq:wj lambda cubic}
        \begin{align}
            \|\td\bfD_{Q}(w_{j+1,e}-\td\bfD_{Q}w_{j,o})\|_{L^2}\lesssim \lambda^{j+2},\label{eq:wj lambda cubic 1}
    		\\
	   	\|\td\bfD_{Q}(w_{j+1,o}-\td\bfD_{Q}w_{j,e})\|_{L^2}\lesssim \lambda^{j+2}. \label{eq:wj lambda cubic 2}
        \end{align} 
    \end{subequations}
    
    \textit{\underline{Step 2a}: Proof of \eqref{eq:wj lambda cubic 1}.}  We first decompose $w=Q+\wt\eps$ in the nonlinear terms,
	\begin{align}
		2(w_{j+1,e}-\td\bfD_{Q}w_{j,o})
		=\,&w_o\calH(\ol{w_e} w_{j,e})+w_e\calH(\ol{w_o} w_{j,e})\label{eq:coer2 step 2 3-1 piece}\\
        &+w_o\calH(\ol{w_o} w_{j,o})+\wt\eps_e\calH(\ol{w_e} w_{j,o})+Q\calH(\ol{\wt \eps_e} w_{j,o}).\label{eq:coer2 step 2 3-2 piece}
	\end{align}
    Since \eqref{eq:coer2 step 2 3-2 piece} is parallel to \eqref{eq:coer2 step 2 3-1 piece}, we only compute \eqref{eq:coer2 step 2 3-1 piece} in detail. 
	We note that
	\begin{align*}
		\|\td\bfD_{Q}[f\calH(gh)]\|_{L^2}
		\lesssim
		\|f\calH(gh)\|_{\dot \calH^1}.
	\end{align*}
    To estimate \eqref{eq:coer2 step 2 3-1 piece}, we observe that the profile part of \eqref{eq:coer2 step 2 3-1 piece} becomes
	\begin{align*}
        T_o\calH(Q P_{j,e})+Q\calH(\ol{T_{o}} P_{j,e})
        =
		i\nu_0\tfrac{y}{2}Q\calH(Q\cdot c_j Q)-Q\calH(i\nu_0\tfrac{y}{2}Q\cdot c_jQ)
		=i\nu_0c_j Q.
	\end{align*}
	Thus, the profile terms for $\td\bfD_{Q}[\eqref{eq:coer2 step 2 3-1 piece}]$ is zero, due to $\td\bfD_{Q}Q=0$. Hence, 
    \begin{align}
       \td\bfD_{Q}[\eqref{eq:coer2 step 2 3-1 piece}]=\,&\td\bfD_{Q}[\eps_o\calH(\ol{w_e} w_{j,e})] \label{eq:coer2 step 2 1st piece 1}
        \\
        &+\td\bfD_{Q}[T_{o}\calH(\ol{\wt\eps_e} w_{j,e}+Q\eps_{j,e})] \label{eq:coer2 step 2 1st piece 2}
        \\
        &+\td\bfD_{Q}[\wt\eps_e\calH(\ol{\wt\eps_o} w_{j,e})]+\td\bfD_{Q}[Q\calH(\ol{\eps_o} w_{j,e}+\ol{T_{o}} \eps_{j,e})].  \label{eq:coer2 step 2 1st piece 3}
    \end{align}
    Here, we used $w_o=\wt\eps_o$ and $\wt\eps_{j,e}=\eps_{j,e}$.
    We estimate \eqref{eq:coer2 step 2 1st piece 1} as follows:
	\begin{align}
		\|\td\bfD_{Q}[\eps_o \calH(\ol{w_e} w_{j,e})]\|_{L^2}
		&\leq
		\|\td\bfD_{Q}[\eps_o \calH(\ol{w_e} w_{j,e}-QP_{j,e})]\|_{L^2}+
		\|\td\bfD_{Q}[\eps_o \calH(Q P_{j,e})]\|_{L^2}\nonumber
		\\
		&\lesssim \|\eps_o \calH(\ol{\wt\eps_e} w_{j,e})\|_{\dot\calH^1}
        \!+\!\|\eps_o \calH(Q\eps_{j,e})\|_{\dot\calH^1}\!+\!\lambda^{j+2}\lesssim \lambda^{j+2}.\label{eq:coer2 step 2 1st piece 1 goal}
	\end{align}
	For $\|\eps_o \calH(\ol{\wt\eps_e} w_{j,e})\|_{\dot\calH^1}$, we used $\wt\eps_o = \nu_0 yQ + \eps_o$ to obtain
    \begin{equation*}
        \begin{aligned}
            \|\eps_o \calH(\ol{\wt\eps_e} w_{j,e})\|_{\dot\calH^1}
    		\lesssim \|\eps_o\|_{\dot \calH^{1,\infty}}\|\wt\eps_e\|_{L^\infty}\|w_{j,e}\|_{L^2}
    		+
    		\|\eps_o\|_{L^\infty}\|\partial_y(\ol{\wt\eps_e} w_{j,e})\|_{L^2}
    		\lesssim \lambda^{j+2}.
        \end{aligned}
    \end{equation*}
    For $\|\eps_o \calH(Q\eps_{j,e})\|_{\dot\calH^1}$, we used \eqref{eq:CommuteHilbertDerivative} with $(1+y^2)\langle y\rangle^{-2}=1$ to obtain
	\begin{align*}
		\|\eps_o \calH(Q\eps_{j,e})\|_{\dot\calH^1}
		\lesssim \|\eps_o\|_{\dot \calH^{1,\infty}}\|Q\eps_{j,e}\|_{L^2}
		+ \|Q \eps_o\|_{L^\infty}\|\eps_{j,e}\|_{\dot\calH^1}
        \lesssim \lambda^{j+2+\frac12}
        \lesssim \lambda^{j+2}.
	\end{align*}

    For \eqref{eq:coer2 step 2 1st piece 2}, similar to \eqref{eq:coer2 step 2 1st piece 1} with $T_o = \nu_0 yQ$ and \eqref{eq:CommuteHilbert}, we have
    \begin{align}
		\|\td\bfD_{Q}[T_{o}\calH(\ol{\wt\eps_e} w_{j,e}+Q\eps_{j,e})]\|_{L^2}
		\lesssim   |\nu_o| (\|\ol{\wt\eps_e} w_{j,e}\|_{\dot{\calH}^1}+ \|Q \eps_{j,e}\|_{L^2})  \lesssim  \lambda^{j+2}.\label{eq:coer2 step 2 1st piece 2 goal}
	\end{align}
    
    For the first term of \eqref{eq:coer2 step 2 1st piece 3}, similar to \eqref{eq:coer2 step 2 1st piece 1} with $T_o = \nu_0 yQ$ and $P_{j,e}=c_jQ$,
	\begin{align}
		\|\wt\eps_e\calH(\ol{\wt\eps_o} w_{j,e})\|_{\dot\calH^1}
        &\!\lesssim\!\|\wt\eps_e\|_{\dot \calH^{1,\infty}}\|\wt\eps_o\|_{L^\infty}\|w_{j,e}\|_{L^2}+
    	\|\wt\eps_e|D|(\ol{\wt\eps_o} w_{j,e})\|_{L^2}\nonumber\\
        &\!\lesssim\!\lambda^{j+2}+ |\nu_0 c_j|\|\wt\eps_e\|_{\dot\calH^2}\!+\!\|\wt\eps_e\|_{L^\infty}\|\ol{\eps_o} w_{j,e}\!+\!\ol{T_o}\eps_{j,e}\|_{\dot H^1}
       \! \lesssim \!\lambda^{j+2}.\label{eq:coer2 step 2 1st piece 3 goal}
	\end{align}

    For the second term of \eqref{eq:coer2 step 2 1st piece 3}, \eqref{eq:CommuteHilbert} and $\td\bfD_{Q}Q=\td\bfD_{Q}(yQ)=0$ imply
    \begin{align}
        \td\bfD_{Q}[Q\calH(\ol{\eps_o} w_{j,e}+\ol{T_{o}} \eps_{j,e})]
        =
        \td\bfD_{Q}[(1+y^2)Q\calH\{\tfrac{1}{1+y^2}(\ol{\eps_o} w_{j,e}+\ol{T_{o}} \eps_{j,e})\}].\label{eq:coer2 step 2 1st piece 3-1}
    \end{align}
    Again using \eqref{eq:CommuteHilbertDerivative}, we obtain the desired bound for the second term of \eqref{eq:coer2 step 2 1st piece 3}:
    \begin{align}
		\|\eqref{eq:coer2 step 2 1st piece 3-1}\|_{L^2}\!
        \lesssim\,&
        \|Q^2(\ol{\eps_o} w_{j,e}+\ol{T_{o}} \eps_{j,e})\|_{L^2}
        +\|(1+y^2)Q|D|\{\tfrac{1}{1+y^2}(\ol{\eps_o} w_{j,e}+\ol{T_{o}} \eps_{j,e})\}\|_{L^2} \nonumber
        \\
        \lesssim\,&\lambda^{j+2}\!+\!
        \|Q^2(\ol{\eps_o} w_{j,e}\!+\!\ol{T_{o}} \eps_{j,e})\|_{\dot{H}^1}\!+\!\|Q(\ol{\eps_o} w_{j,e}\!+\!\ol{T_{o}} \eps_{j,e})\|_{\dot{\calH}^1} \!\lesssim \!\lambda^{j+2}.\label{eq:coer2 step 2 1st piece 4 goal}
	\end{align}
	Therefore, collecting \eqref{eq:coer2 step 2 1st piece 1 goal}--\eqref{eq:coer2 step 2 1st piece 4 goal}, we conclude the desired bound for \eqref{eq:coer2 step 2 3-1 piece}.
    
	We turn to \eqref{eq:coer2 step 2 3-2 piece}. Using $w_o=\wt\eps_o =\eps_o+T_o$, we decompose \eqref{eq:coer2 step 2 3-2 piece} into four parts, $\eps_o\calH(\ol{w_o} w_{j,o})$, $T_o\calH(\ol{w_o} w_{j,o})$, $\wt\eps_e\calH(\ol{w_e} w_{j,o})$ and $Q\calH(\ol{\wt \eps_e} w_{j,o})$. We can easily control them by following the proofs of \eqref{eq:coer2 step 2 1st piece 1 goal}--\eqref{eq:coer2 step 2 1st piece 4 goal}.
    Thus, we finish the proof of \eqref{eq:wj lambda cubic 1}.
    
	\textit{\underline{Step 2b}: Proof of \eqref{eq:wj lambda cubic 2}.}  We will show that
	\begin{align}
		\|2(w_{j+1,o}-\td\bfD_{Q}w_{j,e})-\bm{\beta}_1c_jyQ\|_{\dot\calH^1}\lesssim \lambda^{j+2}. \label{eq:wj+1 profile}
	\end{align}
    Then \eqref{eq:wj lambda cubic 2} comes from $\td\bfD_Q(yQ)=0$. By \eqref{eq:CommuteHilbert} with $(1+y^2)\langle y \rangle^{-2}=1$, we have
	\begin{align}
		2(w_{j+1,o}-\td\bfD_{Q}w_{j,e})
		=\,&\wt\eps_o\calH(Q w_{j,o})\!+\!\wt\eps_e\calH(Q w_{j,e})\!
		+\!Q\calH(\ol{\wt\eps} w_{j})_e\!+\!\{\wt\eps\calH(\ol{\wt\eps} w_{j})\}_o.\nonumber\\
	  =\,&\wt\eps_o\calH(Q w_{j,o})+yQ\calH(\tfrac{y}{1+y^2}\ol{\wt\eps_o} w_{j,o})
		+
		Q\calH(\tfrac{1}{1+y^2}\ol{\wt\eps_o} w_{j,o}) \label{eq:coer2 even 1}
		\\
		&\!+\wt\eps_e\calH(Q w_{j,e})
		\!+\!yQ\calH(\tfrac{y}{1+y^2}\ol{\wt\eps_e} w_{j,e})
		\!+\!
		Q\calH(\tfrac{1}{1+y^2}\ol{\wt\eps_e} w_{j,e}) \label{eq:coer2 even 2}
		\\
		&+\{\wt\eps\calH(\ol{\wt\eps} w_{j})\}_o. \label{eq:coer2 even 3}
	\end{align}
	We claim that the main contribution is \eqref{eq:coer2 even 2}. We further decompose \eqref{eq:coer2 even 1} into:
	\begin{align}
		\eqref{eq:coer2 even 1}=\,
		&T_{o}\calH(Q P_{j,o})+yQ\calH(\tfrac{y}{1+y^2}\ol{T_{o}} P_{j,o})
		+
		Q\calH(\tfrac{1}{1+y^2}\ol{T_{o}} P_{j,o}) \label{eq:coer2 even 1 prof}
		\\
		&
		+\eps_o\calH(Q P_{j,o})+\eps_o\calH(Q \eps_{j,o})+T_{o}\calH(Q \eps_{j,o}) \label{eq:coer2 even 1 rad1}
		\\
		&+yQ\calH(\tfrac{y}{1+y^2}(\ol{\eps_o} w_{j,o}+\ol{T_{o}} \eps_{j,o}))
		+
		Q\calH(\tfrac{1}{1+y^2}(\ol{\eps_o} w_{j,o}+\ol{T_{o}} \eps_{j,o})). \label{eq:coer2 even 1 rad2}
	\end{align}
	Since $T_{o}=\frac{i\nu_0}{2}yQ$, by \eqref{eq:algebraic 1} with $(1+y^2)Q^2=2$ and \eqref{eq:CommuteHilbert}, we can check that
	\begin{align}
		\eqref{eq:coer2 even 1 prof}=
		\tfrac{i\nu_0}{2}\bm{\beta}_j[yQ\calH(yQ^2)
		-yQ\calH(\tfrac{y}{1+y^2}y^2Q^2)
		-Q\calH(\tfrac{1}{1+y^2}y^2Q^2)]=0.\label{eq:coer2 even 1 prof goal}
	\end{align}
    Using \eqref{eq:CommuteHilbertDerivative}, we can estimate each of \eqref{eq:coer2 even 1 rad1} as follows:
        \begin{align}
		\|\eps_o\calH(Q P_{j,o})+\eps_{o}\calH(Q \eps_{j,o})\|_{\dot\calH^1}&\lesssim |\bm{\beta}_j|\|Q^2\eps_o\|_{\dot\calH^1}+\|\eps_o\|_{\dot\calH^{1,\infty}}\|\eps_{j,o}\|_{\dot \calH^1}
		\lesssim \lambda^{j+2+\frac12}, \nonumber \\
        	\|T_{o}\calH(Q \eps_{j,o})\|_{\dot\calH^1}
		&\lesssim \lambda(\|Q\eps_{j,o}\|_{\dot\calH^1}\!+\|\calH(Q\eps_{j,o})\|_{L^\infty})
		\lesssim \lambda^{j+2}.\label{eq:coer2 even 1 rad1 goal}
	\end{align}
	For \eqref{eq:coer2 even 1 rad2}, elementary inequalities with $w_{j,o}=P_{j,o}+\eps_{j,o}$ imply
    \begin{align}
        \|\eqref{eq:coer2 even 1 rad2}\|_{\dot \calH^1}
        \lesssim 
		\|Q\ol{\eps_o} w_{j,o}\|_{H^1}^{\frac{1}{2}}\|Q\ol{\eps_o} w_{j,o}\|_{\dot H^1}^{\frac{1}{2}} +\|Q\ol{T_{o}} \eps_{j,o}\|_{H^1}
		\lesssim \lambda^{j+2}.\label{eq:coer2 even 1 rad2 goal}
	\end{align}
	
	Thus, by \eqref{eq:coer2 even 1 prof goal}, \eqref{eq:coer2 even 1 rad1 goal}, and \eqref{eq:coer2 even 1 rad2 goal}, we obtain
	\begin{align}
		\|\eqref{eq:coer2 even 1}\|_{\dot \calH^1}\lesssim \lambda^{j+2}. \label{eq:coer2 even 1 goal}
	\end{align}
    
	Now, we control \eqref{eq:coer2 even 2}. We further decompose \eqref{eq:coer2 even 2} into
	\begin{align}
		\eqref{eq:coer2 even 2}=\,&
		T_{e}\calH(Q P_{j,e})
		+yQ\calH(\tfrac{y}{1+y^2}\ol{T_{e}} P_{j,e})
		+
		Q\calH(\tfrac{1}{1+y^2}\ol{T_{e}} P_{j,e}) \label{eq:coer2 even 2 prof}
		\\
		&+\eps_e\calH(Q P_{j,e})+\eps_e\calH(Q \eps_{j,e})+T_{e}\calH(Q \eps_{j,e}) \label{eq:coer2 even 2 rad1}
		\\
		&+yQ\calH\{\tfrac{y}{1+y^2}(\ol{\eps_e} P_{j,e}+\ol{\wt\eps_e}\eps_{j,e} )\}
		+
		Q\calH\{\tfrac{1}{1+y^2}(\ol{\eps_e} P_{j,e}+\ol{\wt\eps_e}\eps_{j,e} )\}. \label{eq:coer2 even 2 rad2}
	\end{align}
	Using $T_{e}=-ib_1\frac{y^2}{4}Q-\eta_1\frac{1+y^2}{4}Q$ and $P_{j,e}=c_j Q$, we can check that
	\begin{align}
		\eqref{eq:coer2 even 2 prof}
		 =\bm{\beta}_1c_jyQ. \label{eq:coer2 even 2 prof goal}
	\end{align}
    For \eqref{eq:coer2 even 2 rad1} and \eqref{eq:coer2 even 2 rad2}, using elementary inequalities with \eqref{eq:algebraic 1} and \eqref{eq:CommuteHilbertDerivative}, we have
    \begin{align}
        \|\eqref{eq:coer2 even 2 rad1}\|_{\dot \calH^1}&\lesssim |c_j|\|\eps_e\|_{\dot\calH^2}\!+\!\|\eps_e\!+\!T_e\|_{\dot\calH^{1,\infty}}\|\eps_{j,e}\|_{\dot\calH^1}\lesssim \lambda^{j+2}+\lambda^{j+2+\frac{1}{2}}.\label{eq:coer2 even 2 rad1 goal} \\
        \|\eqref{eq:coer2 even 2 rad2}\|_{\dot \calH^1}&\lesssim
		\|Q(\ol{\eps_e} P_{j,e}+\ol{\wt\eps_e}\eps_{j,e} )\|_{H^1}^{\frac{1}{2}}
		\|Q(\ol{\eps_e} P_{j,e}+\ol{\wt\eps_e}\eps_{j,e} )\|_{\dot{\calH}^1}^{\frac{1}{2}}
		\lesssim 
		\lambda^{j+2}.\label{eq:coer2 even 2 rad2 goal}
	\end{align}
	Therefore, by \eqref{eq:coer2 even 2 prof goal}, \eqref{eq:coer2 even 2 rad1 goal}, and \eqref{eq:coer2 even 2 rad2 goal}, we obtain
	\begin{align}
		\|\eqref{eq:coer2 even 2}-\bm{\beta}_1c_jyQ\|_{\dot\calH^1}\lesssim \lambda^{j+2},\quad \|\td\bfD_{Q}\eqref{eq:coer2 even 2}\|_{L^2}\lesssim \lambda^{j+2}. \label{eq:coer2 even 2 goal}
	\end{align}
	
	For \eqref{eq:coer2 even 3}, using $\ol{\wt\eps} w_{j}-\ol{T_o} P_{j,e} = \ol{\wt\eps}  \wt\eps_{j} + \ol{(\wt\eps-T_{o})} P_{j,e}$ and $\wt\eps-T_{o}=T_{e}+\eps$, we have
    \begin{align}
		\|\wt\eps\calH(\ol{\wt\eps} w_{j})\|_{\dot\calH^1}
        &\lesssim  |\nu_0c_j|\|\wt\eps\|_{\dot\calH^1}\!+\!\|\wt\eps\|_{\dot\calH^1} \|\calH(\ol{\wt\eps} w_{j}\!-\!\ol{T_o} P_{j,e})\|_{L^\infty} \!+\!\|\wt\eps\|_{L^\infty}\|\ol{\wt\eps} w_{j}\!-\!\ol{T_{o}}P_{j,e}\|_{\dot{H}^1}\nonumber\\
        &\lesssim \lambda^{j+2} + \lambda \|\ol{\wt\eps} w_{j}\!-\!\ol{T_{o}}P_{j,e}\|_{L^2}^{\frac{1}{2}} \|\ol{\wt\eps} w_{j}\!-\!\ol{T_{o}}P_{j,e}\|_{\dot{H}^1}^{\frac{1}{2}} +\lambda^{\frac{1}{2}} \|\ol{\wt\eps} w_{j}\!-\!\ol{T_{o}}P_{j,e}\|_{\dot{H}^1}\nonumber\\
        &\lesssim \lambda^{j+2} + \lambda (\lambda^{j+\frac{1}{2}})^{\frac{1}{2}} (\lambda^{j+\frac{3}{2}})^{\frac{1}{2}} + \lambda^{\frac{1}{2}} \lambda^{j+\frac{3}{2}} \lesssim \lambda^{j+2}.\label{eq:coer2 even 3 goal}
	\end{align}
	Therefore, by collecting \eqref{eq:coer2 even 1 goal}, \eqref{eq:coer2 even 2 goal}, and \eqref{eq:coer2 even 3 goal}, we conclude \eqref{eq:wj+1 profile} and so \eqref{eq:wj lambda cubic 2}. 
   
    \underline{\textbf{Step 3.}} We show \eqref{eq:coercivity H2 epsj}, \eqref{eq:coercivity H2 epsj check}, \eqref{eq:betaj and prime relation} and the $\check{\bm{\beta}}_{j}$ part of \eqref{eq:beta j estimate}. We first prove the first line of \eqref{eq:coercivity H2 epsj} and \eqref{eq:coercivity H2 epsj check}. We recall the decomposition
	\begin{align}
		\eps_{j,e}=\check{\eps}_{j,e}+\check{\bm{\beta}}_{j+1} (1+y^2)Q,\quad \check\eps_j\in \{\calZ_c,\calZ_h, \calZ_{c,2}\}^\perp. \label{eq:eps je further decom aux}
	\end{align}
	By \eqref{eq:wj lambda cubic}, \eqref{eq:coercivity H1 epsj}, and the coercivity \eqref{eq:coercivity A2} with $\td\bfD_{Q}^{2}=B_Q^*\calB_2$, we obtain
	\begin{align*}
		\|\check \eps_{j}\|_{\dot\calH^2}\sim \|\td\bfD_{Q}^{2}w_{j}\|_{L^2}
        \leq
        \|\td\bfD_{Q}(w_{j+1}-\td\bfD_{Q}w_{j})\|_{L^2}+\|\td\bfD_{Q}w_{j+1}\|_{L^2}
		\lesssim \lambda^{j+2}.
	\end{align*}
    Thus, we yield the first part of \eqref{eq:coercivity H2 epsj check}. The second part follows by an argument similar to the proof of \eqref{eq:coercivity H1 inf check}. Since $\check \eps_{j,o}=\eps_{j,o}$, we conclude the first line of \eqref{eq:coercivity H2 epsj}.
    
	Next, the $\check{\bm{\beta}}_{j}$ part of \eqref{eq:beta j estimate} comes from \eqref{eq:coercivity H1 epsj} and \eqref{eq:coercivity H2 epsj check},
	\begin{align*}
		|\check{\bm{\beta}}_{j+1} |\lambda^{-\frac32}\sim \|\check{\bm{\beta}}_{j+1} y^2Q\chi_{\lambda^{-1}}\|_{L^2}
		\leq 
		\|\eps_{j,e}\chi_{\lambda^{-1}}\|_{L^2}
		+
		\|\check \eps_{j,e}\chi_{\lambda^{-1}}\|_{L^2} 
		\lesssim  \lambda^{j}.
	\end{align*}
	With this, an argument similar to the proof of \eqref{eq:coercivity H1 inf check} yields the second line of \eqref{eq:coercivity H2 epsj}.
	
	Finally, using \eqref{eq:wj decom 1}, \eqref{eq:CommuteHilbertDerivative} with $(1+y^2)\langle y\rangle^{-2}=1$, and \eqref{eq:eps je further decom aux}, we compute
	\begin{align*}
		\td\bfD_{Q}w_{j,e}=\partial_y\check{\eps}_{j,e}+\tfrac{1}{2}Q\calH(\tfrac{1}{1+y^2}Q \check{\eps}_{j,e})
		+\tfrac{1}{2}yQ\calH(\tfrac{y}{1+y^2}Q \check{\eps}_{j,e})
		+\check{\bm{\beta}}_{j+1} yQ.
	\end{align*}
	That is, the profile part of $\td\bfD_{Q}w_{j,e}$ is $\check{\bm{\beta}}_{j+1} yQ$. Thus, \eqref{eq:wj+1 profile} and \eqref{eq:coercivity H2 epsj check} imply
	\begin{align*}
		|2\check{\bm{\beta}}_{j+1} -2\bm{\beta}_{j+1}+\bm{\beta}_1c_j|\lesssim \lambda^{j+2},
	\end{align*}
    which concludes \eqref{eq:betaj and prime relation}, and we finish the proof of Lemma~\ref{lem:coercivity estimate 2}.
\end{proof}

Now, we prove Lemma~\ref{lem:coercivity H3 epsj}.

\begin{proof}[Proof of Lemma~\ref{lem:coercivity H3 epsj}]
    \underline{\textbf{Step 1.}} We first prove \eqref{eq:coercivity H3}.
    From \eqref{eq:AQw1 w level decom} and \eqref{eq:N12 tilde def}, we have
	\begin{align}
		\td\bfD_{Q}^2 w_{1}=\td\bfD_{Q}^2L_{Q}\wt\eps+\td\bfD_{Q}^2(\td N_1+ \td N_2), \label{eq:pf lemma 5.5 eq1}
	\end{align}
	which yields
	\begin{align*}
		\|\td\bfD_{Q}^2w_{1}-\td\bfD_{Q}^2L_{Q}\wt\eps\|_{L^{2}}
		=\|\partial_y^2 B_{Q} (\td N_1+ \td N_2)\|_{L^{2}}=\|\calB_2 (\td N_1+ \td N_2)\|_{L^{2}}.
	\end{align*}
	Here we recall that $\partial_{yy}B_Q=\calB_2$ and \eqref{eq:N12 tilde def},
	\begin{align*}
		\td N_1=\wt\eps\mathcal{H}(\Re(Q\wt\eps))+\tfrac{1}{2}\wt\eps\mathcal{H}(|\wt\eps|^{2}),
		\quad 
		\td N_2
        =
        \tfrac{1}{2}yQ\mathcal{H}(\tfrac{y}{1+y^2}|\wt\eps|^{2})
		+\tfrac{1}{2}Q\mathcal{H}(\tfrac{1}{1+y^2}|\wt\eps|^{2})
        .
	\end{align*}
	For $\wt\eps\mathcal{H}(\Re(Q\wt\eps))$ in $\td N_1$, $\Re(QT_{o})=0$, $\wt\eps_o-T_o=\eps_o$, \eqref{eq:CommuteHilbert} and \eqref{eq:CommuteHilbertDerivative} imply
    \begin{align*}
        \|\calB_2[\wt\eps\mathcal{H}(\Re(Q\wt\eps))]_o\|_{L^2} &\le \|\calB_2[T_{o}\mathcal{H}(\Re(Q\eps_o))+\eps_o\mathcal{H}(\Re(Q\eps_o))+\wt\eps_e\mathcal{H}(\Re(Q\wt\eps_e))]\|_{L^2}\\
        &\lesssim |\nu_0|\|\eps_o\|_{\dot\calH^2} + \|\eps_o\|_{\dot\calH^{1,\infty}}\|\eps_o\|_{\dot\calH^2}+\beta_1\|\wt\eps_e\|_{\dot\calH^2} \lesssim \lambda^3.
    \end{align*}
	For $\wt\eps\mathcal{H}(|\wt\eps|^{2})$, using $\calB_2(yQ)=\calB_2(y^2Q)=\calB_2((1+y^2)Q)=0$ and \eqref{eq:CommuteHilbert}, we have
	\begin{align*}
		\calB_2(\wt\eps\mathcal{H}(|\wt\eps|^{2}))
		=\,&
		\calB_2 (\eps\mathcal{H}(|\wt\eps|^{2}))
		+\calB_2[(1+y^2)T_{o}\mathcal{H}(\tfrac{1}{1+y^2}|\wt\eps|^{2})]
		\\
		&+\calB_2[yT_{e}\mathcal{H}(\tfrac{y}{1+y^2}|\wt\eps|^{2})]
		+\calB_2[T_{e}\mathcal{H}(\tfrac{1}{1+y^2}|\wt\eps|^{2})].
	\end{align*}
	Similarly to the proof of \eqref{eq:coer2 step 2 1st piece 4 goal}, we obtain the same bound for $\wt\eps\mathcal{H}(|\wt\eps|^{2})$. Thus,
	\begin{align}
		\|\calB_2 (\td N_1)_o\|_{L^{2}}
		\lesssim \lambda^{3}. \label{eq:calB3 N1}
	\end{align}
	Now, we turn to $\td N_2$. We observe that the profile terms in $\td N_2$ are
	\begin{equation*}
		\begin{aligned}
			\tfrac{1}{4}yQ\mathcal{H}(yQ^2|T_{o}|^{2})
			+\tfrac{1}{4}Q\mathcal{H}(Q^2|T_{o}|^{2})
			=\tfrac{\nu_0^2}{16}[yQ\calH(y^3Q^4)+Q\calH(y^2Q^4)]
			=-\tfrac{\nu_0^2}{8}yQ^3.
		\end{aligned} 
	\end{equation*}
    Similarly to above, applying \eqref{eq:CommuteHilbert} and \eqref{eq:CommuteHilbertDerivative} with $(1+y^2)\langle y\rangle^{-2}=1$, we estimate
    \begin{equation}
        \begin{aligned}
            \|\calB_2((\td N_2)_o +\tfrac{\nu_0^2}{8}yQ^3)\|_{L^2}
    		\lesssim\,&
    		\|\calB_2[(1+y^2)yQ\calH(\tfrac{y}{(1+y^2)^2}(|\wt\eps|^2-|T_{o}|^2)_e)]\|_{L^2}
    		\\
    		&+
    		\|\calB_2[(1+y^2)Q\calH(\tfrac{1}{(1+y^2)^2}(|\wt\eps|^2-|T_{o}|^2)_e)]\|_{L^2}
    		\\
    		\lesssim \,&\|Q(|\wt\eps|^2-|T_{o}|^2)_e\|_{\dot\calH^2}
    		\lesssim \lambda^{3}.
        \end{aligned} \label{eq:calB3 N2}
    \end{equation}
	From \eqref{eq:pf lemma 5.5 eq1}, $L_QQ=yQ^3$, and $\td\bfD_{Q}^2=B_Q^*\calB_2$, we have
    \begin{equation}
        \begin{aligned}
            B_Q^*\calB_2 ((\td N_1+ \td N_2)_o+\tfrac{\nu_0^2}{8}yQ^3)
    		&=\td\bfD_{Q}^2w_{1,o}-\td\bfD_{Q}^2L_{Q}(\wt\eps_e-\tfrac{\nu_0^2}{8}Q).
        \end{aligned} \label{eq: calB3 N1 N2 rewrite}
    \end{equation}
	Therefore, \eqref{eq:coercivity H2 epsj}, \eqref{eq:calB3 N1}, \eqref{eq:calB3 N2}, and \eqref{eq: calB3 N1 N2 rewrite} imply
	\begin{align*}
		\|\td\bfD_{Q}^2L_{Q}(\wt\eps_e-\tfrac{\nu_0^2}{8}Q)\|_{L^2}
		\lesssim \|\td\bfD_{Q}^2w_{1,o}\|_{L^2}+\|\calB_2 ((\td N_1+ \td N_2)_o+\tfrac{\nu_0^2}{8}yQ^3)\|_{L^{2}}
		\lesssim \lambda^{3}.
	\end{align*}
    From the decomposition \eqref{eq:check eps decom} together with \eqref{eq:transversality},
	we obtain
	\begin{align*}
		\check{\eps}_e -\tfrac{\nu_0^2}{8}Q\perp \calZ_i \quad \text{for}\quad i=1,2,\cdots,6.
	\end{align*}
    Therefore, applying the coercivity \eqref{eq:coercivity L3}, we conclude \eqref{eq:coercivity H3},
    \begin{equation*}
        \|\check{\eps}_e -\tfrac{\nu_0^2}{8}Q\|_{\dot\calH^3}\sim\|\calL_3(\check\eps_e-\tfrac{\nu_0^2}{8}Q)\|_{L^2}
        =\|\td\bfD_{Q}^2L_{Q}(\wt\eps_e-\tfrac{\nu_0^2}{8}Q)\|_{L^2}
        \lesssim \lambda^{3}. 
    \end{equation*}

    \underline{\textbf{Step 2.}} We next show \eqref{eq:coercivity H3 epsj check}.
    We claim that
    \begin{equation}
        \|2(w_{j+1,o}-\td\bfD_{Q}w_{j,e})+ yQ\cdot \mathcal{P}_j\|_{\dot\calH^2}\lesssim \lambda^{j+2+\frac12},\quad \mathcal{P}_j= O(\lambda^{j+1+\frac12}).\label{eq:wj lambda H3 aux}
    \end{equation}
    With \eqref{eq:wj lambda H3 aux}, $\td\bfD_Q yQ=0$ and $\|\td\bfD_Q^2f\|_{L^2}\lesssim \|f\|_{\dot\calH^2}$ imply
    \begin{equation*}
        \|\td\bfD_{Q}^2w_{j+1,o}-\td\bfD_{Q}^3w_{j,e}\|_{L^2}
        \lesssim 
        \|2(w_{j+1,o}-\td\bfD_{Q}w_{j,e})+ yQ\cdot\mathcal{P}_j\|_{\dot\calH^2}
        \lesssim
        \lambda^{j+2+\frac12}. 
    \end{equation*}
    \eqref{eq:coercivity H2 epsj}, \eqref{eq:eps je further decom aux} with $\td\bfD_{Q}^3w_{j,e}=\td\bfD_{Q}^3\check{\eps}_{j,e}$, and \eqref{eq:coercivity A3} with $\td\bfD_{Q}^3=B_Q^*\calB_3$ yield \eqref{eq:coercivity H3 epsj check}.
    
	To show \eqref{eq:wj lambda H3 aux}, we follow the proof of \eqref{eq:wj lambda cubic}, but we extract the profile terms of size $\lambda^{j+1+\frac12}$ and $\lambda^{j+2}$ from $2(w_{j+1,o}-\td\bfD_{Q}w_{j,e})=\eqref{eq:coer2 even 1}+\eqref{eq:coer2 even 2}+\eqref{eq:coer2 even 3}$, namely, 
    \begin{align}
        	\|\eqref{eq:coer2 even 1}
			+yQ\cdot \mathcal{P}_{j,\eqref{eq:coer2 even 1}}\, +\;\, 0 \cdot yQ^3\|_{\dot \calH^2}
			\lesssim \lambda^{j+2+\frac12},\label{eq:coer2 even 1 add goal}\\
            \|\eqref{eq:coer2 even 2}-yQ\cdot\mathcal{P}_{j,\eqref{eq:coer2 even 2}}
		-\tfrac{\nu_0^2c_j}{4}yQ^3 \|_{\dot\calH^2}\lesssim \lambda^{j+2+\frac12}, \label{eq:coer2 even 2 add goal}\\
        \|\eqref{eq:coer2 even 3}+yQ\cdot \mathcal{P}_{j,\eqref{eq:coer2 even 3}}
		+\tfrac{\nu_0^2c_j}{4}yQ^3\|_{\dot\calH^2}
		\lesssim
		\lambda^{j+2+\frac12},\label{eq:coer2 even 3 add goal}
    \end{align}
    where
	\begin{align*}
        \mathcal{P}_{j,\eqref{eq:coer2 even 1}}&=\tfrac{1}{\pi}{\textstyle \int}\tfrac{1}{1+y^2}(\ol{\wt\eps_o} \wt\eps_{j,o}-\ol{T_{o}} P_{j,o}), \\
		\mathcal{P}_{j,\eqref{eq:coer2 even 2}}&=\bm{\beta}_1c_j
		-\tfrac{ic_j}{2}(\check{b}_1-b_1)
		-\tfrac{c_j}{2}(\check{\eta}_1-\eta_1)
		-\tfrac{1}{\pi}{\textstyle \int}\tfrac{1}{1+y^2}(\ol{\wt\eps_e} \wt\eps_{j,e}), \\
         \mathcal{P}_{j,\eqref{eq:coer2 even 3}} &= \tfrac{i\nu_0}{2\pi}{\textstyle \int}\tfrac{y}{1+y^2}(\ol{\wt\eps} w_{j}-\ol{T_{o}} P_{j,e}).
	\end{align*}
    Adding \eqref{eq:coer2 even 1 add goal}, \eqref{eq:coer2 even 2 add goal}, and \eqref{eq:coer2 even 3 add goal}, we conclude \eqref{eq:wj lambda H3 aux} with $\mathcal{P}_{j}=  \mathcal{P}_{j,\eqref{eq:coer2 even 1}} -  \mathcal{P}_{j,\eqref{eq:coer2 even 2}} + \mathcal{P}_{j,\eqref{eq:coer2 even 3}}= O(\lambda^{j+1+\frac12})$.
   Hence, it suffices to show \eqref{eq:coer2 even 1 add goal}, \eqref{eq:coer2 even 2 add goal}, and \eqref{eq:coer2 even 3 add goal}.
   
    \textit{\underline{Step 2a}: Proof of \eqref{eq:coer2 even 1 add goal}.} We recall $\eqref{eq:coer2 even 1}=\eqref{eq:coer2 even 1 rad1}+\eqref{eq:coer2 even 1 rad2}$. Adapting the proof of \eqref{eq:coer2 even 1 rad1 goal} with \eqref{eq:coercivity H1 epsj}--\eqref{eq:coercivity H2 epsj}, we have 
	\begin{align}
		\|\eqref{eq:coer2 even 1 rad1}\|_{\dot \calH^2}\lesssim \lambda^{j+2+\frac12}. \label{eq:coer2 even 1 rad1 goal add}
	\end{align}
	For \eqref{eq:coer2 even 1 rad2}, applying \eqref{eq:CommuteHilbert} and the fact $\ol{\eps_o} w_{j,o}+\ol{T_{o}} \eps_{j,o}=\ol{\wt\eps_o} \wt\eps_{j,o}-\ol{T_{o}} P_{j,o}$, we obtain
    \begin{align*}
        \eqref{eq:coer2 even 1 rad2} = (1+y^2)Q\calH(\tfrac{1}{1+y^2}(\ol{\eps_o} w_{j,o}+\ol{T_{o}} \eps_{j,o})) - yQ\cdot \tfrac{1}{\pi}{\textstyle \int}\tfrac{1}{1+y^2}(\ol{\wt\eps_o} \wt\eps_{j,o}-\ol{T_{o}} P_{j,o}).
    \end{align*}
    Adapting the proof of \eqref{eq:coer2 even 1 rad2 goal} with \eqref{eq:CommuteHilbertDerivative} and \eqref{eq:coercivity H1 epsj}--\eqref{eq:coercivity H2 epsj}, we deduce
	\begin{align}
		\|\eqref{eq:coer2 even 1 rad2}
			+yQ \tfrac{1}{\pi}{\textstyle \int}\tfrac{1}{1+y^2}(\ol{\wt\eps_o} \wt\eps_{j,o}-\ol{T_{o}} P_{j,o})\|_{\dot \calH^2}\lesssim\|\ol{\eps_o} \eps_{j,o}+\ol{T_{o}} \eps_{j,o}\|_{\dot\calH^2}\lesssim \lambda^{j+2+\frac12}.\label{eq:coer2 even 1 rad2 goal add}
	\end{align}
	Thus, by \eqref{eq:coer2 even 1 rad1 goal add} and \eqref{eq:coer2 even 1 rad2 goal add}, we conclude \eqref{eq:coer2 even 1 add goal}.

    \textit{\underline{Step 2b}: Proof of \eqref{eq:coer2 even 2 add goal}.} We recall $\eqref{eq:coer2 even 2}=\bm{\beta}_1c_jyQ + \eqref{eq:coer2 even 2 rad1}+\eqref{eq:coer2 even 2 rad2}$. We first estimate $\eqref{eq:coer2 even 2 rad1}$.
	Adapting the proof of \eqref{eq:coer2 even 2 rad1 goal} with \eqref{eq:coercivity H1 epsj}--\eqref{eq:coercivity H2 epsj}, we have
	\begin{align}
		\|\eqref{eq:coer2 even 2 rad1}-\eps_e\calH(Q P_{j,e})\|_{\dot\calH^2}=\|\eps_e\calH(Q \eps_{j,e})+T_{e}\calH(Q \eps_{j,e})\|_{\dot\calH^2}
		\lesssim \lambda^{j+3}. \label{eq:coer3 add even 1 goal 1-1}
	\end{align}
	For $\eps_e\calH(Q P_{j,e})$, using $\eps_e=\check{T}_e -T_{e}+\check{\eps}_e$, $\eps_e\calH(Q P_{j,e})$ can be written as follows: 
	\begin{align*}
		\eps_e\calH(Q P_{j,e})
        &=-\tfrac{ic_j}{4}(\check{b}_1-b_1)y^3Q^3
		-\tfrac{c_j}{2}(\check{\eta}_1-\eta_1)yQ+\tfrac{\nu_0^2c_j}{8}yQ^3  + c_jyQ^2(\check{\eps}_e -\tfrac{\nu_0^2}{8}Q).
	\end{align*}
	By \eqref{eq:coercivity H3} with \eqref{eq:coer3 add even 1 goal 1-1}, we deduce
	\begin{align}
		\|\eqref{eq:coer2 even 2 rad1}
		+\tfrac{ic_j}{4}(\check{b}_1-b_1)y^3Q^3
		+\tfrac{c_j}{2}(\check{\eta}_1-\eta_1)yQ-\tfrac{\nu_0^2c_j}{8}yQ^3\|_{\dot\calH^2}
		\lesssim \lambda^{j+3}. \label{eq:coer3 add even 1 goal}
	\end{align}
    Now, we consider $\eqref{eq:coer2 even 2 rad2}$. Using $\wt\eps_e=T_e+\eps_e$, we decompose $\eqref{eq:coer2 even 2 rad2}$ as follows:
	\begin{align}
		\eqref{eq:coer2 even 2 rad2}=\,&
		yQ\calH(\tfrac{y}{1+y^2}\ol{\eps_e} P_{j,e})
		+
		Q\calH(\tfrac{1}{1+y^2}\ol{\eps_e} P_{j,e}) \label{eq:coer2 even 1 rad2 add 1}
		\\
		&+yQ\calH(\tfrac{y}{1+y^2}(\ol{\eps_e} \eps_{j,e}+\ol{T_{e}}\eps_{j,e} ))
		+
		Q\calH(\tfrac{1}{1+y^2}(\ol{\eps_e} \eps_{j,e}+\ol{T_{e}}\eps_{j,e} )). \label{eq:coer2 even 1 rad2 add 2}
	\end{align}
	 From $\eps_e=\check{T}_e -T_{e}+\check{\eps}_e$, we further decompose \eqref{eq:coer2 even 1 rad2 add 1} as follows:
	\begin{align}
		\eqref{eq:coer2 even 1 rad2 add 1}
		=\,&
		c_j[yQ\calH(\tfrac{y}{1+y^2}\ol{(\check{\eps}_e -\tfrac{\nu_0^2}{8}Q)} Q)
		+
		Q\calH(\tfrac{1}{1+y^2}\ol{(\check{\eps}_e -\tfrac{\nu_0^2}{8}Q)} Q)] \nonumber
		\\
		&+\tfrac{\nu_0^2c_j}{8}[yQ\calH(\tfrac{y}{1+y^2}Q^2)
		+
		Q\calH(\tfrac{1}{1+y^2}Q^2)]. \label{eq:coer2 even 1 rad2 add 1-1-2} \\
		&
		+
		i(\check{b}_1-b_1)c_j\tfrac{1}{4}
		[yQ\calH(\tfrac{y}{1+y^2}y^2Q^2)
		+Q\calH(\tfrac{1}{1+y^2}y^2Q^2)]  \label{eq:coer2 even 1 rad2 add 1-2}
		\\
		&-(\check{\eta}_1-\eta_1)c_j\tfrac{1}{2}[
		yQ\calH(\tfrac{y}{1+y^2})
		+
		Q\calH(\tfrac{1}{1+y^2})
		].  \label{eq:coer2 even 1 rad2 add 1-3}
	\end{align}
    We can directly check that $\eqref{eq:coer2 even 1 rad2 add 1-1-2}=\tfrac{\nu_0^2c_j}{8}yQ^3$, $\eqref{eq:coer2 even 1 rad2 add 1-2}=-i(\check{b}_1-b_1)c_j\tfrac{1}{4}yQ^3$ and $\eqref{eq:coer2 even 1 rad2 add 1-3}=0$. Thanks to \eqref{eq:coercivity H3}, we conclude
	\begin{align}
		\|\eqref{eq:coer2 even 1 rad2 add 1}
		-\tfrac{\nu_0^2c_j}{8}yQ^3
		+i\tfrac{c_j}{4}(\check{b}_1-b_1)yQ^3\|_{\dot\calH^2}\lesssim |c_j|\lambda^{2+\frac12}
		\lesssim \lambda^{j+2+\frac12}. \label{eq:coer2 even 1 rad2 add 1 goal}
	\end{align}
	Similarly to \eqref{eq:coer2 even 1 rad2} in \textit{\underline{Step 2a}} with $\eps_{j,e}=\wt\eps_{j,e}$, \eqref{eq:coer2 even 1 rad2 add 2} can be written as follows:
	\begin{align*}
		\eqref{eq:coer2 even 1 rad2 add 2}
		=
		(1+y^2)Q\calH(\tfrac{1}{1+y^2}(\ol{\eps_e} \eps_{j,e}+\ol{T_{e}}\eps_{j,e} ))
		-
		yQ\tfrac{1}{\pi}{\textstyle \int}\tfrac{1}{1+y^2}(\ol{\wt\eps_e} \wt\eps_{j,e}).
	\end{align*}
	Adapting the proof of \eqref{eq:coer2 even 2 rad2 goal} with \eqref{eq:CommuteHilbertDerivative} and \eqref{eq:coercivity H1 epsj}--\eqref{eq:coercivity H2 epsj}, we deduce
	\begin{align*}
		\|\eqref{eq:coer2 even 1 rad2 add 2}+yQ\tfrac{1}{\pi}{\textstyle \int}\tfrac{1}{1+y^2}(\ol{\wt\eps_e} \wt\eps_{j,e})\|_{\dot\calH^2}\lesssim\|Q^{-1}\calH(\tfrac{1}{1+y^2}(\ol{\eps_e} \eps_{j,e}+\ol{T_{e}}\eps_{j,e} ))\|_{\dot\calH^2}\lesssim \lambda^{j+2+\frac{1}{2}}.
	\end{align*}
	By \eqref{eq:coer2 even 1 rad2 add 1 goal} and the above estimate, we obtain
	\begin{align}
		\|\eqref{eq:coer2 even 2 rad2}
		-\tfrac{\nu_0^2c_j}{8}yQ^3
		+i\tfrac{c_j}{4}(\check{b}_1-b_1)yQ^3
		+yQ\tfrac{1}{\pi}{\textstyle \int}\tfrac{1}{1+y^2}(\ol{\wt\eps_e} \wt\eps_{j,e})
		\|_{\dot\calH^2}
		\lesssim \lambda^{j+2+\frac12}. \label{eq:coer3 add even 2 goal}
	\end{align}
	Combining \eqref{eq:coer3 add even 1 goal} and \eqref{eq:coer3 add even 2 goal}, we conclude \eqref{eq:coer2 even 2 add goal}.
	
	\textit{\underline{Step 2c}: Proof of \eqref{eq:coer2 even 3 add goal}.} We decompose \eqref{eq:coer2 even 3} into
    \begin{align}
		\{\wt\eps\calH(\ol{\wt\eps} w_{j})\}_o
		=\,&T_{o}\calH[(\ol{\wt\eps} w_{j})_o-\ol{T_{o}} P_{j,e}] +T_{o}\calH(\ol{T_{o}} P_{j,e})\label{eq:coer2 even 3 add 22}\\
        &+\wt\eps_e\calH(\ol{\wt\eps} w_{j})_e + \eps_o\calH(\ol{\wt\eps} w_{j})_{o}. \label{eq:coer2 even 3 add 11}
	\end{align}
     Thanks to \eqref{eq:CommuteHilbert}, we compute \eqref{eq:coer2 even 3 add 22} as follows:
	\begin{align*}
		\eqref{eq:coer2 even 3 add 22}=\,&\tfrac{i\nu_0}{2}y^2Q\calH\{\tfrac{y}{1+y^2}[(\ol{\wt\eps} w_{j})_o-\ol{T_{o}} P_{j,e}]\}
		+
		\tfrac{i\nu_0}{2}yQ\calH\{\tfrac{1}{1+y^2}[(\ol{\wt\eps} w_{j})_o-\ol{T_{o}} P_{j,e}]\}
		\\
		&-
		\tfrac{i\nu_0}{2\pi}yQ{\textstyle \int}\tfrac{y}{1+y^2}(\ol{\wt\eps} w_{j}-\ol{T_{o}} P_{j,e})-\tfrac{\nu_0^2c_j}{4}yQ^3.
	\end{align*}
   Again using \eqref{eq:CommuteHilbert}, \eqref{eq:CommuteHilbertDerivative} and $(\ol{\wt\eps} w_{j})_o-\ol{T_{o}} P_{j,e}=
		\ol{\eps_o} w_{j,e}+\ol{T_{o}} \eps_{j,e}
		+\ol{\wt\eps_e} w_{j,o}$, we have
   \begin{align*}
        &\|\eqref{eq:coer2 even 3 add 22}
		+\tfrac{i\nu_0}{2\pi}yQ{\textstyle \int}\tfrac{y}{1+y^2}(\ol{\wt\eps} w_{j}-\ol{T_{o}} P_{j,e})+\tfrac{\nu_0^2c_j}{4}yQ^3\|_{\dot\calH^2} \\
       &\lesssim |\nu_0|\|y^2Q\calH\{\tfrac{y}{1+y^2}[(\ol{\wt\eps} w_{j})_o-\ol{T_{o}} P_{j,e}]\} + yQ\calH\{\tfrac{1}{1+y^2}[(\ol{\wt\eps} w_{j})_o-\ol{T_{o}} P_{j,e}]\}\|_{\dot\calH^2} \\
       &\lesssim \lambda\|(\ol{\wt\eps} w_{j})_o-\ol{T_{o}} P_{j,e} \|_{\dot\calH^1}^{\frac{1}{2}}\| (\ol{\wt\eps} w_{j})_o-\ol{T_{o}} P_{j,e}\|_{\dot\calH^2}^{\frac{1}{2}} + \lambda\| (\ol{\wt\eps} w_{j})_o-\ol{T_{o}} P_{j,e}\|_{\dot\calH^2} \lesssim \lambda^{j+2+\frac12}.
   \end{align*}
    Using \eqref{eq:CommuteHilbertDerivative}, we can estimate the first term of \eqref{eq:coer2 even 3 add 11} as follows:
	\begin{align*}
		\|\wt\eps_e\calH(\ol{\wt\eps} w_{j})_e\|_{\dot\calH^2}
		\lesssim\,&
		\|\wt\eps_e\|_{\dot\calH^2}\|\calH(\ol{\wt\eps} w_{j})_e\|_{L^\infty}
		+\|\wt\eps_e\|_{\dot\calH^{1,\infty}}\|(\ol{\wt\eps} w_{j})_e\|_{\dot\calH^1} \\
        &+\|\wt\eps_e\calH\partial_{yy}(\ol{T_{o}} P_{j,o}+\ol{T_{e}} P_{j,e})\|_{L^2}
		\\
		&+
		\|\wt\eps_e\|_{L^\infty}\|\partial_{yy}(\ol{\eps_o} w_{j,o}+ \ol{T_{o}} \eps_{j,o}
		+\ol{\eps_e} P_{j,e}+ \ol{\wt\eps_e} \eps_{j,e} )\|_{L^2}
		\lesssim
		\lambda^{j+2+\frac12}.
	\end{align*}
    The same bound holds for the second term of \eqref{eq:coer2 even 3 add 11}. Hence, we conclude \eqref{eq:coer2 even 3 add goal}.
\end{proof}

\section{Proof of nonlinear estimates}\label{sec:nonlin pf}
In this section, we prove the nonlinear estimates Lemma~\ref{lem:nonlinear esti I} and Lemma~\ref{lem:nonlinear esti II}.

Before proceeding with the proof of Lemma~\ref{lem:nonlinear esti I}, we first clarify the connection between $c_j$ and the radiation. This relation will also play a role in the derivation of \eqref{eq:main nonlinear wt eps} within Lemma~\ref{lem:nonlinear esti I}.
\begin{lem}\label{lem:rad to mod 1} For $1\leq j\leq 2L$, we have
	\begin{align}
		\tfrac{1}{\pi}{\textstyle\int} \tfrac{y}{1+y^2} |\wt\eps|^2
		&=
		-\mu_0
		+O(\lambda^{1+\frac12}), \label{eq:nonli rad c0 relation}
		\\
		\tfrac{1}{\pi}{\textstyle\int} \tfrac{y}{1+y^2}\ol{\wt\eps} \wt\eps_j
		&=
		-2c_{j+1}+i\nu_0c_j
		+O(\lambda^{j+1+\frac12}). \label{eq:nonli rad cj relation}
	\end{align}
\end{lem}
\begin{rem}\label{rem:convert rad to mod}
    Indeed, Lemma~\ref{lem:rad to mod 1} shows that
    \begin{align*}
        \tfrac{1}{\pi}\frakc_{0,0}
        &=\tfrac{1}{\pi}{\textstyle\int} \tfrac{y}{1+y^2} |\wt\eps|^2+O(\lambda^{1+\frac12})
		=-\mu_0
		+O(\lambda^{1+\frac12}),
        \\
        \tfrac{1}{\pi}\frakc_{0,j}
        &=\tfrac{1}{\pi}{\textstyle\int} \tfrac{y}{1+y^2}\ol{\wt\eps} \wt\eps_j
		+O(\lambda^{j+1+\frac12})
        =-2c_{j+1}+i\nu_0c_j
		+O(\lambda^{j+1+\frac12}).
    \end{align*}
    %For general indices $(i,j)$, it remains unclear whether a relation between $\frakc_{i,j}$ and the radiation parameters $c_j$ (or lower-order $\frakc_{i,j}$ terms) exists.
\end{rem}

\begin{proof}[Proof of Lemma~\ref{lem:rad to mod 1}]
	For \eqref{eq:nonli rad c0 relation}, \eqref{eq:def cj} and $w_1=\td \bfD_ww$ yield
	\begin{align}
		c_1{\textstyle\int} Q\calZ_c  =
		{\textstyle\int} \td\bfD_w w \calZ_c = \,&
		{\textstyle\int} \td\bfD_Q w \calZ_c
		+
		{\textstyle\int} \tfrac{1}{2}\wt\eps\calH (\ol{\wt\eps} w) \calZ_c \nonumber
		\\
		&+{\textstyle\int} \tfrac{1}{2}\wt\eps\calH(Qw) \calZ_c 	
		+{\textstyle\int} \tfrac{1}{2}Q\calH (\ol{\wt\eps} w) \calZ_c . \label{eq:identity from c1 to c00}
	\end{align}
    By \eqref{eq:coercivity H2 and H1 inf} and $\td\bfD_QQ=\td\bfD_Q(yQ)=0$, we have
	\begin{align*}
		|{\textstyle\int} \td\bfD_Q w \calZ_c|
		&=
		|{\textstyle\int} \td\bfD_Q \eps_{o} \calZ_c|
		\lesssim \|\eps_{o}\|_{\dot \calH^2}+\|\calH(Q\eps_{o})\|_{L^\infty}
		\lesssim \lambda^{1+\frac12},\\
		|{\textstyle\int} \wt\eps\calH (\ol{\wt\eps} w) \calZ_c|
		&\leq
		|{\textstyle\int} \wt\eps\calH (\ol{\wt\eps} Q) \calZ_c|
		+
		|{\textstyle\int} \wt\eps\calH (|\ol{\wt\eps}|^2) \calZ_c|
		\lesssim \lambda^2.
	\end{align*}
    Now, we decompose \eqref{eq:identity from c1 to c00} as follows:
    \begin{align*}
        \eqref{eq:identity from c1 to c00} =\,& {\textstyle\int} \tfrac{1}{2}T_{o}\calH(Q^2) \calZ_c 
        +{\textstyle\int} \tfrac{1}{2}Q\calH(\ol{T_{o}}Q) \calZ_c + {\textstyle\int} \tfrac{1}{2}Q\calH (|\wt\eps|^2) \calZ_c
		\\
        &+ 
		{\textstyle\int} \tfrac{1}{2}T_{o}\calH( Q\wt\eps) \calZ_c
		+
		{\textstyle\int} \tfrac{1}{2}(\wt\eps-T_{o})\calH (Qw) \calZ_c +
		{\textstyle\int} \tfrac{1}{2}Q\calH (\ol{\wt\eps-T_o} Q) \calZ_c .
    \end{align*}
	From \eqref{eq:coercivity H1}, $\wt\eps-T_o=T_e+\eps$ with \eqref{eq:coercivity H2 and H1 inf}, \eqref{eq:c j estimate}, and \eqref{eq:beta j estimate}, we have
    \[|{\textstyle\int} T_{o}\calH (Q\wt\eps) \calZ_c|
		\lesssim \lambda^2,\quad
		|{\textstyle\int} (\wt\eps-T_{o})\calH (Qw) \calZ_c|
		\lesssim \lambda^{\frac32},\quad |{\textstyle\int} \tfrac{1}{2}Q\calH (\ol{\wt\eps-T_o} Q) \calZ_c|\lesssim \lambda^{\frac32}.\]
	Moreover, using \eqref{eq:CommuteHilbert} with $(1+y^2)\langle y\rangle^{-2}=1$, we can check that
	\begin{align*}
		{\textstyle\int} \tfrac{1}{2}Q\calH (|\wt\eps|^2) \calZ_c
        ={\textstyle\int} \tfrac{1}{2}yQ\calH[\tfrac{y}{1+y^2}|\wt\eps|^2] \calZ_c+
		{\textstyle\int} \tfrac{1}{2}Q\calH[\tfrac{1}{1+y^2} |\wt\eps|^2] \calZ_c-{\textstyle\int} \tfrac{1}{2\pi}Q \calZ_c {\textstyle\int} \tfrac{y}{1+y^2} |\wt\eps|^2,	\\
		|{\textstyle\int} \tfrac{1}{2}yQ\calH[\tfrac{y}{1+y^2}|\wt\eps|^2] \calZ_c|
		\lesssim 
		\|\calH[\tfrac{y}{1+y^2} |\wt\eps|^2]\|_{L^\infty}
		\lesssim \lambda^{1+\frac12},
		\quad
		\|\calH[\tfrac{1}{1+y^2}|\wt\eps|^2]\|_{L^2}
		\lesssim \lambda^{2}.
	\end{align*}
	Thus, we conclude \eqref{eq:nonli rad c0 relation} since $T_{o}\calH(Q^2) +Q\calH(\ol{T_{o}}Q) 
		=i\nu_0 Q $. In the proof of \eqref{eq:nonli rad c0 relation}, we replace $c_1$ to $c_{j+1}$ and use $w_j=c_jQ + \wt\eps_j$, we have \eqref{eq:nonli rad cj relation}:
    \begin{align*}
		c_{j+1}{\textstyle\int} Q\calZ_c
		&={\textstyle\int} \td\bfD_Q w_j \calZ_c+
		{\textstyle\int} \tfrac{1}{2}\wt\eps\calH (\ol{\wt\eps} w_{j}) \calZ_c
		+
		{\textstyle\int} \tfrac{1}{2}\wt\eps\calH (Qw_{j}) \calZ_c 
		+{\textstyle\int} \tfrac{1}{2}Q\calH (\ol{\wt\eps} w_{j}) \calZ_c\\
        &= {\textstyle\int} \tfrac{1}{2} T_{o}\calH(Q\cdot c_jQ) \calZ_c\! +\!{\textstyle\int} \tfrac{1}{2}Q\calH (\ol{T_{o}}\cdot c_jQ)\calZ_c\!+\!{\textstyle\int} \tfrac{1}{2}Q\calH (\ol{\wt\eps} \wt\eps_j) \calZ_c\!+\!O(\lambda^{j+1+\frac12})\\
        &=\tfrac{i\nu_0}{2}c_j {\textstyle\int} Q\calZ_c -{\textstyle\int} \tfrac{1}{2\pi}Q \calZ_c {\textstyle\int} \tfrac{y}{1+y^2} \ol{\wt\eps} \wt\eps_j
		+O(\lambda^{j+1+\frac12}).  \qedhere
	\end{align*}
\end{proof}

Now we prove Lemma~\ref{lem:nonlinear esti I}. We rewrite \eqref{eq:def NL} using $c_1=\frac{i\nu_0+\mu_0}{2}$:
\begin{equation*}
    \begin{aligned}
        \textnormal{NL}_0
    	=\,&(w_1-c_1w)\mathcal{H}(\overline{w}\wt\eps)-w\mathcal{H}\{\overline{(w_1-c_1w)}\wt\eps\},
        \\
        \textnormal{NL}_j
        =\,& (w_1-c_1w)\mathcal{H}(\overline{w}\wt\eps_{j})-w\mathcal{H}\{\overline{(w_1-c_1w)}\wt\eps_{j}\}
        \\
        &-P_{1,o}\mathcal{H}(QP_{j,o})
        +yQ\mathcal{H}\{\tfrac{y}{1+y^2}\overline{P_{1,o}}P_{j,o}\}
        +Q\mathcal{H}\{\tfrac{1}{1+y^2}\overline{P_{1,o}}P_{j,o}\},
        \\
        \textnormal{NL}_{Q}
        =\,&
        (w_1-c_1w)\mathcal{H}(\overline{w}Q)-w\mathcal{H}\{\overline{(w_1-c_1w)}Q\}
        \\
        &-(P_1-c_1Q)\mathcal{H}(Q^2)
        +Q\mathcal{H}\{\overline{(P_1-c_1 Q)}Q\}
        +c_1T_{o}\calH(Q^2)-Q\calH(\ol{c_1 T_{o}} Q).
    \end{aligned} 
\end{equation*}
\begin{proof}[Proof of Lemma~\ref{lem:nonlinear esti I}]
    We proceed in the order \eqref{eq:wH w square}, \eqref{eq:wH w wj}, \eqref{eq:main nonlinear wt eps}, \eqref{eq:main nonlinear Q}, and \eqref{eq:main nonlinear eps j odd}. We will also use Lemmas~\ref{lem:proximity}, \ref{lem:coercivity estimate 2} and the decompositions \eqref{eq:w decomposition}, \eqref{eq:wj decomposition} without further mention.
    
	\underline{\textbf{Step 1.}} We show \eqref{eq:wH w square} and \eqref{eq:wH w wj}. For $1\le j\le 2L$, it suffices to show 
    \begin{align}
         \|Q \{w\calH(|w|^2)-Q\calH(Q^2)\}\|_{L^2} +  \|Q \{w\calH(\ol{w} Q)-Q\calH(Q^2)\}\|_{L^2} &\lesssim \lambda,\label{eq:wH w square goal} \\
        \|Q(w\calH(\ol{w} \wt\eps_{j})+\bm{\beta}_jQ^3)\|_{L^2}&\lesssim \lambda^{j+1}. \label{eq:wH w wj goal}
    \end{align}
    For \eqref{eq:wH w square goal}, we use the following decompositions:
	\begin{align*}
		w\calH(|w|^2)-Q\calH(Q^2)&=\wt\eps \calH(Q^2)+2w\calH\Re(Q\wt\eps)+w\calH( |\wt\eps|^2),
		\\
		w\calH(\ol{w} Q)-Q\calH(Q^2)&=\wt\eps \calH(Q^2)+w\calH(\ol{\wt\eps}Q).
	\end{align*}
	With the fact $\calH(Q^2)=yQ^2$, \eqref{eq:wH w square goal} comes from the following estimates:
	\begin{align*}
		 \|Q \wt\eps \calH(Q^2)\|_{L^2}+\| w\calH\Re(Q\wt\eps)\|_{L^2}+\|w\calH(\ol{\wt\eps}Q)\|_{L^2}
		\lesssim  \|Q^2\wt\eps\|_{L^2}+ \|w\|_{L^\infty}\|Q\wt\eps\|_{L^2}&\lesssim \lambda, \\
        \|w\calH(|\wt\eps|^2)\|_{L^2}
		\lesssim \|\calH(|\wt\eps|^2)\|_{L^{\infty}} \|w\|_{L^2}\lesssim \||\wt\eps|^2\|_{L^2}^{\frac12}\||\wt\eps|^2\|_{\dot H^1}^{\frac12}&\lesssim \lambda.
	\end{align*}
	Since $Q\calH(QP_{j,o})=-\bm{\beta}_jQ^3$, \eqref{eq:wH w wj goal} comes from the following estimates:
    \begin{align*}
		% & \lesssim \|Q\wt\eps\|_{L^\infty}\|\ol{w} \wt\eps_{j}\|_{L^2}+\|\calH(\ol{\wt\eps}\wt\eps_{j}+\ol{\wt\eps}\eps_j)\|_{L^\infty}+\|Q\eps_j\|_{L^2} \lesssim \lambda^{j+1}.\\
		\|Q(w\calH(\ol{w} \wt\eps_{j})+\bm{\beta}_jQ^3)\|_{L^2} &\le  \|Q\wt\eps\calH(\ol{w} \wt\eps_{j})\|_{L^2} + \|Q^2\calH(\ol{\wt\eps}\wt\eps_{j}+\ol{w}\eps_j)\|_{L^2}\\
        & \lesssim \|Q\wt\eps\|_{L^\infty}\!\|\ol{w} \wt\eps_{j}\|_{L^2}\!+\!
		\|\calH(\ol{\wt\eps}\wt\eps_{j}\!+\!\ol{\wt\eps}\eps_j)\|_{L^\infty}\!+\!\|Q\eps_j\|_{L^2} \!\lesssim \lambda^{j+1}.
	\end{align*}
	
	\underline{\textbf{Step 2.}} We prove \eqref{eq:main nonlinear wt eps}. Again, it suffices to show
    \begin{align} 
        \|Q\cdot (\textnormal{NL}_0)_o\|_{L^2}
        &\lesssim \beta_1\lambda+\lambda^3, \label{eq:main nonlinear wt eps odd}\\
     \|Q\cdot \{(\textnormal{NL}_0)_e +(2\ol{c_{2}}+\tfrac{1}{2}(\nu_0^2-\mu_0^2)+i\mu_0\nu_0)Q\} \|_{L^2}
        &\lesssim \lambda^{2+\frac{1}{2}}.\label{eq:main nonlinear wt eps even}
    \end{align}
    
    \textit{\underline{Step 2a}: Proof of \eqref{eq:main nonlinear wt eps odd}.} We rewrite $[\textnormal{NL}_0]_o=\textnormal{NL}_{0,1}^{\text{odd}}-\textnormal{NL}_{0,2}^{\text{odd}}+\textnormal{NL}_{0,3}^{\text{odd}}$, where
	\begin{align*}
		\textnormal{NL}_{0,1}^{\text{odd}}&= [(w_1-c_1 w)\mathcal{H}(\overline{w}\wt\eps)]_o
		- P_{1,o}\mathcal{H}(QT_{o}),
		\\
		\textnormal{NL}_{0,2}^{\text{odd}}&=
        [w\mathcal{H}\{\overline{(w_1-c_1 w)}\wt\eps\}]_o
		- yQ\mathcal{H}\{\tfrac{y}{1+y^2}\overline{P_{1,o}}T_{o}\}
        - Q\mathcal{H}\{\tfrac{1}{1+y^2}\overline{P_{1,o}}T_{o}\},\\
         \textnormal{NL}_{0,3}^{\text{odd}}&=P_{1,o}\mathcal{H}(QT_{o})-yQ\mathcal{H}(\tfrac{y}{1+y^2}\overline{P_{1,o}}T_{o})
        -Q\mathcal{H}(\tfrac{1}{1+y^2}\overline{P_{1,o}}T_{o})
        =-\tfrac{1}{2}b_1\nu_0 yQ^3. 
	\end{align*}
    Hence, $ \textnormal{NL}_{0,3}^{\text{odd}}$ satisfies \eqref{eq:main nonlinear wt eps odd}. For $\textnormal{NL}_{0,1}^{\text{odd}}$, using \eqref{eq:CommuteHilbert}, we can check that
	\begin{align}
		\|Q\cdot \textnormal{NL}_{0,1}^{\text{odd}}\|_{L^2} &\le  \|Q (\eps_1-c_1\wt\eps)\mathcal{H}(\overline{w}\wt\eps)\|_{L^2} +   \|QP_{1,o}\mathcal{H}((\overline{w}\wt\eps)_{o}-QT_{o})\|_{L^2} \nonumber \\
        &\lesssim (\|\eps_1\|_{\dot \calH^1}+|c_1|\|\wt\eps\|_{\dot\calH^1})\|\calH (\overline{w}\wt\eps)\|_{L^\infty} + \beta_1\|Q\mathcal{H}((|\wt\eps|^2)_{o}+Q\eps_{o})\|_{L^2} \nonumber \\
        &\lesssim \lambda^{3} + \beta_1[\|Q(|\wt\eps|^2\!+\!Q\eps)\|_{L^2}
    		\!+\!|{\textstyle\int}\tfrac{y}{1+y^2}(|\wt\eps|^2\!+\!Q\eps)|] \lesssim \lambda^{3}+\beta_1\lambda. \label{eq:NL01 odd goal}
	\end{align}
    
	 Using $w=Q+\wt\eps$ and \eqref{eq:CommuteHilbert}, we decompose $\textnormal{NL}_{0,2}^{\text{odd}}
        =\textnormal{NL}^{\text{odd}}_{0,2,1}
        +\textnormal{NL}^{\text{odd}}_{0,2,2}$, where
	\begin{align*}
		\textnormal{NL}^{\text{odd}}_{0,2,1}
		&= [\wt\eps\mathcal{H}\{\overline{(w_1-c_1 w)}\wt\eps\}]_o,
		\\
		\textnormal{NL}^{\text{odd}}_{0,2,2}
        &=yQ\mathcal{H}[\tfrac{y}{1+y^2}\{\overline{(w_1-c_1 w)}\wt\eps-\overline{P_{1,o}}T_{o}\}_e]
        \!+\!
		Q\mathcal{H}[\tfrac{1}{1+y^2}\{\overline{(w_1-c_1 w)}\wt\eps-\overline{P_{1,o}}T_{o}\}_e].
	\end{align*}
	By Gagliardo--Nirenberg inequality and $ w_1-c_1 w
        =\eps_1-c_1\wt\eps+P_{1,o}$, we have
	\begin{align}
		\|Q\cdot \textnormal{NL}^{\text{odd}}_{0,2,1}\|_{L^2}
        &\lesssim \|Q\wt\eps\|_{L^2}\|\ol{(w_1-c_1 w)}\wt\eps\|_{L^2}^{\frac{1}{2}}\|\ol{(w_1-c_1 w)}\wt\eps\|_{\dot H^1}^{\frac{1}{2}} \lesssim \lambda^{3}. \label{eq:NL021 odd goal} \\
     \|Q \cdot \textnormal{NL}^{\text{odd}}_{0,2,2}\|_{L^2}  &\lesssim \|Q(\{\ol{(\eps_1-c_1 \wt\eps)}\wt\eps\}_e +\ol{P_{1,o}}\eps_{o})\|_{L^2}^{\frac{1}{2}}\|Q(\{\ol{(\eps_1-c_1 \wt\eps)}\wt\eps\}_e +\ol{P_{1,o}}\eps_{o})\|_{\dot H^1}^{\frac{1}{2}}\nonumber \\
        &\lesssim {\textstyle\prod}_{m=0}^1 (\|\eps_1-c_1 \wt\eps\|_{\dot\calH^1}\|\wt\eps\|_{\dot{\calH}^{m,\infty}}
        +\beta_1\|\eps\|_{\dot\calH^{m+1}})^{\frac{1}{2}}\lesssim  \lambda^{3}.\label{eq:NL022 odd goal}
    \end{align}
	Collecting \eqref{eq:NL01 odd goal}, \eqref{eq:NL021 odd goal} and \eqref{eq:NL022 odd goal}, we finally obtain \eqref{eq:main nonlinear wt eps odd}.
	
    \textit{\underline{Step 2b}: Proof of \eqref{eq:main nonlinear wt eps even}.}
	We decompose $\text{NL}_0=\text{NL}_{0,1}-\text{NL}_{0,2}-\text{NL}_{0,3}$, where
    \[\text{NL}_{0,1}=(\wt\eps_1-c_1 \wt\eps)\calH(\ol{w}\wt\eps),\quad\text{NL}_{0,2}=Q\calH[\overline{(w_1-c_1 w)}\wt\eps],\quad\text{NL}_{0,3}=\wt\eps\calH[\overline{(w_1-c_1 w)}\wt\eps].\]
   $\text{NL}_{0,1}$ and $\text{NL}_{0,3}$ can be easily estimated:
    \begin{align}
        \|Q\text{NL}_{0,1}\|_{L^2} &\le \|Q (\wt\eps_1-c_1 \wt\eps)\calH(Q\wt\eps)\|_{L^2} + \|Q (\wt\eps_1-c_1 \wt\eps)\calH(|\wt\eps|^2)\|_{L^2} \nonumber \\
        &\lesssim  \|\wt\eps_1 -c_1\wt\eps\|_{\dot{\calH}^{1,\infty}}\|\wt\eps\|_{\dot{\calH}^1} +  \|\wt\eps_1 -c_1\wt\eps\|_{\dot{\calH}^1} \|\calH(|\wt\eps|^2)\|_{L^\infty} \lesssim  \lambda^{2+\frac12}.\label{eq:NL01e goal}\\
		\|Q\text{NL}_{0,3}\|_{L^2}
		&\lesssim \|Q\wt\eps\|_{L^2}(\|w_1\|_{L^2}+|c_1|\|w\|_{L^2})\|\wt\eps\|_{L^\infty}
		\lesssim \lambda^{2+\frac12}.
        \label{eq:NL03e goal}
	\end{align}
    
    Now, we estimate $\textnormal{NL}_{0,2}$. Using \eqref{eq:CommuteHilbert} twice with the oddness, we obtain
    \begin{align*}
        (\textnormal{NL}_{0,2})_e
        =(1+y^2)Q\mathcal{H}[\tfrac{1}{1+y^2}\{\overline{(\wt\eps_1-c_1 \wt\eps)}\wt\eps\}_{o}]
        -\tfrac{1}{\pi}Q{\textstyle \int}\tfrac{y}{1+y^2}\{\overline{(\wt\eps_1-c_1 \wt\eps)}\wt\eps\}_{o}.
    \end{align*}
    Using \eqref{eq:nonli rad c0 relation}, \eqref{eq:nonli rad cj relation}, and $c_1=\frac{i\nu_0+\mu_0}{2}$, we have
   \begin{align*}
        -\tfrac{1}{\pi}{\textstyle \int}\tfrac{y}{1+y^2}\{\overline{(\wt\eps_1-c_1 \wt\eps)}\wt\eps\}_{o} 
        &=2\ol{c_{2}}+\tfrac{1}{2}(\nu_0^2-\mu_0^2)+i\mu_0\nu_0+O(\lambda^{2+\frac12}).
	\end{align*}
    The first part of $(\textnormal{NL}_{0,2})_e$ has the same estimate of $\textnormal{NL}_{0,1}$. Thus, we arrive at
	\begin{align}
    \|Q \{(\text{NL}_{0,2})_e
		- Q(2\ol{c_{2}}+\tfrac{1}{2}(\nu_0^2-\mu_0^2)+i\mu_0\nu_0)
		\}\|_{L^2}&\lesssim \lambda^{2+\frac12}. \label{eq:NL02e goal}
	\end{align}
	Combining \eqref{eq:NL01e goal}, \eqref{eq:NL03e goal}, and \eqref{eq:NL02e goal}, we obtain \eqref{eq:main nonlinear wt eps even}, which concludes \eqref{eq:main nonlinear wt eps}.
	
	\underline{\textbf{Step 3.}} We prove \eqref{eq:main nonlinear Q}. Similarly to \underline{\textbf{Step 2}}, it suffices to show
    \begin{align}
		\|Q\cdot (\textnormal{NL}_{Q}-2\bm{\beta}_2 yQ )_o \|_{L^2}
		&\lesssim \beta_1\lambda+\lambda^{3}, \label{eq:main nonlinear Q odd}\\
        \|Q\cdot (\textnormal{NL}_{Q})_e \|_{L^2}
		&\lesssim \lambda^{2+\frac{1}{2}}.\label{eq:main nonlinear Q even}
	\end{align}
    Before that, we decompose $\textnormal{NL}_{Q}=\textnormal{NL}_{Q,1}+\textnormal{NL}_{Q,2}+\textnormal{NL}_{Q,3}$, where
	\begin{align*}
		\textnormal{NL}_{Q,1}
		&=(P_1-c_1Q)\mathcal{H}(\ol{\wt\eps}Q)
		-\wt\eps\mathcal{H}[\overline{(w_1-c_1w)}Q],
		\\
		\textnormal{NL}_{Q,2}
		&=\eps_1\mathcal{H}(Q^2)+\eps_1\mathcal{H}(\overline{\wt\eps} Q)
		-Q\mathcal{H}[\overline{\eps_1}Q]
		-c_1\wt\eps\mathcal{H}(\overline{\wt\eps} Q),
		\\
		\textnormal{NL}_{Q,3}
		&=-c_1(\wt\eps-T_{o})\mathcal{H}(Q^2)
		+Q\mathcal{H}[\overline{c_1(\wt\eps-T_{o})}Q].
	\end{align*}
    We first show \eqref{eq:main nonlinear Q odd}. The $\textnormal{NL}_{Q,1}$ terms can be easily checked as follows:
	\begin{align}
		\|Q\textnormal{NL}_{Q,1}\|_{L^2}
		\lesssim \beta_1\|\wt\eps\|_{\dot{\calH}^1}
		+\|\wt\eps\|_{\dot{\calH}^1}(\beta_1+\|\wt\eps_1\|_{\dot{\calH}^1}+|c_1|\|\wt\eps\|_{\dot{\calH}^1})
		\lesssim \beta_1\lambda+\lambda^{3}. \label{eq:NLQ 1 esti}
	\end{align}
    For $\textnormal{NL}_{Q,2}$, we further decompose $(\textnormal{NL}_{Q,2})_o$ by using \eqref{eq:eps je further decom}:
    \begin{align*}
        (\textnormal{NL}_{Q,2})_o
        &=
        2\check{\bm{\beta}}_2  yQ
        +\check\eps_{1,e}\mathcal{H}(Q^2)
        +[\eps_1\mathcal{H}(\overline{\wt\eps} Q)]_o
		-Q\mathcal{H}[\overline{\eps_{1,e}}Q]
		-c_1[\wt\eps\mathcal{H}(\overline{\wt\eps} Q)]_o.
    \end{align*}
    Hence, using \eqref{eq:CommuteHilbert} with $(1+y^2)\langle y \rangle^{-2}=1$ and \eqref{eq:eps je further decom}, we obtain
    \begin{align}
		\|Q(\textnormal{NL}_{Q,2}-2\check{\bm{\beta}}_2  yQ)_o\|_{L^2}
		&\lesssim  \|\check\eps_{1,e}\|_{\dot{\calH}^2}\!+\!\|\eps_1\|_{\dot{\calH}^{1,\infty}}\!\|\wt\eps\|_{\dot{\calH}^1}
        \!+\!\|Q^2\calH(Q \eps_{1,e})\|_{L^2}
        \!+\!|c_1|\|\wt\eps\|_{\dot{\calH}^1}^2 \nonumber
		\\
		&\lesssim \lambda^{3}+\|Q^2\calH(Q \eps_{1,e})\|_{L^2} \lesssim  \lambda^{3}+ \|\check\eps_{1,e}\|_{\dot\calH^2}\lesssim \lambda^3.\label{eq:NLQ 2o esti}
	\end{align}
    Again using \eqref{eq:CommuteHilbert} with $(1+y^2)\langle y \rangle^{-2}=1$, we estimate $\textnormal{NL}_{Q,3}$ as follows:
    \begin{align}
        \|Q(\textnormal{NL}_{Q,3})_o\|_{L^2}
        \lesssim |c_1|\big[\|yQ^3\wt\eps_e\|_{L^2}\!+\!\|Q^2\calH(Q\wt\eps_e)\|_{L^2}\big]\lesssim \lambda\|\wt\eps_e\|_{\dot{\calH}^2} \!\lesssim \!\beta_1\lambda\!+\!\lambda^{3}.\label{eq:NLQ 3o esti}
    \end{align}
    Combining \eqref{eq:NLQ 1 esti}, \eqref{eq:NLQ 2o esti}, \eqref{eq:NLQ 3o esti} and \eqref{eq:betaj and prime relation}, we conclude \eqref{eq:main nonlinear Q odd}. 
    
    Now, we prove \eqref{eq:main nonlinear Q even}. For $\textnormal{NL}_{Q,1}$, we use \eqref{eq:NLQ 1 esti}. Thus, we conclude \eqref{eq:main nonlinear Q} by
	\begin{align*}
		\|Q(\textnormal{NL}_{Q,2})_e\|_{L^2}
		&\lesssim  \|\eps_1\|_{\dot{\calH}^2}+\|\eps_1\|_{\dot{\calH}^{1,\infty}}\|\wt\eps\|_{\dot{\calH}^1}
		+\|\calH(Q \eps_{1,o})\|_{L^\infty}
		+|c_1|\|\wt\eps\|_{\dot{\calH}^1}^2
		\lesssim \lambda^{2+\frac12},\\
		\|Q(\textnormal{NL}_{Q,3})_e\|_{L^2}
		&\lesssim |c_1|(\|Q^2\eps_o\|_{L^2}+\|\calH(Q\eps_o)\|_{L^\infty})
		\lesssim \lambda^{2+\frac12}.
	\end{align*}
	
	\underline{\textbf{Step 4.}} We prove \eqref{eq:main nonlinear eps j odd}.
    It suffices to show that
    \begin{equation}\label{eq:main nonlinear eps j odd Q}
        \|Q[\textnormal{NL}_j-\tfrac{i}{2}\mu_0\nu_0 \bm{\beta}_j yQ^3]_o\|_{L^2}
		\lesssim \beta_1\lambda^{j+1}+\lambda^{j+3}.
    \end{equation}
	We decompose $\textnormal{NL}_j=\calN_{j,1}-\calN_{j,2}$, where
	\begin{align*}
		\calN_{j,1}&\!=\!
		(w_1-c_1w)\mathcal{H}(\overline{w}\wt\eps_{j})
		-P_{1,o}\mathcal{H}(QP_{j,o})+c_1T_{o}\mathcal{H}(QP_{j,o}),
		\\
		\calN_{j,2}&\!=\!
		w\mathcal{H}\{\overline{(w_1-c_1w)}\wt\eps_{j}\!\}
		\!-\!yQ\mathcal{H}\{\!\tfrac{y}{1+y^2}\overline{P_{1,o}}P_{j,o}\!\}
		\!-\!Q\mathcal{H}\{\!\tfrac{1}{1+y^2}\overline{P_{1,o}}P_{j,o}\!\}
		\!+\!c_1T_{o}\mathcal{H}(QP_{j,o}).
	\end{align*}
    It suffices to estimate $\calN_{j,2}$ since $\overline{w}\wt\eps_{j}-QP_{j,o}=\overline{\wt\eps}\wt\eps_{j}+Q\eps_j$, we have 
    \begin{align*}
         \|Q[\calN_{j,1}]_o\|_{L^2} &=\|Q[(\wt\eps_1- c_1\wt\eps)\mathcal{H}(\overline{\wt\eps}\wt\eps_{j}+Q\eps_{j})
		+\{\eps_1-c_1(\wt\eps-T_{o})\}\mathcal{H}(QP_{j,o})]_o\|_{L^2} \nonumber \\
         &\lesssim (\beta_1 \!+\!\lambda^2)\|\mathcal{H}(\overline{\wt\eps}\wt\eps_{j}\!+\!Q\eps_{j})\|_{L^\infty} \!+\! \beta_j \|Q^3(\eps_{1,o}\!-\!c_1\eps_o)\|_{L^2}\! \lesssim \!\beta_1\lambda^{j+1} \!+\! \lambda^{j+3}.%\label{eq:Nj51 goal}
    \end{align*}
    
	We now estimate $\calN_{j,2}$. Using \eqref{eq:CommuteHilbert}, we decompose $[\calN_{j,2}]_o=\calN_{j,2,1}+\calN_{j,2,2}$, where
	\begin{align*}
		\calN_{j,2,1}=\,&[\wt\eps\mathcal{H}\{\overline{(w_1-c_1w)}\wt\eps_{j}\}]_o,
		\\
		\calN_{j,2,2}=\,&yQ\mathcal{H}\{\tfrac{y}{1+y^2}(\overline{(w_1-c_1w)}\wt\eps_{j}-\overline{P_{1,o}}P_{j,o})_e\}\\
		&+Q\mathcal{H}\{\tfrac{1}{1+y^2}(\overline{(w_1-c_1w)}\wt\eps_{j}-\overline{P_{1,o}}P_{j,o})_e\}+c_1T_{o}\mathcal{H}(QP_{j,o}).
	\end{align*}
	For $\calN_{j,2,1}$, Gagliardo--Nirenberg inequality with $\wt\eps_{j}=P_{j,o}+\eps_j$ implies
    \begin{align}
        \|Q \calN_{j,2,1}\|_{L^2}&\lesssim   \|\wt\eps \|_{\dot{\calH}^{1}}\|\calH\{\overline{(\wt\eps_1-c_1\wt\eps)} \wt\eps_{j}\}\|_{L^\infty}\lesssim \lambda^{j+3}+\beta_1 \lambda^{j+1}.\label{eq:Nj521 goal} 
    \end{align}

	Now, we estimate $\calN_{j,2,2}$. We observe that
	\begin{align*}
		\overline{(w_1-c_1w)}\wt\eps_{j}-\overline{P_{1,o}}P_{j,o}
		=
		\overline{P_{1,o}-c_1 T_o}\eps_j
		+\overline{(\eps_1-c_1(\wt\eps-T_{o}))}\wt\eps_{j}
        -\overline{c_1T_{o}}P_{j,o}.
	\end{align*}
	Using this, we further decompose $\calN_{j,2,2}=\calN_{j,2,2,1}+\calN_{j,2,2,2}+\calN_{j,2,2,3}$, where
	\begin{align}
		\calN_{j,2,2,1}&=yQ\mathcal{H}\{\tfrac{y}{1+y^2}\overline{P_{1,o}-c_1 T_o}\eps_{j,o}\}
		+Q\mathcal{H}\{\tfrac{1}{1+y^2}\overline{P_{1,o}-c_1 T_o}\eps_{j,o}\}, \nonumber
		\\
		\calN_{j,2,2,2}&=yQ\mathcal{H}\{\tfrac{y}{1+y^2}\overline{(\eps_1-c_1(\wt\eps-T_{o}))}\wt\eps_{j}\}_o
		+Q\mathcal{H}\{\tfrac{1}{1+y^2}\overline{(\eps_1-c_1(\wt\eps-T_{o}))}\wt\eps_{j}\}_e, \nonumber
		\\
		\calN_{j,2,2,3}&=
		-yQ\mathcal{H}\{\tfrac{y}{1+y^2}\overline{c_1T_{o}}P_{j,o}\}
		-Q\mathcal{H}\{\tfrac{1}{1+y^2}\overline{c_1T_{o}}P_{j,o}\}
		+c_1T_{o}\mathcal{H}(QP_{j,o}) \nonumber \\
        &=-\tfrac{i}{2}\mu_0\nu_0 \bm{\beta}_j yQ^3.\label{eq:Nj5223 goal}
	\end{align}
   Hence, it suffices to check $\calN_{j,2,2,1}$ and $\calN_{j,2,2,2}$. By elementary inequalities, we derive
	\begin{align}
		\|Q\calN_{j,2,2,1}\|_{L^2}&\lesssim \| \overline{P_{1,o}-c_1T_{o}}\eps_j\|_{\dot{\calH}^1}
		\lesssim |\beta_1 + c_1 \nu_0|
		\|\eps_j\|_{\dot{\calH}^1}\lesssim\beta_1\lambda^{j+1}+\lambda^{j+3}. \label{eq:Nj5221 goal} \\
        \|Q\calN_{j,2,2,2}\|_{L^2}  &\lesssim {\textstyle\prod_{m=0}^1} \| Q\{\overline{(\eps_1-c_1(\wt\eps-T_{o}))}\wt\eps_{j}\}_e\|_{\dot{\calH}^m}^{\frac{1}{2}} \lesssim \lambda^{j+3}.\label{eq:Nj5222 goal}
    \end{align}  
    Here, we used $\{\overline{(\eps_1-c_1(\wt\eps-T_{o}))}\wt\eps_{j}\}_e = \overline{\eps_{1,o}}P_{j,o} +(\overline{\eps_1}\eps_{j})_e-\overline{c_1}(\overline{\wt\eps_e} \eps_{j,e}+\overline{\eps_o} \eps_{j,o})$ and 
    \begin{align*}
      \|Q\{\overline{(\eps_1-c_1(\wt\eps-T_{o}))}\wt\eps_{j}\}_e\|_{\dot{\calH}^m}\lesssim \,& \beta_j\|\eps_{1,o}\|_{\dot\calH^{m+1}}+\|\eps_1\|_{\dot \calH^1}\|\eps_j\|_{\dot\calH^{m,\infty}}\\
        &\!+ \lambda(\|\wt\eps_e\|_{\dot\calH^1}\!+\!\|\eps_o\|_{\dot\calH^1})\|\eps_{j}\|_{\dot\calH^{m,\infty}}\! \lesssim \lambda^{j+m+2+\frac12},
	\end{align*}
	for $m=0,1$. By \eqref{eq:Nj521 goal}--\eqref{eq:Nj5222 goal}, we yield \eqref{eq:main nonlinear eps j odd Q}, which concludes \eqref{eq:main nonlinear eps j odd}. \qedhere
\end{proof}

Finally, we finish the proof of Lemma~\ref{lem:nonlinear esti II}.

\begin{proof}[Proof of Lemma~\ref{lem:nonlinear esti II}]
 It suffices to show
    \begin{equation}
        \|Q\cdot(\textnormal{NL}_j-\tfrac{1}{\pi}(\frakc_{1,j}- \ol{c_1}\frakc_{0,j})Q)_e \|_{L^2}\lesssim \lambda^{j+2+\frac12}. \label{eq:main nonlinear eps j even pf}
    \end{equation}
    We decompose $ [\textnormal{NL}_j]_e=[\textnormal{NL}_{j,1}]_e-[\textnormal{NL}_{j,2}]_e -[\textnormal{NL}_{j,3}]_e$, where
    \begin{align*}
        \textnormal{NL}_{j,1}\!=\!(\wt\eps_1-c_1 \wt\eps)\mathcal{H}(\overline{w}\wt\eps_{j}),
        \,
        \textnormal{NL}_{j,2}\!=\!Q\mathcal{H}\{\overline{(w_1-c_1 w)}\wt\eps_{j}\},
       \,
        \textnormal{NL}_{j,3}\!=\!
        \wt\eps\mathcal{H}\{\overline{(w_1-c_1 w)}\wt\eps_{j}\}.
    \end{align*}
    Adapting the proof of \eqref{eq:NL01e goal} and \eqref{eq:NL03e goal} with \eqref{eq:coercivity H1 epsj}--\eqref{eq:coercivity H2 epsj}, we can check that
    \[ \|Q[\textnormal{NL}_{j,1}]_e\|_{L^2} +  \|Q[\textnormal{NL}_{j,3}]_e\|_{L^2}  \lesssim \lambda^{j+2+\frac12}.\]
    Using \eqref{eq:CommuteHilbert} twice with the oddness, we rewrite $[\textnormal{NL}_{j,2}]_e$ as follows:
    \begin{align*}
        [\textnormal{NL}_{j,2}]_e
        =(1+y^2)Q\mathcal{H}[\tfrac{1}{1+y^2}\{\overline{(w_1-c_1 w)}\wt\eps_{j}\}_{o}]
        -\tfrac{1}{\pi}Q{\textstyle \int}\tfrac{y}{1+y^2}\{\overline{(w_1-c_1 w)}\wt\eps_{j}\}_{o}.
    \end{align*}
    For the last term, by $w_1-c_1 w = \wt\eps_1-c_1 \wt\eps$ and \eqref{eq:frakc radiatin relation}, we observe that
    \begin{align*}
        {\textstyle \int}\tfrac{y}{1+y^2}\{\overline{(w_1-c_1 w)}\wt\eps_{j}\}_{o}
        =
        \frakc_{1,j}-\ol{c_1}\frakc_{0,j}+O(\lambda^{j+2+\frac12}).
    \end{align*}
    Thus, we conclude \eqref{eq:main nonlinear eps j even pf} by using $w_1-c_1 w = \eps_1-c_1 \wt\eps+P_{1,o}$, obtaining
    \begin{align*}
        \|Q\cdot([\textnormal{NL}_{j,2}]_e+\tfrac{1}{\pi}(\frakc_{1,j}- \ol{c_1}\frakc_{0,j})Q) \|_{L^2}
        &\lesssim \|Q(\wt\eps_1-c_1 \wt\eps) \|_{L^\infty}\|Q\wt\eps_{j}\|_{L^2}
        \lesssim \lambda^{j+2+\frac12}. \qedhere
    \end{align*}
\end{proof}

\appendix

\section{Subcoercivity}
In the appendix, we collect some facts regarding the functional settings and subcoercivity. The proof is rather standard or found in \cite{KimKwon2020blowup,KimKimKwon2024arxiv}. We assume that $f\in \mathcal{S}$, unless otherwise specified. 

We mainly use the following weighted Hardy inequalities: 
\begin{lem}[{\cite[Lemma~A.1]{KimKwon2020blowup}}] \label{lem:weight hardy} 
    Let $0<r_{1}<r_{2}<\infty$. In addition, let $\varphi:[r_{1},r_{2}]\to\mathbb{R}_{+}$ be a $C^{1}$ weight function such that $\partial_{r}\varphi$ is non-vanishing and $\varphi\lesssim|r\partial_{r}\varphi|$. Then, for smooth $f:[r_{1},r_{2}]\to\mathbb{C}$, we have 
    \begin{align*}
        \int_{r_{1}}^{r_{2}}\left|f\right|^{2}
        |\partial_{r}\varphi|dr
        \lesssim\int_{r_{1}}^{r_{2}}|\partial_{r}f|^{2} (r\varphi) dr+
        \begin{cases}
            \varphi(r_{2})|f(r_{2})|^{2} & if\ \partial_{r}\varphi>0,\\
            \varphi(r_{1})|f(r_{1})|^{2} & if\ \partial_{r}\varphi<0.
        \end{cases}
    \end{align*}
\end{lem}

\begin{lem}[Comparison of $\dot{\mathcal{H}}^{2}$ and $\dot{H}^{2}$]
	\label{lem:H2comparison} We have 
	\begin{align*}
		\|\partial_{xx}f\|_{L^{2}}+\|{\bf 1}_{|x|\lesssim1}(|\partial_{x}f|+|f|)\|_{L^{2}}\sim\|f\|_{\dot{\mathcal{H}}^{2}}.
	\end{align*}
\end{lem}
\begin{proof}
	See \cite[Lemma~A.5]{KimKimKwon2024arxiv}.
\end{proof}

\begin{lem}[Subcoercivity for $\calB_2$ on $\dot\calH^2$]\label{lem:subcoer calB2}
	For $f\in \mathcal{S}$, we have
	\begin{align*}
		\|\calB_2f\|_{L^2}+\|{\bf1}_{|x|\lesssim 1}f\|_{L^2}+\|{\bf1}_{|x|\lesssim 1}\partial_x f\|_{L^2}\sim \|f\|_{\dot\calH^2}.
	\end{align*}
\end{lem}
\begin{proof}
	$(\lesssim)$: This is obvious. $(\gtrsim)$: We have
	\begin{align*}
		\|\calB_2f\|_{L^2}=\|B_Q^*\calB_2f\|_{L^2}=\|A_Q^*A_Qf\|_{L^2}
	\end{align*}
	with $A_Q^*A_Q=-\partial_{xx}+\tfrac{1}{4}Q^4-Q|D|Q$. We observe that $Q^4$ term is perturbative since
	\begin{align*}
		\|Q^4f\|_{L^2}\lesssim \|{\bf1}_{|y|\lesssim r_0}f\|_{L^2}+r_0^{-2}\|f\|_{\dot\calH^2}.
	\end{align*}
	For $Q|D|Q$ part, using \eqref{eq:CommuteHilbert} with $(1+y^2)\langle y \rangle^{-2}=1$, we have
	\begin{align*}
		\|Q|D|(Qf)\|_{L^2}&=\|xQ\calH[\tfrac{x}{1+x^2}\partial_x(Qf)]
		+Q\calH[\tfrac{1}{1+x^2}\partial_x(Qf)]
		-Q\cdot \tfrac{1}{\pi}{\textstyle\int}\tfrac{x}{1+x^2}\partial_x(Qf)\|_{L^2}\\
        &\lesssim
		\|Q^2\partial_x f\|_{L^2}
		+\|Q^3f\|_{L^2}+\|Q^{\frac12+}\|_{L^2}(\|Q^{\frac32-}\partial_xf\|_{L^2}+\|Q^{\frac52-}f\|_{L^2})
		\\
		&\lesssim 
		\|{\bf1}_{|x|\lesssim r_0}f\|_{L^2}+\|{\bf1}_{|x|\lesssim r_0}\partial_xf\|_{L^2}+r_0^{-\frac12+}\|f\|_{\dot\calH^2}.
	\end{align*}
	Then, by Lemma \ref{lem:H2comparison}, taking $r_0\gg 1$ sufficiently large, we conclude the proof.
\end{proof}

\begin{lem}[Comparison of $\dot{\mathcal{H}}^{3}$ and $\dot{H}^{3}$]
	\label{lem:H3comparison} We have 
	\begin{align*}
		\|\partial_{xxx}f\|_{L^{2}}
        +\|{\bf 1}_{|x|\lesssim1}(|\partial_{xx}f|+|\partial_{x}f|+|f|)\|_{L^{2}}
        \sim\|f\|_{\dot{\mathcal{H}}^{3}}.
	\end{align*}
\end{lem}
\begin{proof}
    $(\lesssim)$ directly comes from the definition of $\dot\calH^3$. $(\gtrsim)$: We first observe
    \begin{align*}
        \left\Vert \langle x\rangle^{-j}\partial_{x}^{3-j}f\right\Vert_{L^{2}}
        \lesssim
        \left\Vert {\bf 1}_{|x|\lesssim 1}\partial_{x}^{3-j}f\right\Vert_{L^{2}}
        +
        \left\Vert {\bf 1}_{|x|\gtrsim 1} |x|^{-j}\partial_{x}^{3-j} f\right\Vert_{L^{2}}.
    \end{align*}
    Thus, the proof reduces to, for $j=1,2,3$,
    \begin{align*}
        \left\Vert {\bf 1}_{|x|\gtrsim 1} | x|^{-j}\partial_{x}^{3-j} f\right\Vert_{L^{2}}
        \lesssim
        \|\partial_{xxx}f\|_{L^{2}}+\|{\bf 1}_{|x|\lesssim1}(|\partial_{xx}f|+|\partial_{x}f|+|f|)\|_{L^{2}}.
    \end{align*}
    Applying Lemma~\ref{lem:weight hardy} with $\varphi=r^{-2j+1}$, we obtain
    \begin{align*}
        \left\Vert {\bf 1}_{|x|\geq r_0} | x|^{-j}\partial_{x}^{3-j} f\right\Vert_{L^{2}}^2
        \lesssim 
        \left\Vert {\bf 1}_{|x|\geq r_0} | x|^{-j+1}\partial_{x}^{4-j} f\right\Vert_{L^{2}}^2
        +|r_0^{-2j+1}\partial_{x}^{3-j} f(r_0)|^2.
    \end{align*}
    Iterating this, we arrive at
    \begin{align*}
        \left\Vert {\bf 1}_{|x|\geq r_0} | x|^{-j}\partial_{x}^{3-j} f\right\Vert_{L^{2}}^2
        \lesssim 
        \|\partial_{xxx}f\|_{L^{2}}
        +\textstyle{\sum_{i=1}^3}|r_0^{-2i+1}\partial_{x}^{3-i} f(r_0)|^2.
    \end{align*}
    Averaging the boundary term, we conclude the proof by obtaining
    \[  \left\Vert {\bf 1}_{|x|\gtrsim 1} | x|^{-j}\partial_{x}^{3-j} f\right\Vert_{L^{2}}^2
        \lesssim \|\partial_{xxx}f\|_{L^{2}}
        +\|{\bf 1}_{|x|\sim 1}(|\partial_{xx}f|+|\partial_{x}f|+|f|)\|_{L^{2}}. \qedhere\]
\end{proof}

\begin{lem}[Subcoercivity for $\calL_3$ on $\dot\calH^3_e$]\label{lem:calL3 subcoer}
	For $f\in \mathcal{S}_e$, we have
	\begin{align*}
		\|\calL_3f\|_{L^2}
        +\|{\bf 1}_{|x|\lesssim1}(|\partial_{xx}f|+|\partial_{x}f|+|f|)\|_{L^{2}}
        \sim \|f\|_{\dot\calH^3}.
	\end{align*}
\end{lem}
\begin{proof}
    We recall $\calL_{3}\coloneqq (\calL^{1}_3,\calL^{2}_3)$ with \eqref{eq:newdef BQLQ}. 
    By Lemma~\ref{lem:H3comparison}, it suffices to show
    \begin{align}
        \|\partial_{xxx}f\|_{L^{2}}\lesssim \|\calL_3f\|_{L^2}
        +\|{\bf 1}_{|x|\lesssim1}(|\partial_{xx}f|+|\partial_{x}f|+|f|)\|_{L^{2}}. \label{eq:calL3 subcoer goal}
    \end{align}
    We observe that
    \begin{align*}
        \calL^{1}_3f
        =\,&B_Q\partial_x^{3}f+\partial_x^{2}(\tfrac{1}{\langle x\rangle}f)
        +[\partial_x^3(\tfrac{x}{\langle x\rangle}f)-\tfrac{x}{\langle x\rangle}\partial_x^3f]
        -\calH[\partial_x^{3}(\tfrac{1}{\langle x\rangle}f)-\tfrac{1}{\langle x\rangle}\partial_x^{3}f].
    \end{align*}
    Since $\partial_x(\frac{x}{\langle x\rangle})=\frac{1}{\langle x\rangle^3}$, we derive
    \begin{align*}
        &\|\partial_x^3(\tfrac{x}{\langle x\rangle}f)-\tfrac{x}{\langle x\rangle}\partial_x^3f-\calH[\partial_x^{3}(\tfrac{1}{\langle x\rangle}f)-\tfrac{1}{\langle x\rangle}\partial_x^{3}f]\|_{L^2}\\
        &\lesssim \|{\bf 1}_{|x|\lesssim r_0}(|\partial_{xx}f|+|\partial_{x}f|+|f|)\|_{L^{2}}+r_0^{-1}\|f\|_{\dot{\mathcal{H}}^{3}}.
    \end{align*}
    Thus, the real part of \eqref{eq:calL3 subcoer goal} is reduced to showing
    \begin{align*}
        \|\partial_{xxx}f\|_{L^{2}}
        \lesssim \|B_Q\partial_x^{3}f+\partial_x^{2}(\tfrac{1}{\langle x\rangle}f)\|_{L^2}
        +\|{\bf 1}_{|x|\lesssim1}(|\partial_{xx}f|+|\partial_{x}f|+|f|)\|_{L^{2}}.
    \end{align*}
    Since $f$ is even (so $\partial_x^{3}f$ is odd), \eqref{eq:BQBQstar equal I} yields
    \begin{align}
        \|B_Q\partial_x^{3}f\|_{L^2}=\|\partial_x^{3}f\|_{L^2}. \label{eq:BQ H3 appendix}
    \end{align}
    Hence, we can conclude the real part of \eqref{eq:calL3 subcoer goal} by showing
    \begin{align}
        \|\partial_x^{2}(\tfrac{1}{\langle x\rangle}f)\|_{L^2}
        \lesssim 
        \|\partial_x^{3}f\|_{L^2}
        +\|{\bf 1}_{|x|\lesssim1}(|\partial_{xx}f|+|\partial_{x}f|+|f|)\|_{L^{2}}. \label{eq:calL3 subcoer real 1}
    \end{align}
    As the proof of Lemma~\ref{lem:H3comparison}, we obtain \eqref{eq:calL3 subcoer real 1} by using Lemma~\ref{lem:weight hardy}. 
    
    For the imaginary part of \eqref{eq:calL3 subcoer goal}, $\partial_x(\langle x\rangle g)=\tfrac{x}{\langle x\rangle} g+\langle x\rangle \partial_x g$ implies
    \begin{align*}
        \calL^{2}_3g&=B_Q\partial_x^3 g+\partial_x^{2}(\tfrac{1}{\langle x\rangle}g) +G_{\text{pert}},\\
        G_{\text{pert}}&=
        -\partial_x^{2}(\tfrac{1}{\langle x\rangle^3}g) 
        +[\partial_x^{2}(\tfrac{x}{\langle x\rangle} \partial_x g)- \tfrac{x}{\langle x\rangle} \partial_x^3 g]-\calH[\partial_x^{2}\{\tfrac{1}{\langle x\rangle^2}\partial_x(\langle x\rangle g)\}-\tfrac{1}{\langle x\rangle}\partial_x^{3}g].
    \end{align*}
    Again using $\partial_x(\tfrac{x}{\langle x\rangle})=\tfrac{1}{\langle x\rangle^3}$, we obtain
    \begin{align*}
        \|G_{\text{pert}}\|_{L^2}
        \lesssim
        \|{\bf 1}_{|x|\lesssim r_0}(|\partial_{xx}g|+|\partial_{x}g|+|g|)\|_{L^{2}}+r_0^{-1}\|g\|_{\dot{\mathcal{H}}^{3}}.
    \end{align*}
    Thus, the proof is reduced to
    \begin{align*}
        \|\partial_{xxx}g\|_{L^{2}}
        \lesssim \|B_Q\partial_x^{3}g+\partial_x^{2}(\tfrac{1}{\langle x\rangle}g)\|_{L^2}
        +\|{\bf 1}_{|x|\lesssim1}(|\partial_{xx}g|+|\partial_{x}g|+|g|)\|_{L^{2}}.
    \end{align*}
    By \eqref{eq:BQ H3 appendix} with $g$ in place of $f$, we deduce \eqref{eq:calL3 subcoer real 1}, which completes the proof.
\end{proof}

\begin{lem}[Subcoercivity for $\calB_3$ on $\dot\calH^3_e$]\label{lem:subcoer calB3}
	For $f\in \mathcal{S}_e$, we have
	\begin{align*}
		\|\calB_3f\|_{L^2}
        +\|{\bf 1}_{|x|\lesssim1}(|\partial_{xx}f|+|\partial_{x}f|+|f|)\|_{L^{2}}
        \sim \|f\|_{\dot\calH^3}.
	\end{align*}
\end{lem}
\begin{proof}
    From the definition of $\calB_3$ \eqref{eq:def Bj}, we have
    \begin{align*}
        \calB_3f=B_Q\partial_{x}^3f
        +
        [\partial_{x}^3(\tfrac{x}{\langle x\rangle}f)-\tfrac{x}{\langle x\rangle}\partial_{x}^3f]
        -\calH[\partial_{x}^3(\tfrac{1}{\langle x\rangle}f)-\tfrac{1}{\langle x\rangle}\partial_{x}^3f].
    \end{align*}
    As in the proof of Lemma~\ref{lem:calL3 subcoer}, the identity $\partial_x(\tfrac{x}{\langle x\rangle})=\tfrac{1}{\langle x\rangle^3}$ implies that the last two terms of $ \calB_3f$ is perturbative. Thus, the proof is done by \eqref{eq:BQ H3 appendix}.
\end{proof}

\section{Sketch of proof for \texorpdfstring{\eqref{eq:frac 32 bound}}{(3/2 bound)}} \label{appendix B}
%We provide a sketch of proof for the case $L=0$, that is, we show $\lambda(t)\lesssim (T-t)^{\frac{3}{2}}$ for $H^1$ finite-time blow-up solutions.
%\begin{prop}
%    Let $u(t)$ be a solution to \eqref{CMdnls} satisfying Assumption~\ref{assum:one soliton}, with blow-up time $0<T<+\infty$. Then, the solution $u(t)$ satisfies a universal blow-up bound for energy $H^1$ solutions, i.e., if $u(t)\in H^1$, the scale parameter $\lambda(t)$ in the decomposition \eqref{eq:one soliton cmdnls} satisfies
%    \begin{align*}
%        \lambda(t)\lesssim (T-t)^{\frac32}.
%    \end{align*}
%\end{prop}
\begin{proof}[Proof of \eqref{eq:frac 32 bound}]
    Thanks to Section~\ref{sec:gague correspond}, it suffices to show the corresponding result for \eqref{CMdnls-gauged}. By Proposition~\ref{prop:one soliton}, the solution $v(t)$ to \eqref{CMdnls-gauged} satisfies
    \begin{align*}
       v=[Q+\wt\eps]_{\lambda,\gmm,x},\quad \lambda^2E_1\gtrsim \|\wt\eps\|_{\dot\calH^1}^2.
    \end{align*}
    Then, the evolution equation of $\wt\eps$ is written as follows:
    \begin{align}
		\partial_s\wt\eps=\tfrac{\lambda_s}{\lambda}\Lambda(Q+\wt\eps)-\gamma_si(Q+\wt\eps)+\tfrac{x_{s}}{\lambda}\partial_y(Q+\wt\eps) -i(\calL_Q\wt\eps +R_Q(\wt\eps)) \label{eq:wt eps appendix b}
	\end{align}
	where $R_Q(\eps)$ is defined by
	%\begin{align*}
	%	R_Q(\wt\eps)=\,&Q^{2}\Re(Q\wt\eps)\wt\eps+\tfrac{1}{2}Q^3|\wt\eps|^2+Q(\Re(Q\wt\eps))^{2}	+\tfrac{1}{2}Q^2|\wt\eps|^2\wt\eps+(\Re(Q\wt\eps))^{2}\wt\eps\\
		%&+Q|\wt\eps|^{2}\Re(Q\wt\eps)+|\wt\eps|^{2}\Re(Q\wt\eps)\wt\eps+\tfrac{1}{4}Q|\wt\eps|^4+\tfrac{1}{4}|\wt\eps|^4\wt\eps-Q|D|(|\wt\eps|^2)\\
        %&-2\wt\eps|D|\Re(Q\wt\eps)-\wt\eps|D|(|\wt\eps|^2).
	%\end{align*}
    \begin{align*}
		R_Q(\wt\eps)=\Re(Q\wt\eps)( Q^{2}\wt\eps +w\Re(\ol{w}\wt\eps))
		+\tfrac{w}{2}|\wt\eps|^2(Q^2+\tfrac{|\wt\eps|^2}{2})
		-w|D|(|\wt\eps|^2)-2\wt\eps|D|\Re(Q\wt\eps)
	\end{align*}
    with $w=Q+\wt\eps$. Computing $\partial_s(\wt\eps, \calZ_k)_r$ with a standard argument, we obtain
    \begin{align*}
        |\tfrac{\lambda_s}{\lambda}|+|\gamma_s|+|\tfrac{x_s}{\lambda}|
        \lesssim_M \lambda \sqrt{E_1}.
    \end{align*}
    %Now, we claim $\lambda(t)\lesssim_{M,E_1} (T-t)^{\frac32}$.
    Let $\bm{\eps}=[\wt\eps]_{\lambda,\gamma,0}$ and $f_{\lambda,\gamma}=[f]_{\lambda,\gamma}\coloneqq [f]_{\lambda,\gamma,0}$. We rewrite \eqref{eq:wt eps appendix b} as
    \begin{align*}
        (\tfrac{1}{2}\tfrac{\lambda_t}{\lambda}+i\gamma_t)Q_{\lambda,\gamma}+\partial_t\bm{\eps}
        =\tfrac{\lambda_t}{\lambda}[\Lambda_{-\frac{1}{2}}Q]_{\lambda,\gamma}+x_t\partial_y(Q_{\lambda,\gamma}+\bm{\eps}) -i\lambda^{-2}[\calL_Q\wt\eps +R_Q(\wt\eps)]_{\lambda,\gamma}.
    \end{align*}
    Then, we have
    \begin{align*}
        \|\tfrac{1}{\lambda+|y|}\bm{\eps}\|_{L^2}+\|\bm{\eps}\|_{\dot H^1}\lesssim 1,
        \quad |\lambda_t|+\lambda|\gamma_t|+|x_t|\lesssim_{E_1,M} 1,\quad \lmb(t)\aleq T-t.
    \end{align*}
    We take the inner product with $\psi=\lambda (y^2Q)(\frac{\cdot}{\lambda})\chi_{R(t)}$ where $R(t)=\sqrt{T-t}$. Then,
    \[
(\tfrac{1}{2}\tfrac{\lmb_{t}}{\lmb}+i\gmm_{t})(Q_{\lmb,\gmm},\lmb[y^{2}Q](\tfrac{\cdot}{\lmb})\chi_{R})_{c}=\tfrac{d}{dt}(\sqrt{\lmb}e^{i\gmm})\cdot\textstyle\int[y^{2}Q^{2}](\frac{\cdot}{\lmb})\chi_{R}.
\]
For $\rd_t \bm{\eps}$ term, we differentiate by parts this with $A(t)\coloneqq\int[y^{2}Q^{2}](\frac{\cdot}{\lmb})\chi_{R}\sim R$:
\begin{align*}
(\rd_{t}\bm{\eps},\psi)_{c} &= A(t)\tfrac{d}{dt}[A(t)^{-1}(\bm{\eps},\psi)_{c}]- (\bm{\eps},\rd_{t}\psi)_{c}+\tfrac{\rd_{t}A(t)}{A(t)}(\bm{\eps},\psi)_{c}.
\end{align*}
Using $\psi = \lambda (y^2Q)(\frac{\cdot}{\lambda})\chi_{R(t)}$ with $R(t)=\sqrt{T-t}$, we can check that
\begin{align*}
(\bm{\eps},\psi)_{c} & \lesssim  \|\tfrac{1}{\lambda+|y|}\bm{\eps}\|_{L^2} \| (\lambda^2+|y|^2)\chi_R \|_{L^2} \lesssim R^{\frac52},\\
(\bm{\eps},\rd_{t}\psi)_{c} & \lesssim  \|\tfrac{1}{\lambda+|y|}\bm{\eps}\|_{L^2} \|(\lambda+|y|) (|\lmb_{t}|\chi_R+\tfrac{|y|}{T-t}\chf_{|y|\sim R})  \|_{L^2}\lesssim R^{\frac12},\\ %\textstyle\int_0^{2R} |\bm{\eps}|(|\lmb_{t}|+\tfrac{|y|}{T-t}\chf_{|y|\ge R})
\rd_{t}A(t) &\lesssim \textstyle\int_0^{2R} \frac{\lmb|\lmb_{t}||y|^2}{(\lmb^{2}+|y|^{2})^2}+\frac{|R_t||y|^{2}}{R(\lambda^2+|y|^2)}\chf_{|y|\ge R} \lesssim R^{-1}.
\end{align*}
Due to the decay cancellation from $\Lambda_{-\frac{1}{2}}Q=(1+y\partial_y)Q=\frac{1}{2}Q^3$, we have
\begin{align*}
    \tfrac{\lmb_{t}}{\lmb}([(1+y\rd_{y})Q]_{\lmb,\gmm},\psi)_{c}&=\tfrac{\lmb_{t}}{2\sqrt{\lmb}}e^{i\gmm}\!\textstyle\int [y^2Q^4](\frac{\cdot}{\lmb})\chi_{R} \lesssim \lmb^{\frac{1}{2}}|\lmb_{t}| \textstyle\int_0^{2R} \tfrac{(y/
    \lmb)^2}{(1+(y/\lmb)^2)^2} d(\tfrac{y}{\lmb}) \lesssim \lmb^{\frac{1}{2}}.
\end{align*}
Moreover, $\rd_{y}(Q_{\lmb,\gmm})=0$ implies 
\begin{align*}
    (x_{t}\rd_{y}(Q_{\lmb,\gmm}+\bm{\eps}),\psi)_{c}&=(x_{t}\rd_{y}\bm{\eps},\lmb[y^{2}Q](\tfrac{\cdot}{\lmb})\chi_{R})_{c}\lesssim |x_t|  \|\bm{\eps}\|_{\dot H^1} \| |y| \chi_R \|_{L^2}\lesssim R^{\frac{3}{2}}.
\end{align*}
For the linear term, using $\calL_Q=L_Q^*L_Q$, we have %\textcolor{red}{Key point} 2. decay cancellation form the choice of $\psi$ in $(g, \partial_t(\psi))_r$. Here,
\begin{align*}
      (i\lambda^{-2}[\calL_Q\wt\eps ]_{\lambda,\gamma},\psi)_c
        &= \lambda^{-\frac12}e^{i\gamma}(i\calL_Q\wt\eps,y^2Q\chi_{R/\lambda})_c \lesssim \lambda^{-\frac12}\|\wt\eps\|_{\dot\calH^1}\|y^2Q\chi_{R/\lambda}\|_{\dot\calH^1} \lesssim  R^{\frac12}.
\end{align*}
    For the nonlinear term, we can easily check that
    \begin{align*}
        (i\lambda^{-2}[R_Q(\wt\eps)]_{\lambda,\gamma},\psi)_c=\lambda^{-\frac12}e^{i\gamma}(iR_Q(\wt\eps),y^2Q\chi_{R/\lambda})_c\lesssim \lambda^{-\frac12}(\lambda+\lambda^{\frac12}R^{\frac32}).
    \end{align*}
    Here, we used the following bounds to estimate each term of $R_Q(\wt\eps)$:
    \begin{align*}
        \| Q\wt\eps \|_{L^2}+ \| \wt\eps_y \|_{L^2}\lesssim \lambda,\quad  \| \wt\eps \|_{L^{\infty}}\lesssim \lambda^{\frac{1}{2}},\quad \|w y^2Q \chi_{R/\lmb} \|_{L^2}\lesssim (\tfrac{R}{\lambda})^{\frac12}+\tfrac{R^2}{\lambda}.
    \end{align*}
    %\begin{align*}
     %   (\wt\eps|D|\Re(Q\wt\eps),y^2Q\chi_{R/\lambda})_c&=
      %  (\Re(Q\wt\eps),|D|(\wt\eps y^2Q\chi_{R/\lambda}))_c \lesssim \lambda\| (\wt\eps y^2Q\chi_{R/\lambda})_y \|_{L^2}\lesssim \lambda,
    %\end{align*}
    %and other terms are standard by using H\"older and Gagliardo--Nirenberg inequalities.
Gathering the above and dividing by $A(t)$ and using $R(t)=\sqrt{T-t}$, we conclude
\[
\tfrac{d}{dt}(\sqrt{\lmb}e^{i\gmm}+O((T-t)^{3/4}))=O((T-t)^{-1/4}). \qedhere
\]

\end{proof}

\bibliographystyle{abbrv} 
\bibliography{referenceCM}

@article{JendrejLawrieRodriguez2022ASENS,
	author = {Jendrej, Jacek and Lawrie, Andrew and Rodriguez, Casey},
	doi = {10.24033/asens.2514},
	fjournal = {Annales Scientifiques de l'\'{E}cole Normale Sup\'{e}rieure. Quatri\`eme S\'{e}rie},
	issn = {0012-9593},
	journal = {Ann. Sci. \'{E}c. Norm. Sup\'{e}r. (4)},
	mrclass = {58E20 (35A21 37J51 58E05)},
	mrnumber = {4468859},
	number = {4},
	pages = {1135--1198},
	title = {Dynamics of bubbling wave maps with prescribed radiation},
	volume = {55},
	year = {2022}}

@article{KriegerSchlagTataru2008Invent,
	author = {Krieger, Joachim and Schlag, Wilhelm and Tataru, Daniel},
	doi = {10.1007/s00222-007-0089-3},
	fjournal = {Inventiones Mathematicae},
	issn = {0020-9910},
	journal = {Invent. Math.},
	mrclass = {58J45 (35B40 35L70)},
	mrnumber = {2372807},
	mrreviewer = {Michael Ruzhansky},
	number = {3},
	pages = {543--615},
	title = {Renormalization and blow up for charge one equivariant critical wave maps},
	url = {https://doi.org/10.1007/s00222-007-0089-3},
	volume = {171},
	year = {2008}}

@article{KriegerSchlagTataru2009Duke,
	author = {Krieger, Joachim and Schlag, Wilhelm and Tataru, Daniel},
	doi = {10.1215/00127094-2009-005},
	fjournal = {Duke Mathematical Journal},
	issn = {0012-7094},
	journal = {Duke Math. J.},
	mrclass = {58J45 (35B40 35L70)},
	mrnumber = {2494455},
	mrreviewer = {Michael Ruzhansky},
	number = {1},
	pages = {1--53},
	title = {Slow blow-up solutions for the {$H^1(\mathbb{R}^3)$} critical focusing semilinear wave equation},
	url = {https://doi.org/10.1215/00127094-2009-005},
	volume = {147},
	year = {2009}}

@article{OrtolevaPerelman2013,
	author = {Ortoleva, C. and Perelman, G.},
	doi = {10.1090/S1061-0022-2014-01290-3},
	fjournal = {Rossi\u{\i}skaya Akademiya Nauk. Algebra i Analiz},
	issn = {0234-0852},
	journal = {Algebra i Analiz},
	mrclass = {35Q55 (35B44)},
	mrnumber = {3114854},
	mrreviewer = {Amin\ Esfahani},
	number = {2},
	pages = {162--192},
	title = {Nondispersive vanishing and blow up at infinity for the energy critical nonlinear {S}chr\"{o}dinger equation in {$\mathbb{R}^3$}},
	url = {https://doi.org/10.1090/S1061-0022-2014-01290-3},
	volume = {25},
	year = {2013}}

@article{Perelman2014CMP,
	author = {Perelman, Galina},
	doi = {10.1007/s00220-014-1916-1},
	fjournal = {Communications in Mathematical Physics},
	issn = {0010-3616},
	journal = {Comm. Math. Phys.},
	mrclass = {58J35 (35B44 35K45)},
	mrnumber = {3215578},
	mrreviewer = {Atanas G. Stefanov},
	number = {1},
	pages = {69--105},
	title = {Blow up dynamics for equivariant critical {S}chr\"{o}dinger maps},
	url = {https://doi.org/10.1007/s00220-014-1916-1},
	volume = {330},
	year = {2014}}

@article{Schmid2023arXiv,
	author = {Schmid, Tobias},
	journal = {preprint arXiv:2308.01883},
	title = {Blow up dynamics for the 3D energy-critical Nonlinear {S}chr{\"o}dinger equation},
	year = {2023}}

@article{JeongKim2024arXiv,
      title={Quantized blow-up dynamics for {C}alogero--{M}oser derivative nonlinear {S}chr\"odinger equation}, 
      author={Jeong, Uihyeon and Kim, Taegyu},
    journal = {preprint arXiv:2412.12518},
      year={2024}
}

@article{KimTKwon2024arxivSolResol,
      title={Soliton resolution for {C}alogero--{M}oser derivative nonlinear {S}chr\"{o}dinger equation}, 
      journal = {preprint arXiv:2408.12843},
      author={Kim, Taegyu and Kwon, Soonsik},
      year={2024}
}

@article{KimKwonOh2024arXiv,
	author = {Kim, Kihyun and Kwon, Soonsik and Oh, Sung-Jin},
	journal = {preprint arXiv:2409.13274},
	title = {Blow-up dynamics for radial self-dual {C}hern-{S}imons-{S}chr\"{o}dinger equation with prescribed asymptotic profile},
	year = {2024}}

@article{CFL2025arXiv,
	author = {Chapouto, Andreia and Forlano, Justin and Laurens, Thierry},
	journal = {preprint arXiv:2511.00302},
	title = {On the well-posedness of the intermediate nonlinear {S}chr\"{o}dinger equation on the line},
	year = {2025}}

@article{FrankRead2025arXiv,
	author = {Frank, Rupert L. and Read, Larry},
	journal = {preprint arXiv:2510.11403},
	title = {Jost solutions and direct scattering for the continuum {C}alogero-{M}oser equation},
	year = {2025}}

@article{Chen2025arXivNonzero,
	author = {Chen, Xi},
	journal = {preprint arXiv:2502.17968, to appear in SIAM J. Math. Anal.},
	title = {The defocusing {C}alogero-{M}oser derivative nonlinear {S}chr\"{o}dinger equation with a nonvanishing condition at infinity},
	year = {2025}}

@article{Chen2025arXivScatter,
	author = {Chen, Xi},
	journal = {preprint arXiv:2511.06432},
	title = {Scattering of the defocusing {C}alogero-{M}oser derivative nonlinear {S}chr\"{o}dinger equation},
	year = {2025}}

@article{ABIK2512arXiv,
      title={Global well-posedness for the generalized intermediate {NLS} with a nonvanishing condition at infinity}, 
      author={Akahori, Takafumi and Badreddine, Rana and Ibrahim, Slim and Kishimoto, Nobu},
      year={2025},
      journal = {preprint arXiv:2512.18998}
}

@article{HerreroVelazquez1992,
  title={A blow up result for semilinear heat equations in the supercritical case},
  author={Herrero, Miguel A. and Vel\'{a}zquez, Juan J. L.},
  journal={preprint},
  year={1992}
}

@article {HerreroVelazquez1994CRASPSI,
    AUTHOR = {Herrero, Miguel A. and Vel\'{a}zquez, Juan J. L.},
     TITLE = {Explosion de solutions d'\'{e}quations paraboliques
              semilin\'{e}aires supercritiques},
   JOURNAL = {C. R. Acad. Sci. Paris S\'{e}r. I Math.},
  FJOURNAL = {Comptes Rendus de l'Acad\'{e}mie des Sciences. S\'{e}rie I.
              Math\'{e}matique},
    VOLUME = {319},
      YEAR = {1994},
    NUMBER = {2},
     PAGES = {141--145},
      ISSN = {0764-4442},
   MRCLASS = {35B40 (35K57)},
  MRNUMBER = {1288393},
MRREVIEWER = {Sergey\ A.\ Vakulenko},
}

@article{FilippasHerreroVelazquez2000,
	author = {Filippas, Stathis and Herrero, Miguel A. and Vel\'{a}zquez, Juan J. L.},
	doi = {10.1098/rspa.2000.0648},
	fjournal = {The Royal Society of London. Proceedings. Series A. Mathematical, Physical and Engineering Sciences},
	issn = {1364-5021},
	journal = {R. Soc. Lond. Proc. Ser. A Math. Phys. Eng. Sci.},
	mrclass = {35K55 (35B33 35B40 35C20)},
	mrnumber = {1843848},
	mrreviewer = {Juli\'{a}n Aguirre},
	number = {2004},
	pages = {2957--2982},
	title = {Fast blow-up mechanisms for sign-changing solutions of a semilinear parabolic equation with critical nonlinearity},
	url = {https://doi.org/10.1098/rspa.2000.0648},
	volume = {456},
	year = {2000}}

@article {Mizoguchi2007Math.Ann.,
    AUTHOR = {Mizoguchi, Noriko},
     TITLE = {Rate of type {II} blowup for a semilinear heat equation},
   JOURNAL = {Math. Ann.},
  FJOURNAL = {Mathematische Annalen},
    VOLUME = {339},
      YEAR = {2007},
    NUMBER = {4},
     PAGES = {839--877},
      ISSN = {0025-5831,1432-1807},
   MRCLASS = {35K55 (35B40 35K15)},
  MRNUMBER = {2341904},
MRREVIEWER = {Juli\'{a}n\ Aguirre},
       DOI = {10.1007/s00208-007-0133-z},
       URL = {https://doi.org/10.1007/s00208-007-0133-z},
}

@article {Mizoguchi2011TranAMS,
    AUTHOR = {Mizoguchi, Noriko},
     TITLE = {Blow-up rate of type {II} and the braid group theory},
   JOURNAL = {Trans. Amer. Math. Soc.},
  FJOURNAL = {Transactions of the American Mathematical Society},
    VOLUME = {363},
      YEAR = {2011},
    NUMBER = {3},
     PAGES = {1419--1443},
      ISSN = {0002-9947,1088-6850},
   MRCLASS = {35K91 (35B07 35B44 35K20 57M07)},
  MRNUMBER = {2737271},
MRREVIEWER = {Christian\ Stinner},
       DOI = {10.1090/S0002-9947-2010-04784-1},
       URL = {https://doi.org/10.1090/S0002-9947-2010-04784-1},
}

@article{HadzicRaphael2019JEMS,
	author = {Had\v{z}i\'{c}, Mahir and Rapha\"{e}l, Pierre},
	doi = {10.4171/JEMS/904},
	fjournal = {Journal of the European Mathematical Society (JEMS)},
	issn = {1435-9855},
	journal = {J. Eur. Math. Soc. (JEMS)},
	mrclass = {35R35 (35A20 35K05)},
	mrnumber = {4012340},
	mrreviewer = {Rodica Luca},
	number = {11},
	pages = {3259--3341},
	title = {On melting and freezing for the 2{D} radial {S}tefan problem},
	url = {https://doi.org/10.4171/JEMS/904},
	volume = {21},
	year = {2019}}

@article {CollotGhoulMasmoudiNguyen2022CPAM,
    AUTHOR = {Collot, Charles and Ghoul, Tej-Eddine and Masmoudi, Nader and
              Nguyen, Van Tien},
     TITLE = {Refined description and stability for singular solutions of
              the 2{D} {K}eller-{S}egel system},
   JOURNAL = {Comm. Pure Appl. Math.},
  FJOURNAL = {Communications on Pure and Applied Mathematics},
    VOLUME = {75},
      YEAR = {2022},
    NUMBER = {7},
     PAGES = {1419--1516},
      ISSN = {0010-3640,1097-0312},
   MRCLASS = {35Q92 (92C17)},
  MRNUMBER = {4438587},
}

@article {GhoulIbrahimNguyen2018JDE,
    AUTHOR = {Ghoul, T. and Ibrahim, S. and Nguyen, V. T.},
     TITLE = {Construction of type {II} blowup solutions for the
              1-corotational energy supercritical wave maps},
   JOURNAL = {J. Differential Equations},
  FJOURNAL = {Journal of Differential Equations},
    VOLUME = {265},
      YEAR = {2018},
    NUMBER = {7},
     PAGES = {2968--3047},
      ISSN = {0022-0396,1090-2732},
   MRCLASS = {35L71 (35B40 35B44 35L51 35R01)},
  MRNUMBER = {3812220},
       DOI = {10.1016/j.jde.2018.04.058},
       URL = {https://doi.org/10.1016/j.jde.2018.04.058},
}

@article {MerleRaphaelRodnianskiSzeftel2022Invention,
    AUTHOR = {Merle, Frank and Rapha\"el, Pierre and Rodnianski, Igor and
              Szeftel, Jeremie},
     TITLE = {On blow up for the energy super critical defocusing nonlinear
              {S}chr\"odinger equations},
   JOURNAL = {Invent. Math.},
  FJOURNAL = {Inventiones Mathematicae},
    VOLUME = {227},
      YEAR = {2022},
    NUMBER = {1},
     PAGES = {247--413},
      ISSN = {0020-9910,1432-1297},
   MRCLASS = {35Q55},
  MRNUMBER = {4359478},
       DOI = {10.1007/s00222-021-01067-9},
       URL = {https://doi.org/10.1007/s00222-021-01067-9},
}

@article {RS2011JAMSinhomo,
	AUTHOR = {Rapha\"el, Pierre and Szeftel, Jeremie},
	TITLE = {Existence and uniqueness of minimal blow-up solutions to an
	inhomogeneous mass critical {NLS}},
	JOURNAL = {J. Amer. Math. Soc.},
	FJOURNAL = {Journal of the American Mathematical Society},
	VOLUME = {24},
	YEAR = {2011},
	NUMBER = {2},
	PAGES = {471--546},
	ISSN = {0894-0347,1088-6834},
	MRCLASS = {35Q55 (35B35 35B44)},
	MRNUMBER = {2748399},
	MRREVIEWER = {Masahito\ Ohta},
	DOI = {10.1090/S0894-0347-2010-00688-1},
	URL = {https://doi.org/10.1090/S0894-0347-2010-00688-1},

}

@article{KimMerle2025CPAM,
author = {Kim, Kihyun and Merle, Frank},
title = {On classification of global dynamics for energy-critical equivariant harmonic map heat flows and radial nonlinear heat equation},
journal = {Communications on Pure and Applied Mathematics},
volume = {78},
number = {9},
pages = {1783-1842},
doi = {https://doi.org/10.1002/cpa.22253},
url = {https://onlinelibrary.wiley.com/doi/abs/10.1002/cpa.22253},
eprint = {https://onlinelibrary.wiley.com/doi/pdf/10.1002/cpa.22253},
year = {2025}
}

@article{Jeong2025Anal.PDE,
     TITLE = {Quantized slow blow-up dynamics for the energy-critical corotational wave map problem}, 
    AUTHOR = {Uihyeon Jeong},   
   JOURNAL = {Anal. PDE},
  FJOURNAL = {Analysis \& PDE},
    VOLUME = {18},
      YEAR = {2025},
    NUMBER = {10},
     PAGES = {2415--2480},
       DOI = {10.2140/apde.2025.18.2415},
       URL = {https://doi.org/10.2140/apde.2025.18.2415},
}

@article {GerardLenzmann2024CPAM,
    AUTHOR = {G\'erard, Patrick and Lenzmann, Enno},
     TITLE = {The {C}alogero-{M}oser derivative nonlinear {S}chr\"odinger
              equation},
   JOURNAL = {Comm. Pure Appl. Math.},
  FJOURNAL = {Communications on Pure and Applied Mathematics},
    VOLUME = {77},
      YEAR = {2024},
    NUMBER = {10},
     PAGES = {4008--4062},
      ISSN = {0010-3640,1097-0312},
   MRCLASS = {35Q55 (35C08 37K15 47B35)},
  MRNUMBER = {4814915},
MRREVIEWER = {Felipe\ Linares},
       DOI = {10.1002/cpa.22203},
       URL = {https://doi.org/10.1002/cpa.22203},
}

@article{KimKwon2019,
	author = {Kim, Kihyun and Kwon, Soonsik},
	doi = {10.1090/memo/1409},
	fjournal = {Memoirs of the American Mathematical Society},
	isbn = {978-1-4704-6120-1; 978-1-4704-7447-8},
	issn = {0065-9266,1947-6221},
	journal = {Mem. Amer. Math. Soc.},
	mrclass = {35Q55},
	mrnumber = {4574850},
	number = {1409},
	pages = {vi+128},
	title = {On pseudoconformal blow-up solutions to the self-dual {C}hern-{S}imons-{S}chr\"{o}dinger equation: existence, uniqueness, and instability},
	url = {https://doi.org/10.1090/memo/1409},
	volume = {284},
	year = {2023}}

@article{KimKwon2020blowup,
	author = {Kim, Kihyun and Kwon, Soonsik},
	doi = {10.1007/s40818-023-00147-8},
	fjournal = {Annals of PDE. Journal Dedicated to the Analysis of Problems from Physical Sciences},
	issn = {2524-5317,2199-2576},
	journal = {Ann. PDE},
	mrclass = {35B44 (35Q55)},
	mrnumber = {4564552},
	mrreviewer = {Van\ Duong\ Dinh},
	number = {1},
	pages = {Paper No. 6, 129},
	title = {Construction of blow-up manifolds to the equivariant self-dual {C}hern-{S}imons-{S}chr\"{o}dinger equation},
	url = {https://doi.org/10.1007/s40818-023-00147-8},
	volume = {9},
	year = {2023}}

@article {KimKwonOh2020blowup,
    AUTHOR = {Kim, Kihyun and Kwon, Soonsik and Oh, Sung-Jin},
     TITLE = {Blow-up dynamics for smooth finite energy radial data
              solutions to the self-dual {C}hern-{S}imons-{S}chr\"odinger
              equation},
   JOURNAL = {Ann. Sci. \'Ec. Norm. Sup\'er. (4)},
  FJOURNAL = {Annales Scientifiques de l'\'Ecole Normale Sup\'erieure.
              Quatri\`eme S\'erie},
    VOLUME = {57},
      YEAR = {2024},
    NUMBER = {6},
     PAGES = {1599--1691},
      ISSN = {0012-9593,1873-2151},
   MRCLASS = {35B44 (35Q55)},
  MRNUMBER = {4862502},
}

@article{AbanovBettelheimWiegmann2009FormalContinuum,
	author = {Abanov, Alexander G. and Bettelheim, Eldad and Wiegmann, Paul},
	doi = {10.1088/1751-8113/42/13/135201},
	fjournal = {Journal of Physics. A. Mathematical and Theoretical},
	issn = {1751-8113,1751-8121},
	journal = {J. Phys. A},
	mrclass = {37K10 (81R12)},
	mrnumber = {2485800},
	number = {13},
	pages = {135201, 24},
	title = {Integrable hydrodynamics of {C}alogero-{S}utherland model: bidirectional {B}enjamin-{O}no equation},
	url = {https://doi.org/10.1088/1751-8113/42/13/135201},
	volume = {42},
	year = {2009}}

@article{PelinovskyGrimshaw1995IntermediateDefocusingCMDNLSLaxPair,
	author = {Pelinovsky, Dmitry E. and Grimshaw, Roger H. J.},
	doi = {10.1063/1.530956},
	fjournal = {Journal of Mathematical Physics},
	issn = {0022-2488,1089-7658},
	journal = {J. Math. Phys.},
	mrclass = {35Q55},
	mrnumber = {1341985},
	mrreviewer = {Stanislav\ Z.\ Pakuliak},
	number = {8},
	pages = {4203--4219},
	title = {A spectral transform for the intermediate nonlinear {S}chr\"{o}dinger equation},
	url = {https://doi.org/10.1063/1.530956},
	volume = {36},
	year = {1995}}

@article{OlshanetskyPerelomov1976Invent,
	author = {Olshanetsky, M. A. and Perelomov, A. M.},
	doi = {10.1007/BF01418964},
	fjournal = {Inventiones Mathematicae},
	issn = {0020-9910,1432-1297},
	journal = {Invent. Math.},
	mrclass = {58F05},
	mrnumber = {426053},
	mrreviewer = {Mark\ Adler},
	number = {2},
	pages = {93--108},
	title = {Completely integrable {H}amiltonian systems connected with semisimple {L}ie algebras},
	url = {https://doi.org/10.1007/BF01418964},
	volume = {37},
	year = {1976}}

@article{Moser1975ClassicCalogerMoser,
	author = {Moser, J.},
	doi = {10.1016/0001-8708(75)90151-6},
	fjournal = {Advances in Mathematics},
	issn = {0001-8708},
	journal = {Advances in Math.},
	mrclass = {70.34},
	mrnumber = {375869},
	mrreviewer = {Dieter\ S.\ Schmidt},
	pages = {197--220},
	title = {Three integrable {H}amiltonian systems connected with isospectral deformations},
	url = {https://doi.org/10.1016/0001-8708(75)90151-6},
	volume = {16},
	year = {1975}}

@article{Calogero1971ClassicCalogerMoser,
	author = {Calogero, F.},
	doi = {10.1063/1.1665604},
	fjournal = {Journal of Mathematical Physics},
	issn = {0022-2488,1089-7658},
	journal = {J. Mathematical Phys.},
	mrclass = {81.34},
	mrnumber = {280103},
	mrreviewer = {Z.\ Jankovi\'{c}},
	pages = {419--436},
	title = {Solution of the one-dimensional {$N$}-body problems with quadratic and/or inversely quadratic pair potentials},
	url = {https://doi.org/10.1063/1.1665604},
	volume = {12},
	year = {1971}}

@article{CalogeroMarchioro1974ClassicCalogerMoser,
	author = {Calogero, F. and Marchioro, C.},
	doi = {10.1063/1.1666827},
	fjournal = {Journal of Mathematical Physics},
	issn = {0022-2488,1089-7658},
	journal = {J. Mathematical Phys.},
	mrclass = {81.35},
	mrnumber = {353869},
	mrreviewer = {H.\ Flaschka},
	pages = {1425--1430},
	title = {Exact solution of a one-dimensional three-body scattering problem with two-body and/or three-body inverse-square potentials},
	url = {https://doi.org/10.1063/1.1666827},
	volume = {15},
	year = {1974}}

@article {Badreddine2024PAA,
	AUTHOR = {Badreddine, Rana},
	TITLE = {On the global well-posedness of the {C}alogero-{S}utherland
	derivative nonlinear {S}chr\"odinger equation},
	JOURNAL = {Pure Appl. Anal.},
	FJOURNAL = {Pure and Applied Analysis},
	VOLUME = {6},
	YEAR = {2024},
	NUMBER = {2},
	PAGES = {379--414},
	ISSN = {2578-5885,2578-5893},
	MRCLASS = {35Q55 (37K10)},
	MRNUMBER = {4746420},
	DOI = {10.2140/paa.2024.6.379},
	URL = {https://doi.org/10.2140/paa.2024.6.379},
}

@article {Badreddine2024AIHPC,
    AUTHOR = {Badreddine, Rana},
     TITLE = {Traveling waves and finite gap potentials for the
              {C}alogero-{S}utherland derivative nonlinear {S}chr\"odinger
              equation},
   JOURNAL = {Ann. Inst. H. Poincar\'e{} C Anal. Non Lin\'eaire},
  FJOURNAL = {Annales de l'Institut Henri Poincar\'e{} C. Analyse Non
              Lin\'eaire},
    VOLUME = {42},
      YEAR = {2025},
    NUMBER = {4},
     PAGES = {1037--1092},
      ISSN = {0294-1449,1873-1430},
   MRCLASS = {35C07 (35Q55 37J38)},
  MRNUMBER = {4930525},
       DOI = {10.4171/aihpc/124},
       URL = {https://doi.org/10.4171/aihpc/124},
}

@article {Badreddine2024SIAM,
    AUTHOR = {Badreddine, Rana},
     TITLE = {Zero-dispersion limit of the {C}alogero-{M}oser derivative
              {NLS} equation},
   JOURNAL = {SIAM J. Math. Anal.},
  FJOURNAL = {SIAM Journal on Mathematical Analysis},
    VOLUME = {56},
      YEAR = {2024},
    NUMBER = {6},
     PAGES = {7228--7249},
      ISSN = {0036-1410,1095-7154},
   MRCLASS = {35Q55 (30H10 35B50)},
  MRNUMBER = {4816606},
       DOI = {10.1137/24M1646935},
       URL = {https://doi.org/10.1137/24M1646935},
}

@article{Gerard2023BOequexplicitFormula,
	author = {G\'{e}rard, Patrick},
	doi = {10.2140/tunis.2023.5.593},
	fjournal = {Tunisian Journal of Mathematics},
	issn = {2576-7658,2576-7666},
	journal = {Tunis. J. Math.},
	mrclass = {35Q53 (35C05 37K15 47B35)},
	mrnumber = {4662323},
	number = {3},
	pages = {593--603},
	title = {An explicit formula for the {B}enjamin-{O}no equation},
	url = {https://doi.org/10.2140/tunis.2023.5.593},
	volume = {5},
	year = {2023}}

@article {KLV2025CAMS,
    AUTHOR = {Killip, Rowan and Laurens, Thierry and Vi\c san, Monica},
     TITLE = {Scaling-critical well-posedness for continuum
              {C}alogero--{M}oser models on the line},
   JOURNAL = {Commun. Am. Math. Soc.},
  FJOURNAL = {Communications of the American Mathematical Society},
    VOLUME = {5},
      YEAR = {2025},
     PAGES = {284--320},
      ISSN = {2692-3688},
   MRCLASS = {35Q51 (37K10)},
  MRNUMBER = {4922705},
       DOI = {10.1090/cams/48},
       URL = {https://doi.org/10.1090/cams/48},
}

@article{MouraPilod2010CMDNLSLocal,
	author = {de Moura, Roger Peres and Pilod, Didier},
	fjournal = {Advances in Differential Equations},
	issn = {1079-9389},
	journal = {Adv. Differential Equations},
	mrclass = {35Q55 (35B30 35B65 35R09 76B03)},
	mrnumber = {2677424},
	mrreviewer = {Makoto\ Nakamura},
	number = {9-10},
	pages = {925--952},
	title = {Local well posedness for the nonlocal nonlinear {S}chr\"{o}dinger equation below the energy space},
	volume = {15},
	year = {2010}}

@article {HoganKowalski2024PAA,
    AUTHOR = {Hogan, James and Kowalski, Matthew},
     TITLE = {Turbulent threshold for continuum {C}alogero-{M}oser models},
   JOURNAL = {Pure Appl. Anal.},
  FJOURNAL = {Pure and Applied Analysis},
    VOLUME = {6},
      YEAR = {2024},
    NUMBER = {4},
     PAGES = {941--954},
      ISSN = {2578-5885,2578-5893},
   MRCLASS = {35Q55 (35B44 35C08 37K40)},
  MRNUMBER = {4844677},
MRREVIEWER = {Linjie\ Song},
       DOI = {10.2140/paa.2024.6.941},
       URL = {https://doi.org/10.2140/paa.2024.6.941},
}

@article{RodnianskiSterbenz2010,
	author = {Rodnianski, Igor and Sterbenz, Jacob},
	doi = {10.4007/annals.2010.172.187},
	fjournal = {Annals of Mathematics. Second Series},
	issn = {0003-486X,1939-8980},
	journal = {Ann. of Math. (2)},
	mrclass = {58E20 (35B44 35L70 35Q40 81T10)},
	mrnumber = {2680419},
	mrreviewer = {Andreas\ Gastel},
	number = {1},
	pages = {187--242},
	title = {On the formation of singularities in the critical {${\rm O}(3)$} {$\sigma$}-model},
	url = {https://doi.org/10.4007/annals.2010.172.187},
	volume = {172},
	year = {2010}}

@article{RaphaelRodnianski2012,
	author = {Rapha\"{e}l, Pierre and Rodnianski, Igor},
	doi = {10.1007/s10240-011-0037-z},
	fjournal = {Publications Math\'{e}matiques. Institut de Hautes \'{E}tudes Scientifiques},
	issn = {0073-8301,1618-1913},
	journal = {Publ. Math. Inst. Hautes \'{E}tudes Sci.},
	mrclass = {58E20 (35A20 35B44 35L70 58E15 81T13)},
	mrnumber = {2929728},
	mrreviewer = {Andreas\ Gastel},
	pages = {1--122},
	title = {Stable blow up dynamics for the critical co-rotational wave maps and equivariant {Y}ang-{M}ills problems},
	url = {https://doi.org/10.1007/s10240-011-0037-z},
	volume = {115},
	year = {2012}}

@article{GerardGrellier2015ExplictFormulaCubicSzego1,
	author = {G\'{e}rard, Patrick and Grellier, Sandrine},
	doi = {10.1090/S0002-9947-2014-06310-1},
	fjournal = {Transactions of the American Mathematical Society},
	issn = {0002-9947,1088-6850},
	journal = {Trans. Amer. Math. Soc.},
	mrclass = {37K15 (47B35)},
	mrnumber = {3301889},
	mrreviewer = {Julia\ B.\ Prada},
	number = {4},
	pages = {2979--2995},
	title = {An explicit formula for the cubic {S}zeg{\H{o}} equation},
	url = {https://doi.org/10.1090/S0002-9947-2014-06310-1},
	volume = {367},
	year = {2015}}

@article {GerardPushnitski2024CMP,
	AUTHOR = {G\'{e}rard, Patrick and Pushnitski, Alexander},
	TITLE = {The {C}ubic {S}zeg{\H{o}} {E}quation on the {R}eal {L}ine:
	{E}xplicit {F}ormula and {W}ell-{P}osedness on the {H}ardy
	{C}lass},
	JOURNAL = {Comm. Math. Phys.},
	FJOURNAL = {Communications in Mathematical Physics},
	VOLUME = {405},
	YEAR = {2024},
	NUMBER = {7},
	PAGES = {Paper No. 167},
	ISSN = {0010-3616,1432-0916},
	MRCLASS = {99-06},
	MRNUMBER = {4768537},
	DOI = {10.1007/s00220-024-05040-4},
	URL = {https://doi.org/10.1007/s00220-024-05040-4},
}

@article{MerleRaphaelRodnianski2013Invention,
	author = {Merle, Frank and Rapha\"{e}l, Pierre and Rodnianski, Igor},
	doi = {10.1007/s00222-012-0427-y},
	fjournal = {Inventiones Mathematicae},
	issn = {0020-9910,1432-1297},
	journal = {Invent. Math.},
	mrclass = {35K45 (35B44 35Q55 35Q60 58J99)},
	mrnumber = {3090180},
	mrreviewer = {Dimitrios\ E.\ Tzanetis},
	number = {2},
	pages = {249--365},
	title = {Blowup dynamics for smooth data equivariant solutions to the critical {S}chr\"{o}dinger map problem},
	url = {https://doi.org/10.1007/s00222-012-0427-y},
	volume = {193},
	year = {2013}}

@article{MerleRaphaelRodnianski2015CambJMath,
	author = {Merle, Frank and Rapha\"{e}l, Pierre and Rodnianski, Igor},
	doi = {10.4310/CJM.2015.v3.n4.a1},
	fjournal = {Cambridge Journal of Mathematics},
	issn = {2168-0930,2168-0949},
	journal = {Camb. J. Math.},
	mrclass = {35Q55 (35B44 35C06)},
	mrnumber = {3435273},
	mrreviewer = {Thomas\ Duyckaerts},
	number = {4},
	pages = {439--617},
	title = {Type {II} blow up for the energy supercritical {NLS}},
	url = {https://doi.org/10.4310/CJM.2015.v3.n4.a1},
	volume = {3},
	year = {2015}}

@article{RaphaelSchweyer2013CPAMHeat,
	author = {Rapha\"{e}l, Pierre and Schweyer, Remi},
	doi = {10.1002/cpa.21435},
	fjournal = {Communications on Pure and Applied Mathematics},
	issn = {0010-3640,1097-0312},
	journal = {Comm. Pure Appl. Math.},
	mrclass = {35R01 (35B44 35K91 53Cxx 58J35)},
	mrnumber = {3008229},
	mrreviewer = {Shu-Yu\ Hsu},
	number = {3},
	pages = {414--480},
	title = {Stable blowup dynamics for the 1-corotational energy critical harmonic heat flow},
	url = {https://doi.org/10.1002/cpa.21435},
	volume = {66},
	year = {2013}}

@article{RaphaelSchweyer2014AnalPDEHeatQuantized,
	author = {Rapha\"{e}l, Pierre and Schweyer, Remi},
	doi = {10.2140/apde.2014.7.1713},
	fjournal = {Analysis \& PDE},
	issn = {2157-5045,1948-206X},
	journal = {Anal. PDE},
	mrclass = {35K91 (35B44)},
	mrnumber = {3318739},
	mrreviewer = {Christopher\ P.\ Grant},
	number = {8},
	pages = {1713--1805},
	title = {Quantized slow blow-up dynamics for the corotational energy-critical harmonic heat flow},
	url = {https://doi.org/10.2140/apde.2014.7.1713},
	volume = {7},
	year = {2014}}

@article{Collot2017AnalPDE,
	author = {Collot, Charles},
	doi = {10.2140/apde.2017.10.127},
	fjournal = {Analysis \& PDE},
	issn = {2157-5045,1948-206X},
	journal = {Anal. PDE},
	mrclass = {35K91 (35B20 35B44 35C06 35K58)},
	mrnumber = {3611015},
	mrreviewer = {Hongwei\ Chen},
	number = {1},
	pages = {127--252},
	title = {Nonradial type {II} blow up for the energy-supercritical semilinear heat equation},
	url = {https://doi.org/10.2140/apde.2017.10.127},
	volume = {10},
	year = {2017}}

@article{Collot2018MemAmer,
	author = {Collot, Charles},
	doi = {10.1090/memo/1205},
	fjournal = {Memoirs of the American Mathematical Society},
	isbn = {978-1-4704-2813-6; 978-1-4704-4379-5},
	issn = {0065-9266,1947-6221},
	journal = {Mem. Amer. Math. Soc.},
	mrclass = {35L71 (35B44 58B99)},
	mrnumber = {3778126},
	mrreviewer = {Herbert\ Koch},
	number = {1205},
	pages = {v+163},
	title = {Type {II} blow up manifolds for the energy supercritical semilinear wave equation},
	url = {https://doi.org/10.1090/memo/1205},
	volume = {252},
	year = {2018}}

@article{HillairetRaphael2012AnalPDE,
	author = {Hillairet, Matthieu and Rapha\"{e}l, Pierre},
	doi = {10.2140/apde.2012.5.777},
	fjournal = {Analysis \& PDE},
	issn = {2157-5045,1948-206X},
	journal = {Anal. PDE},
	mrclass = {35L71 (35B44 35L15)},
	mrnumber = {3006642},
	mrreviewer = {Andr\'{e}\ Vicente},
	number = {4},
	pages = {777--829},
	title = {Smooth type {II} blow-up solutions to the four-dimensional energy-critical wave equation},
	url = {https://doi.org/10.2140/apde.2012.5.777},
	volume = {5},
	year = {2012}}

@article{delPinoMussoWeiZhang2020HeatQuantized,
	author = {del Pino, Manuel and Musso, Monica and Wei, Juncheng and Zhang, Qidi and Zhang, Yifu},
	journal = {preprint arXiv:2002.05765},
	title = {Type II Finite time blow-up for the three dimensional energy critical heat equation},
	year = {2020}}

@article{Bogomolnyi1976,
	author = {Bogomol'ny\u{\i}, E. B.},
	fjournal = {Jadernaja Fizika},
	issn = {0044-0027},
	journal = {Jadernaja Fiz.},
	mrclass = {81.49},
	mrnumber = {443719},
	number = {4},
	pages = {861--870},
	title = {The stability of classical solutions},
	volume = {24},
	year = {1976}}

@article{RaphaelSchweyer2014MathAnn,
	author = {Rapha\"{e}l, Pierre and Schweyer, R\'{e}mi},
	doi = {10.1007/s00208-013-1002-6},
	fjournal = {Mathematische Annalen},
	issn = {0025-5831,1432-1807},
	journal = {Math. Ann.},
	mrclass = {35M31 (35B20 35B35 35B44 92C17)},
	mrnumber = {3201901},
	mrreviewer = {Adrien\ Blanchet},
	number = {1-2},
	pages = {267--377},
	title = {On the stability of critical chemotactic aggregation},
	url = {https://doi.org/10.1007/s00208-013-1002-6},
	volume = {359},
	year = {2014}}

@article {MerleZaag2005Math.Ann.,
    AUTHOR = {Merle, Frank and Zaag, Hatem},
     TITLE = {Determination of the blow-up rate for a critical semilinear
              wave equation},
   JOURNAL = {Math. Ann.},
  FJOURNAL = {Mathematische Annalen},
    VOLUME = {331},
      YEAR = {2005},
    NUMBER = {2},
     PAGES = {395--416},
      ISSN = {0025-5831,1432-1807},
   MRCLASS = {35L70 (35B40)},
  MRNUMBER = {2115461},
MRREVIEWER = {Masahito\ Ohta},
       DOI = {10.1007/s00208-004-0587-1},
       URL = {https://doi.org/10.1007/s00208-004-0587-1},
}

@article {MerleZaag2003AJM,
    AUTHOR = {Merle, Frank and Zaag, Hatem},
     TITLE = {Determination of the blow-up rate for the semilinear wave
              equation},
   JOURNAL = {Amer. J. Math.},
  FJOURNAL = {American Journal of Mathematics},
    VOLUME = {125},
      YEAR = {2003},
    NUMBER = {5},
     PAGES = {1147--1164},
      ISSN = {0002-9327,1080-6377},
   MRCLASS = {35L70 (35B40 35L15)},
  MRNUMBER = {2004432},
MRREVIEWER = {Masahito\ Ohta},
       URL ={http://muse.jhu.edu/journals/american_journal_of_mathematics/v125/125.5merle.pdf},
}

@article{MerleRaphaelRodnianskiSzeftel2022-1Ann.Math.,
	author = {Merle, Frank and Rapha\"{e}l, Pierre and Rodnianski, Igor and Szeftel, Jeremie},
	doi = {10.4007/annals.2022.196.2.3},
	fjournal = {Annals of Mathematics. Second Series},
	issn = {0003-486X},
	journal = {Ann. of Math. (2)},
	mrclass = {35Q35 (34C37)},
	mrnumber = {4445442},
	number = {2},
	pages = {567--778},
	title = {On the implosion of a compressible fluid {I}: smooth self-similar inviscid profiles},
	url = {https://doi.org/10.4007/annals.2022.196.2.3},
	volume = {196},
	year = {2022}}

@article{MerleRaphaelRodnianskiSzeftel2022-2Ann.Math.,
	author = {Merle, Frank and Rapha\"{e}l, Pierre and Rodnianski, Igor and Szeftel, Jeremie},
	doi = {10.4007/annals.2022.196.2.4},
	fjournal = {Annals of Mathematics. Second Series},
	issn = {0003-486X},
	journal = {Ann. of Math. (2)},
	mrclass = {35Q35 (35B44)},
	mrnumber = {4445443},
	number = {2},
	pages = {779--889},
	title = {On the implosion of a compressible fluid {II}: singularity formation},
	url = {https://doi.org/10.4007/annals.2022.196.2.4},
	volume = {196},
	year = {2022}}

@article {ErdoganGurelTzirakis2018AIHPC,
    AUTHOR = {Erdo\u{g}an, M. B. and G\"{u}rel, T. B. and Tzirakis, N.},
     TITLE = {The derivative nonlinear {S}chr\"odinger equation on the half
              line},
   JOURNAL = {Ann. Inst. H. Poincar\'e{} C Anal. Non Lin\'eaire},
  FJOURNAL = {Annales de l'Institut Henri Poincar\'e{} C. Analyse Non
              Lin\'eaire},
    VOLUME = {35},
      YEAR = {2018},
    NUMBER = {7},
     PAGES = {1947--1973},
      ISSN = {0294-1449,1873-1430},
   MRCLASS = {35Q55 (35B30)},
  MRNUMBER = {3906860},
MRREVIEWER = {Luigi\ Forcella},
       DOI = {10.1016/j.anihpc.2018.03.006},
       URL = {https://doi.org/10.1016/j.anihpc.2018.03.006},
}

@article {CollianderKeelStaffilaniTakaokaTao2001SIAM,
    AUTHOR = {Colliander, J. and Keel, M. and Staffilani, G. and Takaoka, H.
              and Tao, T.},
     TITLE = {Global well-posedness for {S}chr\"odinger equations with
              derivative},
   JOURNAL = {SIAM J. Math. Anal.},
  FJOURNAL = {SIAM Journal on Mathematical Analysis},
    VOLUME = {33},
      YEAR = {2001},
    NUMBER = {3},
     PAGES = {649--669},
      ISSN = {0036-1410,1095-7154},
   MRCLASS = {35Q55 (35B35)},
  MRNUMBER = {1871414},
MRREVIEWER = {Olivier\ J.\ Goubet},
       DOI = {10.1137/S0036141001384387},
       URL = {https://doi.org/10.1137/S0036141001384387},
}

@article {ShaoWeiZhang2025Pi,
    AUTHOR = {Shao, Feng and Wei, Dongyi and Zhang, Zhifei},
     TITLE = {On blow-up for the supercritical defocusing nonlinear wave
              equation},
   JOURNAL = {Forum Math. Pi},
  FJOURNAL = {Forum of Mathematics. Pi},
    VOLUME = {13},
      YEAR = {2025},
     PAGES = {Paper No. e15, 59},
      ISSN = {2050-5086},
   MRCLASS = {35L71 (35B44 35C06 35L51 35Q35)},
  MRNUMBER = {4884918},
MRREVIEWER = {Pham\ Trieu\ Duong},
       DOI = {10.1017/fmp.2025.7},
       URL = {https://doi.org/10.1017/fmp.2025.7},
}

@article {ShaoWeiZhang2024arXiv,
      title={Self-similar imploding solutions of the relativistic Euler equations}, 
      author={Feng Shao and Dongyi Wei and Zhifei Zhang},
    journal={preprint arXiv:2403.11471},
      year={2024},
}

@article {BuckmasterChen2024arXiv,
      title={Blowup for the defocusing septic complex-valued nonlinear wave equation in $\mathbb{R}^{4+1}$}, 
      author={Tristan Buckmaster and Jiajie Chen},
      year={2024},
      journal={preprint arXiv:2410.15619},
}

@article {CollotMerleRaphael2020JAMS,
    AUTHOR = {Collot, Charles and Merle, Frank and Rapha\"el, Pierre},
     TITLE = {Strongly anisotropic type {II} blow up at an isolated point},
   JOURNAL = {J. Amer. Math. Soc.},
  FJOURNAL = {Journal of the American Mathematical Society},
    VOLUME = {33},
      YEAR = {2020},
    NUMBER = {2},
     PAGES = {527--607},
      ISSN = {0894-0347,1088-6834},
   MRCLASS = {35K91 (35B44 35K15 35K58)},
  MRNUMBER = {4073868},
MRREVIEWER = {Alain\ Brillard},
       DOI = {10.1090/jams/941},
       URL = {https://doi.org/10.1090/jams/941},
}

@article {Harada2020AIHPC,
    AUTHOR = {Harada, Junichi},
     TITLE = {A higher speed type {II} blowup for the five dimensional
              energy critical heat equation},
   JOURNAL = {Ann. Inst. H. Poincar\'e{} C Anal. Non Lin\'eaire},
  FJOURNAL = {Annales de l'Institut Henri Poincar\'e{} C. Analyse Non
              Lin\'eaire},
    VOLUME = {37},
      YEAR = {2020},
    NUMBER = {2},
     PAGES = {309--341},
      ISSN = {0294-1449,1873-1430},
   MRCLASS = {35K91 (35B44)},
  MRNUMBER = {4072807},
MRREVIEWER = {Daniele\ Andreucci},
       DOI = {10.1016/j.anihpc.2019.09.006},
       URL = {https://doi.org/10.1016/j.anihpc.2019.09.006},
}

@article {MartelPilod2024pisa,
    AUTHOR = {Martel, Yvan and Pilod, Didier},
     TITLE = {Finite point blowup for the critical generalized
              {K}orteweg--de {V}ries equation},
   JOURNAL = {Ann. Sc. Norm. Super. Pisa Cl. Sci. (5)},
  FJOURNAL = {Annali della Scuola Normale Superiore di Pisa. Classe di
              Scienze. Serie V},
    VOLUME = {25},
      YEAR = {2024},
    NUMBER = {1},
     PAGES = {371--425},
      ISSN = {0391-173X,2036-2145},
   MRCLASS = {35Q53 (35B44 35C08)},
  MRNUMBER = {4732642},
MRREVIEWER = {Rajkumar\ Roychoudhury},
}

@article{Kim2025JEMS,
	author = {Kim, Kihyun},
	journal = {J. Eur. Math. Soc. (JEMS), published online first},
	title = {Rigidity of smooth finite-time blow-up for equivariant self-dual {C}hern-{S}imons-{S}chr\"{o}dinger equation},
  url = {https://doi.org/10.4171/jems/1569},
	year = {2025}}

@article {KimKwonOh2025AJM,
    AUTHOR = {Kim, Kihyun and Kwon, Soonsik and Oh, Sung-Jin},
     TITLE = {Soliton resolution for equivariant self-dual
              {C}hern-{S}imons-{S}chr\"odinger equation in weighted
              {S}obolev class},
   JOURNAL = {Amer. J. Math.},
  FJOURNAL = {American Journal of Mathematics},
    VOLUME = {147},
      YEAR = {2025},
    NUMBER = {6},
     PAGES = {1475--1508},
      ISSN = {0002-9327,1080-6377},
   MRCLASS = {35},
  MRNUMBER = {4995127},
}

@article {MerleZaag2012AJM,
    AUTHOR = {Merle, Frank and Zaag, Hatem},
     TITLE = {Existence and classification of characteristic points at
              blow-up for a semilinear wave equation in one space dimension},
   JOURNAL = {Amer. J. Math.},
  FJOURNAL = {American Journal of Mathematics},
    VOLUME = {134},
      YEAR = {2012},
    NUMBER = {3},
     PAGES = {581--648},
      ISSN = {0002-9327,1080-6377},
   MRCLASS = {35L71 (35B44 35L15)},
  MRNUMBER = {2931219},
MRREVIEWER = {Mohammad\ A.\ Rammaha},
       DOI = {10.1353/ajm.2012.0021},
       URL = {https://doi.org/10.1353/ajm.2012.0021},
}

@article {CoteZaag2013CPAM,
    AUTHOR = {C\^ote, Rapha\"el and Zaag, Hatem},
     TITLE = {Construction of a multisoliton blowup solution to the
              semilinear wave equation in one space dimension},
   JOURNAL = {Comm. Pure Appl. Math.},
  FJOURNAL = {Communications on Pure and Applied Mathematics},
    VOLUME = {66},
      YEAR = {2013},
    NUMBER = {10},
     PAGES = {1541--1581},
      ISSN = {0010-3640,1097-0312},
   MRCLASS = {35L71 (35B44 35L15)},
  MRNUMBER = {3084698},
       DOI = {10.1002/cpa.21452},
       URL = {https://doi.org/10.1002/cpa.21452},
}

@article{MartelMerleRaphael2015pisa,
	author = {Martel, Yvan and Merle, Frank and Rapha\"{e}l, Pierre},
	fjournal = {Annali della Scuola Normale Superiore di Pisa. Classe di Scienze. Serie V},
	issn = {0391-173X},
	journal = {Ann. Sc. Norm. Super. Pisa Cl. Sci. (5)},
	mrclass = {35Q53 (35B33 35B40 35B44 35C06 35C08)},
	mrnumber = {3410473},
	mrreviewer = {Peter E. Zhidkov},
	number = {2},
	pages = {575--631},
	title = {Blow up for the critical g{K}d{V} equation {III}: exotic regimes},
	volume = {14},
	year = {2015}}

@article{MartelMerleRaphael2015JEMS,
	author = {Martel, Yvan and Merle, Frank and Rapha\"{e}l, Pierre},
	doi = {10.4171/JEMS/547},
	fjournal = {Journal of the European Mathematical Society (JEMS)},
	issn = {1435-9855},
	journal = {J. Eur. Math. Soc. (JEMS)},
	mrclass = {35Q53 (35B41 35B44 35C08 35Q51)},
	mrnumber = {3372073},
	mrreviewer = {Olivier J. Goubet},
	number = {8},
	pages = {1855--1925},
	title = {Blow up for the critical g{K}d{V} equation. {II}: {M}inimal mass dynamics},
	url = {https://doi.org/10.4171/JEMS/547},
	volume = {17},
	year = {2015}}

@article{MartelMerleRaphael2014Acta,
	author = {Martel, Yvan and Merle, Frank and Rapha\"{e}l, Pierre},
	doi = {10.1007/s11511-014-0109-2},
	fjournal = {Acta Mathematica},
	issn = {0001-5962,1871-2509},
	journal = {Acta Math.},
	mrclass = {35Q53 (35B33 35B44 35C08)},
	mrnumber = {3179608},
	mrreviewer = {Masayoshi\ Tsutsumi},
	number = {1},
	pages = {59--140},
	title = {Blow up for the critical generalized {K}orteweg--de {V}ries equation. {I}: {D}ynamics near the soliton},
	url = {https://doi.org/10.1007/s11511-014-0109-2},
	volume = {212},
	year = {2014}}

@article {KimKimKwon2024arxiv,
	author={Kim, Kihyun and Kim, Taegyu and Kwon, Soonsik},
	journal = {preprint arXiv:2404.09603, to appear in Mem. Amer. Math. Soc.},
	title={Construction of smooth chiral finite-time blow-up solutions to {C}alogero--{M}oser derivative nonlinear {S}chr\"{o}dinger equation}, 
	year={2024}
}

@article {Raphael2005MathAnnalen,
	AUTHOR = {Raphael, Pierre},
	TITLE = {Stability of the log-log bound for blow up solutions to the
	critical non linear {S}chr\"odinger equation},
	JOURNAL = {Math. Ann.},
	FJOURNAL = {Mathematische Annalen},
	VOLUME = {331},
	YEAR = {2005},
	NUMBER = {3},
	PAGES = {577--609},
	ISSN = {0025-5831,1432-1807},
	MRCLASS = {35Q55 (35B35 35B40)},
	MRNUMBER = {2122541},
	MRREVIEWER = {Masahito\ Ohta},
	DOI = {10.1007/s00208-004-0596-0},
	URL = {https://doi.org/10.1007/s00208-004-0596-0},
}

@article {MerleRaphael2005CMP,
	AUTHOR = {Merle, Frank and Raphael, Pierre},
	TITLE = {Profiles and quantization of the blow up mass for critical
	nonlinear {S}chr\"odinger equation},
	JOURNAL = {Comm. Math. Phys.},
	FJOURNAL = {Communications in Mathematical Physics},
	VOLUME = {253},
	YEAR = {2005},
	NUMBER = {3},
	PAGES = {675--704},
	ISSN = {0010-3616,1432-0916},
	MRCLASS = {35Q55 (35B40)},
	MRNUMBER = {2116733},
	MRREVIEWER = {Justin\ A.\ Holmer},
	DOI = {10.1007/s00220-004-1198-0},
	URL = {https://doi.org/10.1007/s00220-004-1198-0},
}

@article {MerleRaphael2005AnnMath,
	AUTHOR = {Merle, Frank and Raphael, Pierre},
	TITLE = {The blow-up dynamic and upper bound on the blow-up rate for
	critical nonlinear {S}chr\"odinger equation},
	JOURNAL = {Ann. of Math. (2)},
	FJOURNAL = {Annals of Mathematics. Second Series},
	VOLUME = {161},
	YEAR = {2005},
	NUMBER = {1},
	PAGES = {157--222},
	ISSN = {0003-486X,1939-8980},
	MRCLASS = {35Q55 (35B40)},
	MRNUMBER = {2150386},
	MRREVIEWER = {John\ Albert},
	DOI = {10.4007/annals.2005.161.157},
	URL = {https://doi.org/10.4007/annals.2005.161.157},
}

@article {MerleRaphael2003GAFA,
	AUTHOR = {Merle, F. and Raphael, P.},
	TITLE = {Sharp upper bound on the blow-up rate for the critical
	nonlinear {S}chr\"odinger equation},
	JOURNAL = {Geom. Funct. Anal.},
	FJOURNAL = {Geometric and Functional Analysis},
	VOLUME = {13},
	YEAR = {2003},
	NUMBER = {3},
	PAGES = {591--642},
	ISSN = {1016-443X,1420-8970},
	MRCLASS = {35Q55 (35B40)},
	MRNUMBER = {1995801},
	DOI = {10.1007/s00039-003-0424-9},
	URL = {https://doi.org/10.1007/s00039-003-0424-9},
}

@article {MerleRaphael2004Invent,
	AUTHOR = {Merle, Frank and Raphael, Pierre},
	TITLE = {On universality of blow-up profile for {$L^2$} critical
	nonlinear {S}chr\"odinger equation},
	JOURNAL = {Invent. Math.},
	FJOURNAL = {Inventiones Mathematicae},
	VOLUME = {156},
	YEAR = {2004},
	NUMBER = {3},
	PAGES = {565--672},
	ISSN = {0020-9910,1432-1297},
	MRCLASS = {35Q55 (35B33 35B40 35B65 35Q51)},
	MRNUMBER = {2061329},
	MRREVIEWER = {Masahito\ Ohta},
	DOI = {10.1007/s00222-003-0346-z},
	URL = {https://doi.org/10.1007/s00222-003-0346-z},
}

@article {MerleRaphael2006JAMS,
	AUTHOR = {Merle, Frank and Raphael, Pierre},
	TITLE = {On a sharp lower bound on the blow-up rate for the {$L^2$}
	critical nonlinear {S}chr\"odinger equation},
	JOURNAL = {J. Amer. Math. Soc.},
	FJOURNAL = {Journal of the American Mathematical Society},
	VOLUME = {19},
	YEAR = {2006},
	NUMBER = {1},
	PAGES = {37--90},
	ISSN = {0894-0347,1088-6834},
	MRCLASS = {35Q55 (35B40)},
	MRNUMBER = {2169042},
	MRREVIEWER = {Masahito\ Ohta},
	DOI = {10.1090/S0894-0347-05-00499-6},
	URL = {https://doi.org/10.1090/S0894-0347-05-00499-6},
}

@article {Matsuno2023,
	AUTHOR = {Matsuno, Yoshimasa},
	TITLE = {Multiphase solutions and their reductions for a nonlocal
	nonlinear {S}chr\"odinger equation with focusing nonlinearity},
	JOURNAL = {Stud. Appl. Math.},
	FJOURNAL = {Studies in Applied Mathematics},
	VOLUME = {151},
	YEAR = {2023},
	NUMBER = {3},
	PAGES = {883--922},
	ISSN = {0022-2526,1467-9590},
	MRCLASS = {35Q55},
	MRNUMBER = {4655451},
	MRREVIEWER = {Luigi\ Forcella},
	DOI = {10.1111/sapm.12610},
	URL = {https://doi.org/10.1111/sapm.12610},
}

@article {ChenPelinovsky2025,
    AUTHOR = {Chen, Jinbing and Pelinovsky, Dmitry E.},
     TITLE = {Traveling periodic waves and breathers in the nonlocal
              derivative {NLS} equation},
   JOURNAL = {Nonlinearity},
  FJOURNAL = {Nonlinearity},
    VOLUME = {38},
      YEAR = {2025},
    NUMBER = {7},
     PAGES = {Paper No. 075016, 39},
      ISSN = {0951-7715,1361-6544},
   MRCLASS = {35C07 (35Q55)},
  MRNUMBER = {4920349},
       DOI = {10.1088/1361-6544/addfa2},
       URL = {https://doi.org/10.1088/1361-6544/addfa2},
}

@article {BF2023PhysD,
	AUTHOR = {Berntson, Bjorn K. and Fagerlund, Alexander},
	TITLE = {A focusing-defocusing intermediate nonlinear {S}chr\"odinger
	system},
	JOURNAL = {Phys. D},
	FJOURNAL = {Physica D. Nonlinear Phenomena},
	VOLUME = {451},
	YEAR = {2023},
	PAGES = {Paper No. 133762, 26},
	ISSN = {0167-2789,1872-8022},
	MRCLASS = {35Q55},
	MRNUMBER = {4591360},
	DOI = {10.1016/j.physd.2023.133762},
	URL = {https://doi.org/10.1016/j.physd.2023.133762},
}

@article {Sun2024LMPsystem,
    AUTHOR = {Sun, Ruoci},
     TITLE = {The intertwined derivative {S}chr\"odinger system of
              {C}alogero-{M}oser-{S}utherland type},
   JOURNAL = {Lett. Math. Phys.},
  FJOURNAL = {Letters in Mathematical Physics},
    VOLUME = {114},
      YEAR = {2024},
    NUMBER = {3},
     PAGES = {Paper No. 74, 32},
      ISSN = {0377-9017,1573-0530},
   MRCLASS = {37K10 (47B35)},
  MRNUMBER = {4754345},
       DOI = {10.1007/s11005-024-01815-x},
       URL = {https://doi.org/10.1007/s11005-024-01815-x},
}

@article{BCD2024arxivNumerical,
	title={Spectrally accurate fully discrete schemes for some nonlocal and nonlinear integrable {PDE}s via explicit formulas}, 
	author={Bronsard, Yvonne Alama and Chen, Xi and Dolbeault, Matthieu},
	year={2024},
	journal = {preprint arXiv:2412.13480, to appear in Found. Comput. Math}
}

\end{document}